\documentclass[12pt,oneside,openany,headsepline,halfparskip,normalheadings,leqno]{scrbook}
\usepackage[sort,numbers]{natbib}
\usepackage[a4paper,text={14cm,22cm},top=3.5cm,left=3.5cm]{geometry}
\usepackage{scrpage2,setspace}
\usepackage[latin1]{inputenc}
\usepackage[T1]{fontenc}
\usepackage{mathptmx}
\usepackage[scaled=.90]{helvet}
\usepackage{courier}
\usepackage{amssymb,amsmath,bbm}
\usepackage{euro}
\usepackage{graphicx}
\usepackage{epic}
\usepackage{pict2e}
\usepackage{longtable}
\usepackage{afterpage}
\newtheorem{theo}{Theorem}
\newtheorem{pro}{Proposition}
\newtheorem{defi}{Definition}
\newtheorem{lem}{Lemma}
\newtheorem{cor}{Corollary}
\newtheorem{rem}{Remark}
\newcounter{fussall}

\usepackage{Diplomarbeit} % MYSTYLE
\titlehead{
  \begin{center}
  {\Large }
    \end{center}
  \begin{center}    
      { \LARGE Philipp Hungerländer\\}
  \end{center}
  \begin{center}
  {\Large }
    \end{center}
  \begin{center}
      { \LARGE Discrete-Time Dynamic\\ Noncooperative Game Theory}\\
  \end{center}
}
\subject{\LARGE \sc Diplomarbeit}
\author{\large Zur Erlangung des akamdemischen Grades}
\title{}
\date{\large
 Magister der Sozial- und Wirtschaftswissenschaften\\
 \large Angewandte Betriebswirtschaft\\
 \large Alpen-Adria-Universität Klagenfurt\\
 \large Fakultät für Wirtschaftswissenschaften}
\publishers{\large\begin{flushleft} 
  Begutachter: O.Univ.-Prof. Mag. Dr. Reinhard Neck\\
  Institut: Volkswirtschaftslehre \end{flushleft}
  \begin{flushright} 04/2008 \end{flushright}   }

\begin{document}

%
% Vorspann
%

\frontmatter % Roman-Nummerierung, S. 1
\maketitle
Ich erkläre ehrenwörtlich, dass ich die vorliegende wissenschaftliche Arbeit selbst\-ständig
angefertigt und die mit ihr unmittelbar verbundenen Tätigkeiten selbst erbracht habe. Ich
erkläre weiters, dass ich keine anderen als die angegebenen Hilfsmittel benutzt habe.
Alle aus gedruckten, ungedruckten oder dem Internet im Wortlaut oder im wesentlichen
Inhalt übernommenen Formulierungen und Konzepte sind gemäß den Regeln für
wissenschaftliche Arbeiten zitiert und durch Fußnoten bzw. durch andere genaue
Quellenangaben gekennzeichnet.\\

Die während des Arbeitsvorganges gewährte Unterstützung einschließlich signifikanter
Betreuungshinweise ist vollständig angegeben.\\

Die wissenschaftliche Arbeit ist noch keiner anderen Prüfungsbehörde vorgelegt worden.
Diese Arbeit wurde in gedruckter und elektronischer Form abgegeben. Ich bestätige, dass
der Inhalt der digitalen Version vollständig mit dem der gedruckten Version
übereinstimmt.\\

Ich bin mir bewusst, dass eine falsche Erklärung rechtliche Folgen haben wird.\\\\

(Unterschrift)  \;\;\;\;\;\;\;\;\;\;\;\;\;\;\;\;\;\;\;\;\;\;\;\;\;\;\;\;\;\;\;\;\;\;\;\;\;\;\;\;\;\;\;\;\;\;\;\;\;\;\;\;\;\;\;\;\;\;\;\;\;\;\;\;\;\;\;\;\;\;\;\;\;\;\;(Ort, Datum) 

\setcounter{page}{1}
\tableofcontents
%\listoffigures
%\addcontentsline{toc}{chapter}{Abbildungsverzeichnis}
%\listoftables
%\addcontentsline{toc}{chapter}{Tabellenverzeichnis}
%\input{dp_bez}

%
% Hauptteil
%

\mainmatter % Arabic-Nummerierung, S. 1
\setcounter{fussall}{\thefootnote}
\chapter {Introduction}
\setcounter{footnote}{\thefussall}

Dynamic noncooperative game theory is a field of mathematics and economics in which a lot of research is being carried out at present featuring a great number of applications in many different areas of economics and management science like:

\begin{itemize}
\item capital accumulation and investments,
\item R\&D and technological innovations,
\item macroeconomics,
\item microeconomics,
\item pricing and advertising decisions in marketing,
\item natural resource extraction
\item pollution control \footnote{For studies in the above-mentioned areas cf. e.g. Dockner et al. (2000) \cite{doc}.}\\\\
\end{itemize}

The aim of this diploma thesis is

\begin{itemize}
\item to correct some results for discrete-time affine-quadratic dynamic games of prespecified fixed duration with open-loop and feedback information patterns (cf. Subsections \eqref{ssaqcp}, \eqref{lqfsg} and \eqref{sslqols}).
\item to give shorter and more convenient proofs for some results already stated in the literature (cf. Subsection \eqref{ssaqoln}).
\item to present some extensions for the open-loop and feedback Stackelberg equilibrium solutions of discrete-time affine-quadratic dynamic games of prespecified fixed duration, concerning the number of followers, the structure of the cost und state functions and the possibility of an algorithmic disintegration (cf. Subsections \eqref{ssaqfsg}, \eqref{OLS_ss2i} and \eqref{OLS_1i}).
\end{itemize}

\clearpage
\setcounter{chapter}{1}
\setcounter{fussall}{\thefootnote}
\chapter{Basic Definitions and Basic Insights}
\setcounter{footnote}{\thefussall}
%\section{Introduction}
%
%
%
%\clearpage
\section{Basic Definitions}

In this section the central notions are defined for discrete-time dynamic noncooperative games that will be used permanently throughout the next chapters.\footnote{ The definitions are geared to the ones given in Ba\c{s}ar and Olsder (1999)\cite{baol}. They were modified insofar as it was helpful to keep the diploma thesis consistent.}

\subsection{Game structure}

In this subsection definitions are given for the kinds of games, information structures and cost functionals examined in this paper.

\begin{defi}
An \emph {\textit{n}-person discrete-time deterministic infinite dynamic game (also known as an \textit{n}-person deterministic multi-stage game) of prespecified fixed duration} involves
\begin{enumerate}
\item An index set $N :=\{1,\ldots,n\}$ called the \emph {players' set}.
\item An index set $K:=\{1,\ldots,T\}$ denoting the \emph {stages} of the game, where T is the maximum possible number of moves a player is allowed to make in the game
\item An infinite set X with some topological structure, called the \emph {state set (space)} of the game, to which the state of the game ($x_{k-1}$) belongs for all $k \in K $.
\item An infinite set $U_k^i$ with some topological structure, defined for each $k \in K $ and $i \in N$, which is called the \emph {action (control) set} of player i (\textbf{P}i) at stage k. Its elements are the permissible actions $u_k^i$ of \textbf{P}i at stage k.
\item A function $f_k : X \times U_k^1 \times \ldots \times U_k^n \to X$, defined by
\begin{align}
\nonumber x_k = f_{k-1}(x_{k-1},u_k^1,\ldots,u_k^n) \; , \; x_0 \in X (initial\; state) \; , \; k \in K
\end{align}
which is called the \emph {state equation} of the dynamic game. It describes the evolution of the underlying decision process.
\item A set $Y_k^i$ with some topological structure (defined for: $k \in K$ , $i \in N$) called the \emph {observation set} of \textbf{P}i at stage k, to which the oberservation $y_k^i$ of \textbf{P}i belongs at stage k.
\item A function $h_k^i : X \to Y_k^i$ (defined for: $k \in K$ , $i \in N$) given by
\begin{align}
\nonumber y_k^i = h_k^i(x_{k-1}) \; , \; i \in N \; , \; k \in K
\end{align}
which is the \emph {state-measurement (-observation) equation} of \textbf{P}i concerning the value of $x_{k-1}$.
\item A finite set $\eta_k^i$ (defined for $k \in K$ , $i \in N$) as a subset of $\{y_1^1 , \ldots , y_k^1; \ldots ;\\ y_1^n , \ldots , y_k^n ; u_1^1, \ldots , u_{k-1}^1 ; \ldots ; u_1^n , \ldots , u_{k-1}^n\}$, which determines the information gained and recalled by \textbf{P}i at stage k of the game. The specification of $\eta_k^i$ characterizes the \emph {information structure (pattern)} of \textbf{P}i, and the collection of these information structures for all $i \in N$ is the information structure of the game.
\item A set $N_k^i$ (defined for $k \in K$ , $i \in N$) as a subset of $\{(Y_1^1 \times \ldots \times Y_k^1) \times \ldots \times (Y_1^n \times \ldots \times Y_k^n) \times (U_1^1 \times \ldots \times U_{k-1}^1) \times \ldots \times (U_1^n \times \ldots \times U_{k-1}^n)\}$ designed to be compatible with $\eta_k^i$. $N_k^i$ is called the \emph {information space} of \textbf{P}i at stage k, induced by his information $\eta_k^i$.  
\item A prespecified class $\Gamma_k^i$ (defined for $k \in K$ , $i \in N$) of mappings  $\gamma_k^i: N_k^i \to U_k^i$ which are the \emph {permissible strategies} of \textbf{P}i at stage k of the game. The aggregate mapping $\gamma^i  = \{\gamma_1^i,\ldots,\gamma_T^i\}$ is a strategy of \textbf{P}i in the game. Furthermore the class $\Gamma^i$ of all mappings $\gamma^i$ is the \emph {strategy set (space)} of \textbf{P}i.
\item A functional $L^i : (X \times U_1^1 \times \ldots \times U_1^n) \times (X \times U_2^1 \times \ldots \times U_2^n) \times \ldots \times (X \times U_T^1 \times \ldots \times U_T^n) \to \textbf{R}$ (defined for $i \in N$) called the \emph {cost functional} of \textbf{P}i in the game of fixed duration.
\end{enumerate}
\label{defspiel}
\end{defi}

\begin{defi}
In an n-person discrete-time deterministic infinite dynamic game of prespecified fixed duration (cf. Def. \eqref{defspiel}), \textbf{P}i's  $(i \in N)$ information structure is called a(n)
\begin{enumerate}
\item \emph {open-loop} (OL) patern if $\eta_k^i = \{x_0\}, (k \in K),$
\item \emph {closed-loop perfect state information} (CLPS) pattern if $\eta_k^i = \{x_0,\ldots,x_{k-1}\}, \\(k \in K),$
\item \emph {closed-loop imperfect state information} (CLIS) pattern if $\eta_k^i = \{y_1^i,\ldots,y_{k}^i\},\\ (k \in K),$
\item \emph {memoryless perfect state information} (MPS) pattern if $\eta_k^i = \{x_0,x_{k-1}\}, \\(k \in K),$
\item \emph {feedback (perfect state) information} (FB) pattern if $\eta_k^i = \{x_{k-1}\}, \\(k \in K),$
\end{enumerate}
\label{definfo}
\end{defi}

\begin{defi}
For an n-person discrete-time deterministic infinite dynamic game of prespecified fixed duration (cf. Def. \eqref{defspiel}), \textbf {P}i's $(i \in N)$ cost functional is said to be \emph {stage-additive} if there exist $g_k^i: X \times X \times U_k^1 \times \ldots \times U_k^n \to \textbf R , (k \in K),$ so that ($i,j \in N$)
\begin{align}
\nonumber L^i(x_0,u^1,\ldots ,u^n) = \sum_{k=1}^T g_k^i(x_{k},u_k^1,\ldots ,u_k^n,x_{k-1}), where\;\; u^j \; \widehat{=} \; (u_1^j,\ldots ,u_T^j).
\end{align}
Furthermore, if $L^i(u^1,\ldots ,u^n)$ depends only on $x_{T}$ (the terminal state), then it is called \emph {terminal cost function}.\label{defcost}
\end{defi}

\begin{defi}
An n-person discrete-time deterministic infinite dynamic game of prespecified fixed duration (cf. Def. \eqref{defspiel}) is of \emph {affine-quadratic} type if

\begin{multline}
\nonumber f_{k-1}(x_{k-1},u_k^1,\ldots,u_k^n) = A_kx_{k-1} + \sum_{j\in N}B_k^ju_k^j+s_k \\
\end{multline}

\vspace*{-5ex}

\begin{multline}
\nonumber L^i(u^{1},\ldots,u^{n})= \sum_{k = 1}^T g_k^i(x_k,u_k^{1},\ldots,u_k^{n},x_{k-1}) \\
\end{multline}

\vspace*{-5ex}

\begin{multline}
 \nonumber g_k^i(x_k,u_k^1,\ldots,u_k^n,x_{k-1}) = \frac{1}{2} (x_k^{'}Q_k^ix_k + \sum_{j\in N}u_k^{j'}R_k^{ij}u_k^j)\; \\ +
 \frac{1}{2} (\tilde{x}_k^{i'}Q_k^i\tilde{x}_k^i + \sum_{j\in N}\tilde{u}_k^{ij'}R_k^{ij}\tilde{u}_k^{ij}) - \tilde{x}_k^{i'}Q_k^ix_k - \sum_{j\in N}\tilde{u}_k^{ij'}R_k^{ij}u_k^j \\
\end{multline}

where $U_k^i = \textbf{R}^{m_i}$ and $s_k \in \textbf{R}^p$. $A_k, B_k^i, Q_k^i, R_k^{ij}, \tilde{u}_k^{ij}, \tilde{x}_k^i$ (defined for $k \in K$ , $i \in N$, $j \in N$) are fixed sequences of matrices or vectors of appropriate dimensions. Furthermore $Q_k^i$ and $R_k^{ij}$ are symmetric. An affine-quadratic game is of the \emph {linear-quadratic} type if $s_k \equiv 0$.
\label{defspielaq}
\end{defi}

\begin{rem}
The cost function of player i at stage k $[g_k^i(x_k,u_k^1,\ldots,u_k^n,x_{k-1})]$ can also be written in the following way

\begin{align}
\nonumber \frac{1}{2} ((x_k^{'}-\tilde{x}_k^{i'})Q_k^i(x_k-\tilde{x}_k^{i'}) + \sum_{j\in N}(u_k^{j'}-\tilde{u}_k^{ij'})R_k^{ij}(u_k^j-\tilde{u}_k^{ij'})) 
\end{align}

Therefore $\tilde{x}_k^{i'}$ and $\tilde{u}_k^{ij'}$ can be interpreted as desired (target) values of each player for all variables of the game.
\end{rem}
\clearpage

\subsection{Solution concepts}

In this subsection the Nash and Stackelberg equilibrium solution concepts are introduced for an \textit{n}-person discrete-time deterministic infinite dynamic game of prespecified fixed duration. The Nash equilibrium solution concept provides a reasonable noncooperative equilibrium solution when no single player dominates the decision making process and therefore the roles of the players are symmetric. However, there are other types of noncooperative decision problems in which one of the players is a so-called leader and has the ability to enforce his strategy on the other players, the so-called followers. For that kind of decision problems a hierarchical equilibrium solution concept, the Stackelberg equilibrium solution, is introduced.

\begin{defi}
An n-tuple of strategies $\{ \gamma^{1*},\gamma^{2*},\ldots ,\gamma^{n*}\}$ with $\gamma^{i*} \in \Gamma^i, i \in N,$ is said to constitute a \emph {noncooperative Nash equilibrium solution} for an n-person discrete-time deterministic infinite dynamic game of prespecified fixed duration (cf. Def. \eqref{defspiel}) if the following \textit{n} inequalities are satisfied for all $\gamma^i \in \Gamma^i, i \in N:$\\
\begin{align}
\nonumber L^{1*}\; \widehat = \; L^1(\gamma^{1*}, \gamma^{2*}, \gamma^{3*},\ldots,\gamma^{n*}) \leq L^1(\gamma^{1},\gamma^{2*}, \gamma^{3*},\ldots,\gamma^{n*}) 
\end{align}
\begin{align}
\nonumber L^{2*}\; \widehat = \;L^2(\gamma^{1*},\gamma^{2*}, \gamma^{3*},\ldots,\gamma^{n*}) \leq  L^2(\gamma^{1*},\gamma^{2}, \gamma^{3*},\ldots,\gamma^{n*}) 
\end{align}
\begin{align}
\nonumber \cdots
\end{align}
\begin{align}
\nonumber \cdots
\end{align}
\begin{align}
\nonumber L^{n*}\; \widehat = \;L^n(\gamma^{1*},\gamma^{2*},\ldots,\gamma^{n-1*},\gamma^{n*}) \leq  L^n(\gamma^{1*},\gamma^{2},\ldots,\gamma^{n-1*},\gamma^{n}) 
\end{align}
\begin{align}
\nonumber where \; \; \gamma^i \; \widehat{=} \; (\gamma_1^i,\ldots ,\gamma_T^i) 
\end{align}
The n-tuple of quantities $\{L^{1*},\ldots,L^{n*}\}$ is known as a \emph {Nash equilibrium outcome} of the discrete-time deterministic infinite dynamic game of prespecified fixed duration.
\end{defi}

\begin{defi}
In an n-person discrete-time deterministic infinite dynamic game of prespecified fixed duration (cf. Def. \eqref{defspiel}) with \textbf{P}1 as the leader the unique element $r^i(\gamma^1) \in \Gamma^i (i \in N, i \not = 1)$ defined for each $\gamma^1 \in \Gamma^1$ by
\begin{multline}
\nonumber r^i(\gamma^1)= \min_{\gamma^i \in \Gamma^i} L^i(\gamma^1,r^2(\gamma^1),\ldots,r^{i-1}(\gamma^1),\gamma^i,r^{i+1}(\gamma^1),\ldots,r^{n}(\gamma^1))\\
\nonumber where \; \; \gamma^l \; \widehat{=} \; (\gamma_1^l,\ldots ,\gamma_T^l) \; and \; r^i \; \widehat{=} \; (r_1^i,\ldots ,r_T^i) \; ; \; l \in N \\
\end{multline}

is called the \emph {unique optimal response (rational reaction)} of \textbf{P}i to the strategy $\gamma^1 \in \Gamma^1$ of  \textbf{P}1.
\end{defi}

\begin{defi}
In an n-person discrete-time deterministic infinite dynamic game of prespecified fixed duration (cf. Def. \eqref{defspiel}) with \textbf{P}1 as the leader, a strategy  $\gamma^{1*} \in \Gamma^1$ is called a \emph{Stackelberg equilibrium strategy} for the leader if
\begin{multline}
\nonumber L^{1*}\; \widehat = \; L^1(\gamma^{1*},r^2(\gamma^{1*}),\ldots,r^n(\gamma^{1*})) \leq L^1(\gamma^1,r^2(\gamma^1),\ldots r^{n}(\gamma^1)) , \\ \forall \gamma^1 \in \Gamma^1, \forall r^i(\cdot) \in R^i(\cdot), i \in N, i \not = 1\} \\
\nonumber where \; \; \gamma^1 \; \widehat{=} \; (\gamma_1^1,\ldots ,\gamma_T^1) \; and \; r^i \; \widehat{=} \; (r_1^i,\ldots ,r_T^i) \\
\end{multline}
The quantity $L^{1*}$ is called the \emph{Stackelberg cost} of the leader.
\end{defi}

\begin{defi}
In an n-person discrete-time deterministic infinite dynamic game of prespecified fixed duration (cf. Def. \eqref{defspiel}) with \textbf{P}1 as the leader, the element $\gamma^{i*} \in R^i(\gamma^{1*}) (i \in N, i \not = 1)$ is called an \emph{unique optimal strategy} for follower i that is \emph{in equilibrium} with $\gamma^{1*}$. The n-tuple $(\gamma^{1*}, \gamma^{2*},\ldots,\gamma^{n*})$ is a \emph{Stackelberg solution} for the n-person discrete-time deterministic infinite dynamic game of prespecified fixed duration with \textbf{P}1 as the leader, and the n-tuple of quantities $\{L^{1*}, L^{2*},\ldots,L^{n*}\}$ is known as the corresponding \emph {Stackelberg equilibrium outcome}.
\end{defi}

\clearpage

\subsection{Time consistency}

The issue of time consistency has pervaded the economics literature during the past three decades, following the important paper by Kydland and Prescott (1977) \cite{kypr}. Based on Ba\c{s}ar's paper (1989) \cite{bas} we impose further refinements on the class of equilibrium strategies drawing a distinction between time inconsistent, weakly and strongly time consistent optimal strategies that will be used later on to characterize the quality of different equilibrium solution concepts under different information patterns.

%\begin{defi}
%For an n-person discrete-time deterministic infinite dynamic game of prespecified fixed duration (cf. Def. \eqref{defspiel}) let the strategies of all players, except \textbf{P}i $(i \in N)$ , be fixed at $\gamma^j \in \Gamma^j,j \in N, j \not = i$. Then, a strategy $\gamma^i \in \Gamma^i$ for \textbf{P}i is a \emph {representation} of another strategy $\tilde{\gamma}^i \in \Gamma^i$, if
%\begin{enumerate}
%\item the n-tuples $\{\gamma^i,\gamma^j\}$ and $\{\tilde{\gamma}^i,\gamma^j\}$ generate the same unique state trajectory, and
%\item $\gamma^i$ and $\tilde{\gamma}^i$ have the same open-loop value on this trajectory.
%\end{enumerate}
%\end{defi}

\begin{defi}
An n-tuple of policies (equilibrium strategies) $(\gamma^{1*},\ldots,\gamma^{n*}) \in \Gamma$ solving the dynamic game defined in Definition \eqref{defspiel} for any particular information pattern defined in Definition \eqref{definfo} is \emph{weakly time consistent} (WTC) if its trunctation of stages to the interval [s,T] (for $s \in [1,\ldots,T]$), $(\gamma_{[s,T]}^{1*},\ldots,\gamma_{[s,T]}^{n*})$ solves the truncated game. If an n-tuple of policies $(\gamma^{1*},\ldots,\gamma^{n*}) \in \Gamma$ is not WTC, then it is called \emph {time inconsistent}.
\end{defi}

\begin{defi}
An n-tuple of policies (optimal strategies) $(\gamma^{1*},\ldots,\gamma^{n*}) \in \Gamma$ solving the dynamic game defined in Definition \eqref{defspiel} for any particular information pattern defined in Definition \eqref{definfo} is \emph{strongly time consistent} (STC) if its trunctation of stages to the interval [s,T] (for $s \in [1,\ldots,T]$), $(\gamma_{[s,T]}^{1*},\ldots,\gamma_{[s,T]}^{n*})$ solves the truncated game for every permissible n-tuple of strategies played in the interval [0,s) [$\widehat {=} (\gamma_{[0,s)}^{1},\ldots,\gamma_{[0,s)}^{n}) \in \Gamma$].
\end{defi}

\clearpage
\section{Basic insights}

This section presents the most important optimization tools for discrete-time dynamic noncooperative game theory and some results about matrix identities and the definiteness of matrices that will be used throughout the next chapters.

\subsection{Dynamic programming}

The method of \emph {dynamic programming} was developed by Richard Bellman (cf. Bellman (1957) \cite{bel}) and is a tool for solving games with feedback information pattern. It is based on the \emph{principle of optimality}, which states that an optimal strategy has the property that, whatever the initial state and time are, all remaining decisions (from that particular initial state and time onwards) must also constitute an optimal strategy. To make use of this principle in a mathematical dynamic (game theoretic) framework, we have to work backwards in time, starting at all possible final states with the corresponding final times.\\
\\
To apply the \emph{principle of optimality} to our game theoretic framework (cf. Def. \eqref{defspiel}), we have to consider a stage-additive cost functional (cf. Def. \eqref{defcost}) and feedback (perfect state) information pattern (cf. Def. \eqref{definfo}):

\begin{align}
\nonumber L^i(u^1,\ldots ,u^n) = \sum_{k=1}^T g_k^i(x_{k},u_k^1,\ldots ,u_k^n,x_{k-1}), 
\end{align}

\vspace*{-5ex}

\begin{align}
\nonumber where\;\; u^j \widehat{=} (u_1^j,\ldots ,u_T^j)\;\;and\;\;u_k^i = \gamma_k^i(x_{k-1})\;;\;j \in N
\end{align} 

\vspace*{-5ex}

\begin{align}
\nonumber x_k = f_{k-1}(x_{k-1},u_k^1,\ldots,u_k^n)\; , \; k \in K
\end{align}

On this basis we can define an expression for the minimal cost of \textbf {P}i $(i \in N)$ for any starting point and any corresponding intial time.

\begin{defi}
For an n-person discrete-time deterministic infinite dynamic game of prespecified fixed duration (cf. Def. \eqref{defspiel}), \textbf {P}i's $(i \in N)$ \emph{value function} is defined as:

\begin{align}
\nonumber V^i(k,x_{k-1}) = \min_{u^i_k \ldots u_T^i} \sum_{j=k}^T g_j^i(x_{j},u_j^1,\ldots ,u_j^n,x_{j-1}) 
\end{align}

\end{defi}

Because of the principle of optimality, the value function is equivalent to the following recursive relation:

\vspace*{-5ex}

\begin{multline}
\nonumber V^i(k,x_{k-1}) = \min_{u_k^i \in U_k^i}[g_k^i(f_{k-1}(x_{k-1},u_k^1,\ldots,u_k^n),u_k^1,\ldots ,u_k^n,x_{k-1})\\ +  V^i(k+1,f_{k-1}(x_{k-1},u_k^1,\ldots,u_k^n)] \\
\end{multline}

\vspace*{-5ex}

Solving dynamic problems by using this recursive relation is known as \emph {dynamic programming}.

\clearpage

\subsection{The minimum principle}

\emph{The minimum principle} was developed (in continuous time) by Lew Semjonowitsch Pontryagin (cf. Pontryagin et al. (1962) \cite{pon}) and is a tool for solving games with open-loop (and feedback) information pattern. In the version presented below sufficient conditions for the existence of optimal solutions are given because the following Theorem \eqref{minprin} will be applied (via Theorem \eqref{OLN_opt}) to affine-quadratic games which fulfill the stated sufficient conditions. The derivation of Theorem \eqref{minprin} can be found in either Canon et al. (1970) \cite{can} or Boltyanski (1978) \cite{bol}.

%Theorem 5.5 Seite 246, auch schaun, ob Proposition 5.1 als feedback rep. der ol lösung ...

\begin{theo}
For the discrete-time optimal control problem (cf. Def. \eqref{defspiel} with N = \{1\}) let
\begin{itemize}
\item $f_k(\cdot,\cdot)$ be continuously differentiable on $\textbf{R}^p \times \textbf{R}^{m}$ (defined for: $k \in K$)
\item $g_k(\cdot,\cdot, \cdot)$ be continuously differentiable on $\textbf{R}^p \times \textbf{R}^{m} \times \textbf{R}^p$ (defined for: $k \in K$)
\item $f_k(\cdot,\cdot)$ be convex on $\textbf{R}^p \times \textbf{R}^{m}$ (defined for: $k \in K$)
\item $g_k(\cdot,\cdot, \cdot)$ be convex on $\textbf{R}^p \times \textbf{R}^{m} \times \textbf{R}^p$ (defined for: $k \in K$)
\item the cost function be stage-additive (cf. Def. \eqref{defcost}).
\end{itemize}
Then $\{\gamma_k^{*}(x_0) = u_k^{*}\}$ denotes an \emph {optimal control sequence} and $\{x_k^*; k \in K\}$ is the corresponding state trajectory and there exists a finite sequence of p-dimensional costate vectors $\{p_1,\ldots,p_T\}$ so that the following relations are satisfied:\\

\vspace*{-5ex}

\begin{align}
x_k^* = f_{k-1}(x_{k-1}^*,u_k^{*}) \; , \; x_0^* = x_0 
\end{align}

\vspace*{-5ex}

\begin{align}
\nabla_{u_k} H_k(p_k, u_k^*, x_{k-1}^*) = 0
\end{align}

\vspace*{-5ex}

\begin{multline}
p_k = \frac {\partial}{\partial x_k} f_{k}(x_{k}^*,u_{k+1}^{*})^{'}[p_{k+1} + (\frac {\partial}{\partial x_{k+1}}g_{k+1}(x_{k+1}^*, u_{k+1}^{*}, x_{k}^*))^{'}]\\ + [\frac {\partial}{\partial x_{k}}g_{k+1}(x_{k+1}^*, u_{k+1}^{*}, x_{k}^*)]^{'} \; ; \; p_T = 0 \\
\end{multline}

\vspace*{-5ex}

where 

\vspace*{-5ex}

\begin{multline}
H_k(p_k, u_k, x_{k-1}) \; \widehat{=} \; g_k(f_{k-1}(x_{k-1},u_{k}), u_{k}, x_{k-1}) + p_{k}^{'}f_{k-1}(x_{k-1},u_{k}) \\
\end{multline}

\label{minprin}

\end{theo}

\clearpage

\subsection{Some results about the definiteness of matrices}

This subsection is devoted to the derivation of Corollary \eqref{pdc}, which is needed for convexity analyses later on.

\begin{lem}
$A+B > 0$ if

\begin{itemize}
\item $A > 0$, $B \geq 0$.\\
\end{itemize}

\label{AB}

\end{lem}

\textsc {Proof:}

\vspace*{-5ex}

\begin{multline}
x'(A + B)x = x'Ax + x'Bx > 0 \;\;\;\;\;\;\boxed{}\\\\
\end{multline}

\vspace*{-5ex}

\begin{lem}

$B'AB \geq 0$ if

\begin{itemize}
\item $A \geq 0$.\\
\end{itemize}

\label{BAB}
\end{lem}

\textsc {Proof:}

For showing

\vspace*{-5ex}

\begin{multline}
x'B'ABx \geq 0\;\;\; \forall x\\
\end{multline}

\vspace*{-5ex}

we define a linear mapping L

\vspace*{-5ex}

\begin{multline}
L: \tilde{x} = Bx \;\;\; \forall x\\
\end{multline}

\vspace*{-5ex}

If \textit{B} has full rank, the vectors $\tilde{x}$ cover the whole space and $\tilde{x}$ is only zero if the corresponding vector \textit{x} is also zero. If \textit{B} has not full rank, some vectors $\tilde{x}$ are mapped to zero for vectors \textit{x} unequal to zero. Therefore the vectors $\tilde{x}$ cover only a part of the space. In both cases

\vspace*{-5ex}

\begin{multline}
\tilde{x}'A\tilde{x} \geq 0 \;\;\; \forall \tilde{x}\\
\end{multline}

\vspace*{-5ex}

holds true because $A \geq 0$.\;\;\;\;\;\;\boxed{}\\\\

\begin{cor}
$A + C'BC > 0$ if
\begin{itemize}
\item $A > 0$, $B \geq 0$.\\
\end{itemize}
\label{pdc}
\end{cor}

\textsc {Proof:}

\textit{D}:=\textit{C'BC} has to be $\geq 0$ because of Lemma \eqref{BAB}. Now we can apply Lemma \eqref{AB} to \textit{A} + \textit{D}.\;\;\;\;\;\;\boxed{}\\

\clearpage

\subsection{Some results about matrix identities}

In this subsection two matrix identities are deduced. These are needed later on in some propositions which relate my findings with the results stated in Ba\c{s}ar and Olsder (1999) \cite{baol}.\\

First a general result is given concerning the eigenvalues of a product of two matrices with certain properties, whose derivation can be found in Horn and Johnson (1991, p. 465)\cite{hojo}.

\begin{lem}

The product of a positive definite matrix A and a Hermitian matrix B is a diagonalizable matrix, all of whose eigenvalues are real. The matrix AB has the same number of positive, negative and zero eigenvalues as B.
\label{Horn}
\end{lem}

Now Lemma \eqref{Horn} is applied to deduce two matrix identities, which will be used several times in Proposition \eqref{bofs} and Proposition \eqref{bsolsäa}.

\begin{lem}

Let A be positive definite and B a be matrix of appropriate dimension, then the following two matrix identities hold:

\begin{enumerate}
			
\item $[I - AB(I + B'AB)^{-1}B'] \equiv (I + ABB')^{-1}$ \label{mi1}
			
\item $[I - B(I + B'AB)^{-1}B'A] \equiv (I + BB'A)^{-1}$\label{mi2}\\

\end{enumerate}
\label{maid}
\end{lem}

\textsc {Proof:}

We will prove matrix identity 1 and the proof of 2 can be done in essentially the same way.\\

At first note that if \textit{A} is positive definite, \textit{B'AB} is positive semidefinite with Lemma \eqref{BAB} and \textit{ABB}' is positive semidefinite with Lemma \eqref{Horn}.\footnote{Marcus and Minc (1964, p. 24) \cite{mami} show that both terms have the same eigenvalues.} Therefore the matrices $(I + B'AB)^{-1}$ and $(I + ABB')^{-1}$ exist.\\

Now we have to show that $[I - AB(I + B'AB)^{-1}B']$ multiplied by $(I + ABB')$ gives the identity matrix \textit{I}.

\vspace*{-5ex}

\begin{multline}
\nonumber [I - AB(I + B'AB)^{-1}B'](I + ABB') \\
\end{multline}

\vspace*{-5ex}

Expanding and prescinding yields

\vspace*{-5ex}

\begin{multline}
\nonumber I + ABB' - AB(I + B'AB)^{-1}B' - AB(I + B'AB)^{-1}B'ABB' \\
\end{multline}

\vspace*{-5ex}
\vspace*{-5ex}

\begin{multline}
\nonumber I  + AB[I - (I + B'AB)^{-1} - (I + B'AB)^{-1}B'AB]B' \\
\end{multline}

\vspace*{-5ex}
\vspace*{-5ex}

\begin{multline}
\nonumber I  + AB[I - (I + B'AB)^{-1}(I + B'AB)]B' = I + AB[I - I]B' = I \;\;\;\;\;\;\boxed{}\\
\end{multline}

\vspace*{-5ex}

\clearpage

\setcounter{chapter}{2}
\setcounter{fussall}{\thefootnote}
\chapter {Discrete-Time Infinite Dynamic Games with Feedback Information Pattern}
\setcounter{footnote}{\thefussall}
%\section{Introduction}
%
%\clearpage
\section {Feedback Nash Equilibrium Solutions}

This section is devoted to the derivation of the so-called feedback Nash equilibrium solution for affine-quadratic games. First a general result is stated about the existence and uniqueness of a feedback Nash equilibrium solution in \textit{n}-person discrete-time deterministic infinite dynamic games of prespecified fixed duration (cf. Def. \eqref{defspiel}) with feedback information pattern. Then this result is applied to affine-quadratic games and finally the solution for the affine-quadratic control problem is deduced as a special case.

\subsection{Optimality conditions}

In this subsection a theorem is stated which gives necessary and sufficient conditions for the existence of a feedback Nash equilibrium solution. Results about feedback Nash equilibria in infinite dynamic games first appeared in continuous time in the works of Starr and Ho (1969) \cite{sthob}, \cite{sthoa} and Case (1969) \cite{cas}. A proof of Theorem \eqref{FBN_Opt}, which is the counterpart of the above-mentioned results in discrete time, can be found in Ba\c{s}ar and Olsder (1999, pp. 278-279)\cite{baol}.
 
\begin{theo}
For an n-person discrete-time deterministic infinite dynamic game of prespecified fixed duration (cf. Def. \eqref{defspiel}) with feedback information pattern, the set of strategies $\{\gamma_k^{i*}(x_k);k \in  K , i \in N \}$ provides a \emph {feedback Nash equilibrium solution} if, and only if, functions $V^i(k,\cdot) : R^n \to R (k \in K, i \in N$) exist such that the following recursive relations are satisfied:\\

\vspace*{-5ex}

\begin{multline}
V^i(k,x_{k-1}) = \min_{u_k^i \in U_k^i}[g_k^i(\tilde{f}_{k-1}^{i*}(x_{k-1},u_k^i),\gamma_k^{1*}(x_{k-1}),\ldots,\gamma_k^{i-1*}(x_{k-1}),u_k^i,\\\gamma_k^{i+1*}(x_{k-1}),\ldots,\gamma_k^{n*}(x_{k-1}),x_{k-1}) +  V^i(k+1,\tilde{f}_{k-1}^{i*}(x_{k-1},u_k^i))]  =\\ 
g_k^i[\tilde{f}_{k-1}^{i*}(x_{k-1},\gamma_k^{i*}(x_{k-1})),\gamma_k^{1*}(x_{k-1}),\ldots,\gamma_k^{n*}(x_{k-1}),x_{k-1}] \\ + V^i[k+1,\tilde{f}_{k-1}^{i*}(x_{k-1},\gamma_k^{i*}(x_{k-1}))] \; ;\;V^i(T+1,x_T) = 0 \\\label{vfn}
\end{multline} \vspace*{-5ex}

\vspace*{-5ex}

\begin{multline}
\nonumber where \;\; \tilde{f}_{k-1}^{i*}(x_{k-1},z) \; \widehat{=} \\ f_{k-1}(x_{k-1},\gamma_k^{1*}(x_{k-1}),\ldots,\gamma_k^{i-1*}(x_{k-1}),z,\gamma_k^{i+1*}(x_{k-1}),\ldots,\gamma_k^{n*}(x_{k-1}) \\
\end{multline} \vspace*{-5ex}

Every such equilibrium solution is \emph {strongly time consistent}, and the corresponding Nash equilibrium cost for \textbf {P}i is $V^i(1,x_0)$.
\label{FBN_Opt}
\end{theo}

\clearpage

\subsection{Results for affine-quadratic games with arbitrarily many players}

In the following, the results of Theorem \eqref{FBN_Opt} are applied to an affine-quadratic dynamic game with arbitrarily many players. Theorem \eqref{FBNTheo}, which is an extension (concerning the cost functionals) of Corollary 6.1 in Ba\c{s}ar and Olsder (1999, pp. 279-281)\cite{baol}, presents equilibirum equations that can easily be used for an algorithmic disintegration of the given Nash game.\\

Furthermore in Proposition \eqref{profbn} the equivalence of the equations in Theorem \eqref{FBNTheo} with terminologically different equations is shown.

\begin{theo}
An n-person affine-quadratic dynamic game (cf. Def. \eqref{defspielaq}) admits a \emph{unique feedback Nash equilibrium solution} if  
\begin{itemize}
\item $Q_k^i \geq 0$, $R_k^{ii} > 0$ (defined for $k \in K$ , $i \in N$).
\item \eqref{P} and \eqref{alpha} admit unique optimal solution sets {$P_k^{i*}$} and {$\alpha_k^{i*}$} (defined for: $k \in K$ , $i \in N$)
\end{itemize}
If these conditions are satisfied, the unique equilibrium strategies are given by \eqref{gamman} and the corresponding feedback Nash equilibrium cost for each player is stated in \eqref{costn}. \footnote{ For all equations belonging to this theorem and its proof, $i \in N$ and $k \in K$ if nothing different is stated.}

\vspace*{-5ex}

\begin{align}
f_{k-1}(x_{k-1},u_k^1,\ldots,u_k^n) = A_kx_{k-1} + \sum_{j\in N}B_k^ju_k^j+s_k \label{fn}
\end{align}

\vspace*{-5ex} 

\begin{align}
L^i(x_0,u^{1},\ldots,u^{n})= \sum_{k = 1}^T g_k^i(x_k,u_k^{1},\ldots,u_k^{n},x_{k-1}) 
\end{align}

\vspace*{-5ex} \begin{multline}
 g_k^i(x_k,u_k^1,\ldots,u_k^n,x_{k-1}) = \frac{1}{2} (x_k^{'}Q_k^ix_k + \sum_{j\in N}u_k^{j'}R_k^{ij}u_k^j)\; \\ +
 \frac{1}{2} (\tilde{x}_k^{i'}Q_k^i\tilde{x}_k^i + \sum_{j\in N}\tilde{u}_k^{ij'}R_k^{ij}\tilde{u}_k^{ij}) - \tilde{x}_k^{i'}Q_k^ix_k - \sum_{j\in N}\tilde{u}_k^{ij'}R_k^{ij}u_k^j \\\label{gn}
\end{multline} \vspace*{-5ex}

\vspace*{-5ex} 

\begin{align}
\gamma_k^{i*}(x_{k-1}) = u_k^{i*} = -P_k^{i*}x_{k-1}-\alpha_k^{i*} \label{gamman}
\end{align}

\vspace*{-5ex} 

\begin{align}
P_k^{i*} = (R_k^{ii} + B_k^{i'}Z_k^iB_k^i)^{-1}[B_k^{i'}Z_k^i(A_k - \sum_{j\in N \;,\; j \not = i}B_k^jP_k^{j*})]\label{P} 
\end{align}

\vspace*{-5ex} \begin{multline}
\alpha_k^{i*} = (R_k^{ii} + B_k^{i'}Z_k^iB_k^i)^{-1} [B_k^{i'}(Z_k^i(s_k - \sum_{j\in N \;,\; j \not = i}B_k^j\alpha_k^{j*})\\+ \zeta_{k}^i -Q_k^i\tilde{x}_k^i)-R_k^{ii}\tilde{u}_k^{ii}]\\\label{alpha}
\end{multline} \vspace*{-5ex}

\vspace*{-5ex} \begin{multline}
Z_{k-1}^i = (A_k - \sum_{j \in N}B_k^jP_k^{j*})^{'}Z_k^i(A_k - \sum_{j \in N}B_k^jP_k^{j*}) \\ + \sum_{j\in N}P_k^{j*'}R_k^{ij}P_k^{j*} + Q_{k-1}^i\; ;\; Z_T^i = Q_T^i \\\label{Zfn}
\end{multline} \vspace*{-5ex}

\vspace*{-5ex} \begin{multline}
\zeta_{k-1}^i = (A_k - \sum_{j \in N}B_k^jP_k^{j*})^{'}[\zeta_{k}^i + Z_k^i(s_k-\sum_{j\in N}B_k^j\alpha_k^{j*})-Q_k^i\tilde{x}_k^i] \\ + \sum_{j\in N}P_k^{j*'}R_k^{ij}(\alpha_k^{j*} + \tilde{u}_k^{ij})\; ;\; \zeta_{T}^i = 0 \\\label{zetafn}
\end{multline} \vspace*{-5ex}

\vspace*{-5ex} 

\begin{align}
V^i(1,x_0) = \frac{1}{2} x_0^{'}Z_0^ix_0 + \zeta_0^{i'}x_0 + n_0^i \label{costn}
\end{align}

\vspace*{-5ex} \begin{multline}
n_{k-1}^i = n_{k}^i + \frac{1}{2} (s_k-\sum_{j\in N}B_k^j\alpha_k^{j*})^{'}Z_k^i(s_k-\sum_{j\in N}B_k^j\alpha_k^{j*}) \\+ \zeta_{k}^{i'}(s_k-\sum_{j\in N}B_k^j\alpha_k^{j*}) + \frac{1}{2}\sum_{j\in N }\alpha_k^{j*'}R_k^{ij}\alpha_k^{j*} -  \tilde{x}_k^{i'}Q_k^i(s_k-\sum_{j\in N}B_k^j\alpha_k^{j*}) \\ + \sum_{j\in N}\tilde{u}_k^{ij'}R_k^{ij}\alpha_k^{j*} + \frac{1}{2}(\tilde{x}_k^{i'}Q_k^i\tilde{x}_k^i + \sum_{j\in N}\tilde{u}_k^{ij'}R_k^{ij}\tilde{u}_k^{ij})\; ;\; n_{T}^i = 0 \\\\\label{nn}
\end{multline} \vspace*{-5ex}

\label{FBNTheo}
\end{theo}

\textsc {Proof:}

The proof is carried out by using an induction argument to show that the strategies $\gamma_k^{i*}(x_{k-1})$, given by \eqref{gamman}, minimize the strictly convex functionals \eqref{Vdn} at each stage of the game. But the minimization of \eqref{Vdn} is exactly \eqref{vfn} applied to the specific state equation and cost functionals of the above game and therefore, considering Theorem \eqref{FBN_Opt}, it follows that the $\gamma_k^{i*}(x_{k-1})$ are the unique equilibrium strategies.

\vspace*{-5ex} \begin{multline}
\frac{1}{2} ((A_kx_{k-1} + B_k^iu_k^i + \sum_{j\in N\;,\; j \not = i}B_k^ju_k^{j*}+s_k)^{'}Q_k^i(A_kx_{k-1} + B_k^iu_k^i \\+ \sum_{j\in N\;,\; j \not = i}B_k^ju_k^{j*}+s_k) + u_k^{i'}R_k^{ii}u_k^i + \sum_{j\in N\;,\; j \not = i}u_k^{j*'}R_k^{ij}u_k^{j*})\; +
\frac{1}{2} (\tilde{x}_k^{i'}Q_k^i\tilde{x}_k^i \\+ \sum_{j\in N \;,\; j \not = i}\tilde{u}_k^{ij'}R_k^{ij}\tilde{u}_k^{ij}) - \tilde{x}_k^{i'}Q_k^i(A_kx_{k-1} + B_k^iu_k^i + \sum_{j\in N\;,\; j \not = i}B_k^ju_k^{j*}+s_k) - \tilde{u}_k^{ii'}R_k^{ii}u_k^{i} \\- \sum_{j\in N\;,\; j \not = i}\tilde{u}_k^{ij'}R_k^{ij}u_k^{j*} + V^i(k+1,\tilde{f}_{k-1}^{i*}(x_{k-1},u_k^i)) \; ;\; V^i(T+1,x_T) = 0 \\
\\\label{Vdn}
\end{multline} \vspace*{-5ex}

The induction argument runs from \textit {T}+1 to 1 and proves that the value function of player \textit {i} at stage \textit {k} can be written as stated below in \eqref{Vgn}:

\vspace*{-5ex} \begin{multline}
\frac{1}{2} x_{k-1}^{'}(Z_{k-1}^i-Q_{k-1}^i)x_{k-1} + \zeta_{k-1}^{i'}x_{k-1} + n_{k-1}^i \\\label{Vgn}
\end{multline} \vspace*{-5ex}

\textbf{Basis:}

The induction starts at \textit {k} = \textit {T} + 1. First we make use of the general optimality conditions for $V_k^i$ at stage \textit {T}+1.

\vspace*{-5ex} \begin{multline}
V^i(T+1,x_T) = 0\\\label{optVlsn}
\end{multline} \vspace*{-5ex}

Now we show that $\eqref{Vgn}_{k=T+1}$ fulfills \eqref{optVlsn}.

\vspace*{-5ex} \begin{multline}
\nonumber\eqref{Vgn}_{k=T+1} \; \; \frac{1}{2} x_{T}^{'}(Z_{T}^i-Q_{T}^i)x_{T} + \zeta_{T}^{i'}x_{T} + n_{T}^i\\
\end{multline} \vspace*{-5ex}

Making use of $\eqref{Zfn}_{k=T+1}$, $\eqref{zetafn}_{k=T+1}$ and $\eqref{nn}_{k=T+1}$ yields

\vspace*{-5ex} \begin{multline}
\nonumber \frac{1}{2} x_{T}^{'}(Q_{T}^i-Q_{T}^i)x_{T} + 0^{'}x_{T} + 0 = 0 = V^i(T+1,x_T)\\
\end{multline} \vspace*{-5ex}

\textbf{Inductive step:}

As an induction hypothesis, the equations $\eqref{Vgn}_{k=l+1}$ are assumed to be equal to $V^i(l+1,x_l)$ respectively. Now we have to prove that the relation also holds at stage \textit {l}. In other words: we have to show $\eqref{Vgn}_{k=l} = V^i(l,x_{l-1})$.\\

\vspace*{-5ex} \begin{multline}
\nonumber\eqref{Vgn}_{k=l+1} \; \; \frac{1}{2} x_{l}^{'}(Z_{l}^i-Q_{l}^i)x_{l} + \zeta_{l}^{i'}x_{l} + n_{l}^i \\
\end{multline} \vspace*{-5ex}

Using $\eqref{fn}_{k=l}$ and considering that the optimal control vector for player \textit {i} remains to be deduced in the further induction argument, gives

\vspace*{-5ex} \begin{multline}
\frac{1}{2} (A_lx_{l-1} + B_l^iu_l^i + \sum_{j\in N\;,\; j \not = i}B_l^ju_l^{j*}+s_l)^{'}(Z_{l}^i-Q_{l}^i)(A_lx_{l-1} + B_l^iu_l^i \\+ \sum_{j\in N\;,\; j \not = i}B_l^ju_l^{j*}+s_l) + \zeta_{l}^{i'}(A_lx_{l-1} + B_l^iu_l^i + \sum_{j\in N\;,\; j \not = i}B_l^ju_l^{j*}+s_l) + n_{l}^i \\\label{Vgis}
\end{multline} \vspace*{-5ex}

%\vspace*{-5ex} \begin{multline}
%\nonumber\eqref{P}_{k=l+1}\;\; P_{l+1}^{i*} = (R_{l+1}^{ii} + B_{l+1}^{i'}Z_{l+1}^iB_{l+1}^i)^{-1}[B_{l+1}^{i'}Z_{l+1}^i(A_{l+1} - \sum_{j\in N \;,\; j \not = i}B_{l+1}^jP_{l+1}^{j*})]\\
%\end{multline} \vspace*{-5ex}

%\vspace*{-5ex} \begin{multline}
%\nonumber\eqref{alpha}_{k=l+1}\;\; \alpha_{l+1}^{i*} = (R_{l+1}^{ii} + B_{l+1}^{i'}Z_{l+1}^iB_{l+1}^i)^{-1} %[B_{l+1}^{i'}(Z_{l+1}^i(s_{l+1} - \sum_{j\in N \;,\; j \not = i}B_{l+1}^j\alpha_{l+1}^{j*})\\+ \zeta_{l+2}^i %-Q_{l+1}^i\tilde{x}_{l+1}^i)-R_{l+1}^{ii}\tilde{u}_{l+1}^{ii}]\\
%\end{multline} \vspace*{-5ex}

First we prove that $\eqref{gamman}_{k=l}$ minimizes $\eqref{Vdn}_{k=l}$. To do this, we substitute $V^i(l+1,\tilde{f}_{l-1}^{i*}(x_{l-1},u_{l}^i))$ in $\eqref{Vdn}_{k=l}$ with the help of the induction hypothesis in \eqref{Vgis}:

\vspace*{-5ex} \begin{multline}
\nonumber \eqref{Vdn}_{k=l} \; \; \frac{1}{2} ((A_lx_{l-1} + B_l^iu_l^i + \sum_{j\in N\;,\; j \not = i}B_l^ju_l^{j*}+s_l)^{'}Q_l^i(A_lx_{l-1} + B_l^iu_l^i \\+ \sum_{j\in N\;,\; j \not = i}B_l^ju_l^{j*}+s_l) + u_l^{i'}R_l^{ii}u_l^i + \sum_{j\in N\;,\; j \not = i}u_l^{j*'}R_l^{ij}u_l^{j*})\; +
\frac{1}{2} (\tilde{x}_l^{i'}Q_l^i\tilde{x}_l^i \\+ \sum_{j\in N \;,\; j \not = i}\tilde{u}_l^{ij'}R_l^{ij}\tilde{u}_l^{ij}) - \tilde{x}_l^{i'}Q_l^i(A_lx_{l-1} + B_l^iu_l^i + \sum_{j\in N\;,\; j \not = i}B_l^ju_l^{j*}+s_l) - \tilde{u}_l^{ii'}R_l^{ii}u_l^{i} \\- \sum_{j\in N\;,\; j \not = i}\tilde{u}_l^{ij'}R_l^{ij}u_l^{j*} + V^i(l+1,\tilde{f}_{l-1}^{i*}(x_{l-1},u_l^i)) \\
\end{multline} \vspace*{-5ex}

\vspace*{-5ex} \begin{multline}
\frac{1}{2} ((A_lx_{l-1} + B_l^iu_l^i + \sum_{j\in N\;,\; j \not = i}B_l^ju_l^{j*}+s_l)^{'}Q_l^i (A_lx_{l-1} + B_l^iu_l^i \\+ \sum_{j\in N\;,\; j \not = i}B_l^ju_l^{j*}+s_l) + u_l^{i'}R_l^{ii}u_l^i + \sum_{j\in N\;,\; j \not = i}u_l^{j*'}R_l^{ij}u_l^{j*})\; +
\\ \frac{1}{2} (\tilde{x}_l^{i'}Q_l^i\tilde{x}_l^i + \sum_{j\in N}\tilde{u}_l^{ij'}R_l^{ij}\tilde{u}_l^{ij}) - \tilde{x}_l^{i'}Q_l^i(A_lx_{l-1} + B_l^iu_l^i + \sum_{j\in N\;,\; j \not = i}B_l^ju_l^{j*}+s_l) \\- \tilde{u}_l^{ii'}R_l^{ii}u_l^{i} - \sum_{j\in N\;,\; j \not = i}\tilde{u}_l^{ij'}R_l^{ij}u_l^{j*} + \frac{1}{2} (A_lx_{l-1} + B_l^iu_l^i + \sum_{j\in N\;,\; j \not = i}B_l^ju_l^{j*}+s_l)^{'}(Z_{l}^i-Q_{l}^i)\\(A_lx_{l-1} + B_l^iu_l^i + \sum_{j\in N\;,\; j \not = i}B_l^ju_l^{j*}+s_l) + \zeta_{l}^{i'}(A_lx_{l-1} + B_l^iu_l^i + \sum_{j\in N\;,\; j \not = i}B_l^ju_l^{j*}+s_l) + n_{l}^i \\\label{Vln} 
\end{multline} \vspace*{-5ex}

$\eqref{Vdn}_{k=l}$ is strictly convex in $u_l^i$. This can be seen by applying Corollary \eqref{pdc} to \eqref{connfb}. Therefore there has to be a unique equilibrium strategy for player \textit {i} at stage \textit {l}.

\vspace*{-5ex} \begin{multline}
\frac {\partial^2} {\partial u_l^{i^2}} \eqref{Vln} \;\;\ R_l^{ii} + B_l^{i'}Z_l^iB_l^i \\\label{connfb}
\end{multline} \vspace*{-5ex}

This unique optimal strategy can be found by using the first-order necessary and sufficient (because of strict convexity of $\eqref{Vdn}_{k=l}$) conditions for minimization

\vspace*{-5ex} \begin{multline}
\nonumber\frac {\partial} {\partial u_l^{i}} \eqref{Vln} = 0 \Rightarrow \;\;\ -(R_l^{ii} + B_l^{i'}Z_l^iB_l^i)u_l^{i*} - B_l^{i'}Z_l^i \\ \sum_{j\in N \;,\; j \not = i}B_l^ju_l^{j*} = B_l^{i'}[Z_l^i(A_lx_{l-1} + s_l) + \zeta_{l}^i  - Q_l^i\tilde{x}_l^i] - R_l^{ii}\tilde{u}_l^{ii} \\
\end{multline} \vspace*{-5ex}
 
As the right hand side of the above equation is affine in $x_{l-1}$, the left hand side also has to be affine in $x_{l-1}$. Therefore, the substitution

\vspace*{-5ex} \begin{multline}
u_l^{i*} = -P_l^{i*}x_{l-1}-\alpha_l^{i*} \\\label{subnfbu}
\end{multline} \vspace*{-5ex}
%hier ist kritscher punkt!!!!!!!!!!!!!!!!1

is allowed and leads to

\vspace*{-5ex} \begin{multline}
-(R_l^{ii} + B_l^{i'}Z_l^iB_l^i)(-P_l^{i*}x_{l-1}-\alpha_l^{i*}) - B_l^{i'}Z_l^i \\ \sum_{j\in N \;,\; j \not =  i}B_l^j(-P_l^{j*}x_{l-1}-\alpha_l^{j*}) = B_l^{i'}[Z_l^i(A_lx_{l-1} + s_l) + \zeta_{l}^i  - Q_l^i\tilde{x}_l^i] - R_l^{ii}\tilde{u}_l^{ii} \\\label{gammanl}
\end{multline} \vspace*{-5ex}

By comparison of coefficients follows

\vspace*{-5ex} \begin{multline}
\nonumber\eqref{gammanl}_{x_{l-1}} \; \; (R_l^{ii} + B_l^{i'}Z_l^iB_l^i)P_l^{i*} + B_l^{i'}Z_l^i \sum_{j\in N \;,\; j \not = i}B_l^jP_l^{j*} = B_l^{i'}Z_l^iA_l \\
\end{multline} \vspace*{-5ex}

\vspace*{-5ex} \begin{multline}
\nonumber\eqref{gammanl}_{const.} \; \; (R_l^{ii} + B_l^{i'}Z_l^iB_l^i)\alpha_l^{i*} + B_l^{i'}Z_l^i  \sum_{j\in N \;,\; j \not = i}B_l^j\alpha_l^{j*} = B_l^{i'}[Z_l^is_l + \zeta_{l}^i  - Q_l^i\tilde{x}_l^i] - R_l^{ii}\tilde{u}_l^{ii}\\
\end{multline} \vspace*{-5ex}

Making $P_l^{i*}$ and $\alpha_l^{i*}$ explicit yields

\vspace*{-5ex} \begin{multline}
\nonumber\eqref{P}_{k=l} \;\;P_l^{i*} = (R_l^{ii} + B_l^{i'}Z_l^iB_l^i)^{-1}[B_l^iZ_l^i(A_l - \sum_{j\in N \;,\; j \not = i}B_l^jP_l^{j*})]\\
\end{multline} \vspace*{-5ex}

\vspace*{-5ex} \begin{multline}
\nonumber\eqref{alpha}_{k=l}\;\;\alpha_l^{i*} = (R_l^{ii} + B_l^{i'}Z_l^iB_l^i)^{-1} [B_l^{i'}(Z_l^i(s_l  - \sum_{j\in N \;,\; j \not = i}B_l^j\alpha_l^{j*}) + \zeta_{l}^i - Q_l^i\tilde{x}_l^i)-R_l^{ii}\tilde{u}_l^{ii}] \\
\end{multline} \vspace*{-5ex}

Now, after finding the optimal strategies for the players, we are able to rewrite \eqref{Vln} as

\vspace*{-5ex} \begin{multline}
\nonumber V^i(l,x_{l-1}) = \frac{1}{2} [(A_lx_{l-1} + \sum_{j\in N}B_l^ju_l^{j*}+s_l)^{'}Z_l^i(A_lx_{l-1} + \sum_{j\in N}B_l^ju_l^{j*}+s_l) \\+ \sum_{j\in N}u_l^{j*'}R_l^{ij}u_l^{j*}] +
 \frac{1}{2} (\tilde{x}_l^{i'}Q_l^i\tilde{x}_l^i + \sum_{j\in N}\tilde{u}_l^{ij'}R_l^{ij}\tilde{u}_l^{ij}) - \tilde{x}_l^{i'}
Q_l^i(A_lx_{l-1} + \sum_{j\in N}B_l^ju_l^{j*}+s_l) \\- \sum_{j\in N}\tilde{u}_l^{ij'}R_l^{ij}u_l^{j*} + \zeta_{l}^{i'}(A_lx_{l-1} + \sum_{j\in N}B_l^ju_l^{j*}+s_l) + n_{l} \\
\end{multline} \vspace*{-5ex}

Making use of \eqref{subnfbu} yields

\vspace*{-5ex} \begin{multline}
\nonumber V^i(l,x_{l-1}) = \frac{1}{2} [(A_lx_{l-1} + \sum_{j\in N}B_l^j(-P_l^{j*}x_{l-1}-\alpha_l^{j*}) +s_l)^{'}Z_l^i(A_lx_{l-1} \\+ \sum_{j\in N}B_l^j(-P_l^{j*}x_{l-1}-\alpha_l^{j*})+s_l) + \sum_{j\in N}(-P_l^{j*}x_{l-1}-\alpha_l^{j*}) ^{'}R_l^{ij}(-P_l^{j*}x_{l-1}-\alpha_l^{j*}) ]\\ +
 \frac{1}{2} (\tilde{x}_l^{i'}Q_l^i\tilde{x}_l^i + \sum_{j\in N}\tilde{u}_l^{ij'}R_l^{ij}\tilde{u}_l^{ij}) - \tilde{x}_l^{i'}Q_l^i(A_lx_{l-1} + \sum_{j\in N}B_l^j(-P_l^{j*}x_{l-1}-\alpha_l^{j*})+s_l) \\- \sum_{j\in N}\tilde{u}_l^{ij'}R_l^{ij}(-P_l^{j*}x_{l-1}-\alpha_l^{j*}) + \zeta_{l}^{i'}(A_lx_{l-1} + \sum_{j\in N}B_l^j(-P_l^{j*}x_{l-1}-\alpha_l^{j*})+s_l) + n_{l}  \\
\end{multline} \vspace*{-5ex}

Rewriting the above equation to the power of $x_{l-1}$ gives

\vspace*{-5ex} \begin{multline}
\nonumber V^i(l,x_{l-1}) = \frac{1}{2} x_{l-1}^{'}[(A_l - \sum_{j \in N}B_l^jP_l^{j*})^{'}Q_l^i(A_l - \sum_{j \in N}B_l^jP_l^{j*}) + \sum_{j\in N}P_l^{j*'}R_l^{ij}P_l^{j*}]x_{l-1} \\ + [(A_l - \sum_{j \in N}B_l^jP_l^{j*})^{'}[\zeta_{l}^{i}Z_l^i(s_l-\sum_{j\in N}B_l^j\alpha_l^{j*}) - Q_l^i\tilde{x}_k^i] + \sum_{j\in N}P_l^{j*'}R_l^{ij}(\alpha_l^{j*} + \tilde{u}_l^{ij}) ]^{'}x_{l-1} \\+ \frac{1}{2} (s_l-\sum_{j\in N}B_l^j\alpha_l^{j*})^{'}Z_l^i(s_l-\sum_{j\in N}B_l^j\alpha_l^{j*})  + \frac{1}{2}\sum_{j\in N }\alpha_l^{j*'}R_l^{ij}\alpha_l^{j*} \\-  \tilde{x}_l^{i'}Q_l^i(s_l-\sum_{j\in N}B_l^j\alpha_l^{j*}) + \sum_{j\in N}\tilde{u}_l^{ij'}R_l^{ij}\alpha_l^{j*} + \frac{1}{2}(\tilde{x}_l^{i'}Q_l^i\tilde{x}_l^i + \sum_{j\in N}\tilde{u}_l^{ij'}R_l^{ij}\tilde{u}_l^{ij}) \\+ \zeta_{l}^{i'}(s_l - \sum_{j\in N}B_l^j\alpha_l^{j*}) + n_{l} \\
\end{multline} \vspace*{-5ex}

To finish off the inductive step and consequently the induction argument, we use the recursive equations $\eqref{Zfn}_{k=l}$, $\eqref{zetafn}_{k=l}$ and $\eqref{nn}_{k=l}$ in the above equation. This leads to

\vspace*{-5ex} \begin{multline}
\nonumber\eqref{Vgn}_{k=l} \;\;V^i(l,x_{l-1}) = \frac{1}{2} x_{l-1}^{'}(Z_{l-1}^i-Q_{l-1}^i)x_{l-1} + \zeta_{l-1}^{i'}x_{l-1} + n_{l-1}^i \\
\end{multline} \vspace*{-5ex}

The expression for the total costs of the game for player \textit {i} given by \eqref{costn} is equal to the function the induction argument was based on at stage \textit {l}. In other words, \eqref{costn} is equal to $\eqref{Vgn}_{k=1}$.\;\;\;\;\;\;\;\;\boxed{}

\begin{rem}
The proof of Theorem \eqref{FBNTheo} is a formalization of the heuristic argumentation presented in Ba\c{s}ar and Olsder (1999, pp. 280-281)\cite{baol}.
\end{rem}

\begin{rem}
$Z_k^i$, given by \eqref{Zfn}, is positive definite for all $k \in \{0, \ldots, T\}$. This can be proven by a straightforward induction argument, starting at stage T and using Corollary \eqref{pdc} in the inductive step.
\end{rem}

\begin{rem}
Special attention should be paid to the observation that using an affine state equation together with a quadratic cost functional yields affine equilibrium strategies. This also holds true for feedback Stackelberg, open-loop Nash and open-loop Stackelberg games.
\end{rem}

\begin{rem}
To solve the Nash game algorithmically, the following order of application of the equations of Theorem \eqref{FBNTheo} is advisable ($i \in N$):
\begin{enumerate}
\item For k running backward from T to 1
\begin{enumerate}
\item $Z_k^i$, $\zeta_k^i$ and $n_k^i$
\item $P_{k}^{i*}$, $\alpha_{k}^{i*}$ 
\end{enumerate}
\item $Z_0^i$, $\zeta_0^i$ and $n_0^i$
\item $V^i(1,x_0)$
\item For k running forward from 1 to T\\
$\gamma_k^{i*}(x_{k-1})$\\\\
\end{enumerate}
\end{rem}

\begin{pro}\footnote{In this proposition we rewrite the equilibrium equations in a notation that was used at our department in the past to enable comparison.}
The systems of equations defining the unique equilibrium strategies $\gamma_{k}^{i*}(x_{k-1})$ in Theorem \eqref{FBNTheo} can also be written in the following way:\footnote{For all equations belonging to this proposition and its proof, $k \in K$ and $i \in N$ if nothing different is stated.}

\vspace*{-5ex} 

\begin{align}
\gamma_k^{i*}(x_{k-1}) = G_k^{i*}x_{k-1} + g_k^{i*} \label{gammanä}
\end{align}

\vspace*{-5ex} 

\begin{align}
G_k^{i*} = -(D_k^i)^{-1}[B_k^{i'}H_k^i(A_k + \sum_{j\in N \;,\; j \not = i}B_k^jG_k^{j*})] \label{Pä} 
\end{align}

\vspace*{-5ex} 

\begin{align}
g_k^{i*} = -(D_k^i)^{-1} [B_k^{i'}H_k^i\sum_{j\in N \;,\; j \not = i}B_k^jg_k^{j*} + v_k^i] \label{alphaä}
\end{align}

\vspace*{-5ex} 

\begin{align}
H_{k-1}^i = K_k^{'}Z_k^iK_k + \sum_{j\in N}G_k^{j'}R_k^{ij}G_k^j + Q_{k-1}^i\; ;\; H_T^i = Q_T^i \label{Hknä}
\end{align}

\vspace*{-5ex} 

\begin{align}
h_{k-1}^i = Q_{k-1}^i\tilde{x}_{k-1}^i - K_k^{'}[H_k^ik_k - h_{k+1}^i] + \sum_{j\in N}G_k^{j'}R_k^{ij}(\tilde{u}_k^{ij} - g_k^j)\; ;\; h_{T}^i = Q_T^i\tilde{x}_T^i  \label{zetanä}
\end{align}

\vspace*{-5ex}

\begin{align}
K_k = A_k + \sum_{j \in N}B_k^iG_k^j \label{Ksubn}
\end{align}

\vspace*{-5ex} 

\begin{align}
k_k = s_k + \sum_{j\in N}B_k^jg_k^j \label{ksubn}
\end{align}

\vspace*{-5ex}

\begin{align}
D_k^i = R_k^{ii} + B_k^{i'}H_k^iB_k^i \label{Dsubn}
\end{align}

\vspace*{-5ex} 

\begin{align}
v_k^i = B_k^{i'}(H_k^is_k - h_{k}^i) - R_k^{ii}\tilde{u}_k^{ii} \; ; \; i \in N \label{vsubn}
\end{align}\\\\

\label{profbn}
\end{pro}

\textsc {Proof:}

The proof is carried out by renaming some matrices and then showing that the relations for the equilibrium strategies $\gamma_k^{i*}(x_{k-1})$ of Theorem \eqref{FBSTheo} can be rewritten in the way stated above.\\

Let us start by renaming the feedback matrices $P_{k}^{i*}$ and $\alpha_k^{i*}$ and the matrices $Z_k^i$ and $\zeta_{k+1}^i$.

\vspace*{-5ex} \begin{multline}
P_k^i \; \widehat {=} \; -G_k^i \; ; \; \alpha_k^i \; \widehat {=} \; -g_k^i  \; ; \; Z_k^i \; \widehat {=} \; H_k^i \; ; \; \zeta_{k}^i \; \widehat {=} \; - h_k^i + Q_k^i\tilde{x}_k^i \\\label{äqusubnfb}
\end{multline} \vspace*{-5ex}

Next we prove that $Z_k^i$ and $\zeta_{k}^i$ fulfill \eqref{Hknä} and \eqref{zetanä} respectively.\\

Taking the renaming \eqref{Zfn} into account gives

\vspace*{-5ex} \begin{multline}
\nonumber H_{k-1}^i = (A_k + \sum_{j \in N}B_k^iG_k^j)^{'}H_k^i(A_k + \sum_{j \in N}B_k^jG_k^j) + \sum_{j\in N}G_k^{j'}R_k^{ij}G_k^j + Q_{k-1}^i \\
\end{multline} \vspace*{-5ex}

Making use of \eqref{Ksubn} yields

\vspace*{-5ex} \begin{multline}
\nonumber\eqref{Hknä} \; \; H_{k-1}^i = K_k^{'}H_k^iK_k + \sum_{j\in N}G_k^{j'}R_k^{ij}G_k^j + Q_{k-1}^i \\
\end{multline} \vspace*{-5ex}

Now we show the correctness of equation \eqref{zetanä}. To do so we start with stage \textit {T} and use \eqref{zetafn} and \eqref{äqusubnfb} to get

\vspace*{-5ex} \begin{multline}
\nonumber 0 = \zeta_{T}^i = - h_T^i + Q_T^i\tilde{x}_T^i \\
\end{multline} \vspace*{-5ex}

\vspace*{-5ex} \begin{multline}
\nonumber\eqref{zetanä}_{k=T} \; \; h_T^i = Q_T^i\tilde{x}_T^i \\
\end{multline} \vspace*{-5ex}

For the general stage \textit{k} rewrite \eqref{zetafn} taking consideration of \eqref{äqusubnfb}

\vspace*{-5ex} \begin{multline}
\nonumber -h_{k-1}^i + Q_{k-1}^i\tilde{x}_{k-1}^i = (A_k + \sum_{j \in N}B_k^jG_k^j)^{'}[-h_{k}^i + H_k^i(s_k + \sum_{j\in N}B_k^jg_k^j)] - \sum_{j\in N}G_k^{j'}R_k^{ij}(\tilde{u}_k^{ij} - g_k^j) \\
\end{multline} \vspace*{-5ex}

Using \eqref{Ksubn} and \eqref{ksubn} and making $h_{k-1}^i$ explicit yields

\vspace*{-5ex} \begin{multline}
\nonumber\eqref{zetanä} \; \; h_{k-1}^i  = Q_{k-1}^i\tilde{x}_{k-1}^i - K_k^{'}[H_k^ik_k - h_{k}^i] + \sum_{j\in N}G_k^{j'}R_k^{ij}(\tilde{u}_k^{ij} - g_k^j) \\
\end{multline} \vspace*{-5ex}

Eventually the correctness of the rewritten equilibrium strategies $\gamma_k^{i*}(x_{k-1})$ given by \eqref{gammanä} - \eqref{alphaä} has to be shown.\\

First substituting $P_{k}^{i*}$ and $\alpha_{k}^{i*}$ in \eqref{gamman} with the help of \eqref{äqusubnfb} leads to

\vspace*{-5ex} \begin{multline}
\nonumber\eqref{gammanä} \; \; \gamma_k^{i*}(x_{k-1}) = G_k^{i*}x_{k-1} + g_k^{i*} \\
\end{multline} \vspace*{-5ex}

Taking \eqref{äqusubnfb} into consideration, the feedback matrices given by \eqref{P} and \eqref{alpha} can be rewritten as

\vspace*{-5ex} \begin{multline}
\nonumber G_k^{i*} = -(R_k^{ii} + B_k^{i'}H_k^iB_k^i)^{-1}[B_k^{i'}H_k^i(A_k + \sum_{j\in N \;,\; j \not = i}B_k^jG_k^{j*})] \\
\end{multline} \vspace*{-5ex}

\vspace*{-5ex} \begin{multline}
g_k^{i*} = -(R_k^{ii} + B_k^{i'}H_k^iB_k^i)^{-1} [B_k^{i'}(H_k^i(s_k + \sum_{j\in N \;,\; j \not = i}B_k^jg_k^{j*}) - h_k^i) -R_k^{ii}\tilde{u}_k^{ii}] \\\label{ginco}
\end{multline} \vspace*{-5ex}

Finally using \eqref{Dsubn} in the two equations above and additionally using \eqref{vsubn} in \eqref{ginco} gives

\vspace*{-5ex} \begin{multline}
\nonumber\eqref{Pä} \; \; G_k^{i*} = -(D_k^i)^{-1}[B_k^{i'}H_k^i(A_k + \sum_{j\in N \;,\; j \not = i}B_k^jG_k^{j*})] \\
\end{multline} \vspace*{-5ex}

\vspace*{-5ex} \begin{multline}
\nonumber\eqref{alphaä} \; \; g_k^{i*} = -(D_k^i)^{-1} [B_k^{i'}H_k^i\sum_{j\in N \;,\; j \not = i}B_k^jg_k^{j*} + v_k^i] \;\;\;\;\;\;\;\;\;\;\;\;\;\;\;\;\boxed{}\\
\end{multline} \vspace*{-5ex}

\clearpage

\subsection{Special case: The affine-quadratic control problem}

In this subsection the results of Theorem \eqref{FBNTheo} are first specialized in Corollary \eqref{cthb} by reducing the number of players from \textit {n} to one and then in Proposition \eqref{procä} the specialized results are transformed into the terminology used in Proposition 5.1 in Ba\c{s}ar and Olsder (1999, pp. 234-235)\cite{baol} to point out some serious mistakes stated there.

\begin{cor}
An affine-quadratic control problem (cf. Def. \eqref{defspiel} with N = \{1\} and \eqref{fc} - \eqref{gc}) admits the \emph{unique control solution} if
\begin{itemize}
\item $Q_k \geq 0$, $R_k > 0$ (defined for: $k \in K$).
\end{itemize}
If these conditions are satisfied,  the unique optimal strategies are given by \eqref{gammac}. The corresponding \emph {minimum value} is stated in \eqref{costc}. \footnote {For all equations belonging to this corollary, its proof and its equivalence analysis, $k \in K$ if nothing different is stated.}

\vspace*{-5ex} \begin{align}
f_{k}(x_{k},u_k) = A_kx_{k} + B_ku_k + c_k \label{fc}
\end{align}

\vspace*{-5ex} \begin{align}
L(u)= \sum_{k = 1}^T g_k(x_{k+1},u_k, x_{k}) 
\end{align}

\vspace*{-5ex} \begin{align}
 g_k(x_{k+1}, u_k, x_{k}) = \frac{1}{2} (x_{k+1}^{'}Q_{k+1}x_{k+1} + u_k^{'}R_ku_k) \label{gc}
\end{align}

\vspace*{-5ex} \begin{align}
\gamma_k(x_{k}) = -P_k^{*}x_{k}-\alpha_k^{*} \label{gammac}
\end{align}

\vspace*{-5ex} \begin{align}
P_k^{*} = (R_k + B_kZ_{k+1}B_k)^{-1}B_k^{'}Z_{k+1}A_k\label{Pc} 
\end{align}

\vspace*{-5ex} \begin{align}
\alpha_k^{*}= (R_k + B_k^{'}Z_{k+1}B_k)^{-1} B_k^{'}(Z_{k+1}c_k + \zeta_{k+1})\label{alphac}
\end{align}

\vspace*{-5ex} \begin{align}
Z_{k} = (A_k - B_kP_k^{*})^{'}Z_{k+1}(A_k - B_kP_k^{*}) + P_k^{*'}R_kP_k^{*} + Q_{k}\; ;\; Z_T = Q_T \label{Zc}
\end{align}

\vspace*{-5ex} \begin{align}
\zeta_k = (A_k - B_kP_k^{*})^{'}[\zeta_{k+1} + Z_{k+1}(c_k - B_k\alpha_k^{*})] + P_k^{*'}R_k\alpha_k^{*} \; ;\; \zeta_{T+1} = 0 \label{zetac}
\end{align}

\vspace*{-5ex} \begin{align}
V(1,x_1) = \frac{1}{2} x_1^{'}Z_1x_1 + \zeta_1^{'}x_1 + n_1 \label{costc}
\end{align}

\vspace*{-5ex} \begin{multline}
n_k = n_{k+1} + \frac{1}{2} (c_k - B_k\alpha_k^{*})^{'}Z_{k+1}(c_k - B_k\alpha_k^{*}) \\+ \zeta_{k+1}^{'}(c_k - B_k\alpha_k^{*}) + \frac{1}{2}\alpha_k^{*'}R_k\alpha_k^{*} \\\label{nc}
\end{multline} \vspace*{-5ex}

\label{cthb}
\end{cor}

\textsc {Proof:}

Corollary \eqref{cthb} is proven in the same way as Theorem \eqref{FBNTheo} taking into consideration simplifications resulting from the reduction in the number of players to one and the modified state equation and cost functionals.\;\;\;\;\;\;\;\;\;\;\;\;\;\;\;\;\boxed{}\\

\begin{pro}
The systems of equations defining the unique equilibrium strategies $\gamma_{k+1}^{i*}(x_{k})$ in Corollary \eqref{cthb} can also be written in the following way:\footnote{For all equations belonging to this proposition and its proof, $k \in K$ and $i \in N$ if nothing different is stated. Equations \eqref{Zcä}, \eqref{zetacä} and \eqref{ncä} are wrong in Ba\c{s}ar and Olsder.}

\vspace*{-5ex} \begin{align}
\gamma_k^{*}(x_{k}) = -P_kS_{k+1}A_k x_{k} - P_k(s_{k+1} + S_{k+1}c_k) \label{gammacä}
\end{align}

\vspace*{-5ex} \begin{align}
P_k = (R_k + B_k^{'}S_{k+1}B_k)^{-1}B_k^{'} \label{Psubc}
\end{align}

\vspace*{-5ex} \begin{align}
S_{k} = (A_k - B_kP_k^{*})^{'}S_{k+1}(A_k - B_kP_k^{*}) + P_k^{*'}R_kP_k^{*} + Q_{k}\; ;\; S_{T+1} = Q_{T+1} \label{Zcä}
\end{align}

\vspace*{-5ex} \begin{align}
s_k = (A_k - B_kP_k^{*})^{'}[s_{k+1} + S_{k+1}(c_k - B_k\alpha_k^{*})] + P_k^{*'}R_k\alpha_k^{*} \; ;\; s_{T+1} = 0 \label{zetacä}
\end{align}

\vspace*{-5ex} \begin{align}
L(1,x_1) = \frac{1}{2} x_1^{'}S_1x_1 + s_1^{'}x_1 + q_1 \label{costcä}
\end{align}

\vspace*{-5ex} \begin{multline}
q_k = q_{k+1} + \frac{1}{2} (c_k - B_kP_k(S_{k+1}c_k + s_{k+1}))^{'}S_{k+1}(c_k - B_kP_k(S_{k+1}c_k + s_{k+1})) \\+ s_{k+1}^{'}(c_k - B_kP_k(S_{k+1}c_k + s_{k+1})) + \frac{1}{2}(P_k(S_{k+1}c_k + s_{k+1}))^{'}R_k(P_k(S_{k+1}c_k + s_{k+1})) \\\\\label{ncä}
\end{multline} \vspace*{-5ex}
\label{procä}

\end{pro}

\textsc {Proof:}\\

The proof is carried out by renaming some matrices and then showing that the relations for the optimal strategies $\gamma_k^{*}(x_{k})$ of Corollary \eqref{cthb} can be rewritten in the way stated above.\\

Let us start by renaming the matrices $Z_k$, $\zeta_{k}$ and $n_k$.

\vspace*{-5ex} \begin{multline}
Z_k \; \widehat {=} \; S_k \; ; \; \zeta_{k} \; \widehat {=} \; s_k  \; ; \; n_k \; \widehat {=} \; q_k \\\label{äqusubcfb}
\end{multline} \vspace*{-5ex}

Considering \eqref{äqusubcfb}, the feedback matrices given by \eqref{Pc} and \eqref{alphac} can be rewritten as

\vspace*{-5ex} \begin{multline}
\nonumber P_k^{*} = (R_k + B_kS_{k+1}B_k)^{-1}B_k^{'}S_{k+1}A_k \\
\end{multline} \vspace*{-5ex}

\vspace*{-5ex} \begin{multline}
\nonumber \alpha_k^{*}= (R_k + B_k^{'}S_{k+1}B_k)^{-1}B_k^{'}(S_{k+1}c_k + s_{k+1}) \\
\end{multline} \vspace*{-5ex}

Using \eqref{Psubc} in the two equations above gives

\vspace*{-5ex} \begin{multline}
P_k^{*} = P_kS_{k+1}A_k \\\label{Pcä}
\end{multline} \vspace*{-5ex}

\vspace*{-5ex} \begin{multline}
\alpha_k^{*} = P_k(S_{k+1}c_k + s_{k+1}) \\\label{alphacä}
\end{multline} \vspace*{-5ex}

Substituting $P_{k}^{*}$ and $\alpha_{k}^{*}$ in \eqref{gammac} with the help of \eqref{Pcä} and \eqref{alphacä} yields

\vspace*{-5ex} \begin{multline}
\nonumber\eqref{gammacä} \; \; \gamma_k^{*}(x_{k}) = -P_kS_{k+1}A_k x_{k} - P_k(s_{k+1} + S_{k+1}c_k)\\
\end{multline} \vspace*{-5ex}

Next we prove that $S_k$ and $s_{k}$ fulfill \eqref{Zcä} and \eqref{zetacä} respectively.\\

Taking the renaming and equation \eqref{Pcä} into consideration, \eqref{Zc} gives

\vspace*{-5ex} \begin{multline}
\nonumber\eqref{Zcä} \; \; S_{k} = (A_k - B_kP_kS_{k+1}A_k)^{'}S_{k+1}(A_k - B_kP_kS_{k+1}A_k) \\+ (P_kS_{k+1}A_k)^{'}R_k(P_kS_{k+1}A_k) + Q_{k}\; ;\; S_T = Q_T \\
\end{multline} \vspace*{-5ex}

Now we show the correctness of equation \eqref{zetacä}. To do so we rewrite \eqref{zetac} taking consideration of \eqref{äqusubcfb}, \eqref{Pcä} and \eqref{alphacä}

\vspace*{-5ex} \begin{multline}
\nonumber\eqref{zetacä} \; \; s_k = (A_k - B_kP_kS_{k+1}A_k)^{'}[s_{k+1} + S_{k+1}(c_k - B_kP_k(s_{k+1} + S_{k+1}c_k))] \\+ (P_kS_{k+1}A_k)^{'}R_kP_k(s_{k+1} + S_{k+1}c_k) \; ;\; s_{T+1} = 0 \\
\end{multline} \vspace*{-5ex}

Eventually the correctness of the rewritten value function $L(1,x_1)$ given by \eqref{costcä} has to be shown. Making use of \eqref{äqusubcfb} in \eqref{costc} and \eqref{nc} respectively leads to

\vspace*{-5ex} \begin{multline}
\nonumber\eqref{costcä} \; \; L(1,x_1) = \frac{1}{2} x_1^{'}S_1x_1 + s_1^{'}x_1 + q_1 \\
\end{multline} \vspace*{-5ex}

\vspace*{-5ex} \begin{multline}
q_k = q_{k+1} + \frac{1}{2} (c_k - B_k\alpha_k^{*})^{'}S_{k+1}(c_k - B_k\alpha_k^{*}) \\+ s_{k+1}^{'}(c_k - B_k\alpha_k^{*}) + \frac{1}{2}\alpha_k^{*'}R_k\alpha_k^{*} \\\label{qincom}
\end{multline} \vspace*{-5ex}

Finally substituting $\alpha_k^{*}$ in \eqref{qincom} with the help of \eqref{alphacä} yields

\vspace*{-5ex} \begin{multline}
\nonumber\eqref{ncä} \; \; q_k = q_{k+1} + \frac{1}{2} (c_k - B_kP_k(S_{k+1}c_k + s_{k+1}))^{'}S_{k+1}(c_k - B_kP_k(S_{k+1}c_k + s_{k+1})) \\+ s_{k+1}^{'}(c_k - B_kP_k(S_{k+1}c_k + s_{k+1})) + \frac{1}{2}(P_k(S_{k+1}c_k + s_{k+1}))^{'}R_k(P_k(S_{k+1}c_k + s_{k+1})) \;\;\;\;\;\;\;\;\;\;\;\;\;\;\;\;\boxed{}\\
\end{multline} \vspace*{-5ex}

\label{ssaqcp}

\clearpage
\section {Feedback Stackelberg Equilibrium Solutions}

This section is devoted to the derivation of the feedback Stackelberg equilibrium solution with one leader and arbitrarily many followers for affine-quadratic games. First a general result is stated about the existence of a Stackelberg equilibrium solution with one leader and arbitrarily many followers in \textit{n}-person discrete-time deterministic infinite dynamic games of prespecified fixed duration (cf. Def. \eqref{defspiel}) with feedback information pattern. Then this result is applied to affine-quadratic games and finally the feedback Stackelberg equilibrium solutions with one leader and one follower for affine-quadratic and linear-quadratic games are deduced as special cases.

\subsection{Optimality conditions}

In this subsection a theorem is stated which gives necessary and sufficient conditions for the existence of a feedback Stackelberg equilibrium solution with one leader and arbitrarily many followers.  Results about feedback Stackelberg equilibria in infinite dynamic games first appeared in discrete time in the works of Simaan and Cruz (1973) \cite{sicra}, \cite{sicrb}.

\begin{theo}
For an n-person discrete-time deterministic infinite dynamic game of prespecified fixed duration (cf. Def. \eqref{defspiel}) with feedback information pattern, the set of strategies $\{\gamma_k^{i*}(x_k); k \in K , i \in  N \}$ provides a \emph {feedback Stackelberg equilibrium solution} with \textbf{P}1 as the leader and \textbf{P}2 $\ldots$ \textbf{P}n as followers if, and only if, functions $V^i(k,\cdot) : R^n \to R,k \in  K, i \in N,$ exist such that the following recursive relations are satisfied:

\vspace*{-5ex} \begin{multline}
V^1(k,x_{k-1}) = \min_{u_k^1 \in U_k^1 , u_k^2 \in R_k^2(u_k^1) , \ldots, u_k^n \in R_k^n(u_k^1)}[g_k^1(f_{k-1}(x_{k-1},\\u_k^{1},u_k^2,\ldots,u_k^n),u_k^{1},u_k^2,\ldots,u_k^{n}, x_{k-1}) +  V^1(k+1,f_{k-1}(x_{k-1},u_k^{1},u_k^2,\ldots,u_k^n))]  =\\ 
g_k^1[\tilde{f}_{k-1}^{1}(x_{k-1},\gamma_k^{1*}(x_{k-1})),\gamma_k^{1*}(x_{k-1}),\ldots,\gamma_k^{n*}(x_{k-1}),x_{k-1}] \\ + V^1[k+1,\tilde{f}_{k-1}^{1}(x_{k-1},\gamma_k^{1*}(x_{k-1}))] \; ;\;V^1(T+1,x_T) = 0 \\\\
V^i(k,x_{k-1}) =  \min_{u_k^i \in U_k^i}[g_k^i(\tilde{f}_{k-1}^{i}(x_{k-1},u_k^i),\gamma_k^{1*}(x_{k-1}),\ldots,\gamma_k^{i-1*}(x_{k-1}),u_k^i,\\\gamma_k^{i+1*}(x_{k-1}),\ldots,\gamma_k^{n*}(x_{k-1}),x_{k-1}) +  V^i(k+1,\tilde{f}_{k-1}^{i}(x_{k-1},u_k^i))] = \\
g_k^i[\tilde{f}_{k-1}^{i}(x_{k-1},\gamma_k^{i*}(x_{k-1})),\gamma_k^{1*}(x_{k-1}),\ldots,\gamma_k^{n*}(x_{k-1}),x_{k-1}]\\ + V^i[k+1,\tilde{f}_{k-1}^{i}(x_{k-1},\gamma_k^{i*}(x_{k-1}))] \; ;\;V^i(T+1,x_T) = 0 \; ; \; i \in \{2\ldots n\} \\\\
where \;\; R_k^i \; (i \in \{2\ldots n\}) \;is\;a\;singleton\;set\;defined\;by\\
R_k^i(u_k^1)=\{r_k^i(u_k^1) \in \Gamma_k^i\; \: : \: g_k^i(f_k(x_{k-1},u_k^1,r_k^2(u_k^1),\ldots,r_k^n(u_k^1)),u_k^1,r_k^2(u_k^1),\ldots,r_k^n(u_k^1),x_{k-1}) \\+ V^i[k+1,f(x_{k-1},u_k^1,r_k^2(u_k^1),\ldots,r_k^n(u_k^1)) = \\\min_{u_k^i} g_k^i(f_k(x_{k-1},u_k^1,r_k^2(u_k^1),\ldots,r_k^{i-1}(u_k^1),u_k^i,r_k^{i+1}(u_k^1),\ldots,r_k^n(u_k^1)),u_k^1,r_k^2(u_k^1),\ldots,r_k^{i-1}(u_k^1),\\u_k^i,r_k^{i+1}(u_k^1),\ldots,r_k^n(u_k^1),x_{k-1}) + V^i[k+1,f(x_{k-1},u_k^1,r_k^2(u_k^1),\ldots,r_k^{i-1}(u_k^1),u_k^i,r_k^{i+1}(u_k^1),\ldots,r_k^n(u_k^1))]\}   \\\\
R_k^i(\gamma_k^{1*}) = \gamma_k^{i*} \\\\
\tilde{f}_{k-1}^{1}(x_{k-1},z) \; \widehat{=} f_{k-1}(x_{k-1},z,r_k^{2}(z),\ldots,r_k^{n}(z)) \\ \\
\tilde{f}_{k-1}^{i}(x_{k-1},z) \; \widehat{=} f_{k-1}(x_{k-1},\gamma_k^{1*}(x_{k-1}),\ldots,\gamma_k^{i-1*}(x_{k-1}),z,\gamma_k^{i+1*}(x_{k-1}),\ldots,\gamma_k^{n*}(x_{k-1})) \\\\\label{vfsl}
\end{multline}

\vspace*{-5ex}

Every such equilibrium solution is \emph {strongly time consistent}, and the corresponding Stackelberg equilibrium cost for \textbf {P}i is $V^i(1,x_0)$.\\\\

\label{FBS_opt}
\end{theo}

\textsc {Proof:}

Theorem \eqref{FBS_opt} can be proven in the same way as Theorem \eqref{FBN_Opt}, bearing in mind that the leader additionally accounts for the influence of his strategy on the followers' strategies when minimizing his cost functional.\;\;\;\;\;\;\;\;\;\;\boxed{}

\begin{rem}
For n-person affine-quadratic dynamic games (cf. Def. \eqref{defspielaq}) the assumption of a "unique follower response" ($R_k^i \; (i \in \{2\ldots n\})$ is a singleton set) is met if the  followers' cost functions are strictly convex over $R^{m_i}$.
\end{rem}
\clearpage
\subsection{Results for affine-quadratic games with one leader and arbitrarily many followers}

In the following, the results of Theorem \eqref{FBS_opt} are applied to an affine-quadratic dynamic game with one leader and arbitrarily many followers. Theorem \eqref{FBSTheo} is a generalization of Corollary 7.2 in Ba\c{s}ar and Olsder (1999, pp. 374-375)\cite{baol}. On the one hand a more general state equation and more general cost functionals are considered and on the other hand the number of followers is extended from one to arbitrarily many.\\

Furthermore in Proposition \eqref{profbsnb} the equivalence of the equations in Theorem \eqref{FBSTheo} with terminologically different equations is shown.

\begin{theo}
An n-person affine-quadratic dynamic game (cf. Def. \eqref{defspielaq}) admits a \emph{unique feedback Stackelberg equilibrium solution with one leader and arbitrarily many followers} if 
\begin{enumerate}
\item $Q_k^i \geq 0$, $R_k^{ii} > 0$ and $R_k^{ij} \geq 0$ (defined for $k \in K$ , $i,j \in N$, $j \not = i$)
\item \eqref{Pf}, \eqref{alphaf}, \eqref{rbar}, \eqref{W} and \eqref{w} admit unique optimal solutions $P_k^{i*}$, $\alpha_k^{i*}$,  $\bar r_k^{i}$, $W_k^{i}$ and $w_k^{i}$ (defined for $k \in K$ , $i \in \{2,\ldots,n\}$)
\end{enumerate}
If these conditions are satisfied, the unique equilibrium strategies are given by \eqref{gamma} and the corresponding feedback Stackelberg equilibrium cost for each player is stated in \eqref{cost}.\footnote {For all equations belonging to this theorem and its proof, $i \in N$ and $k \in K$ if nothing different is stated.}

\vspace*{-5ex} \begin{align}
f_{k-1}(x_{k-1},u_k^1,\ldots,u_k^n) = A_kx_{k-1} + \sum_{j\in N}B_k^ju_k^j+s_k \label{f}
\end{align}

\vspace*{-5ex} \begin{align}
L^i(x_0,u^{1},\ldots,u^{n})= \sum_{k = 1}^T g_k^i(x_k,u_k^{1},\ldots,u_k^{n},x_{k-1}) 
\end{align}

\vspace*{-5ex} \begin{multline}
 g_k^i(x_k,u_k^1,\ldots,u_k^n,x_{k-1}) = \frac{1}{2} (x_k^{'}Q_k^ix_k + \sum_{j\in N}u_k^{j'}R_k^{ij}u_k^j)\; + \\
 \frac{1}{2} (\tilde{x}_k^{i'}Q_k^i\tilde{x}_k^i + \sum_{j\in N}\tilde{u}_k^{ij'}R_k^{ij}\tilde{u}_k^{ij}) - \tilde{x}_k^{i'}Q_k^ix_k - \sum_{j\in N}\tilde{u}_k^{ij'}R_k^{ij}u_k^j  \\\label{g}
\end{multline} \vspace*{-5ex}

\vspace*{-5ex} \begin{align}
\gamma_k^{i*}(x_{k-1}) = u_k^{i*} = -P_k^{i*}x_{k-1}-\alpha_k^{i*} \label{gamma}
\end{align}

\vspace*{-5ex} \begin{multline}
P_k^{1*} = [(B_k^1 + \sum_{j\in 2\ldots n}B_k^j\bar {r}_k^{j})^{'}Z_k^1(B_k^1 + \sum_{j\in 2\ldots n}B_k^j\bar {r}_k^{j}) + R_k^{11} + \sum_{j\in 2\ldots n}\bar {r}_k^{j'}R_k^{1j}\bar {r}_k^{j}]^{-1} \\ [(B_k^1 + \sum_{j\in 2\ldots n}B_k^j\bar {r}_k^{j})^{'}Z_k^1(A_k + \sum_{j\in 2\ldots n}B_k^jW_k^{j}) + \sum_{j\in 2\ldots n}\bar {r}_k^{j'}R_k^{1j} W_k^{j}]\\\label{P1}
\end{multline} \vspace*{-5ex}

\vspace*{-5ex} \begin{multline}
\alpha_k^{1*} = [(B_k^1 + \sum_{j\in 2\ldots n}B_k^j\bar {r}_k^{j})^{'}Z_k^1(B_k^1 + \sum_{j\in 2\ldots n}B_k^j\bar {r}_k^{j}) + R_k^{11}  \\ + \sum_{j\in 2\ldots n}\bar {r}_k^{j'}R_k^{1j}\bar {r}_k^{j}]^{-1} [ (B_k^1 + \sum_{j\in 2\ldots n}B_k^j\bar {r}_k^{j})^{'}Z_k^1(s_k - \sum_{j\in 2\ldots n}B_k^jw_k^{j})\\ + \sum_{j\in 2\ldots n}\bar {r}_k^{j'}R_k^{1j} w_k^{j}  - R_k^{11}\tilde{u}_k^{11} \\- \sum_{j\in 2\ldots n}\bar {r}_k^{j'}R_k^{1j} \tilde{u}_k^{1j} + (B_k^1 + \sum_{j\in 2\ldots n}B_k^j\bar{r}_k^{j})^{'}(\zeta_{k}^1 - Q_k^1\tilde{x}_k^1)]\\\label{alpha1}
\end{multline} \vspace*{-5ex}

\vspace*{-5ex} \begin{multline}
P_k^{i*} = (R_k^{ii} + B_k^{i'}Z_k^iB_k^i)^{-1}[B_k^{i'}Z_k^i(A_k - B_k^1P_k^{1*} \\- \sum_{j\in 2\ldots n \;,\; j \not = i}B_k^jP_k^{j*})] = - W_T^{i} + \bar {r}_T^{i}P^{1*}_{T} \; ; \; i \in \{2\ldots n\} \\\label{Pf} 
\end{multline} \vspace*{-5ex}

\vspace*{-5ex} \begin{multline}
\alpha_k^{i*} = (R_k^{ii} + B_k^{i'}Z_k^iB_k^i)^{-1} [B_k^{i'}(Z_k^i(s_k - B_k^1\alpha_k^{1*} - \sum_{j\in 2\ldots n \;,\; j \not = i}B_k^j\alpha_k^{j*})\\+ \zeta_{k}^i -Q_k^i\tilde{x}_k^i)-R_k^{ii}\tilde{u}_k^{ii}] = -w_k^i + \bar {r}_k^{i}\alpha_{k}^{1*} \; ; \; i \in \{2\ldots n\}\\\label{alphaf}
\end{multline} \vspace*{-5ex}

\vspace*{-5ex} \begin{multline}
Z_{k-1}^i = (A_k - \sum_{j \in N}B_k^jP_k^{j*})^{'}Z_k^i(A_k - \sum_{j \in N}B_k^jP_k^{j*}) \\ + \sum_{j\in N}P_k^{j*'}R_k^{ij}P_k^{j*} + Q_{k-1}^i\; ;\; Z_T^i = Q_T^i \\\label{Zfs}
\end{multline} \vspace*{-5ex}

\vspace*{-5ex} \begin{align}
\bar r_k^{i} = -(R_k^{ii} + B_k^{i'}Z_k^iB_k^i)^{-1}[B_k^{i'}Z_k^i(B_k^1 + \sum_{j\in 2\ldots n \;,\; j \not = i}B_k^j\bar r_k^{j})] \; ; \; i \in \{2\ldots n\} \label{rbar}
\end{align}

\vspace*{-5ex} \begin{multline}
\zeta_{k-1}^i = (A_k - \sum_{j \in N}B_k^jP_k^{j*})^{'}[\zeta_{k}^i + Z_k^i(s_k-\sum_{j\in N}B_k^j\alpha_k^{j*})-Q_k^i\tilde{x}_k^i] \\ + \sum_{j\in N}P_k^{j*'}R_k^{ij}(\alpha_k^{j*} + \tilde{u}_k^{ij})\; ;\; \zeta_{T}^i = 0 \\\label{zetaki}
\end{multline} \vspace*{-5ex}

\vspace*{-5ex} \begin{align}
W_k^{i} = -(R_k^{ii} + B_k^{i'}Z_k^iB_k^i)^{-1}[B_k^{i'}Z_k^i(\sum_{j\in 2\ldots n\;,\; j \not = i}B_k^jW_k^{j} + A_k)] \; ; \; i \in \{2\ldots n\}\label{W}
\end{align}

\vspace*{-5ex} \begin{multline}
w_k^{i} = -(R_k^{ii} + B_k^{i'}Z_k^iB_k^i)^{-1}[B_k^{i'}(Z_k^i(\sum_{j\in 2\ldots n\;,\; j \not = i}B_k^jw_k^{j} + s_k) \\+ \zeta_{k}^i - Q_k^i\tilde{x}_k^i)-R_k^{ii}\tilde{u}_k^{ii}] \; ; \; i \in \{2\ldots n\}\\\label{w}
\end{multline} \vspace*{-5ex}

\vspace*{-5ex} \begin{align}
V^i(1,x_0) = \frac{1}{2} x_0^{'}Z_0^ix_0 + \zeta_0^{i'}x_0 + n_0^i \label{cost}
\end{align}

\vspace*{-5ex} \begin{multline}
n_{k-1}^i = n_{k}^i + \frac{1}{2} (s_k-\sum_{j\in N}B_k^j\alpha_k^{j*})^{'}Z_k^i(s_k-\sum_{j\in N}B_k^j\alpha_k^{j*}) \\+ \zeta_{k}^{i'}(s_k-\sum_{j\in N}B_k^j\alpha_k^{j*}) + \frac{1}{2}\sum_{j\in N }\alpha_k^{j*'}R_k^{ij}\alpha_k^{j*} -  \tilde{x}_k^{i'}Q_k^i(s_k-\sum_{j\in N}B_k^j\alpha_k^{j*}) \\ + \sum_{j\in N}\tilde{u}_k^{ij'}R_k^{ij}\alpha_k^{j*} + \frac{1}{2}(\tilde{x}_k^{i'}Q_k^i\tilde{x}_k^i + \sum_{j\in N}\tilde{u}_k^{ij'}R_k^{ij}\tilde{u}_k^{ij}) \; ;\; n_{T}^i = 0 \\\\\label{n}
\end{multline} \vspace*{-5ex}

\label{FBSTheo}
\end{theo}

\textsc {Proof:}\footnote{The basis and the last part of the inductive step (after finding the optimal strategies for the leader and the followers) are proven in the same way as in Theorem \eqref{FBNTheo}, because in these parts the distinction between leader and followers is not essential.}

The proof is done using an induction argument to show that the strategy of the leader and the strategies of the followers $\gamma_k^{i*}(x_{k-1})$, given by \eqref{gamma}, minimize\footnote{Taking into consideration the structural advantage of the leader, of course.} the strictly convex functionals \eqref{Vdls} and \eqref{Vdfs} at each stage of the game. But the minimization of \eqref{Vdls} and \eqref{Vdfs} is exactly \eqref{vfsl} applied to the specific state equation and cost functionals of the above game and therefore, considering Theorem \eqref{FBS_opt}, it follows that the $\gamma_k^{i*}(x_{k-1})$ are the unique equilibrium strategies.

\vspace*{-5ex} \begin{multline}
\frac{1}{2} ((A_kx_{k-1} + B_k^1u_k^1 + \sum_{j\in \{2\ldots n\}}B_k^jr_k^j+s_k)^{'}Q_k^1(A_kx_{k-1} + B_k^1u_k^1 \\+ \sum_{j\in \{2\ldots n\}}B_k^jr_k^j+s_k) + u_k^{1'}R_k^{11}u_k^1 + \sum_{j\in \{2\ldots n\}}r_k^{j'}R_k^{1j}r_k^j) +
\frac{1}{2} (\tilde{x}_k^{1'}Q_k^1\tilde{x}_k^1 \\+ \sum_{j\in \{2\ldots n\}}\tilde{u}_k^{1j'}R_k^{1j}\tilde{u}_k^{1j}) - \tilde{x}_k^{1'}Q_k^1(A_kx_{k-1} + B_k^1u_k^1 + \sum_{j\in \{2\ldots n\}}B_k^jr_k^j+s_k) - \tilde{u}_k^{11'}R_k^{11}u_k^{1} \\- \sum_{j\in \{2\ldots n\}}\tilde{u}_k^{1j'}R_k^{1j}r_k^j + V^1(k+1,\tilde{f}_{k-1}^{1}(x_{k-1},u_k^1)) \; ;\; V^1(T+1,x_T) = 0 \\\label{Vdls}
\end{multline} \vspace*{-5ex}

\vspace*{-5ex} \begin{multline}
\frac{1}{2} ((A_kx_{k-1} + B_k^iu_k^i + \sum_{j\in N\;,\; j \not = i}B_k^ju_k^{j*}+s_k)^{'}Q_k^i(A_kx_{k-1} + B_k^iu_k^i \\+ \sum_{j\in N\;,\; j \not = i}B_k^ju_k^{j*}+s_k) + u_k^{i'}R_k^{ii}u_k^i + \sum_{j\in N\;,\; j \not = i}u_k^{j*'}R_k^{ij}u_k^{j*}) +
\frac{1}{2} (\tilde{x}_k^{i'}Q_k^i\tilde{x}_k^i \\+ \sum_{j\in N \;,\; j \not = i}\tilde{u}_k^{ij'}R_k^{ij}\tilde{u}_k^{ij}) - \tilde{x}_k^{i'}Q_k^i(A_kx_{k-1} + B_k^iu_k^i + \sum_{j\in N\;,\; j \not = i}B_k^ju_k^{j*}+s_k) - \tilde{u}_k^{ii'}R_k^{ii}u_k^{i} \\- \sum_{j\in N\;,\; j \not = i}\tilde{u}_k^{ij'}R_k^{ij}u_k^{j*} + V^i(k+1,\tilde{f}_{k-1}^{i}(x_{k-1},u_k^i)) \; ;\; V^i(T+1,x_T) = 0 \; ; \; i \in \{2\ldots n\} \\\label{Vdfs}
\end{multline} \vspace*{-5ex}

Furthermore the optimal reactions $r_k^j$ of the followers to an arbitrary strategy $u_k^1$ of the leader at stage \textit {k}, subject to the assumption that all strategies for all players from stage \textit {k}+1 to stage \textit {T} are optimal, are derived by minimization over $u_l^i$ of the below equations

\vspace*{-5ex} \begin{multline}
\frac{1}{2} ((A_lx_{l-1} + B_l^1u_l^1 + B_l^iu_l^i + \sum_{j \in \{2\ldots n\}\;,\; j \not = i}B_l^jr_l^{j}+s_l)^{'}Q_l^i (A_lx_{l-1} + B_l^1u_l^1 + B_l^iu_l^i \\+ \sum_{j\in \{2\ldots n\}\;,\; j \not = i}B_l^jr_l^{j}+s_l) + u_l^{1'}R_l^{i1}u_l^1 + u_l^{i'}R_l^{ii}u_l^i + \sum_{j\in \{2\ldots n\}\;,\; j \not = i}r_l^{j'}R_l^{ij}r_l^{j})\; \\+ \frac{1}{2} (\tilde{x}_l^{i'}Q_l^i\tilde{x}_l^i + \sum_{j\in N}\tilde{u}_l^{ij'}R_l^{ij}\tilde{u}_l^{ij}) - \tilde{x}_l^{i'}Q_l^i(A_lx_{l-1} + B_l^1u_l^1 + B_l^iu_l^i \\+ \sum_{j\in \{2\ldots n\}\;,\; j \not = i}B_l^jr_l^{j}+s_l) - \tilde{u}_l^{i1'}R_l^{i1}u_l^{1} - \tilde{u}_l^{ii'}R_l^{ii}u_l^{i} - \sum_{j\in \{2\ldots n\}\;,\; j \not = i}\tilde{u}_l^{ij'}R_l^{ij}r_l^{j} \\+ V^i(k+1,f_{k-1}^{i}(x_{k-1},u_k^1, r_k^{2}, \ldots, r_k^{i-1}, u_k^i, r_k^{i+1}, \ldots, r_k^{n})) \; ; \; i \in \{2\ldots n\}\\\label{Vdfr}
\end{multline} \vspace*{-5ex}
%\vspace*{-5ex} \begin{multline}
%r_k^{i} = -(R_k^{ii} + B_k^{i'}Z_k^iB_k^i)^{-1} [B_k^{i'}(Z_k^i(B_k^1u_k^1 + \sum_{j \in 2 \ldots n \; , \; j \not = i}B_k^jr_k^j\\ + A_kx_{k-1} + s_k) +\zeta_{k}^i - Q_k^i\tilde{x}_k^i) - R_k^{ii}\tilde{u}_k^{ii}] \; ; \; i \in \{2\ldots n\}\\\label{r}
%\end{multline} \vspace*{-5ex}

The induction argument runs from \textit {T}+1 to 1 and proves that the value function for player \textit {i} at stage \textit {k} can be written as stated below in \eqref{Vg}.

\vspace*{-5ex} \begin{multline}
\frac{1}{2} x_{k-1}^{'}(Z_{k-1}^i-Q_{k-1}^i)x_{k-1} + \zeta_{k-1}^{i'}x_{k-1} + n_{k-1}^i \\\label{Vg}
\end{multline} \vspace*{-5ex}

%\vspace*{-5ex} \begin{multline}
%g_k^1(f_{k-1}(x_{k-1},u_k^{1},r_k^2,\ldots,r_k^n),u_k^{1},r_k^2,\ldots,r_k^{n}, x_{k-1}) \\+  V^1(k+1,f_{k-1}(x_{k-1},u_k^{1},r_k^2,\ldots,r_k^n)\; ;\;V^1(T+1,x_T) = 0 \\\label{Vl}
%\end{multline} \vspace*{-5ex}
%
%\vspace*{-5ex} \begin{multline}
%g_k^i(f_k(x_{k-1},u_k^1,r_k^2,\ldots,r_k^{i-1},u_k^i,r_k^{i+1},\ldots,r_k^n),u_k^1,r_k^2,\ldots,r_k^{i-1},u_k^i,\\r_k^{i+1},\ldots,r_k^n,x_{k-1}) + V^i[k+1,f(x_{k-1},u_k^1,r_k^2,\ldots,r_k^{i-1},u_k^i,r_k^{i+1},\ldots,r_k^n)] \\V^i(T+1,x_T) = 0 \; ; \; i \in \{2\ldots n\}\\\label{Vf}
%\end{multline} \vspace*{-5ex}

\textbf{Basis:}

The induction starts at \textit {k} = \textit {T}. First we make use of the general optimality conditions for $V_k^i$ at stage \textit {T}+1.

\vspace*{-5ex} \begin{multline}
V^i(T+1,x_T) = 0\\\label{optVlss}
\end{multline} \vspace*{-5ex}

Now we show that $\eqref{Vg}_{k=T+1}$ fulfills \eqref{optVlss}.

\vspace*{-5ex} \begin{multline}
\nonumber\eqref{Vg}_{k=T+1} \; \; \frac{1}{2} x_{T}^{'}(Z_{T}^i-Q_{T}^i)x_{T} + \zeta_{T}^{i'}x_{T} + n_{T}^i\\
\end{multline} \vspace*{-5ex}

Making use of $\eqref{Zfs}_{k=T+1}$, $\eqref{zetaki}_{k=T+1}$ and $\eqref{n}_{k=T+1}$ yields

\vspace*{-5ex} \begin{multline}
\nonumber \frac{1}{2} x_{T}^{'}(Q_{T}^i-Q_{T}^i)x_{T} + 0^{'}x_{T} + 0 = 0 = V^i(T+1,x_T)\\
\end{multline} \vspace*{-5ex}

\textbf{Inductive step:}

As an induction hypothesis, the equations $\eqref{Vg}_{k=l+1}$ are assumed to be equal to $V^i(l+1,x_l)$ respectively. Now we have to prove that the relation also holds at stage \textit {l}. In other words: we have to show $\eqref{Vg}_{k=l} = V^i(l,x_{l-1})$.\\

\vspace*{-5ex} \begin{multline}
\nonumber\eqref{Vg}_{k=l+1} \;\; \frac{1}{2} x_{l}^{'}(Z_{l}^i-Q_{l}^i)x_{l} + \zeta_{l}^{i'}x_{l} + n_{l-1}^i \\
\end{multline} \vspace*{-5ex}

The inductive step is done by  first showing what the optimal response of the followers to an arbitrary strategy of the leader looks like. In other words: first we deduce the optimal reactions $r_l^{j}$ ($j \in \{2, \ldots, n\}$) of the followers. Then we derive the optimal strategy of the leader by minimizing $\eqref{Vdls}_{k=l}$ over $u_l^1$ considering the optimal reactions of the followers, and finally we derive the optimal strategies of the followers as the optimal reactions of the followers to the optimal strategy of the leader.\\

Using $\eqref{f}_{k=l}$ in $\eqref{Vg}_{k=l+1}$, bearing in mind that the optimal control vector of the leader remains to be deduced in the further induction argument, gives

\vspace*{-5ex} \begin{multline}
\frac{1}{2} (A_lx_{l-1} + B_l^1u_l^1 + \sum_{j\in \{2\ldots n\}}B_l^jr_l^{j}+s_l)^{'}(Z_{l}^1-Q_{l}^1)(A_lx_{l-1} \\+ B_l^1u_l^1 + \sum_{j\in \{2\ldots n\}}B_l^jr_l^{j}+s_l) + \zeta_{l}^{1'}(A_lx_{l-1} + B_l^1u_l^1 + \sum_{j\in \{2\ldots n\}}B_l^jr_l^{j}+s_l) + n_{l}^1 \\\label{Vgil}
\end{multline} \vspace*{-5ex}

for the leader and

\vspace*{-5ex} \begin{multline}
\frac{1}{2} (A_lx_{l-1} + B_l^1u_l^1 + B_l^iu_l^i + \sum_{j\in \{2\ldots n\}\;,\; j \not = i}B_l^jr_l^{j}+s_l)^{'}(Z_{l}^i-Q_{l}^i)\\(A_lx_{l-1} + B_l^1u_l^1 + B_l^iu_l^i + \sum_{j \in \{2\ldots n\}\;,\; j \not = i}B_l^jr_l^{j}+s_l) + \zeta_{l}^{i'}(A_lx_{l-1} + B_l^1u_l^1 \\+ B_l^iu_l^i + \sum_{j\in \{2\ldots n\}\;,\; j \not = i}B_l^jr_l^{j*}+s_l) + n_{l}^i \; ; \; i \in \{2\ldots n\}\\\label{Vgif}
\end{multline} \vspace*{-5ex}

for the followers.\\

We start with the derivation of the optimal reactions of the followers $r_l^{i}$ ($i \in \{2, \ldots, n\}$) to an arbitrary strategy $u_l^1$ of the leader. To do this we substitute $V^i(l+1,\tilde{f}_{k-1}^{i}(x_{k-1},u_k^1, r_k^{2}, \ldots, r_k^{i-1}, u_k^i, r_k^{i+1}, \ldots, r_k^{n}))$ ($i \in \{2, \ldots, n\}$) in $\eqref{Vdfr}_{k=l}$ with the help of the induction hypothesis given by \eqref{Vgif}:

\vspace*{-5ex} \begin{multline}
\frac{1}{2} ((A_lx_{l-1} + B_l^1u_l^1 + B_l^iu_l^i + \sum_{j\in N\;,\; j \not = i}B_l^jr_l^{j}+s_l)^{'}Q_l^i (A_lx_{l-1} + B_l^1u_l^1 + B_l^iu_l^i \\+ \sum_{j\in \{2\ldots n\}\;,\; j \not = i}B_l^jr_l^{j}+s_l) + u_l^{1'}R_l^{i1}u_l^1 + u_l^{i'}R_l^{ii}u_l^i + \sum_{j\in \{2\ldots n\}\;,\; j \not = i}r_l^{j'}R_l^{ij}r_l^{j})\; + \frac{1}{2} (\tilde{x}_l^{i'}Q_l^i\tilde{x}_l^i \\+ \sum_{j\in N}\tilde{u}_l^{ij'}R_l^{ij}\tilde{u}_l^{ij}) - \tilde{x}_l^{i'}Q_l^i(A_lx_{l-1} + B_l^1u_l^1 + B_l^iu_l^i + \sum_{j\in \{2\ldots n\}\;,\; j \not = i}B_l^jr_l^{j}+s_l) - \\\tilde{u}_l^{i1'}R_l^{i1}u_l^{1} - \tilde{u}_l^{ii'}R_l^{ii}u_l^{i} - \sum_{j\in \{2\ldots n\}\;,\; j \not = i}\tilde{u}_l^{ij'}R_l^{ij}r_l^{j} + \frac{1}{2} (A_lx_{l-1} + B_l^1u_l^1 + B_l^iu_l^i + \\\sum_{j\in \{2\ldots n\}\;,\; j \not = i}B_l^jr_l^{j}+s_l)^{'}(Z_{l}^i-Q_{l}^i)(A_lx_{l-1} + B_l^1u_l^1 + B_l^iu_l^i + \sum_{j\in \{2\ldots n\}\;,\; j \not = i}B_l^jr_l^{j}+s_l) \\+ \zeta_{l}^{i'}(A_lx_{l-1} + B_l^1u_l^1 + B_l^iu_l^i + \sum_{j\in \{2\ldots n\}\;,\; j \not = i}B_l^jr_l^{j}+s_l) + n_{l}^i\; ; \; i \in \{2\ldots n\}\\\label{Vlf}
\end{multline} \vspace*{-5ex}

$\eqref{Vlf}$ is strictly convex in $u_l^i$ . This can be seen by applying Corollary \eqref{pdc} to \eqref{consfbr}. Therefore there has to be a unique optimal response of follower \textit{i} ($i \in \{2, \ldots, n\}$) to an arbitrary strategy $u_l^1$ of the leader.

\vspace*{-5ex} \begin{multline}
\frac {\partial^2} {\partial u_l^{i^2}} \eqref{Vlf} \;\;\ R_l^{ii} + B_l^{i'}Z_l^iB_l^i \; ; \; i \in \{2\ldots n\}\\\label{consfbr}
\end{multline} \vspace*{-5ex}

This unique optimal response can be found by using the first-order necessary and sufficient (because of strict convexity of \eqref{Vlf}) conditions for minimization

\vspace*{-5ex} \begin{multline}
\nonumber\frac {\partial} {\partial u_l^{i}} \eqref{Vlf} = 0 \Rightarrow \;\;\ -(R_l^{ii} + B_l^{i'}Z_l^iB_l^i)r_l^{i} - B_l^{i'}Z_l^i \sum_{j\in \{2\ldots n\} \;,\; j \not = i}B_l^jr_l^{j} = \\B_l^{i'}[Z_l^i(A_lx_{l-1} + s_l) + \zeta_{l}^i  - Q_l^i\tilde{x}_l^i] - R_l^{ii}\tilde{u}_l^{ii} + B_l^{i'}Z_l^iB_l^1u_l^1 \; ; \; i \in \{2\ldots n\}\\
\end{multline} \vspace*{-5ex}

Making $r_l^{i}$ explicit leads to

\vspace*{-5ex} \begin{multline}
r_l^i = -(R_l^{ii} + B_l^{i'}Z_l^iB_l^i)^{-1} [ B_l^{i'}(Z_l^i(B_l^1u_l^1 + \sum_{j \in 2 \ldots n \; , \; j \not = i}B_l^jr_l^j \\+ A_lx_{l-1} + s_l) +\zeta_{l}^i  - Q_l^i\tilde{x}_l^i) - R_l^{ii}\tilde{u}_l^{ii}] \; ; \; i \in \{2\ldots n \} \\\label{rli}
\end{multline} \vspace*{-5ex}

As the right hand side of the above equation is affine in $x_{l-1}$ and $u_l^1$, the left hand side also has to be affine in $x_{l-1}$ and $u_l^1$. Therefore the substitution

\vspace*{-5ex} \begin{multline}
r_l^{i} = W_l^{i}x_{l-1} + \bar {r}_l^{i}u_l^1 + w_l^i \; ; \; i \in \{2\ldots n\}\\\label{rlia}
\end{multline} \vspace*{-5ex}

is allowed and leads to

\vspace*{-5ex} \begin{multline}
W_l^{i}x_{l-1} + \bar {r}_l^{i}u_l^1 + w_l^i = -(R_l^{ii} + B_l^{i'}Z_l^iB_l^i)^{-1} [ B_l^{i'}(Z_l^i(B_l^1u_l^1 \\+ \sum_{j \in 2 \ldots n \; , \; j \not = i}B_l^j(W_l^{j}x_{l-1} + \bar {r}_l^{j}u_l^1 + w_l^j) + A_lx_{l-1} + s_l) +\zeta_{l}^i  - Q_l^i\tilde{x}_l^i) - R_l^{ii}\tilde{u}_l^{ii}] \; ; \; i \in \{2\ldots n\}\\\label{puli}
\end{multline} \vspace*{-5ex}

Differentiating \eqref{rlia} gives

\vspace*{-5ex} \begin{multline}
\nonumber \frac {\partial} {\partial u_l^1}\eqref{rlia} \;\;\ \frac {\partial} {\partial u_l^1}r_l^{i} = \bar {r}_l^{i}\; ; \; i \in \{2\ldots n\}\ \\
\end{multline} \vspace*{-5ex}

Differentiating \eqref{puli} and using the above equation yields

\vspace*{-5ex} \begin{multline}
\frac {\partial} {\partial u_l^1}\eqref{puli} \;\; \bar r_l^{i} = -(R_l^{ii} + B_l^{i'}Z_l^iB_l^i)^{-1} \\ [B_l^{i'}Z_l^i(B_l^1 + \sum_{j\in 2\ldots n \;,\; j \not = i}B_l^j\bar r_l^{j})] \; ; \; i \in \{2\ldots n\} \\\label{rbarl}
\end{multline} \vspace*{-5ex}

Rearranging the terms in \eqref{rlia} gives

\vspace*{-5ex} \begin{multline}
\nonumber W_l^{i}x_{l-1} + w_l^i = r_l^{i} - \bar {r}_l^{i}u_l^1 \; ; \; i \in \{2\ldots n\}\\
\end{multline} \vspace*{-5ex}

Now we substitute $r_l^{i}$ and $\bar {r}_l^{i}$ with the help of \eqref{rli} and \eqref{rbarl}

\vspace*{-5ex} \begin{multline}
\nonumber W_l^{i}x_{l-1} + w_l^i = -(R_l^{ii} + B_l^{i'}Z_l^iB_l^i)^{-1} \\ [ B_l^{i'}(Z_l^i(B_l^1u_l^1 + \sum_{j \in 2 \ldots n \; , \; j \not = i}B_l^jr_l^j + A_lx_{l-1} + s_l) + \zeta_{l}^i - Q_l^i\tilde{x}_l^i) - R_l^{ii}\tilde{u}_l^{ii}] \\+ (R_l^{ii} + B_l^{i'}Q_l^iB_l^i)^{-1}[-B_l^{i'}Z_l^i(B_l^1 + \sum_{j\in 2\ldots n \;,\; j \not = i}B_l^j\bar r_l^{j})]u_l^1 \; ; \; i \in \{2\ldots n\}\\
\end{multline} \vspace*{-5ex}

Making use of \eqref{rlia} leads to

\vspace*{-5ex} \begin{multline}
W_l^{i}x_{l-1} + w_l^i = -(R_l^{ii} + B_l^{i'}Z_l^iB_l^i)^{-1}  [ B_l^{i'}(Z_l^i(\sum_{j \in 2 \ldots n \; , \; j \not = i}B_l^j\\(W_l^{j}x_{l-1} + w_l^j) + A_lx_{l-1} + s_l) + \zeta_{l}^i - Q_l^i\tilde{x}_l^i) - R_l^{ii}\tilde{u}_l^{ii}] \; ; \; i \in \{2\ldots n\}\\\label{Wil}
\end{multline} \vspace*{-5ex}

Comparing the coefficients it follows that

\vspace*{-5ex} \begin{multline}
\nonumber\eqref{Wil}_{x_{l-1}}=\eqref{W}_{k=l} \;\;W_l^{i} = (R_k^{ii} + B_l^{i'}Z_l^iB_l^i)^{-1} [B_l^{i'}Z_l^i(\sum_{j\in 2\ldots n\;,\; j \not = i}B_l^jW_l^{j} + A_l)] \; ; \; i \in \{2\ldots n\}\\
\end{multline} \vspace*{-5ex}

\vspace*{-5ex} \begin{multline}
\nonumber\eqref{Wil}_{const.}=\eqref{w}_{k=l} \;\; w_l^{i} = -(R_l^{ii} + B_l^{i'}Z_l^iB_l^i)^{-1} [B_l^{i'}(Z_l^i(\sum_{j\in 2\ldots n\;,\; j \not = i}B_l^jw_l^{j} \\+ s_l) + \zeta_{l}^i - Q_l^i\tilde{x}_l^i) - R_l^{ii}\tilde{u}_l^{ii}] \; ; \; i \in \{2\ldots n\}\\
\end{multline} \vspace*{-5ex}

After deducing the unique optimal response of follower \textit{i} ($i \in \{2, \ldots, n\}$) to an arbitrary strategy $u_l^1$ of the leader, we can substitute $V^1(k+1,\tilde{f}_{k-1}^{1}(x_{k-1},u_k^1))$ in $\eqref{Vdls}_{k=l}$ ($i \in \{2, \ldots, n\}$) with the help of the induction hypothesis for the leader, given by \eqref{Vgil}.

\vspace*{-5ex} \begin{multline}
\nonumber \eqref{Vdls}_{k=l} \; \; \frac{1}{2} ((A_lx_{l-1} + B_l^1u_l^1 + \sum_{j\in \{2\ldots n\}}B_l^jr_l^j+s_l)^{'}Q_l^1(A_lx_{l-1} + B_l^1u_l^1 \\+ \sum_{j\in \{2\ldots n\}}B_l^jr_l^j+s_l) + u_l^{1'}R_l^{11}u_l^1 + \sum_{j\in \{2\ldots n\}}r_l^{j'}R_l^{1j}r_l^j) +
\frac{1}{2} (\tilde{x}_l^{1'}Q_l^1\tilde{x}_l^1 \\+ \sum_{j\in \{2\ldots n\}}\tilde{u}_l^{1j'}R_l^{1j}\tilde{u}_l^{1j}) - \tilde{x}_l^{1'}Q_l^1(A_lx_{l-1} + B_l^1u_l^1 + \sum_{j\in \{2\ldots n\}}B_l^jr_l^j+s_l) - \tilde{u}_l^{11'}R_l^{11}u_l^{1} \\- \sum_{j\in \{2\ldots n\}}\tilde{u}_l^{1j'}R_l^{1j}r_l^j + V^1(l+1,\tilde{f}_{l-1}^{1}(x_{l-1},u_l^1)) \\
\\
\end{multline} \vspace*{-5ex}

\vspace*{-5ex} \begin{multline}
\nonumber \frac{1}{2} ((A_lx_{l-1} + B_l^1u_l^1 + \sum_{j\in \{2\ldots n\}}B_l^jr_l^j+s_l)^{'}Q_l^1(A_lx_{l-1} + B_l^1u_l^1 + \sum_{j\in \{2\ldots n\}}B_l^jr_l^j+s_l) \\+ u_l^{1'}R_l^{11}u_l^1 + \sum_{j\in \{2\ldots n\}}r_l^{j'}R_l^{1j}r_l^j)\; +
\frac{1}{2} (\tilde{x}_l^{1'}Q_l^1\tilde{x}_l^1 + \sum_{j\in \{2\ldots n\}}\tilde{u}_l^{1j'}R_l^{1j}\tilde{u}_l^{1j}) \\- \tilde{x}_l^{1'}Q_l^1(A_lx_{l-1} + B_l^1u_l^1 + \sum_{j\in \{2\ldots n\}}B_l^jr_l^j+s_l) - \tilde{u}_l^{11'}R_l^{11}u_l^{1} - \sum_{j\in \{2\ldots n\}}\tilde{u}_l^{1j'}R_l^{1j}r_l^j \\+ \frac{1}{2} (A_lx_{l-1} + B_l^1u_l^1 + \sum_{j\in \{2\ldots n\}}B_l^jr_l^{j}+s_l)^{'}(Z_{l}^1-Q_{l}^1)(A_lx_{l-1} + B_l^1u_l^1 + \sum_{j\in \{2\ldots n\}}B_l^jr_l^{j}+s_l) \\+ \zeta_{l}^{1'}(A_lx_{l-1} + B_l^1u_l^1 + \sum_{j\in \{2\ldots n\}}B_l^jr_l^{j}+s_l) + n_{l}^1  \\
\end{multline} \vspace*{-5ex}

Next we substitute $r_k^i$ using \eqref{rlia}

\vspace*{-5ex} \begin{multline}
\frac{1}{2} ((A_lx_{l-1} + B_l^1u_l^1 + \sum_{j \in 2 \ldots n}B_l^j(W_l^{i}x_{l-1} + w_l^i + \bar {r}_l^{i}u_l^1)+s_l)^{'}\\Z_l^1 (A_lx_{l-1} + B_l^1u_l^1 + \sum_{j \in 2 \ldots n}B_l^j(W_l^{i}x_{l-1} + w_l^i + \bar {r}_l^{i}u_l^1)+s_l) + u_l^{1'}R_l^{11}u_l^1 + \\\sum_{j \in 2 \ldots n}r_l^{j'}R_l^{1j}(W_l^{i}x_{l-1} + w_l^i + \bar {r}_l^{i}u_l^1))\; +
 \frac{1}{2} (\tilde{x}_l^{1'}Q_l^1\tilde{x}_l^1 + \sum_{j\in N}\tilde{u}_l^{1j'}R_l^{1j}\tilde{u}_l^{1j}) - \tilde{x}_l^{1'}Q_l^1(A_lx_{l-1} + B_l^1u_l^1 + \\\sum_{j \in 2 \ldots n}B_l^j(W_l^{i}x_{l-1} + w_l^i + \bar {r}_l^{i}u_l^1)+s_l) - \tilde{u}_l^{11'}R_l^{11}u_l^{1} - \sum_{j \in 2 \ldots n}\tilde{u}_l^{1j'}R_l^{1j}(W_l^{i}x_{l-1} + w_l^i + \bar {r}_l^{i}u_l^1) + \\\frac{1}{2} (A_lx_{l-1} + B_l^1u_l^1 + \sum_{j \in 2 \ldots n}B_l^j(W_l^{i}x_{l-1} + w_l^i + \bar {r}_l^{i}u_l^1)+s_l)^{'}(Z_{l}^1-Q_{l}^1)\\(A_lx_{l-1} + B_l^1u_l^1 + \sum_{j \in 2 \ldots n}B_l^j(W_l^{i}x_{l-1} + w_l^i + \bar {r}_l^{i}u_l^1)+s_l) + \\\zeta_{l}^{1'}(A_lx_{l-1} + B_l^1u_l^1 + \sum_{j \in 2 \ldots n}B_l^j(W_l^{i}x_{l-1} + w_l^i + \bar {r}_l^{i}u_l^1)+s_l) + n_{l}^1\\\\\label{Vll}
\end{multline} \vspace*{-5ex}

$\eqref{Vdls}_{k=l}$ is strictly convex in $u_l^1$. This can be seen by applying Corollary \eqref{pdc} to \eqref{consll}. Therefore there has to be a unique optimal strategy of the leader at stage \textit{l}.

\vspace*{-5ex} \begin{multline}
\frac {\partial^2} {\partial u_l^{1^2}} \eqref{Vll} \;\;\ R_l^{11} + \sum_{j\in 2\ldots n}\bar {r}_l^{j'}R_l^{1j}\bar {r}_l^{j} + (B_l^{1} + \sum_{j\in 2\ldots n}B_l^j\bar {r}_l^{j})^{'}Z_l^i(B_l^{1} + \sum_{j\in 2\ldots n}B_l^j\bar {r}_l^{j}) \\\label{consll}
\end{multline} \vspace*{-5ex}

This unique optimal strategy for the leader can be found by using the first-order necessary and sufficient (because of strict convexity of $\eqref{Vdls}_{k=l}$) conditions for minimization

\vspace*{-5ex} \begin{multline}
\nonumber\frac {\partial} {\partial u_l^{1}} \eqref{Vll} \; = 0 \;\Rightarrow \; (B_l^{1} + \sum_{j\in 2\ldots n}B_l^j\bar {r}_l^{j})^{'}Z_l^1(A_lx_{l-1} + B_l^1u_l^{1*} + \sum_{j\in 2\ldots n}B_l^j\\(W_l^{j}x_{l-1}+w_l^j + \bar {r}_l^{j}u_l^{1*}) + s_l) +R_l^{11}u_l^{1*} + \sum_{j\in 2\ldots n}\bar {r}_l^{j'}R_l^{1j}(W_l^{j}x_{l-1}+w_l^j+ \\ \bar {r}_l^{j}u_l^{1*}) - (B_l^{1} + \sum_{j\in 2\ldots n}B_l^j\bar {r}_l^{j})^{'}Q_l^1\tilde{x}_l^{1} - R_l^{11}\tilde{u}_l^{11} - \sum_{j\in 2\ldots n}\bar {r}_l^{j'}R_l^{1j}\tilde{u}_l^{1j} + (B_l^{1} + \sum_{j\in 2\ldots n}B_l^j\bar {r}_l^{j})^{'}\zeta_{l}^{1} = 0 \\
\end{multline} \vspace*{-5ex}

Rearranging the terms in the above equation yields

\vspace*{-5ex} \begin{multline}
\nonumber -(B_l^{1} + \sum_{j\in 2\ldots n}B_l^j\bar {r}_l^{j})^{'}Z_l^1(B_l^1 + \sum_{j\in 2\ldots n}B_l^j\bar {r}_l^{j})u_l^{1*} - (R_l^{11} + \sum_{j\in 2\ldots n}\bar {r}_l^{j}R_l^{1j}\bar {r}_l^{j}) u_l^{1*} \\= (B_l^{1} + \sum_{j\in 2\ldots n}B_l^j\bar {r}_l^{j})^{'}Z_l^1(A_lx_{l-1} + \sum_{j\in 2\ldots n}B_l^j(W_l^{j}x_{l-1}+w_l^j) + s_l) + \sum_{j\in 2\ldots n}\bar {r}_l^{j'}R_l^{1j}(W_l^{j}x_{l-1}+w_l^j) \\- (B_l^{1} + \bar {r}_l^{j})^{'}Q_l^1\tilde{x}_l^1 - R_l^{11}\tilde{u}_l^{11} - \sum_{j\in 2\ldots n}\bar {r}_l^{j'}R_l^{1j}\tilde{u}_l^{1j} + (B_l^{1} + \sum_{j\in 2\ldots n}B_l^j\bar {r}_l^{j})^{'}\zeta_{l}^{1} \\
\end{multline} \vspace*{-5ex}

The structure of the above equation justifies the following substitution

\vspace*{-5ex} \begin{multline}
\eqref{gamma}_{k=l}^{i=1} \;\; u_{l}^{1*} = -P_l^{1*}x_{l-1}-\alpha_l^{1*} \\\label{u1subns}
\end{multline} \vspace*{-5ex}

\vspace*{-5ex} \begin{multline}
(B_l^{1} + \sum_{j\in 2\ldots n}B_l^j\bar {r}_l^{j})^{'}Z_l^1(B_l^1 + \sum_{j\in 2\ldots n}B_l^j\bar {r}_l^{j})(P_k^{1*}x_{k-1} + \alpha_k^{1*}) \\+ (R_l^{11} + \sum_{j\in 2\ldots n}\bar {r}_l^{j}R_l^{1j}\bar {r}_l^{j})(P_k^{1*}x_{k-1} + \alpha_k^{1*})   = (B_l^{1} + \sum_{j\in 2\ldots n}B_l^j\bar {r}_l^{j})^{'}\\Z_l^1(A_lx_{l-1} + \sum_{j\in 2\ldots n}B_l^j(W_l^{j}x_{l-1}+w_l^j) + s_l) + \sum_{j\in 2\ldots n}\bar {r}_l^{j'}R_l^{1j}(W_l^{j}x_{l-1}+w_l^j) \\- (B_l^{1} + \bar {r}_l^{j})^{'}Q_l^1\tilde{x}_l^1 - R_l^{11}\tilde{u}_l^{11} - \sum_{j\in 2\ldots n}\bar {r}_l^{j'}R_l^{1j}\tilde{u}_l^{1j} + (B_l^{1} + \sum_{j\in 2\ldots n}B_l^j\bar {r}_l^{j})^{'}\zeta_{l}^{1} \\\label{gamma1l}
\end{multline} \vspace*{-5ex}

Comparing coefficients gives

\vspace*{-5ex} \begin{multline}
\nonumber\eqref{gamma1l}_{x_{l}} \;\;\ [(B_l^{1} + \sum_{j\in 2\ldots n}B_l^j\bar {r}_l^{j})^{'}Z_l^1(B_l^1 + \sum_{j\in 2\ldots n}B_l^j\bar {r}_l^{j}) + R_l^{11} + \sum_{j\in 2\ldots n}\bar {r}_l^{j}R_l^{1j}\bar {r}_l^{j}]P_k^{1*}  =\\ (B_l^{1} + \sum_{j\in 2\ldots n}B_l^j\bar {r}_l^{j})^{'}Z_l^1(A_l + \sum_{j\in 2\ldots n}B_l^jW_l^{j}) + \sum_{j\in 2\ldots n}\bar {r}_l^{j'}R_l^{1j}W_l^{j}\\
\end{multline} \vspace*{-5ex}

\vspace*{-5ex} \begin{multline}
\nonumber\eqref{gamma1l}_{const.}\;\;[(B_l^{1} + \sum_{j\in 2\ldots n}B_l^j\bar {r}_l^{j})^{'}Z_l^1(B_l^1 + \sum_{j\in 2\ldots n}B_l^j\bar {r}_l^{j}) + R_l^{11} + \sum_{j\in 2\ldots n}\bar {r}_l^{j}R_l^{1j}\bar {r}_l^{j}]\alpha_k^{1*}   =\\ (B_l^{1} + \bar {r}_l^{j})^{'}Z_l^1(\sum_{j\in 2\ldots n}B_l^jw_l^{j} + s_l) + \sum_{j\in 2\ldots n}\bar {r}_l^{j'}R_l^{1j}w_l^{j} \\- (B_l^{1} + \bar {r}_l^{j})^{'}Q_l^1\tilde{x}_l^1 - R_l^{11}\tilde{u}_l^{11} - \sum_{j\in 2\ldots n}\bar {r}_l^{j'}R_l^{1j}\tilde{u}_l^{1j} + (B_l^{1} + \sum_{j\in 2\ldots n}B_l^j\bar {r}_l^{j})^{'}\zeta_{l}^{1} \\
\end{multline} \vspace*{-5ex}

Making $P_l^{1*}$ and $\alpha_l^{1*}$ explicit finally leads to

\vspace*{-5ex} \begin{multline}
\nonumber\eqref{P1}_{k = l} \;\;P_l^{1*} = [(B_l^1 + \sum_{j\in 2\ldots n}B_l^j\bar {r}_l^{j})^{'}Z_l^1(B_l^1 + \sum_{j\in 2\ldots n}B_l^j\bar {r}_l^{j}) + R_l^{11} + \sum_{j\in 2\ldots n}\bar {P}_l^{j'}R_l^{1j}\bar {r}_l^{j}]^{-1} \\ [(B_l^1 + \sum_{j\in 2\ldots n}B_l^j\bar {r}_l^{j})^{'}Z_l^1(A_l - \sum_{j\in 2\ldots n}B_l^jW_l^{j}) + \sum_{j\in 2\ldots n}\bar {r}_l^{j'}R_l^{1j} W_l^{j}]\\
\end{multline} \vspace*{-5ex}

\vspace*{-5ex} \begin{multline}
\nonumber\eqref{alpha1}_{k = l} \;\; \alpha_l^{1*} = [(B_l^1 + \sum_{j\in 2\ldots n}B_l^j\bar {r}_l^{j})^{'}Z_l^1(B_l^1 + \sum_{j\in 2\ldots n}B_l^j\bar {r}_l^{j}) + R_l^{11} + \sum_{j\in 2\ldots n}\bar {r}_l^{j'}R_l^{1j}\bar {r}_l^{j}]^{-1} \\  [(B_l^1 + \sum_{j\in 2\ldots n}B_l^j\bar {r}_l^{j})^{'}Z_l^1(\sum_{j\in 2\ldots n}B_l^jw_l^{j} + s_l) + \sum_{j\in 2\ldots n}\bar {r}_l^{j'}R_l^{1j} w_l^{j} \\- R_l^{11}\tilde{u}_l^{11} - \sum_{j\in 2\ldots n}\bar {r}_l^{j'}R_l^{1j} \tilde{u}_l^{1j} + (B_l^{1} + (\sum_{j\in 2\ldots n}B_l^j\bar {r}_l^{j})^{'}(\zeta_{l}^{1} - Q_l^1\tilde{x}_l^1)]\\
\end{multline} \vspace*{-5ex}

By making use of the unique optimal strategy $u_l^{1*}$ of the leader in the optimal response functions of the followers given by \eqref{rli} and, in a different presentation, by \eqref{rlia} we get the optimal strategies of the followers

\vspace*{-5ex} \begin{multline}
\nonumber r_l^{i*} = u_l^{i*} = -(R_l^{ii} + B_l^{i'}Z_l^iB_l^i)^{-1} [ B_l^{i'}(Z_l^i(B_l^1u_l^{1*} + \sum_{j \in 2 \ldots n \; , \; j \not = i}B_l^ju_l^{j*} \\+ A_lx_{l-1} + s_l) +\zeta_{l}^i  - Q_l^i\tilde{x}_l^i) - R_l^{ii}\tilde{u}_l^{ii}] \; ; \; i \in \{2\ldots n \}\\
\end{multline} \vspace*{-5ex}

\vspace*{-5ex} \begin{multline}
\nonumber r_l^{i*} = u_l^{i*} = W_l^{i}x_{l-1} + w_l^i + \bar {r}_l^{i}u_l^{1*} \; ; \; i \in \{2\ldots n \} \\
\end{multline} \vspace*{-5ex}

Making use of \eqref{u1subns} in the above two systems of equations yields

\vspace*{-5ex} \begin{multline}
\nonumber u_l^{i*} = -(R_l^{ii} + B_l^{i'}Z_l^iB_l^i)^{-1} [ B_l^{i'}(Z_l^i(B_l^1(P^{1*}_{l}x_{l-1} + \alpha_{l}^{1*}) + \sum_{j \in 2 \ldots n \; , \; j \not = i}B_l^ju_l^{j*} \\+ A_lx_{l-1} + s_l) +\zeta_{l}^i  - Q_l^i\tilde{x}_l^i) - R_l^{ii}\tilde{u}_l^{ii}] \; ; \; i \in \{2\ldots n \}\\
\end{multline} \vspace*{-5ex}

\vspace*{-5ex} \begin{multline}
\nonumber u_l^{i*} = W_l^{i} + w_l^i - \bar {r}_l^{i}(P^{1*}_{l}x_{l-1} + \alpha_{l}^{1*}) \; ; \; i \in \{2\ldots n \} \\
\end{multline} \vspace*{-5ex}

The structure of the above two systems of equations justifiies the following substitution

\vspace*{-5ex} \begin{multline}
\eqref{gamma}_{k=l} \;\; u_{l}^{i*} = -P_k^{i*}x_{k-1}-\alpha_k^{i*} \; ; \; i \in \{2\ldots n \}\\\label{uisubns}
\end{multline} \vspace*{-5ex}

\vspace*{-5ex} \begin{multline}
P_k^{i*}x_{k-1} + \alpha_k^{i*} = (R_l^{ii} + B_l^{i'}Z_l^iB_l^i)^{-1} [ B_l^{i'}(Z_l^i(B_l^1(P^{1*}_{l}x_{l-1} + \alpha_{l}^{1*}) \\+ \sum_{j \in 2 \ldots n \; , \; j \not = i}B_l^j(-P_k^{i*}x_{k-1}-\alpha_k^{i*}) + A_lx_{l-1} + s_l) +\zeta_{l}^i  - Q_l^i\tilde{x}_l^i) - R_l^{ii}\tilde{u}_l^{ii}] \; ; \; i \in \{2\ldots n \}\\\label{gammail}
\end{multline} \vspace*{-5ex}

\vspace*{-5ex} \begin{multline}
P_k^{i*}x_{k-1} + \alpha_k^{i*} = -W_l^{i} - w_l^i + \bar {r}_l^{i}(P^{1*}_{l}x_{l-1} + \alpha_{l}^{1*}) \; ; \; i \in \{2\ldots n \} \\\label{folad}
\end{multline} \vspace*{-5ex}

By comparing coefficients it follows that

\vspace*{-5ex} \begin{multline}
\nonumber\eqref{gammail}_{x_{l}}=\eqref{folad}_{x_{l}}=\eqref{Pf}_{k=l} \; \; P_l^{i*} = (R_l^{ii} + B_l^{i'}Z_l^iB_l^i)^{-1}[B_l^iZ_l^i(A_l - B_l^1P_l^{1*} \\- \sum_{j\in 2\ldots n \;,\; j \not = i}B_l^jP_l^{j*})] = - W_l^{i} + \bar {r}_l^{i}P^{1*}_{l}\; ; \; i \in \{2\ldots n\}\\
\end{multline} \vspace*{-5ex}

\vspace*{-5ex} \begin{multline}
\nonumber \eqref{gammail}_{const.}=\eqref{folad}_{const.}=\eqref{alphaf}_{k=l} \; \; \alpha_l^{i*} = (R_l^{ii} + B_l^{i'}Z_l^iB_l^i)^{-1} [B_l^{i'}(Z_l^is_l - B_l^1\alpha_l^{1*}\\ - \sum_{j\in 2\ldots n \;,\; j \not = i}B_l^j\alpha_l^{j*} + \zeta_{l}^{i} - Q_l^i\tilde{x}_l^i)-R_l^{ii}\tilde{u}_l^{ii}] =  -w_l^i + \bar {r}_l^{i}\alpha_{l}^{1*} \; ; \; i \in \{2\ldots n\}\\
\end{multline} \vspace*{-5ex}

Now, after finding the optimal strategies for the players, we are able to rewrite \eqref{Vgil} and \eqref{Vgif} as

\vspace*{-5ex} \begin{multline}
\nonumber V^i(l,x_{l-1}) = \frac{1}{2} [(A_lx_{l-1} + \sum_{j\in N}B_l^ju_l^{j*}+s_l)^{'}Z_l^i(A_lx_{l-1} + \sum_{j\in N}B_l^ju_l^{j*}+s_l) \\+ \sum_{j\in N}u_l^{j*'}R_l^{ij}u_l^{j*}] +
 \frac{1}{2} (\tilde{x}_l^{i'}Q_l^i\tilde{x}_l^i + \sum_{j\in N}\tilde{u}_l^{ij'}R_l^{ij}\tilde{u}_l^{ij}) - \tilde{x}_l^{i'}
Q_l^i(A_lx_{l-1} + \sum_{j\in N}B_l^ju_l^{j*}+s_l) \\- \sum_{j\in N}\tilde{u}_l^{ij'}R_l^{ij}u_l^{j*} + \zeta_{l}^{i'}(A_lx_{l-1} + \sum_{j\in N}B_l^ju_l^{j*}+s_l) + n_{l} \\
\end{multline} \vspace*{-5ex}

Making use of \eqref{uisubns} and \eqref{u1subns} respectively yields

\vspace*{-5ex} \begin{multline}
\nonumber V^i(l,x_{l-1}) = \frac{1}{2} [(A_lx_{l-1} + \sum_{j\in N}B_l^j(-P_l^{j*}x_{l-1}-\alpha_l^{j*}) +s_l)^{'}Z_l^i(A_lx_{l-1} \\+ \sum_{j\in N}B_l^j(-P_l^{j*}x_{l-1}-\alpha_l^{j*})+s_l) + \sum_{j\in N}(-P_l^{j*}x_{l-1}-\alpha_l^{j*}) ^{'}R_l^{ij}(-P_l^{j*}x_{l-1}-\alpha_l^{j*}) ]\\ +
 \frac{1}{2} (\tilde{x}_l^{i'}Q_l^i\tilde{x}_l^i + \sum_{j\in N}\tilde{u}_l^{ij'}R_l^{ij}\tilde{u}_l^{ij}) - \tilde{x}_l^{i'}Q_l^i(A_lx_{l-1} + \sum_{j\in N}B_l^j(-P_l^{j*}x_{l-1}-\alpha_l^{j*})+s_l) \\- \sum_{j\in N}\tilde{u}_l^{ij'}R_l^{ij}(-P_l^{j*}x_{l-1}-\alpha_l^{j*}) + \zeta_{l}^{i'}(A_lx_{l-1} + \sum_{j\in N}B_l^j(-P_l^{j*}x_{l-1}-\alpha_l^{j*})+s_l) + n_{l}  \\
\end{multline} \vspace*{-5ex}

Rewriting the above equation to the power of $x_{l-1}$ gives

\vspace*{-5ex} \begin{multline}
\nonumber V^i(l,x_{l-1}) = \frac{1}{2} x_{l-1}^{'}[(A_l - \sum_{j \in N}B_l^jP_l^{j*})^{'}Q_l^i(A_l - \sum_{j \in N}B_l^jP_l^{j*}) + \sum_{j\in N}P_l^{j*'}R_l^{ij}P_l^{j*}]x_{l-1} \\ + [(A_l - \sum_{j \in N}B_l^jP_l^{j*})^{'}[\zeta_{l}^{i}Z_l^i(s_l-\sum_{j\in N}B_l^j\alpha_l^{j*}) - Q_l^i\tilde{x}_k^i] + \sum_{j\in N}P_l^{j*'}R_l^{ij}(\alpha_l^{j*} + \tilde{u}_l^{ij}) ]^{'}x_{l-1} \\+ \frac{1}{2} (s_l-\sum_{j\in N}B_l^j\alpha_l^{j*})^{'}Z_l^i(s_l-\sum_{j\in N}B_l^j\alpha_l^{j*})  + \frac{1}{2}\sum_{j\in N }\alpha_l^{j*'}R_l^{ij}\alpha_l^{j*} \\-  \tilde{x}_l^{i'}Q_l^i(s_l-\sum_{j\in N}B_l^j\alpha_l^{j*}) + \sum_{j\in N}\tilde{u}_l^{ij'}R_l^{ij}\alpha_l^{j*} + \frac{1}{2}(\tilde{x}_l^{i'}Q_l^i\tilde{x}_l^i + \sum_{j\in N}\tilde{u}_l^{ij'}R_l^{ij}\tilde{u}_l^{ij}) \\+ \zeta_{l}^{i'}(s_l - \sum_{j\in N}B_l^j\alpha_l^{j*}) + n_{l} \\
\end{multline} \vspace*{-5ex}

To finish off the inductive step and consequently the induction argument, we use the recursive equations $\eqref{Zfs}_{k=l}$, $\eqref{zetaki}_{k=l}$ and $\eqref{n}_{k=l}$ in the above equation. This leads to

\vspace*{-5ex} \begin{multline}
\nonumber\eqref{Vg}_{k=l} \;\;V^i(l,x_{l-1}) = \frac{1}{2} x_{l-1}^{'}(Z_{l-1}^i-Q_{l-1}^i)x_{l-1} + \zeta_{l-1}^{i'}x_{l-1} + n_{l-1}^i \\
\end{multline} \vspace*{-5ex}

The expression for the total costs of the game for player \textit {i} given by \eqref{costn} is equal to the function the induction argument was based on at stage 1. In other words, \eqref{costn} is equal to $\eqref{Vg}_{k=1}$.\;\;\;\;\;\;\;\;\boxed{}

\begin{rem}
The proof of Theorem \eqref{FBSTheo} is a formalization and generalization of the heuristic argumentation presented in Ba\c{s}ar and Olsder (1999, pp. 274-275)\cite{baol}.
\end{rem}

\begin{rem}
$Z_k^i$ given by \eqref{Zfs} is positive definite for all $k \in \{0, \ldots, T\}$. This can be proven by a straightforward induction argument, starting at stage T and using Corollary \eqref{pdc} in the inductive step.
\end{rem}

\begin{rem}
To solve the Stackelberg game algorithmically, the following order of application of the equations of Theorem \eqref{FBSTheo} is advisable ($i \in N \; , \; j \in \{2, \ldots, n\}$):
\begin{enumerate}
\item For k running backward from T to 1
\begin{enumerate}
\item $Z_k^i$, $\zeta_k^i$ and $n_k^i$
\item $\bar{r}_k^j$, $W_k^j$ and $w_k^j$
\item $P_{k}^{1*}$, $\alpha_{k}^{1*}$ 
\item $P_{k}^{j*}$, $\alpha_{k}^{j*}$ 
\end{enumerate}
\item $Z_0^i$, $\zeta_0^i$ and $n_0^i$
\item $V^i(1,x_0)$
\item For k running forward from 1 to T\\
$\gamma_k^{i*}(x_{k-1})$\\\\
\end{enumerate}
\end{rem}

\begin{pro}\footnote{In this proposition we rewrite the equilibrium equations in a notation that was used at our department in the past to enable comparison.}
The systems of equations defining the unique equilibrium strategies $\gamma_{k}^{i*}(x_{k-1})$ in Theorem \eqref{FBSTheo} can also be written in the following way:\footnote{For all equations belonging to this proposition and its proof, $k \in K$ and $i \in N$ if nothing different is stated.}

\vspace*{-5ex} \begin{align}
\gamma_k^{i*}(x_{k-1}) = G_k^{i*}x_{k-1} + g_k^{i*} \label{gammaä}
\end{align}

\vspace*{-5ex} \begin{align}
G_k^{1*} = -[\Lambda]^{-1} [ \bar{B}_k^{'}H_k^1A_k + \sum_{j\in 2\ldots n}\bar{D}_k^jW_k^{j}] \label{P1ä}
\end{align}

\vspace*{-5ex} \begin{align}
g_k^{1*} = -[\bar {\Lambda}_k]^{-1} [v_k^1 + \bar{v}_k + \sum_{ j \in 2\ldots n}\bar{D}_k^jw_k^j] \label{alpha1ä}
\end{align}

\vspace*{-5ex} \begin{align}
G_k^{i*} = W_k^i + \Psi_k^iG_k^{1*} \; ; \; i \in \{2\ldots n\} \label{Pfä} 
\end{align}

\vspace*{-5ex} \begin{align}
g_k^{i*} = w_k^i + \Psi_k^ig_k^{1*} \; ; \; i \in \{2\ldots n\} \label{alphafä}
\end{align}

\vspace*{-5ex} \begin{align}
\bar \Lambda_k = \bar{B}_k^{'}H_k^1\bar {B}_k + R_k^{11} + \sum_{j\in 2\ldots n}\Psi_k^{j'}R_k^{1j}\Psi_k^{j} \label{lamque}
\end{align}

\vspace*{-5ex} \begin{align}
\bar {B}_k = B_k^1 + \sum_{j\in 2\ldots n}B_k^j\Psi_k^{j} \label{Bqs}
\end{align}

\vspace*{-5ex} \begin{align}
\Psi_k^{i} = -(D_k^i)^{-1}[B_k^{i'}Z_k^i(B_k^1 + \sum_{j\in 2\ldots n \;,\; j \not = i}B_k^j\Psi_k^{j})] \; ; \; i \in \{2\ldots n\} \label{rbarä}
\end{align}

\vspace*{-5ex} \begin{align}
H_{k-1}^i = K_k^{'}Z_k^iK_k + \sum_{j\in N}G_k^{j'}R_k^{ij}G_k^j + Q_{k-1}^i\; ;\; H_T^i = Q_T^i \label{Hksä}
\end{align}

\vspace*{-5ex} \begin{align}
h_{k-1}^i = Q_{k-1}^i\tilde{x}_{k-1}^i - K_k^{'}[H_k^ik_k - h_{k}^i] + \sum_{j\in N}G_k^{j'}R_k^{ij}(\tilde{u}_k^{ij} - g_k^j)\; ;\; h_{T}^i = Q_T^i\tilde{x}_T^i  \label{zetasä}
\end{align}

\vspace*{-5ex} \begin{align}
K_k = A_k + \sum_{j \in N}B_k^iG_k^j \label{Ksubs}
\end{align}

\vspace*{-5ex} \begin{align}
k_k = s_k + \sum_{j\in N}B_k^jg_k^j \label{ksubs}
\end{align}

\vspace*{-5ex} \begin{align}
D_k^i = R_k^{ii} + B_k^{i'}H_k^iB_k^i \label{Dsubs}
\end{align}

\vspace*{-5ex} \begin{align}
\bar{D}_k = \bar{B}_k^{'}H_k^1B_k^i + \Psi_k^{i'}R_k^{1i} \label{Dqi}
\end{align}

\vspace*{-5ex} \begin{align}
W_k^{i} = -(D_k^i)^{-1}[B_k^{i'}H_k^i(\sum_{j\in 2\ldots n\;,\; j \not = i}B_k^jW_k^{j} + A_k)] \; ; \; i \in \{2\ldots n\} \label{Wä}
\end{align}

\vspace*{-5ex} \begin{align}
w_k^{i} = -(D_k^i)^{-1}[B_k^{i'}H_k^i\sum_{j\in 2\ldots n\;,\; j \not = i}B_k^jw_k^{j} + v_k^i] \; ; \; i \in \{2\ldots n\} \label{wä}
\end{align}

\vspace*{-5ex} \begin{align}
v_k^i = B_k^{i'}(H_k^is_k - h_{k}^i) - R_k^{ii}\tilde{u}_k^{ii} \; ; \; i \in N \label{vsubs}
\end{align}

\vspace*{-5ex} \begin{align}
\bar{v}_k = \sum_{j\in 2\ldots n}\Psi_k^{j'}(B_k^{j'}H_k^1s_k - R_k^{1j}\tilde{u}_k^{1j} - B_k^{j'}h_{k}^1) \label{vqs}
\end{align}\\\\

\label{profbsnb}
\end{pro}

\textsc {Proof:}

The proof is carried out by renaming some matrices and then showing that the relations for the equilibrium strategies $\gamma_k^{i*}(x_{k-1})$ of Theorem \eqref{FBSTheo} can be rewritten in the way stated above.\\

Let us start by renaming the feedback matrices $P_{k}^{i*}$ and $\alpha_k^{i*}$, the reaction coefficients $\bar{r}_k^j$ (j $\in \{2, \ldots, n\}$) and the matrices $Z_k^i$ and $\zeta_{k}^i$.

\vspace*{-5ex} \begin{multline}
P_k^i \; \widehat {=} \; -G_k^i \; ; \; \alpha_k^i \; \widehat {=} \; -g_k^i \; ; \;  \bar{r}_k^j \; \widehat {=} \; \Psi_k^j \; ; \; Z_k^i \; \widehat {=} \; H_k^i \; ; \;\\ \zeta_{k}^i \; \widehat {=} \; - h_k^i + Q_k^i\tilde{x}_k^i \;;\; j \in \{2, \ldots, n\}\\\label{äqusubsfb}
\end{multline} \vspace*{-5ex}

Next we prove that $Z_k^i$ and $\zeta_{k}^i$ fulfill \eqref{Hksä} and \eqref{zetasä} respectively.\\

Taking consideration of the renaming \eqref{Zfs} gives

\vspace*{-5ex} \begin{multline}
\nonumber H_{k-1}^i = (A_k + \sum_{j \in N}B_k^iG_k^j)^{'}H_k^i(A_k + \sum_{j \in N}B_k^jG_k^j) + \sum_{j\in N}G_k^{j'}R_k^{ij}G_k^j + Q_{k-1}^i \\
\end{multline} \vspace*{-5ex}

Making use of \eqref{Ksubs} yields

\vspace*{-5ex} \begin{multline}
\nonumber\eqref{Hksä} \; \; H_{k-1}^i = K_k^{'}H_k^iK_k + \sum_{j\in N}G_k^{j'}R_k^{ij}G_k^j + Q_{k-1}^i \\
\end{multline} \vspace*{-5ex}

Now we demonstrate the correctness of equation \eqref{zetasä}. To do so we start at stage \textit{T} and use \eqref{zetaki} and \eqref{äqusubsfb} to get

\vspace*{-5ex} \begin{multline}
\nonumber 0 = \zeta_{T}^i = - h_T^i + Q_T^i\tilde{x}_T^i \\
\end{multline} \vspace*{-5ex}

\vspace*{-5ex} \begin{multline}
\nonumber\eqref{zetasä}_{k=T} \; \; h_T^i = Q_T^i\tilde{x}_T^i \\
\end{multline} \vspace*{-5ex}

For the general stage \textit{k} rewrite \eqref{zetaki} taking consideration of \eqref{äqusubsfb}

\vspace*{-5ex} \begin{multline}
\nonumber -h_{k-1}^i + Q_{k-1}^i\tilde{x}_{k-1}^i = (A_k + \sum_{j \in N}B_k^jG_k^j)^{'}[-h_{k}^i + H_k^i(s_k + \sum_{j\in N}B_k^jg_k^j)] - \sum_{j\in N}G_k^{j'}R_k^{ij}(\tilde{u}_k^{ij} - g_k^j) \\
\end{multline} \vspace*{-5ex}

Using \eqref{Ksubs} and \eqref{ksubs} and making $h_{k-1}^i$ explicit yields

\vspace*{-5ex} \begin{multline}
\nonumber\eqref{zetasä} \; \; h_{k-1}^i  = Q_{k-1}^i\tilde{x}_{k-1}^i - K_k^{'}[H_k^ik_k - h_{k}^i] + \sum_{j\in N}G_k^{j'}R_k^{ij}(\tilde{u}_k^{ij} - g_k^j) \\
\end{multline} \vspace*{-5ex}

Furthermore the relations for $W_k^{i}$, $w_k^{i}$ and $\Psi_k^i$ (i $\in \{2, \ldots, n\}$), given by \eqref{Wä}, \eqref{wä} and \eqref{rbarä}, have to be deduced.\\

Making use of \eqref{äqusubsfb} and \eqref{Dsubs} in \eqref{W} leads to

\vspace*{-5ex} \begin{multline}
\nonumber W_k^{i} = -(R_k^{ii} + B_k^{i'}H_k^iB_k^i)^{-1}[B_k^{i'}H_k^i(\sum_{j\in 2\ldots n\;,\; j \not = i}B_k^jW_k^{j} + A_k)] \; ; \; i \in \{2\ldots n\} \\
\end{multline} \vspace*{-5ex}

\vspace*{-5ex} \begin{multline}
\nonumber\eqref{Wä} \; \; W_k^{i} = -(D_k^i)^{-1}[B_k^{i'}H_k^i(\sum_{j\in 2\ldots n\;,\; j \not = i}B_k^jW_k^{j} + A_k)] \; ; \; i \in \{2\ldots n\} \\
\end{multline} \vspace*{-5ex}

Using \eqref{äqusubsfb} and \eqref{Dsubs} in \eqref{w} yields

\vspace*{-5ex} \begin{multline}
\nonumber w_k^{i} = -(R_k^{ii} + B_k^{i'}H_k^iB_k^i)^{-1}[B_k^{i'}(H_k^i(\sum_{j\in 2\ldots n\;,\; j \not = i}B_k^jw_k^{j} + s_k) + h_k^i)-R_k^{ii}\tilde{u}_k^{ii}] \; ; \; i \in \{2\ldots n\}\\
\end{multline} \vspace*{-5ex}

\vspace*{-5ex} \begin{multline}
\nonumber\eqref{wä} \; \; w_k^{i} = -(D_k^i)^{-1}[B_k^{i'}H_k^i\sum_{j\in 2\ldots n\;,\; j \not = i}B_k^jw_k^{j} + v_k^i] \; ; \; i \in \{2\ldots n\}\\
\end{multline} \vspace*{-5ex}

Considering \eqref{äqusubsfb} and \eqref{Dsubs}, \eqref{rbar} gives

\vspace*{-5ex} \begin{multline}
\nonumber \Psi_k^{i} = -(R_k^{ii} + B_k^{i'}H_k^iB_k^i)^{-1}[B_k^{i'}H_k^i(B_k^1 + \sum_{j\in 2\ldots n \;,\; j \not = i}B_k^j\Psi_k^{j})] \; ; \; i \in \{2\ldots n\}\\
\end{multline} \vspace*{-5ex}

\vspace*{-5ex} \begin{multline}
\nonumber\eqref{rbarä} \; \; \Psi_k^{i} = -(D_k^i)^{-1}[B_k^{i'}H_k^i(B_k^1 + \sum_{j\in 2\ldots n \;,\; j \not = i}B_k^j\Psi_k^{j})] \; ; \; i \in \{2\ldots n\}\\
\end{multline} \vspace*{-5ex}

Eventually the correctness of the rewritten equilibrium strategies $\gamma_k^{i*}(x_{k-1})$ given by \eqref{gammaä} - \eqref{alphafä} has to be shown.\\

First substitute $P_{k}^{i*}$ and $\alpha_{k}^{i*}$ in \eqref{gamma} with the help of \eqref{äqusubsfb}

\vspace*{-5ex} \begin{multline}
\nonumber\eqref{gammaä} \; \; \gamma_k^{i*}(x_{k-1}) = G_k^{i*}x_{k-1} + g_k^{i*} \\
\end{multline} \vspace*{-5ex}

Next we deduce the feedback matrices for the leader. To do so we start by rewriting \eqref{P1} taking consideration of \eqref{äqusubsfb}

\vspace*{-5ex} \begin{multline}
\nonumber G_k^{1*} = -[(B_k^1 + \sum_{j\in 2\ldots n}B_k^j\Psi_k^{j})^{'}H_k^1(B_k^1 + \sum_{j\in 2\ldots n}B_k^j\Psi_k^{j}) + R_k^{11} + \sum_{j\in 2\ldots n}\Psi_k^{j'}R_k^{1j}\Psi_k^{j}]^{-1}\\ [ (B_k^1 + \sum_{j\in 2\ldots n}B_k^j\Psi_k^{j})^{'}H_k^1(A_k + \sum_{j\in 2\ldots n}B_k^jW_k^{j}) + \sum_{j\in 2\ldots n}\Psi_k^{j'}R_k^{1j}W_k^{j}]\\
\end{multline} \vspace*{-5ex}

Using \eqref{lamque} and \eqref{Bqs} yields

\vspace*{-5ex} \begin{multline}
\nonumber G_k^{1*} = -[\Lambda]^{-1} [ \bar{B}_k^{'}H_k^1A_k + \bar{B}_k^{'}H_k^1\sum_{j\in 2\ldots n}B_k^jW_k^{j} + \sum_{j\in 2\ldots n}\Psi_k^{j'}R_k^{1j}W_k^{j}]\\
\end{multline} \vspace*{-5ex}

Making use of \eqref{Dqi} gives

\vspace*{-5ex} \begin{multline}
\nonumber\eqref{P1ä} \; \; G_k^{1*} = -[\Lambda]^{-1} [ \bar{B}_k^{'}H_k^1A_k + \sum_{j\in 2\ldots n}\bar{D}_k^jW_k^{j}] \\
\end{multline} \vspace*{-5ex}

The constant part of the equilibrium strategy of the leader can be rewritten (considering \eqref{äqusubsfb}) in the following way

\vspace*{-5ex} \begin{multline}
\nonumber g_k^{1*} = -[(B_k^1 + \sum_{j\in 2\ldots n}B_k^j\Psi_k^{j})^{'}H_k^1(B_k^1 + \sum_{j\in 2\ldots n}B_k^j\Psi_k^{j}) + R_k^{11} + \sum_{j\in 2\ldots n}\Psi_k^{j'}R_k^{1j}\Psi_k^{j}]^{-1}\\ [ (B_k^1 + \sum_{j\in 2\ldots n}B_k^j\Psi_k^{j})^{'}H_k^1(s_k + \sum_{j\in 2\ldots n}B_k^jw_k^{j}) + \sum_{j\in 2\ldots n}\Psi_k^{j'}R_k^{1j}w_k^{j}  \\- R_k^{11}\tilde{u}_k^{11} - \sum_{j\in 2\ldots n}\Psi_k^{j'}R_k^{1j}\tilde{u}_k^{1j} - (B_k^1 + \sum_{j\in 2\ldots n}B_k^j\Psi_k^{j})^{'}h_{k}^1]\\
\end{multline} \vspace*{-5ex}

Using \eqref{lamque} and \eqref{Bqs} yields

\vspace*{-5ex} \begin{multline}
\nonumber g_k^{1*} = -[\Lambda]^{-1} [(B_k^{1'} + \sum_{j\in 2\ldots n}\Psi_k^{j'}B_k^{j'})H_k^1s_k + \bar{B}_k^{'}H_k^1\sum_{j\in 2\ldots n}B_k^jw_k^{j} + \sum_{j\in 2\ldots n}\Psi_k^{j'}R_k^{1j}w_k^{j}  \\- R_k^{11}\tilde{u}_k^{11} - \sum_{j\in 2\ldots n}\Psi_k^{j'}R_k^{1j}\tilde{u}_k^{1j} - B_k^{1'}h_{k}^1 - \sum_{j\in 2\ldots n}\Psi_k^{j'}B_k^{j'}h_{k}^1]\\
\end{multline} \vspace*{-5ex}

Making use of \eqref{Dqi}, $\eqref{vsubs}_{i=1}$ and \eqref{vqs} gives

\vspace*{-5ex} \begin{multline}
\nonumber\eqref{alpha1ä} \; \;  g_k^{1*} = -[\Lambda]^{-1} [\sum_{j\in 2\ldots n}\bar{D}_k^jw_k^{j} + \bar{v}_k + v_k^1]\\
\end{multline} \vspace*{-5ex}

Finally considering \eqref{äqusubsfb}, the feedback matrices of the followers given by \eqref{Pf} and \eqref{alphaf} can be rewritten as

\vspace*{-5ex} \begin{multline}
\nonumber\eqref{Pfä} \; \; G_k^{i*} = W_k^i + \Psi_k^iG_k^{1*} \; ; \; i \in \{2\ldots n\} \\ 
\end{multline} \vspace*{-5ex}

\vspace*{-5ex} \begin{multline}
\nonumber\eqref{alphafä} \; \; g_k^{i*} = w_k^i + \Psi_k^ig_k^{1*} \; ; \; i \in \{2\ldots n\}  \;\;\;\;\;\;\;\;\;\;\;\;\;\;\;\;\boxed{}\\
\end{multline} \vspace*{-5ex}

\label{ssaqfsg}

\clearpage
\subsection{Special case: Affine-quadratic games with one leader and one follower}

In this subsection, the results of Theorem \eqref{FBSTheo} are specialized by reducing the number of followers from \textit{n} to one.

\begin{cor}
A 2-person affine-quadratic dynamic game (cf. Def. \eqref{defspielaq}) admits a \emph{unique feedback Stackelberg equilibrium solution} if 
\begin{enumerate}
\item $Q_k^i \geq 0$, $R_k^{ii} > 0$ and $R_k^{ij} \geq 0$ (defined for $k \in K$ , $i,j \in \{1,2\}$, $j \not = i$)
\end{enumerate}
If these conditions are satisfied, the unique equilibrium strategies $\gamma_{k}^{i*}(x_{k-1})$ (i $\in \{1,2\}$) are given by \eqref{gamma2} and the corresponding feedback Stackelberg equilibrium costs for the two players are stated in \eqref{cost2}. \footnote{For all equations belonging to this theorem and its proof, $k \in K$ if nothing different is stated.}

\vspace*{-5ex} \begin{align}
 f_{k-1}(x_{k-1},u_k^1,u_k^2) = A_kx_{k-1} + \sum_{j\in \{1,2\}}B_k^ju_k^j+s_k \label{f2}
\end{align}

\vspace*{-5ex} \begin{align}
L^i(u^{1}, u^{2})= \sum_{k = 1}^T g_k^i(x_k,u_k^{1}, u_k^{2},x_{k-1}) 
\end{align}

\vspace*{-5ex} \begin{multline}
 g_k^i(x_k,u_k^1,\ldots,u_k^n,x_{k-1}) = \frac{1}{2} (x_k^{'}Q_k^ix_k + \sum_{j\in \{1,2\}}u_k^{j'}R_k^{ij}u_k^j)\; + \\
 \frac{1}{2} (\tilde{x}_k^{i'}Q_k^i\tilde{x}_k^i + \sum_{j\in \{1,2\}}\tilde{u}_k^{ij'}R_k^{ij}\tilde{u}_k^{ij}) - \tilde{x}_k^{i'}Q_k^ix_k - \sum_{j\in \{1,2\}}\tilde{u}_k^{ij'}R_k^{ij}u_k^j \\\label{g2}
\end{multline} \vspace*{-5ex}

\vspace*{-5ex} \begin{align}
\gamma_k^{i*}(x_{k-1}) = -P_k^{i*}x_{k-1}-\alpha_k^{i*} \label{gamma2}
\end{align}

\vspace*{-5ex} \begin{multline}
P_k^{1*} = [(B_k^1 + B_k^2\bar {r}_k^{2})^{'}Z_k^1(B_k^1 + B_k^2\bar {r}_k^{2}) + R_k^{11} + \bar {r}_k^{2'}R_k^{12}\bar {r}_k^{2}]^{-1} \\ [ (B_k^1 + B_k^2\bar {r}_k^{2})^{'}Z_k^1(A_k + B_k^2W_k^{2}) + \bar {r}_k^{2'}R_k^{1j} W_k^{2}] \\\label{P12}
\end{multline} \vspace*{-5ex} 

\vspace*{-5ex} \begin{multline}
\alpha_k^{1*} = [(B_k^1 + B_k^2\bar {r}_k^{2})^{'}Z_k^1(B_k^1 + B_k^2\bar {r}_k^{2}) + R_k^{11}   + \bar {r}_k^{2'}R_k^{12}\bar {r}_k^{2}]^{-1}  [ (B_k^1 + B_k^2\bar {r}_k^{2})^{'}\\Z_k^1(s_k -B_k^2w_k^{2}) + \bar {r}_k^{2'}R_k^{12} w_k^{2}  - R_k^{11}\tilde{u}_k^{11} - \bar {r}_k^{2'}R_k^{12} \tilde{u}_k^{12} + (B_k^1 + B_k^2\bar{r}_k^{2})^{'}(\zeta_{k}^1) - Q_k^1\tilde{x}_k^1]\\\label{alpha12}
\end{multline} \vspace*{-5ex} 

\vspace*{-5ex} \begin{align}
P_k^{2*} = (R_k^{22} + B_k^{2'}Z_k^2B_k^2)^{-1}[B_k^{2'}Z_k^2(A_k - B_k^1P_k^{1*})] \label{Pf2} 
\end{align}

\vspace*{-5ex} \begin{align}
\alpha_k^{2*} = (R_k^{22} + B_k^{2'}Z_k^i2B_k^2)^{-1} [B_k^{2'}(Z_k^2(s_k - B_k^1\alpha_k^{1*})+ \zeta_{k}^2 -Q_k^2\tilde{x}_k^2)-R_k^{22}\tilde{u}_k^{22}] \label{alphaf2}
\end{align}

\vspace*{-5ex} \begin{multline}
Z_{k-1}^i = (A_k - \sum_{j \in \{1,2\}}B_k^jP_k^{j*})^{'}Z_k^i(A_k - \sum_{j \in \{1,2\}}B_k^jP_k^{j*}) + \\\sum_{j\in \{1,2\}}P_k^{j*'}R_k^{ij}P_k^{j*} + Q_{k-1}^i\; ;\; Z_T^i = Q_T^i\\
\end{multline} \vspace*{-5ex} 

\vspace*{-5ex} \begin{align}
\bar r_k^{2} = -(R_k^{22} + B_k^{2'}Z_k^2B_k^2)^{-1}B_k^{2'}Z_k^2B_k^1 \label{rbar2}
\end{align}

\vspace*{-5ex} \begin{multline}
\zeta_{k-1}^i = (A_k - \sum_{j \in \{1,2\}}B_k^jP_k^{j*})^{'}[\zeta_{k}^i + Z_k^i(s_k-\sum_{j\in \{1,2\}}B_k^j\alpha_k^{j*})-Q_k^i\tilde{x}_k^i] \\ + \sum_{j\in \{1,2\}}P_k^{j*'}R_k^{ij}(\alpha_k^{j*} + \tilde{u}_k^{ij})\; ;\; \zeta_{T}^i = 0 \\
\end{multline} \vspace*{-5ex} 

\vspace*{-5ex} \begin{align}
W_k^{2} = -(R_k^{22} + B_k^{2'}Z_k^2B_k^2)^{-1}B_k^{2'}Z_k^2 A_k \label{W2}
\end{align}

\vspace*{-5ex} \begin{align}
w_k^{2} = -(R_k^{22} + B_k^{2'}Z_k^2B_k^2)^{-1}[B_k^{2'}(Z_k^2s_k + \zeta_{k}^2 - Q_k^2\tilde{x}_k^2)-R_k^{22}\tilde{u}_k^{22}] \label{w2}
\end{align}

\vspace*{-5ex} \begin{align}
V^i(1,x_0) = \frac{1}{2} x_0^{'}Z_0^ix_0 + \zeta_0^{i'}x_0 + n_0^i \label{cost2}
\end{align}

\vspace*{-5ex} \begin{multline}
n_{k-1}^i = n_{k}^i + \frac{1}{2} (s_k-\sum_{j\in \{1,2\}}B_k^j\alpha_k^{j*})^{'}Z_k^i(s_k-\sum_{j\in \{1,2\}}B_k^j\alpha_k^{j*}) + \zeta_{k}^{i'}\\(s_k-\sum_{j\in \{1,2\}}B_k^j\alpha_k^{j*}) + \frac{1}{2}\sum_{j\in \{1,2\} }\alpha_k^{j*'}R_k^{ij}\alpha_k^{j*} -  \tilde{x}_k^{i'}Q_k^i(s_k-\sum_{j\in \{1,2\}}B_k^j\alpha_k^{j*}) \\ + \sum_{j\in \{1,2\}}\tilde{u}_k^{ij'}R_k^{ij}\alpha_k^{j*} + \frac{1}{2}(\tilde{x}_k^{i'}Q_k^i\tilde{x}_k^i + \sum_{j\in \{1,2\}}\tilde{u}_k^{ij'}R_k^{ij}\tilde{u}_k^{ij})\; ;\; n_{T}^i = 0\\\\\\\\\label{n2}
\end{multline} \vspace*{-5ex} 

\label{fbscoaq}
\end{cor}

\textsc {Proof:}

Corollary \eqref{fbscoaq} is proven in the same way as Theorem \eqref{FBSTheo} taking into consideration simplifications resulting from the different number of followers.\;\;\;\;\;\;\;\;\;\;\;\;\;\;\;\;\boxed{}\\

\begin{rem}
Special attention is drawn to the fact that the assumption about the existence of unique optimal solution sets of the systems of equations \eqref{Pf}, \eqref{alphaf}, \eqref{rbar}, \eqref{W} and \eqref{w} in Theorem \eqref{FBSTheo} is not needed in Corollary \eqref{fbscoaq} (and Corollary \eqref{lqfbs2p}), because these systems of equations degenerate to easily solvable equations in the case of only one follower.
\end{rem}

\label{FBS_ssaq}

\clearpage

\subsection{Special case: Linear-quadratic games with one leader and one follower}

In this subsection, first the results of the previous subsection \eqref{FBS_ssaq} are further specialized to a linear-quadratic 2-person game in Corollary \eqref{lqfbs2p} and then in Proposition \eqref{bofs} the specialized results are transformed into the terminology used in Corollary 7.2 in Ba\c{s}ar and Olsder (1999, pp. 374-375)\cite{baol} to point out some serious mistakes stated there.

\begin{cor}
A 2-person linear-quadratic dynamic game (cf. Def. \eqref{defspielaq}) admits a \emph{unique feedback Stackelberg equilibrium solution} if 
\begin{enumerate}
\item $Q_k^i \geq 0$, $R_k^{ii} > 0$ and $R_k^{ij} \geq 0$ (defined for $k \in K$ , $i,j \in \{1,2\}$, $j \not = i$)
\end{enumerate}
If these conditions are satisfied, the unique equilibrium strategies $\gamma_{k}^{i*}(x_{k})$ (i $\in \{1,2\}$) are given by \eqref{gamma2l} and the corresponding feedback Stackelberg equilibrium costs for the two players are stated in \eqref{cost2l}. \footnote {For all equations belonging to this theorem and its proof, $k \in K$ if nothing different is stated.}

\vspace*{-5ex} \begin{align}
 f_{k-1}(x_k,u_k^1,u_k^2) = A_kx_{k} + \sum_{j\in \{1,2\}}B_k^ju_k^j\label{f2l}
\end{align}

\vspace*{-5ex} \begin{align}
L^i(u^{1}, u^{2})= \sum_{k = 1}^T g_k^i(x_{k+1},u_k^{1}(x_k), u_k^{2}(x_k),x_k)
\end{align}

\vspace*{-5ex} \begin{align}
 g_k^i(x_{k+1},u_k^1,u_k^2,x_k) = \frac{1}{2} (x_{k+1}^{'}Q_{k+1}^ix_{k+1} + u_k^{i'}u_k^i+ u_k^{j'}R_k^{ij}u_k^j)\;, \; j \not = i \label{g2l}
\end{align}

\vspace*{-5ex} \begin{align}
\gamma_k^{i*}(x_k) = -P_k^{i*}x_k \label{gamma2l}
\end{align}

\vspace*{-5ex} \begin{multline}
P_k^{1*} = [(B_k^1 + B_k^2\bar {r}_k^{2})^{'}Z_{k+1}^1(B_k^1 + B_k^2\bar {r}_k^{2}) + I  \\ + \bar {r}_k^{2'}R_k^{12}\bar {r}_k^{2}]^{-1} [ (B_k^1 + B_k^2\bar {r}_k^{2})^{'}Z_{k+1}^1(A_k + B_k^2W_k^{2}) + \bar {r}_k^{2'}R_k^{12} W_k^{2}]\\\label{P12l}
\end{multline} \vspace*{-5ex} 

\vspace*{-5ex} \begin{align}
P_k^{2*} = (I + B_k^{2'}Z_{k+1}^2B_k^2)^{-1}[B_k^{2'}Z_{k+1}^2(A_k - B_k^1P_k^{1*})] \label{Pf2l} 
\end{align}

\vspace*{-5ex} \begin{multline}
Z_{k}^i = (A_k - B_k^iP_k^{i*} - B_k^jP_k^{j*})^{'}Z_{k+1}^i(A_k - B_k^iP_k^{i*} - B_k^jP_k^{j*}) \\+ P_k^{i*'}P_k^{i*} + P_k^{j*'}R_k^{ij}P_k^{j*} + Q_{k}^i\; ;\; Z_{T+1}^i = Q_{T+1}^i \; , \; i,j \in \{1,2\} \; , \; j \not = i\\\label{Zkilq}
\end{multline} \vspace*{-5ex} 

\vspace*{-5ex} \begin{align}
\bar r_k^{2} = -(I + B_k^{2'}Z_{k+1}^2B_k^2)^{-1}B_k^{2'}Z_{k+1}^2B_k^1 \label{rbar2l}
\end{align}

\vspace*{-5ex} \begin{align}
\zeta_k^i = (A_k - \sum_{j \in \{1,2\}}B_k^jP_k^{j*})^{'}(\zeta_{k+1}^i) + P_k^{i*} + R_k^{ij}P_k^{j*}\; ;\; \zeta_{T+1}^i = 0\;, \; j \not = i 
\end{align}

\vspace*{-5ex} \begin{align}
W_k^{2} = -(I + B_k^{2'}Z_{k+1}^2B_k^2)^{-1}B_k^{2'}Z_{k+1}^2A_k \label{W2l}
\end{align}

\vspace*{-5ex} \begin{align}
V^i(1,x_1) = \frac{1}{2} x_1^{'}Z_1^ix_1 + \zeta_1^{i'}x_1\label{cost2l}
\end{align}

\label{lqfbs2p}
\end{cor}

\textsc {Proof:}

For the proof of Corollary \eqref{lqfbs2p}, the same arguments are valid as for the proof of Corollary \eqref{fbscoaq}. Additionally simplifications result from the modified state equation and cost functionals.\;\;\;\;\;\;\;\;\;\;\;\;\;\;\;\;\boxed{}\\

\begin{rem}
Special attention should be paid to the observation that using a linear state equation together with a quadratic cost function without a term linearly dependent on $x_{k+1}$ yields linear equilibrium strategies. This also holds for feedback Nash, open-loop Nash and open-loop Stackelberg games.
\end{rem}

\begin{pro}
The systems of equations defining the unique equilibrium strategies $\gamma_{k}^{i*}(x_{k})$ (i $\in \{1,2\}$) in Corollary \eqref{lqfbs2p} can also be written in the following way:\footnote{For all equations belonging to this proposition and its proof, $k \in K$ if nothing different is stated. \eqref{Sk1b} is wrong in Ba\c{s}ar and Olsder.}

\vspace*{-5ex} \begin{align}
\gamma_k^{i*}(x_k) = -S_k^{i}x_k \label{gamma2lb}
\end{align}

\vspace*{-5ex} \begin{multline}
S_k^{1} = [B_k^{1'}(I + Z_{k+1}^2B_k^2B_k^{2'})^{-1}Z_{k+1}^1(I + B_k^2B_k^{2'}Z_{k+1}^2)^{-1}B_k^{1} + B_k^{1'}Z_{k+1}^{2'}B_k^{2}\\(I + B_k^{2'}Z_{k+1}^iB_k^2)^{-1}R_k^{12}(I + B_k^{2'}Z_{k+1}^2B_k^2)^{-1}B_k^{2'}Z_{k+1}^2B_k^1 + I]^{-1} B_k^{1'}[(I + Z_{k+1}^2B_k^2B_k^{2'})^{-1}Z_{k+1}^1\\(I + B_k^2B_k^{2'}Z_{k+1}^2)^{-1} + Z_{k+1}^{2'}B_k^{2'}(I + B_k^{2'}Z_{k+1}^iB_k^2)^{-1}R_k^{12}(I + B_k^{2'}Z_{k+1}^2B_k^2)^{-1}B_k^{2'}Z_{k+1}^2]A_k\\\label{Sk1b}
\end{multline} \vspace*{-5ex} 

\vspace*{-5ex} \begin{align}
S_k^{2} = (I + B_k^{2'}L_{k+1}^2B_k^2)^{-1}[B_k^{2'}L_{k+1}^2(A_k - B_k^1S_k^{1})] \label{Sk2b} 
\end{align}

\vspace*{-5ex} \begin{multline}
L_{k}^i = (A_k - B_k^iS_k^i - B_k^jS_k^j)^{'}L_{k+1}^i(A_k - B_k^iS_k^i - B_k^jS_k^j) \\+ S_k^{i'}S_k^i + S_k^{j'}R_k^{ij}S_k^j + Q_{k}^i\; ;\; Z_{T+1}^i = Q_{T+1}^i \; ,\; i,j \in \{1,2\} \; ,\; j \not = i\\\\\label{Lkb}
\end{multline} \vspace*{-5ex}

\label{bofs}
\end{pro}

\textsc {Proof:}\\

The proof is carried out by renaming some matrices and then showing that the relations for the equilibrium strategies $\gamma_k^{i*}$ of Corollary \eqref{lqfbs2p} can be rewritten in the way stated above.\\

Let us start by renaming the feedback matrices $P_{k}^{i*}$ and the matrices $Z_k^i$ (i $\in \{1,2\}$).

\vspace*{-5ex} \begin{multline}
P_{k}^{i*}\; \widehat {=}\; S_k^i \; ; \; Z_{k}^{i}\;   \widehat {=}\; L_k^i \; ; \; k \in K \; ; \; i \in \{1,2\}\\\label{äqufbsubs2}
\end{multline} \vspace*{-5ex} 

Next we prove that the $L_k^i$ fulfill \eqref{Lkb}. Taking consideration of the renaming \eqref{Zkilq} gives

\vspace*{-5ex} \begin{multline}
\nonumber \eqref{Lkb} \; \; L_{k}^i = (A_k - B_k^iS_k^i - B_k^jS_k^j)^{'}L_{k+1}^i(A_k - B_k^iS_k^i - B_k^jS_k^j) \\+ S_k^{i'}S_k^i + S_k^{j'}R_k^{ij}S_k^j + Q_{k}^i\; ;\; L_{T+1}^i = Q_{T+1}^i \; , \; i,j \in \{1,2\} \; , \; j \not = i\\
\end{multline} \vspace*{-5ex} 

Applying \eqref{äqufbsubs2} to \eqref{Pf2l} yields

\vspace*{-5ex} \begin{multline}
\nonumber \eqref{Sk2b} \; \; S_k^{2} = (I + B_k^{2'}L_{k+1}^2B_k^2)^{-1}[B_k^{2'}L_{k+1}^2(A_k - B_k^1S_k^{1})] \\
\end{multline} \vspace*{-5ex} 

Eventually the correctness of the rewritten feedback matrix $S_k^1$ has to be shown. For this purpose rewrite \eqref{P12l} considering \eqref{äqufbsubs2}.

\vspace*{-5ex} \begin{multline}
S_k^{1} = [(B_k^1 + B_k^2\bar {r}_k^{2})^{'}L_{k+1}^1(B_k^1 + B_k^2\bar {r}_k^{2}) + I  \\ + \bar {r}_k^{2'}R_k^{12}\bar {r}_k^{2}]^{-1} [ (B_k^1 + B_k^2\bar {r}_k^{2})^{'}L_{k+1}^1(A_k + B_k^2W_k^{2}) + \bar {r}_k^{2'}R_k^{12} W_k^{2}]\\
\end{multline} \vspace*{-5ex} 

Now using \eqref{rbar2l} and \eqref{W2l} in \eqref{P12l} (also considering \eqref{äqufbsubs2}) gives

\vspace*{-5ex} \begin{multline}
\nonumber S_k^{1} = [(B_k^1 - B_k^2(I + B_k^{2'}L_{k+1}^2B_k^2)^{-1}B_k^{2'}L_{k+1}^2B_k^1)^{'}L_{k+1}^1(B_k^1 - B_k^2(I + B_k^{2'}L_{k+1}^2B_k^2)^{-1}B_k^{2'}\\L_{k+1}^2B_k^1) + I + B_k^{1'}L_{k+1}^{2'}B_k^{2'}(I + B_k^{2}L_{k+1}^iB_k^2)^{-1}R_k^{12}(I + B_k^{2'}L_{k+1}^2B_k^2)^{-1}B_k^{2'}L_{k+1}^2B_k^1]^{-1}\\ [(B_k^1 - B_k^2(I + B_k^{2'}L_{k+1}^2B_k^2)^{-1}B_k^{2'}L_{k+1}^2B_k^1)^{'}L_{k+1}^1(A_k - B_k^2(I + B_k^{2'}L_{k+1}^2B_k^2)^{-1}B_k^{2'}L_{k+1}^2A_k) \\+ B_k^{1'}L_{k+1}^{2'}B_k^{2'}(I + B_k^{2'}L_{k+1}^iB_k^2)^{-1}R_k^{12}(I + B_k^{2'}L_{k+1}^2B_k^2)^{-1}B_k^{2'}L_{k+1}^2A_k]\\
\end{multline} \vspace*{-5ex}

Finally making some algebraic manipulations and applying Lemma \eqref{maid} to the four particular terms below gives
 
\vspace*{-5ex} \begin{multline}
\nonumber S_k^{1} = [B_k^{1'}(I - B_k^2(I + B_k^{2'}L_{k+1}^2B_k^2)^{-1}B_k^{2'}L_{k+1}^2)^{'}L_{k+1}^1(I - B_k^2(I + B_k^{2'}L_{k+1}^2B_k^2)^{-1}B_k^{2'}\\L_{k+1}^2)B_k^1 + I + B_k^{1'}L_{k+1}^{2'}B_k^{2'}(I + B_k^{2}L_{k+1}^iB_k^2)^{-1}R_k^{12}(I + B_k^{2'}L_{k+1}^2B_k^2)^{-1}B_k^{2'}L_{k+1}^2B_k^1]^{-1}\\ B_k^{1'}[(I - B_k^2(I + B_k^{2'}L_{k+1}^2B_k^2)^{-1}B_k^{2'}L_{k+1}^2)^{'}L_{k+1}^1(I - B_k^2(I + B_k^{2'}L_{k+1}^2B_k^2)^{-1}B_k^{2'}L_{k+1}^2) \\+ L_{k+1}^{2'}B_k^{2'}(I + B_k^{2'}L_{k+1}^iB_k^2)^{-1}R_k^{12}(I + B_k^{2'}L_{k+1}^2B_k^2)^{-1}B_k^{2'}L_{k+1}^2]A_k\\
\end{multline} \vspace*{-5ex} 

\vspace*{-5ex} \begin{multline}
\nonumber \eqref{Sk1b} \;\; S_k^{1} = [B_k^{1'}(I + L_{k+1}^2B_k^2B_k^{2'})^{-1}L_{k+1}^1(I + B_k^2B_k^{2'}L_{k+1}^2)^{-1}B_k^{1} + B_k^{1'}L_{k+1}^{2'}B_k^{2}\\(I + B_k^{2'}L_{k+1}^iB_k^2)^{-1}R_k^{12}(I + B_k^{2'}L_{k+1}^2B_k^2)^{-1}B_k^{2'}L_{k+1}^2B_k^1 + I]^{-1} B_k^{1'}[(I + L_{k+1}^2B_k^2B_k^{2'})^{-1}L_{k+1}^1\\(I + B_k^2B_k^{2'}L_{k+1}^2)^{-1} + L_{k+1}^{2'}B_k^{2'}(I + B_k^{2'}L_{k+1}^iB_k^2)^{-1}R_k^{12}(I + B_k^{2'}L_{k+1}^2B_k^2)^{-1}B_k^{2'}L_{k+1}^2]A_k\;\;\;\;\;\boxed{}\\
\end{multline} \vspace*{-5ex}

\label{lqfsg}

\clearpage
\setcounter{chapter}{3}
\setcounter{fussall}{\thefootnote}
\chapter {Discrete-Time Infinite Dynamic Games with Open-Loop Information Pattern}
\setcounter{footnote}{\thefussall}
%\section{Introduction}
%
%\clearpage
\section{Open-loop Nash Equilibrium Solutions}

This section is devoted to the derivation of the so-called open-loop Nash equilibrium solution for affine-quadratic games. First a general result is stated about the existence of a Nash equilibrium solution in \textit{n}-person discrete-time deterministic infinite dynamic games of prespecified fixed duration (cf. Def. \eqref{defspiel}) with open-loop information pattern. Then this result is applied to affine-quadratic games.

\subsection{Optimality conditions}

This subsection contains a theorem that gives sufficient conditions for the existence of an open-loop Nash equilibrium solution and provides equations for state, control and costate vectors which have to be satisfied on the equilibrium path.  Results about open-loop Nash equilibria in infinite dynamic games first appeared in continuous time in the works of Starr and Ho (1969) \cite{sthob}, \cite{sthoa} and Case (1969) \cite{cas}.

\begin{theo}
For an n-person discrete-time deterministic infinite dynamic game of prespecified fixed duration (cf. Def. \eqref{defspiel}) with open-loop information pattern let
\begin{itemize}
\item $f_k(\cdot,\cdot ,\ldots,\cdot)$ be continuously differentiable on $\textbf{R}^p \times \textbf{R}^{m_1} \times \ldots \times \textbf{R}^{m_n}$ (defined for $k \in K$)
\item $g_k^i(\cdot,\cdot ,\ldots,\cdot,\cdot)$ be continuously differentiable on $\textbf{R}^p \times \textbf{R}^{m_1} \times \ldots \times \textbf{R}^{m_n} \times \textbf{R}^p$ (defined for $k \in K$, $i \in N$)
\item $f_k(\cdot,\cdot ,\ldots,\cdot)$ be convex on $\textbf{R}^p \times \textbf{R}^{m_1} \times \ldots \times \textbf{R}^{m_n}$ (defined for $k \in K$)
\item $g_k^i(\cdot,\cdot ,\ldots,\cdot,\cdot)$ be convex on $\textbf{R}^p \times \textbf{R}^{m_1} \times \ldots \times \textbf{R}^{m_n} \times \textbf{R}^p$ (defined for $k \in K$, $i \in N$)
\item the cost functionals be stage-additive (cf. Def. \eqref{defcost}).
\end{itemize}
Then the set of strategies $\{\gamma_k^{i*}(x_0) = u_k^{i*}; i \in N \}$ provides an \emph {open-loop Nash equilibrium solution}. $\{x_k^*; k \in K\}$ is the corresponding state trajectory and a finite sequence of p-dimensional costate vectors $\{p_1^i,\ldots,p_K^i\}$ (defined for $i \in N$) exists so that the following relations are satisfied:\\

\vspace*{-5ex}  \begin{align}
x_k^* = f_{k-1}(x_{k-1}^*,u_k^{1*},\ldots,u_k^{n*}) \; , \; x_0^* = x_0 \label{OLN_x}
\end{align}

\vspace*{-5ex}  \begin{align}
\nabla_{u_k^i} H_k^i(p_k^i, u_k^{1*}, \ldots, u_k^{i-1*}, u_k^{i}, u_k^{i+1*},\ldots, u_k^{n*}, x_{k-1}^*) = 0
\end{align}

\vspace*{-5ex}  \begin{multline}
p_k^i = \frac {\partial}{\partial x_k} f_{k}(x_{k}^*,u_{k+1}^{1*},\ldots,u_{k+1}^{n*})^{'}[p_{k+1}^i + (\frac {\partial}{\partial x_{k+1}}g_{k+1}^i(x_{k+1}^*,\\ u_{k+1}^{1*}, \ldots, u_{k+1}^{n*}, x_{k}^*))^{'}] + [\frac {\partial}{\partial x_{k}}g_{k+1}^i(x_{k+1}^*, u_{k+1}^{1*}, \ldots, u_{k+1}^{n*}, x_{k}^*)]^{'} \; ; \; p_T^i = 0 \\
\end{multline} \vspace*{-5ex} 

where 

\vspace*{-5ex}  \begin{multline}
H_k^i(p_k^i, u_k^{1}, \ldots, u_k^{n}, x_{k-1}) \; \widehat{=} \; g_k^i(f_{k-1}(x_{k-1},u_{k}^{1},\ldots,u_{k}^{n}), \\ u_{k}^{1},\ldots,u_{k}^{n},x_{k-1}) + p_{k}^{i'}f_{k-1}(x_{k-1},u_{k}^{1},\ldots,u_{k}^{n}) \\\label{OLN_Hi}
\end{multline} \vspace*{-5ex} 

Every such equilibrium solution is \emph {weakly time consistent}.

\label{OLN_opt}
\end{theo}

\textsc {Proof:}

Theorem \eqref{FBS_opt} can be proven in the same way as Theorem 6.1 in Ba\c{s}ar and Olsder (1999, pp. 267-268)\cite{baol}, bearing in mind that in the above theorem we additionally assume (implicitly) the convexity of the cost functionals. That is why we apply Theorem \eqref{minprin} (instead of Theorem 5.5 in Ba\c{s}ar and Olsder (1999, p. 246)\cite{baol}) to the standard optimal control problem for player \textit{i} which emerges in the course of the derivation. That is why not only necessary but necessary and sufficient conditions are stated in the above theorem. \;\;\;\;\;\;\;\;\;\;\;\;\;\;\;\;\boxed{}

\clearpage

\subsection{Results for affine-quadratic games with arbitrarily many players}

In the following the results of Theorem \eqref{OLN_opt} are applied to an affine-quadratic dynamic game with arbitrarily many players. Theorem \eqref{OLNgc}, which is an extension (concerning the cost functionals) of Theorem 6.2 in Ba\c{s}ar and Olsder (1999, pp. 269-271)\cite{baol}, presents equilibrium equations that can easily be used for an algorithmic disintegration of the given Nash game.\\

Furthermore in Proposition \eqref{prooln} the equivalence of the equations of Theorem \eqref{OLNgc} with terminologically different equations is shown.

\begin{theo}
An n-person affine-quadratic dynamic game (cf. Def. \eqref{defspielaq}) admits a \emph{unique open-loop Nash equilibrium solution} if  
\begin{itemize}
\item $Q_k^i \geq 0$, $R_k^{ii} > 0$ (defined for $k \in K$ , $i \in N$).
\item $(I + \sum_{j\in N}B_k^j(R_k^{jj})^{-1}B_k^{j'}M_k^j)^{-1}$ (defined for $k \in K$) exist.
\end{itemize}
If these conditions are satisfied, the unique equilibrium strategies are given by \eqref{gammanol}, where the associated state trajectory $x_{k+1}^*$ is given by \eqref{xopt}. \footnote {For all equations belonging to this theorem and its proof, $i \in N$ and $k \in \{0 , \ldots ,T-1\}$ if nothing different is stated.}

\vspace*{-5ex}  \begin{align}
f_{k-1}(x_{k-1},u_k^1,\ldots,u_k^n) = A_kx_{k-1} + \sum_{j\in N}B_k^ju_k^j+s_k \; ; \; k \in K\label{fnol}
\end{align}

\vspace*{-5ex}  \begin{align}
L^i(x_0,u^{1},\ldots,u^{n})= \sum_{k = 1}^T g_k^i(x_k,u_k^{1},\ldots,u_k^{n},x_{k-1}) \label{lnol}
\end{align}

\vspace*{-5ex}  \begin{multline}
 g_k^i(x_k,u_k^1,\ldots,u_k^n,x_{k-1}) = \frac{1}{2} (x_k^{'}Q_k^ix_k + \sum_{j\in N}u_k^{j'}R_k^{ij}u_k^j)\; \\ +
 \frac{1}{2} (\tilde{x}_k^{i'}Q_k^i\tilde{x}_k^i + \sum_{j\in N}\tilde{u}_k^{ij'}R_k^{ij}\tilde{u}_k^{ij}) - \tilde{x}_k^{i'}Q_k^ix_k - \sum_{j\in N}\tilde{u}_k^{ij'}R_k^{ij}u_k^j \; ; \; k \in K\\\label{gnol}
\end{multline} \vspace*{-5ex}

\vspace*{-5ex}  \begin{align}
x_{k+1}^* = \Phi_kx_{k}^* + \phi_k \label{xopt}
\end{align}

\vspace*{-5ex}  \begin{align}
\Phi_k = (I + \sum_{j\in N}B_{k+1}^j(R_{k+1}^{jj})^{-1}B_{k+1}^{j'}M_{k+1}^j)^{-1}A_{k+1} \label{Phikn}
\end{align}

\vspace*{-5ex}  \begin{multline}
\phi_k = (I + \sum_{j\in N}B_{k+1}^j(R_{k+1}^{jj})^{-1}B_{k+1}^{j'}M_{k+1}^j)^{-1}(s_{k+1} \\- \sum_{j\in N}B_{k+1}^j((R_{k+1}^{jj})^{-1}B_{k+1}^{j'}(m_{k+1}^j - Q_{k+1}^j\tilde{x}_{k+1}^{j'}) - \tilde{u}_{k+1}^{jj})) \\\label{phikn}
\end{multline} \vspace*{-5ex}

\vspace*{-5ex}  \begin{align}
M_k^i = Q_{k}^i +  A_{k+1}^{'}M_{k+1}^i\Phi_k \; ; \; M_T^i = Q_{T}^i\label{M}
\end{align}

\vspace*{-5ex}  \begin{align}
m_{k}^i =  A_{k+1}^{'}[M_{k+1}^i\phi_k + m_{k+1}^i  - Q_{k+1}^i\tilde{x}_{k+1}^i ] \; ; \; m_T^i = 0 \label{m}
\end{align}

\vspace*{-5ex}  \begin{align}
\gamma_{k+1}^{i*}(x_0) = u_{k+1}^{i*} = -P_{k+1}^{i*}x_{k}^* - \alpha_{k+1}^{i*} \label{gammanol}
\end{align}

\vspace*{-5ex}  \begin{align}
P_{k+1}^{i*} = (R_{k+1}^{ii})^{-1}B_{k+1}^{i'}M_{k+1}^i\Phi_k \label{Pnol}
\end{align}

\vspace*{-5ex}  \begin{align}
\alpha_{k+1}^i = (R_{k+1}^{ii})^{-1}B_{k+1}^{i'}(M_{k+1}^i\phi_k + m_{k+1}^i - Q_{k+1}^i\tilde{x}_{k+1}^{i'}) + \tilde{u}_{k+1}^{ii} \label{alphanol}
\end{align}\\\\

\label{OLNgc}
\end{theo}

\textsc {Proof:}

Theorem \eqref{OLN_opt} can be applied to the given affine-quadratic game, since all conditions are satisfied for the given state equation \eqref{fnol} and cost functionals \eqref{lnol}. Furthermore $g_k^i$ is strictly convex in $u_k^i$. This can be seen by applying Corollary \eqref{pdc} to \eqref{conol}. Therefore there has to be a unique optimal equilibrium solution.

\vspace*{-5ex}  \begin{multline}
\frac{\partial^2}{\partial u_i^2} g_k^i(x_k,u_k^1 ,\ldots , u_k^n, x_{k-1}) = B_k^{i'}Q_k^iB_k^i + R_k^{ii} \\\label{conol}
\end{multline} \vspace*{-5ex} 

To obtain relations which satisfy this unique solution, we have to adapt \eqref{OLN_x} - \eqref{OLN_Hi} to the given state equation and cost functionals. This yields

\vspace*{-5ex}  \begin{multline}
H_k^i =   \frac{1}{2} (x_k^{'}Q_k^ix_k + \sum_{j\in N}u_k^{j'}R_k^{ij}u_k^j)\;  +
 \frac{1}{2} (\tilde{x}_k^{i'}Q_k^i\tilde{x}_k^i + \sum_{j\in N}\tilde{u}_k^{ij'}R_k^{ij}\tilde{u}_k^{ij}) \\- \tilde{x}_k^{i'}Q_k^ix_k - \sum_{j\in N}\tilde{u}_k^{ij'}R_k^{ij}u_k^j + p_k^{i'}(A_kx_{k-1} + \sum_{j\in N}B_k^ju_k^j+s_k)\\\label{Hik}
\end{multline} \vspace*{-5ex} 

\vspace*{-5ex}  \begin{multline}
p_k^i = A_{k+1}^{'}[p_{k+1} + Q_{k+1}^i( x_{k+1}^* - \tilde{x}_{k+1}^i )] \; ; \; p_T^i = 0\\\label{pik}
\end{multline} \vspace*{-5ex} 

\vspace*{-5ex}  \begin{multline}
\nonumber \frac{\partial}{\partial u_i} \eqref{Hik} = 0\; \Rightarrow \;\;\ B_k^{i'}Q_k^ix_k^* + R_k^{ii}u_k^{i*} -  B_k^{i'}Q_k^i\tilde{x}_k^{i'} - R_k^{ii}\tilde{u}_k^{ii} + B_k^{i'}p_{k}^i = 0\\
\end{multline} \vspace*{-5ex} 

\vspace*{-5ex}  \begin{multline}
u_k^{i*} =  -(R_k^{ii})^{-1}B_k^{i'}(Q_k^i(x_k^*-\tilde{x}_k^{i'}) + p_{k}^i ) + \tilde{u}_k^{ii}\\\label{ukxk}
\end{multline} \vspace*{-5ex} 

\vspace*{-5ex}  \begin{multline}
x_{k}^* = A_kx_{k-1}^* + \sum_{j\in N}B_k^ju_k^{j*}+s_k \; ; \; x_0^* = x_0\\\label{xk}
\end{multline} \vspace*{-5ex} 

In the following induction argument we will give proof that \eqref{psub} is valid and the recursive relations for $M_{k}^{i}$ and $m_{k}^i$ (stated in the above theorem) are correct.

\vspace*{-5ex}  \begin{multline}
p_k^i = (M_{k}^i-Q_{k}^i)x_{k}^* + m_{k}^i \\\label{psub}
\end{multline} \vspace*{-5ex}

\textbf{Basis:}

The induction starts at \textit{k} = \textit{T}. First we make use of the general optimality conditions for $p_k^i$ at stage \textit{T}.

\vspace*{-5ex}  \begin{multline}
\nonumber\eqref{pik}_{k=T} \;\; p_{T}^i = 0\\
\end{multline} \vspace*{-5ex} 

Now we can substitute the $p_T^i$ with functions affinely dependent on $x_T$.

\vspace*{-5ex}  \begin{multline}
\nonumber\eqref{psub}_{k=T} \;\;p_{T}^i = (M_{T}^i-Q_{T}^i)x_{T}^* + m_{T}^i \\
\end{multline} \vspace*{-5ex} 

\vspace*{-5ex}  \begin{multline}
\nonumber (M_{T}^i-Q_{T}^i)x_{T}^* + m_{T}^i = 0 \\
\end{multline} \vspace*{-5ex} 

Comparing coefficients gives

\vspace*{-5ex}  \begin{multline}
\nonumber\eqref{M}_{k=T} \;\; M_T^i = Q_{T}^i \\
\end{multline} \vspace*{-5ex} 

\vspace*{-5ex}  \begin{multline}
\nonumber\eqref{m}_{k=T} \;\; m_T^i = 0 \\
\end{multline} \vspace*{-5ex}

\textbf{Inductive step:}

As an induction hypothesis, the system of equations \eqref{psub} is assumed to be true at stage \textit{l}+1. Now we have to prove that this system of equations is fulfilled at stage \textit{l} and determines the corresponding recursive relations for $M_{l}^i$ and $m_{l}^i$.

\vspace*{-5ex}  \begin{multline}
\nonumber\eqref{psub}_{k=l+1} \;\; p_{l+1}^i = (M_{l+1}^i-Q_{(l+1)}^i)x_{l+1}^* + m_{l+1}^i \\
\end{multline} \vspace*{-5ex} 

First the induction hypothesis is used in the general optimality conditions for $p_k^i$ at stage \textit{l}.

\vspace*{-5ex}  \begin{multline}
\nonumber\eqref{pik}_{k=l} \;\; p_{l}^i = A_{l+1}^{'}[p_{l+1} + Q_{l+1}^i(x_{l+1}^* - \tilde{x}_{l+1}^i )] \\
\end{multline} \vspace*{-5ex} 

\vspace*{-5ex}  \begin{multline}
\nonumber p_{l}^i = A_{l+1}^{'}[(M_{l+1}^i-Q_{l+1}^i)x_{l+1}^* + m_{l+1}^i  + Q_{l+1}^i(x_{l+1}^* - \tilde{x}_{l+1}^i )] \\
\end{multline} \vspace*{-5ex} 

\vspace*{-5ex}  \begin{multline}
\nonumber p_{l}^i = A_{l+1}^{'}[M_{l+1}^ix_{l+1}^* + m_{l+1}^i  - Q_{l+1}^i\tilde{x}_{l+1}^i ] \\
\end{multline} \vspace*{-5ex}

To complete the inductive step we have to prove that the $p_l^i$ can be written as affine functions of $x_l$. Therefore an interrelation between $x_l$ and $x_{l+1}$, which does not depend on the controls of the players nor on the costate variables $p_{l+1}^i$, has to be deduced. To do so, first we have to substitute $u_{l+1}^{i*}$ in the equation stated below for the evolution of the optimal state vector $x_{l+1}^*$ by terms which are affine in  $x_l$ and $x_{l+1}$ and furthermore only contain $M_{l+1}^i$, $m_{l+1}^i$ and matrices and vectors given by the game definition.

\vspace*{-5ex}  \begin{multline}
\nonumber\eqref{xk}_{k=l+1}\;\; x_{l+1}^* = A_{l+1}x_{l}^* + \sum_{j\in N}B_{l+1}^ju_{l+1}^{j*} + s_{l+1} \\
\end{multline} \vspace*{-5ex}

In the first instance the optimality condition for $u_{l+1}^{i*}$ can be rewritten with the help of $\eqref{psub}_{k=l+1}$, which are the induction hypotheses for $p_{l+1}^i$. 

\vspace*{-5ex}  \begin{multline}
\nonumber\eqref{ukxk}_{k=l+1}\;\; u_{l+1}^{i*} =  -(R_{l+1}^{ii})^{-1}B_{l+1}^{i'}(Q_{l+1}^i(x_{l+1}^*-\tilde{x}_{l+1}^{i'}) + p_{l+1}^i ) + \tilde{u}_{l+1}^{ii}\\
\end{multline} \vspace*{-5ex} 

\vspace*{-5ex}  \begin{multline}
\nonumber u_{l+1}^{i*} =  -(R_{l+1}^{ii})^{-1}B_{l+1}^{i'}(Q_{l+1}^i(x_{l+1}^*-\tilde{x}_{l+1}^{i'}) \\+ (M_{l+1}^i-Q_{l+1}^i)x_{l+1}^* + m_{l+1}^i) + \tilde{u}_{l+1}^{ii}\\
\end{multline} \vspace*{-5ex} 

\vspace*{-5ex}  \begin{multline}
\nonumber u_{l+1}^{i*} =  -(R_{l+1}^{ii})^{-1}B_{l+1}^{i'}(M_{l+1}^ix_{l+1}^* + m_{l+1}^i - Q_{l+1}^i\tilde{x}_{l+1}^{i'}) + \tilde{u}_{l+1}^{ii}\\
\end{multline} \vspace*{-5ex}

Now the control variables can be replaced in the optimal state equation by terms affine in $x_{l+1}$.

\vspace*{-5ex}  \begin{multline}
\nonumber x_{l+1}^* = A_{l+1}x_{l}^* - \sum_{j\in N}B_{l+1}^j((R_{l+1}^{jj})^{-1}[B_{l+1}^{j'}(M_{l+1}^jx_{l+1}^* \\+ m_{l+1}^j - Q_{l+1}^j\tilde{x}_{l+1}^{j'})] + \tilde{u}_{l+1}^{jj}) + s_{l+1} \\
\end{multline} \vspace*{-5ex} 

Making $x_{l+1}$ explicit gives

\vspace*{-5ex}  \begin{multline}
\nonumber x_{l+1}^* = (I + \sum_{j\in N}B_{l+1}^j(R_{l+1}^{jj})^{-1}B_{l+1}^{j'}M_{l+1}^j)^{-1} [A_{l+1}x_{l}^* \\- \sum_{j\in N}B_{l+1}^j((R_{l+1}^{jj})^{-1}B_{l+1}^{j'}(m_{l+1}^j - Q_{l+1}^j\tilde{x}_{l+1}^{j'}) - \tilde{u}_{l+1}^{jj}) + s_{l+1}] \\
\end{multline} \vspace*{-5ex} 

The structure of the above equation justifies the following substitution

\vspace*{-5ex}  \begin{multline}
\nonumber\eqref{xopt}_{k=l}\;\; x_{l+1}^* = \Phi_lx_{l}^* + \phi_l \\
\end{multline} \vspace*{-5ex} 

\vspace*{-5ex}  \begin{multline}
\Phi_lx_{l}^* + \phi_l = (I + \sum_{j\in N}B_{l+1}^j(R_{l+1}^{jj})^{-1}B_{l+1}^{j'}M_{l+1}^j)^{-1} [A_{l+1}x_{l}^* \\- \sum_{j\in N}B_{l+1}^j((R_{l+1}^{jj})^{-1}B_{l+1}^{j'}(m_{l+1}^j - Q_{l+1}^j\tilde{x}_{l+1}^{j'}) - \tilde{u}_{l+1}^{jj}) + s_{l+1}] \\\label{xlgn}
\end{multline} \vspace*{-5ex} 

By comparing coefficients it follows that

\vspace*{-5ex}  \begin{multline}
\nonumber\eqref{xlgn}_{x_{l}^*}=\eqref{Phikn}_{k=l}\;\; \Phi_l = (I + \sum_{j\in N}B_{l+1}^j(R_{l+1}^{jj})^{-1}B_{l+1}^{j'}M_{l+1}^j)^{-1}A_{l+1} \\
\end{multline} \vspace*{-5ex} 

\vspace*{-5ex}  \begin{multline}
\nonumber\eqref{xlgn}_{const.}=\eqref{phikn}_{k=l}\;\; \phi_l = (I + \sum_{j\in N}B_{l+1}^j(R_{l+1}^{jj})^{-1}B_{l+1}^{j'}M_{l+1}^j)^{-1}(s_{l+1} - \sum_{j\in N}B_{l+1}^j((R_{l+1}^{jj})^{-1}\\B_{l+1}^{j'}(m_{l+1}^j - Q_{l+1}^j\tilde{x}_{l+1}^{j'}) - \tilde{u}_{l+1}^{jj})) \\
\end{multline} \vspace*{-5ex}

The above relation between $x_{l+1}$ and $x_l$, given by equation $\eqref{xopt}_{k=l}$, can be used to finish the inductive step for $p_l^i$.

\vspace*{-5ex}  \begin{multline}
\nonumber p_{l}^i = A_{l+1}^{'}[M_{l+1}^ix_{l+1}^* + m_{l+1}^i  - Q_{l+1}^i\tilde{x}_{l+1}^i ] \\
\end{multline} \vspace*{-5ex} 

\vspace*{-5ex}  \begin{multline}
\nonumber p_{l}^i = A_{l+1}^{'}[M_{l+1}^i(\Phi_lx_{l}^* + \phi_l) + m_{l+1}^i  - Q_{l+1}^i\tilde{x}_{l+1}^i ] \\
\end{multline} \vspace*{-5ex} 

The structure of the above equations justifies the following substitutions

\vspace*{-5ex}  \begin{multline}
\nonumber\eqref{psub}_{k=l}\;\; p_l^i = (M_{l}^i-Q_{l}^i)x_{l}^* + m_{l}^i \\
\end{multline} \vspace*{-5ex} 

\vspace*{-5ex}  \begin{multline}
(M_{l}^i-Q_{l}^i)x_{l}^* + m_{l}^i = A_{l+1}^{'}[M_{l+1}^i(\Phi_lx_{l}^* + \phi_l) + m_{l+1}^i  - Q_{l+1}^i\tilde{x}_{l+1}^i ] \\\label{plgn}
\end{multline} \vspace*{-5ex} 

Comparing coefficients gives

\vspace*{-5ex}  \begin{multline}
\nonumber\eqref{plgn}_{x_{l}^*}=\eqref{M}_{k=l} \;\; M_l^i = Q_{l}^i +  A_{l+1}^{'}M_{l+1}^i\Phi_l\\
\end{multline} \vspace*{-5ex} 

\vspace*{-5ex}  \begin{multline}
\nonumber\eqref{plgn}_{const.}=\eqref{m}_{k=l} \;\; m_{l}^i =  A_{l+1}^{'}[M_{l+1}^i\phi_l + m_{l+1}^i  - Q_{l+1}^i\tilde{x}_{l+1}^i ]\\
\end{multline} \vspace*{-5ex}

At this point the inductive step and hence the induction argument is completed, but we try to transform $u_{l+1}^{i*}$ so that their evolution depends affinely on $x_l^*$ and therefore, their algorithmic computation is straightforward. This is done by using $\eqref{xopt}_{k=l}$ in the equation deduced above for $u_{l+1}^{i*}$.

\vspace*{-5ex}  \begin{multline}
\nonumber u_{l+1}^{i*} =  -(R_{l+1}^{ii})^{-1}B_{l+1}^{i'}(M_{l+1}^ix_{l+1}^* + m_{l+1}^i - Q_{l+1}^i\tilde{x}_{l+1}^{i'}) + \tilde{u}_{l+1}^{ii}\\
\end{multline} \vspace*{-5ex} 

\vspace*{-5ex}  \begin{multline}
\nonumber u_{l+1}^{i*} =  -(R_{l+1}^{ii})^{-1}B_{l+1}^{i'}(M_{l+1}^i(\Phi_lx_{l}^* + \phi_l) + m_{l+1}^i - Q_{l+1}^i\tilde{x}_{l+1}^{i'}) + \tilde{u}_{l+1}^{ii}\\
\end{multline} \vspace*{-5ex} 

The structure of the above equations justifies the following substitutions

\vspace*{-5ex}  \begin{multline}
\nonumber\eqref{gammanol}_{k=l} \;\; \gamma_{l+1}^{i*}(x_0) = u_{l+1}^{i*} = -P_{l+1}^{i*}x_{l}^* - \alpha_{l+1}^{i*} \\
\end{multline} \vspace*{-5ex} 

\vspace*{-5ex}  \begin{multline}
-P_{l+1}^{i*}x_{l}^* - \alpha_{l+1}^{i*} =  -(R_{l+1}^{ii})^{-1}B_{l+1}^{i'}\\(M_{l+1}^i(\Phi_lx_{l}^* + \phi_l) + m_{l+1}^i - Q_{l+1}^i\tilde{x}_{l+1}^{i'}) + \tilde{u}_{l+1}^{ii}\\\label{uievn}
\end{multline} \vspace*{-5ex} 

Comparing coefficients it follows that

\vspace*{-5ex}  \begin{multline}
\nonumber\eqref{uievn}_{x_{l}^*}=\eqref{Pnol}_{k=l} \;\; P_{l+1}^{i*} = (R_{l+1}^{ii})^{-1}B_{l+1}^{i'}M_{l+1}^i\Phi_l\\
\end{multline} \vspace*{-5ex} 

\vspace*{-5ex}  \begin{multline}
\nonumber\eqref{uievn}_{const.}=\eqref{alphanol}_{k=l} \;\; \alpha_{l+1}^i = (R_{l+1}^{ii})^{-1}B_{l+1}^{i'}(M_{l+1}^i\phi_l + m_{l+1}^i - Q_{l+1}^i\tilde{x}_{l+1}^{i'}) + \tilde{u}_{l+1}^{ii} \;\;\;\;\boxed{}\\\\
\end{multline} \vspace*{-5ex}

\begin{rem}
By making the induction hypothesis that $p_k^i$ is dependent on $x_k$ (instead of $x_{k+1}$ in Ba\c{s}ar and Olsder (1999, pp.269 - 271)\cite{baol}) the number of algebraic manipulations and substitutions is decreased to less than half and also the complexity and length of the equations is considerably reduced.
\end{rem}

\begin{rem}
Special attention should be paid to the fact that the induction argument only runs over the costate vectors $p_k^i$, and the general relations for $x_k$ and $u_k^i$ given in \eqref{xopt} and \eqref{gammanol} arise from the combination of the general optimality conditions for $x_k$ and $u_k^i$ given in \eqref{xk} and \eqref{ukxk} (which have to hold true at each stage of the game) and the induction hypothesis for $p_k^i$. In contrast the proof in Ba\c{s}ar and Olsder (1999, pp.269 - 271)\cite{baol}) looks as though an induction for $p_k^i$, $x_k$ and $u_k^i$ is being made with the basis at the last stage (which is not allowed because the basis for $x_k$ and $u_k^i$ cannot be verified (in a correct way) at the last stage), but not using the induction hypothesis for $x_k$ and $u_k^i$ in the inductive step.
\end{rem}

\begin{rem}
To solve the Nash game algorithmically, the following order of application of the equations of Theorem \eqref{OLNgc} is advisable ($i \in N$):
\begin{enumerate}
\item $M_{T}^{i}$, $m_{T}^{i}$ 
\item For k running backward from T-1 to 0
\begin{enumerate}
\item $\Phi_k$, $\phi_k$
\item $M_{k}^{i}$, $m_{k}^{i}$  
\end{enumerate}
\item $x_{0}^*$
\item For k running forward from 0 to T-1
\begin{enumerate}
\item $P_{k+1}^{i}$, $\alpha_{k+1}^i$ 
\item $u_{k+1}^{i*}$
\item $x_{k+1}^*$
\item $g_{k+1}^i(x_k,u_k^1,\ldots,u_k^n,x_{k-1})$
\end{enumerate}
\item $L^i(x_0,u^{1},\ldots,u^{n})$
\end{enumerate}
\end{rem}

\begin{pro}\footnote{In this proposition we rewrite the equilibrium equations in a notation that was used at our department in the past to enable comparison.}
The systems of equations defining the unique equilibrium strategies $\gamma_{k+1}^{i*}(x_0)$ and the associated state trajectory $x_{k+1}^*$ in Theorem \eqref{OLNgc} can also be written in the following way:\footnote{For all equations belonging to this proposition and its proof, $k \in K$ and $i \in N$ if nothing different is stated.}

\vspace*{-5ex}  \begin{align}
x_{k+1}^* = (\Lambda_{k+1})^{-1}(A_{k+1}x_{k}^* + \eta_{k+1}) \label{xoptä}
\end{align}

\vspace*{-5ex}  \begin{align}
\Lambda_{k+1} = I + \sum_{j\in N}B_{k+1}^j(R_{k+1}^{jj})^{-1}B_{k+1}^{j'}H_{k+1}^j \label{laminv}
\end{align}

\vspace*{-5ex}  \begin{align}
\eta_{k+1} = s_{k+1} + \sum_{j\in N}B_{k+1}^j(\tilde{u}_{k+1}^{jj} - (R_{k+1}^{jj})^{-1}B_{k+1}^{j'}h_{k+1}^j) \label{etaä}
\end{align}

\vspace*{-5ex}  \begin{align}
H_k^i = Q_{k}^i +  A_{k+1}^{'}H_{k+1}^i(\Lambda_{k+1})^{-1}A_{k+1} \; ; \; H_T^i = Q_{T}^i\label{Mä}
\end{align}

\vspace*{-5ex}  \begin{align}
h_{k}^i =  - Q_k^i\tilde{x}_k^i + A_{k+1}^{'}[H_{k+1}^i(\Lambda_{k+1})^{-1}\eta_{k+1} + h_{k+1}^i]  \; ; \; h_T^i = -Q_T^i\tilde{x}_T^i \label{mä}
\end{align}

\vspace*{-5ex}  \begin{align}
\gamma_{k+1}^{i*}(x_0) = u_{k+1}^{i*} = \tilde{u}_{k+1}^{ii} - (R_{k+1}^{ii})^{-1}B_{k+1}^{i'}(H_{k+1}^ix_{k+1}^* + h_{k+1}^i)\label{gammanolä}
\end{align}\\\\

\label{prooln}
\end{pro}

\textsc {Proof:}

The proof is done by renaming the costate matrices and then showing that the relations for the costate matrices and the optimal state and control vectors of Theorem \eqref{OLNgc} can be rewritten in the way stated above.

Let us start by renaming the costate matrices $M_{k}^{i}$ and $m_{k}^{i}$.

\vspace*{-5ex}  \begin{multline}
M_k^i \; \widehat {=} \; H_k^i \; ; \;   m_k^i \; \widehat {=} \; h_k^i + Q_k^i\tilde{x}_k^i \\\label{äqusubsol}
\end{multline} \vspace*{-5ex} 

Next we prove that the costate matrices fulfill \eqref{Mä} and \eqref{mä} respectively.\\

Taking the renaming into consideration \eqref{M} gives

\vspace*{-5ex}  \begin{multline}
\nonumber H_k^i = Q_{k}^i +  A_{k+1}^{'}H_{k+1}^i\Phi_k \\
\end{multline} \vspace*{-5ex} 

Substituting $\Phi_k$ with the help of \eqref{Phikn} (and considering \eqref{äqusubsol}) leads to

\vspace*{-5ex}  \begin{multline}
\nonumber H_k^i = Q_{k}^i +  A_{k+1}^{'}H_{k+1}^i(I + \sum_{j\in N}B_{k+1}^j(R_{k+1}^{jj})^{-1}B_{k+1}^{j'}H_{k+1}^j)^{-1}A_{k+1} \\
\end{multline} \vspace*{-5ex} 

Making use of \eqref{laminv} yields

\vspace*{-5ex}  \begin{multline}
\nonumber\eqref{Mä} \; \; H_k^i = Q_{k}^i +  A_{k+1}^{'}H_{k+1}^i(\Lambda_{k+1})^{-1}A_{k+1} \\
\end{multline} \vspace*{-5ex} 

Now we show the correctness of \eqref{mä}. To do this we start at stage \textit{T} and use \eqref{m} and \eqref{äqusubsol} to get

\vspace*{-5ex}  \begin{multline}
\nonumber 0 = m_T^i = h_T^i + Q_T^i\tilde{x}_T^i \\
\end{multline} \vspace*{-5ex} 

\vspace*{-5ex}  \begin{multline}
\nonumber\eqref{mä}_{k=T} \; \; h_T^i = - Q_T^i\tilde{x}_T^i \\
\end{multline} \vspace*{-5ex} 

For the general stage k rewrite \eqref{m} taking consideration of \eqref{äqusubsol}

\vspace*{-5ex}  \begin{multline}
\nonumber h_{k}^i + Q_k^i\tilde{x}_k^i =  A_{k+1}^{'}[H_{k+1}^i\phi_k + h_{k+1}^i] \\
\end{multline} \vspace*{-5ex} 

Substituting $\phi_k$ with the help of \eqref{phikn} (and considering \eqref{äqusubsol}) yields

\vspace*{-5ex}  \begin{multline}
\nonumber h_{k}^i =  - Q_k^i\tilde{x}_k^i + A_{k+1}^{'}[H_{k+1}^i(I + \sum_{j\in N}B_{k+1}^j(R_{k+1}^{jj})^{-1}B_{k+1}^{j'}H_{k+1}^j)^{-1}\\(s_{k+1} - \sum_{j\in N}B_{k+1}^j((R_{k+1}^{jj})^{-1}B_{k+1}^{j'}h_{k+1}^j - \tilde{u}_{k+1}^{jj})) + h_{k+1}^i] \\
\end{multline} \vspace*{-5ex} 

Eventually using \eqref{laminv} and \eqref{etaä} gives

\vspace*{-5ex}  \begin{multline}
\nonumber \eqref{mä} \; \; h_{k}^i =  - Q_k^i\tilde{x}_k^i + A_{k+1}^{'}[H_{k+1}^i(\Lambda_{k+1})^{-1}\eta_{k+1} + h_{k+1}^i] \\
\end{multline} \vspace*{-5ex} 

Next we show the correctness of \eqref{xoptä}. To do this we start with equation \eqref{xopt} and substitute $\Phi_k$ and $\phi_k$ (and also consider \eqref{äqusubsol})

\vspace*{-5ex}  \begin{multline}
\nonumber x_{k+1}^* = (I + \sum_{j\in N}B_{k+1}^j(R_{k+1}^{jj})^{-1}B_{k+1}^{j'}H_{k+1}^j)^{-1}A_{k+1}x_{k}^* + (I + \sum_{j\in N}B_{k+1}^j(R_{k+1}^{jj})^{-1}B_{k+1}^{j'}H_{k+1}^j)^{-1}\\(s_{k+1} - \sum_{j\in N}B_{k+1}^j((R_{k+1}^{jj})^{-1}B_{k+1}^{j'}h_{k+1}^j - \tilde{u}_{k+1}^{jj}))
\end{multline} \vspace*{-5ex} 

Making use of \eqref{laminv} and \eqref{etaä} we get

\vspace*{-5ex}  \begin{multline}
\nonumber \eqref{xoptä} \; \; x_{k+1}^* = (\Lambda_{k+1})^{-1}(A_{k+1}x_{k}^* + \eta_{k+1}) \\
\end{multline} \vspace*{-5ex} 

Eventually the correctness of the rewritten equilibrium strategies $u_{k+1}^{i*}$, given by \eqref{gammanolä}, has to be shown.\\

First substituting $P_{k+1}^i$ and $\alpha_{k+1}^i$ in \eqref{gammanol} with the help of \eqref{Pnol} and  \eqref{alphanol} (and also considering \eqref{äqusubsol}) gives

\vspace*{-5ex}  \begin{multline}
\nonumber u_{k+1}^{i*} = -(R_{k+1}^{ii})^{-1}B_{k+1}^{i'}H_{k+1}^i\Phi_kx_{k}^* - (R_{k+1}^{ii})^{-1}B_{k+1}^{i'}(H_{k+1}^i\phi_k + h_{k+1}^i) + \tilde{u}_{k+1}^{ii} \\
\end{multline} \vspace*{-5ex} 

Simplifying the above equation leads to

\vspace*{-5ex}  \begin{multline}
\nonumber u_{k+1}^{i*} = \tilde{u}_{k+1}^{ii} - (R_{k+1}^{ii})^{-1}B_{k+1}^{i'}H_{k+1}^i(\Phi_kx_{k}^* + \phi_k + h_{k+1}^i) \\
\end{multline} \vspace*{-5ex} 

Finally using \eqref{xopt} yields

\vspace*{-5ex}  \begin{multline}
\nonumber \eqref{gammanolä} \; \; u_{k+1}^{i*} = \tilde{u}_{k+1}^{ii} - (R_{k+1}^{ii})^{-1}B_{k+1}^{i'}H_{k+1}^i(x_{k+1}^* + h_{k+1}^i) \;\;\;\;\;\;\;\;\;\;\;\;\;\;\;\;\boxed{}\\
\end{multline} \vspace*{-5ex} 

\label{ssaqoln}

\clearpage
%%% Local Variables: 
%%% mode: latex
%%% TeX-master: "Diplomarbeit"
%%% End: 

\section {Open-loop Stackelberg Equilibrium Solutions}

This section is devoted to the derivation of the so-called open-loop Stackelberg equilibrium solution with one leader and arbitrarily many followers for affine-quadratic games. First a general result is stated about the existence and uniqueness of a Stackelberg equilibrium solution with one leader and arbitrarily many followers in $n$-person discrete-time deterministic infinite dynamic games of prespecified fixed duration (cf. Def. \eqref{defspiel}) with open-loop information pattern. Then this result is applied to affine-quadratic games. First a proof geared to the one indicated in Ba\c{s}ar and Olsder (1999, p. 372)\cite{baol} is presented. But this proof "produces" a hardly algorithmically solvable system of equilibrium equations. Therefore another way of deriving the equilibrium solution for affine-quadratic games is presented giving us a system of equilibrium equations that can easily be used for an algorithmic disintegration of the given Stackelberg game.

\subsection{Optimality conditions}

The following theorem gives sufficient conditions for the existence of an open-loop Stackelberg equilibrium solution with one leader and arbitrarily many followers and provides equations for state, control, costate and cocontrol vectors, which have to be satisfied on the equilibrium path. Results about open-loop Stackelberg equilibria in infinite dynamic games first appeared in continuous time in the works of Chen and Cruz (1972) \cite{chcr} and Simaan and Cruz (1973) \cite{sicra}, \cite{sicrb}.

\begin{theo}
For an n-person discrete-time deterministic infinite dynamic game of prespecified fixed duration (cf. Def. \eqref{defspiel}) with open-loop information pattern let
\begin{itemize}
\item $f_k(x_{k-1}, \cdot, u_k^2, \ldots,u_k^n)$ be continuously differentiable on $\textbf{R}^{m_1}$ (defined for $k \in K$)
\item $f_k(\cdot,u_k^1, \cdot , \ldots, \cdot)$ be twice continuously differentiable on $\textbf{R}^p \times \textbf{R}^{m_2} \times \ldots \times \textbf{R}^{m_n}$ (defined for $k \in K$)
\item $g_k^1(\cdot,\cdot ,\ldots,\cdot,\cdot)$ be continuously differentiable on $\textbf{R}^p \times \textbf{R}^{m_1} \times \ldots \times \textbf{R}^{m_n} \times \textbf{R}^p$ (defined for $k \in K$)
\item $g_k^i(x_k , \cdot, u_k^2, \ldots,u_k^n, x_{k-1})$ be continuously differentiable on $\textbf{R}^{m_1}$ (defined for $k \in K$, $i \in \{2 \ldots n\}$)
\item $g_k^i(\cdot,u_k^1, \cdot ,\ldots,\cdot,\cdot)$ be twice continuously differentiable on $\textbf{R}^p \times \textbf{R}^{m_2} \times \ldots \times \textbf{R}^{m_n} \times \textbf{R}^p$ (defined for $k \in K$, $i \in \{2 \ldots n\}$)
\item $f_k(\cdot,\cdot ,\ldots,\cdot)$ be convex on $\textbf{R}^p \times \textbf{R}^{m_1} \times \ldots \times \textbf{R}^{m_n}$ (defined for $k \in K$)
\item $g_k^i(\cdot,\cdot ,\ldots,\cdot,\cdot)$ be strictly convex on $\textbf{R}^p \times \textbf{R}^{m_1} \times \ldots \times \textbf{R}^{m_n} \times \textbf{R}^p$ (defined for $k \in K$, $i \in N$)
\item the cost functionals be stage-additive (cf. Def. \eqref{defcost}).
\end{itemize}
Then the set of strategies $\{\gamma_k^{i*}(x_0) = u^{i*}; i \in N \}$ is unique and provides an \emph {open-loop Stackelberg equilibrium solution} with \textbf{P}1 as the leader and \textbf{P}2 $\ldots$ \textbf{P}n as followers. Furthermore the corresponding state trajectory $\{x_k^*; k \in K\}$, the m-dimensional cocontrol vectors of the leader $\{\nu_1^i,\ldots,\nu_{T-1}^i$ ; $i \in \{2\ldots n\}\}$ and the p-dimensional costate vectors $\{\lambda_1, \ldots, \lambda_T, \mu_1^i, \ldots, \mu_T^i, p_1^{i*},\ldots,p_T^{i*}\}$ (defined for $i \in \{2\ldots n\}$) exist such that the following relations are satisfied:\\

\vspace*{-5ex}  \begin{align}
x_k^* = f_{k-1}(x_{k-1}^*,u_k^{1*},\ldots,u_k^{n*}) \; , \; x_0^* = x_0 \label{OLS_x}
\end{align}

\vspace*{-5ex}  \begin{multline}
\nabla_{u_k^1} H_k^1(\lambda_k, \mu_{k-1}^2, \ldots, \mu_{k-1}^n, \nu_{k-1}^2, \ldots, \nu_{k-1}^n, \\p_k^{2*}, \ldots, p_k^{n*}, u_k^{1*}, \ldots, u_k^{n*}, x_{k-1}^*) = 0 \\\label{guk1}
\end{multline} \vspace*{-5ex} 

\vspace*{-5ex}  \begin{multline}
\nabla_{u_k^i} H_k^1(\lambda_k, \mu_k^2, \ldots, \mu_{k-1}^2, \ldots, \mu_{k-1}^n, \nu_{k-1}^2, \ldots, \nu_{k-1}^n,\\ p_k^{2*}, \ldots, p_k^{n*}, u_k^{1*}, u_k^{2*}, \ldots, u_k^{n*}, x_{k-1}^*) = 0 \; ; \; i \in \{2\ldots n\} \\
\end{multline} \vspace*{-5ex} 

\vspace*{-5ex}  \begin{multline}
\lambda_{k-1}^{'} = \frac {\partial}{\partial x_{k-1}} H_k^1(\lambda_k, \mu_{k-1}^2, \ldots, \mu_{k-1}^n, \nu_{k-1}^2, \ldots, \nu_{k-1}^n,\\ p_k^{2*}, \ldots, p_k^{n*}, u_k^{1*}, \ldots, u_k^{n*}, x_{k-1}^*) \; ; \; \lambda_T^i = 0 \\
\end{multline} \vspace*{-5ex} 

\vspace*{-5ex}  \begin{multline}
\mu_{k}^{i'} = \frac {\partial}{\partial p_k} H_k^1(\lambda_k, \mu_{k-1}^2, \ldots, \mu_{k-1}^n, \nu_{k-1}^2, \ldots, \nu_{k-1}^n,\\ p_k^{2*}, \ldots, p_k^{n*}, u_k^{1*}, \ldots, u_k^{n*}, x_{k-1}^*)\; ; \; \mu_0^i = 0 \; ; \; i \in \{2\ldots n\}\\\label{gmuki}
\end{multline} \vspace*{-5ex} 

\vspace*{-5ex}  \begin{align}
\nabla_{u_k^i} H_k^i(p_k^{i*}, u_k^{1*}, \ldots, u_k^{n*}, x_{k-1}^*) = 0\; ; \; i \in \{2\ldots n\}
\end{align}

\vspace*{-5ex}  \begin{multline}
p_{k-1}^{i*} =  F_{k-1}^i(x_{k-1}^*,u_{k}^{1*},\ldots, u_{k}^{n*}, p_{k}^i)\; ; \; p_T^{i*} = 0 \; ; \; i \in \{2\ldots n\}\\
\end{multline} \vspace*{-5ex} 

where 

\vspace*{-5ex}  \begin{multline}
H_k^1(\lambda_k, \mu_{k-1}^2, \ldots, \mu_{k-1}^n, \nu_{k-1}^2, \ldots, \nu_{k-1}^n, p_k^2, \ldots, p_k^n, u_k^{1*}, \ldots, u_k^{n*}, x_{k-1}^*) \; \widehat{=} \;  \\g_k^1(f_k(x_{k-1},u_k^1, \ldots, u_k^n),u_k^1, \ldots, u_k^n, x_{k-1})+ \lambda_k^{'}f_k(x_{k-1},u_k^1, \ldots, u_k^n) +  \\\sum_{j \in \{2\ldots n\}}\mu_{k-1}^{j'}F_{k-1}^j(x_{k-1},u_{k}^1,\ldots, u_{k}^n, p_{k}^j) + \sum_{j \in \{2\ldots n\}}\nu_{k-1}^{j'}(\nabla_{u_k^j} H_k^j(p_k^j, u_k^{1*}, \ldots, u_k^{n*}, x_{k-1}^*)) \\
\end{multline} \vspace*{-5ex} 

\vspace*{-5ex}  \begin{multline}
F_k^i(x_{k},u_{k+1}^1,\ldots, u_{k+1}^n, p_{k+1}^i) \; \widehat{=} \; \frac {\partial}{\partial x_k} f_{k}(x_{k},u_{k+1}^{1},\ldots,u_{k+1}^{n})^{'}[p_{k+1}^i \\+ (\frac {\partial}{\partial x_{k+1}}g_{k+1}^i(x_{k+1}, u_{k+1}^{1}, \ldots, u_{k+1}^{n}, x_{k}))^{'}] + [\frac {\partial}{\partial x_{k}}g_{k+1}^i(x_{k+1}, u_{k+1}^{1}, \ldots, u_{k+1}^{n}, x_{k})]^{'} \; ; \; i \in \{2\ldots n\}\\
\end{multline} \vspace*{-5ex} 

\vspace*{-5ex}  \begin{multline}
H_k^i(p_k^i, u_k^{1}, \ldots, u_k^{n}, x_{k-1}) \; \widehat{=} \; g_k^i(f_{k-1}(x_{k-1},u_{k}^{1},\ldots,u_{k}^{n}), \\ u_{k}^{1},\ldots,u_{k}^{n},x_{k-1}) + p_{k}^{i'}f_{k-1}(x_{k-1},u_{k}^{1},\ldots,u_{k}^{n}) \; ; \; i \in \{2\ldots n\}\\\label{OLS_Hi}
\end{multline} \vspace*{-5ex} 
\label{OLS_opt}
\end{theo}

\textsc {Proof:}

Theorem \eqref{OLS_opt} can be proven in the same way as Theorem \eqref{OLN_opt}, bearing in mind that the leader additionally accounts for the influence of his strategy on the followers' strategies when minimizing his cost functional. Thus the minimization problem of the leader is equivalent to a finite dimensional nonlinear programming problem. The solution vectors of this problem are stated in \eqref{guk1} - \eqref{gmuki} and a derivation of this solution can be found e.g. in Canon et al. (1970, p. 51) \cite{can}. \;\;\;\;\;\;\;\;\;\;\;\;\;\;\;\;\boxed{}

\clearpage

\subsection{The "interwoven-inductions-results" for affine-quadratic games with one leader and arbitrarily many followers}

In the following, the results of Theorem \eqref{OLS_opt} are applied to an affine-quadratic dynamic game with one leader and arbitrarily many followers. Theorem \eqref{OLS2i} is a generalization of Corollary 7.1 in Ba\c{s}ar and Olsder (1999, pp. 371-372) \cite{baol}. On the one hand a more general state equation and more general cost functionals are considered and on the other hand the number of followers is extended from one to arbitrarily many.\\ 

The structure of the proof is geared to the one indicated in Ba\c{s}ar and Olsder (1999, p. 372) \cite{baol} where two induction arguments are interwoven. The induction for $\mu^i_k$ ($ i \in \{2 \ldots n\}$ , being the costate vectors of the leader, which are associated with the costate vectors of the followers) runs forward in time from \textit{k}=0 to \textit{k}=\textit{T}-1 and the induction for $p^i_k$ and $\lambda_k$ ($i \in \{2 \ldots n\}$ , being the costate vectors of the followers and the costate vectors of the leader) runs backward in time from \textit{k}=\textit{T} to \textit{k}=1. In the inductive step the induction hypotheses of the two inductions are used together.\\

This fact causes severe problems if we intend to use the obtained equilibrium equations for an algorithmic disintegration of the game because the evolution of $C^i_k$ and $c^i_k$ ($ i \in \{2 \ldots n\}$, $\mu_k^i = C^i_kx_k + c_k^i$) is dependent on $M_{k+1}^i$, $m_{k+1}^i$, $L_{k+1}$ and $l_{k+1}$ ($ i \in \{2 \ldots n\}$, $p_k^i = M^i_kx_k + m_k^i$, $\lambda_k^i = L_kx_k + l_k$) and the evolution of $M_{k+1}^i$, $m_{k+1}^i$, $L_{k+1}$ and $l_{k+1}$ is dependent on $C^i_k$ and $c^i_k$. But $C^i_k$, $c^i_k$, $M_{k+1}^i$, $m_{k+1}^i$, $L_{k+1}$ and $l_{k+1}$ are needed to determine $\Phi_k$ and $\phi_k$, which define the evolution of the state variable and the control variables. This means that to solve the game we have to determine $C^i_k$ (determining $c^i_k$ is even more complicated) as a function of $(L_1, \ldots, L_{k}, M_1^i, \ldots, M_{k}^i)$ for each \textit{k} $\in$ $\{2, \ldots, T\}$ (which is all but impossible even for a very small number of stages) and then substitute $C^i_k$ when determining $L_{k+1}$ and $M_{k+1}^i$ backwards in time. %This yields a system of \textit{n*(T-1)} associated matrix equations. Such a system is very hard to solve and furthermore general conditions for the solvability of this system cannot be stated because they depend on \textit{T}.\\

To overcome these intractabilities, another proof structure has to be used to solve the dynamic, affine-quadratic Stackelberg game. In subsection \eqref{OLS_1i} a theorem is stated which is proven by showing that $p^i_k$ and $\lambda_k$ ($\in \{2 \ldots n\}$) can be determined as functions linearly dependent on $x_k$ and $\mu_k^i$ ($\in \{2 \ldots n\}$). Solving the dynamic, affine-quadratic Stackelberg game in that way yields equilibrium equations that can easily be used for an algorithmic disintegration of that game.

\begin{theo}
An n-person affine-quadratic dynamic game (cf. Def. \eqref{defspielaq}) admits a \emph{unique open-loop Stackelberg equilibrium solution with one leader and arbitrarily many followers} if  
\begin{itemize}
\item $Q_k^i \geq 0$, $R_k^{ii} > 0$ (defined for $k \in K$ , $i \in N$).
\item  $(I- B_{k+1}^1W_{k+1}^{1} - \sum_{j\in \{2 \ldots n\}}B_{k+1}^jT_{k+1}^{j})^{-1}$ and $(A_{k+1} + B_{k+1}^1W_{k+1}^{0})^{-1}$ (defined for $k \in K$) exist.
\item \eqref{Nk1na}, \eqref{Nk0na} and \eqref{nkna} admit unique solutions $N_k^{1i}$, $N_k^{0i}$ and $n_k^{i}$ (defined for $k \in K$ , $i \in \{2,\ldots,n\}$)
\end{itemize}
If these conditions are satisfied, the unique equilibrium strategies are given by \eqref{gammaolsna}, where the associated state trajectory $x_{k+1}^*$ is given by \eqref{xoptsna}.\footnote{For all equations belonging to this theorem and its proof, $i \in N$ and $k \in \{0 , \ldots ,T-1\}$ if nothing different is stated.}

\vspace*{-5ex}  \begin{align}
f_{k-1}(x_{k-1},u_k^1,\ldots,u_k^n) = A_kx_{k-1} + \sum_{j\in N}B_k^ju_k^j + s_k \; ; \; k \in K \label{folsna}
\end{align}

\vspace*{-5ex}  \begin{align}
L^i(x_0,u^{1},\ldots,u^{n})= \sum_{k = 1}^T g_k^i(x_k,u_k^{1},\ldots,u_k^{n},x_{k-1}) \label{lolsna} 
\end{align}

\vspace*{-5ex}  \begin{multline}
 g_k^i(x_k,u_k^1,\ldots,u_k^n,x_{k-1}) = \frac{1}{2} (x_k^{'}Q_k^ix_k + \sum_{j\in N}u_k^{j'}R_k^{ij}u_k^j)\; \\ +
 \frac{1}{2} (\tilde{x}_k^{i'}Q_k^i\tilde{x}_k^i + \sum_{j\in N}\tilde{u}_k^{ij'}R_k^{ij}\tilde{u}_k^{ij}) - \tilde{x}_k^{i'}Q_k^ix_k - \sum_{j\in N}\tilde{u}_k^{ij'}R_k^{ij}u_k^j \; ; \; k \in K \\\label{golsna}
\end{multline} \vspace*{-5ex}

\vspace*{-5ex}  \begin{align}
 x_{k+1}^* = \Phi_{k}x_{k}^* + \phi_k \; ; \; x_0^* = x_0\label{xoptsna}
\end{align}

\vspace*{-5ex}  \begin{align}
\Phi_{k} = (I- B_{k+1}^1W_{k+1}^{1} - \sum_{j\in \{2 \ldots n\}}B_{k+1}^jT_{k+1}^{j})^{-1}(A_{k+1} + B_{k+1}^1W_{k+1}^{0})\label{Phikna}
\end{align}

\vspace*{-5ex}  \begin{multline}
\phi_{k} = (I- B_{k+1}^1W_{k+1}^{1} - \sum_{j\in \{2 \ldots n\}}B_{k+1}^jT_{k+1}^{j}){-1}\\(B_{k+1}^1w_{k+1} + \sum_{j\in \{2 \ldots n\}}B_{k+1}^jt_{k+1}^j + s_{k+1}) \\\label{phikna}
\end{multline} \vspace*{-5ex}

\vspace*{-5ex}  \begin{align}
W_{k+1}^{1} =   -(R_{k+1}^{11})^{-1} (B_{k+1}^{1'}L_{k+1} + \sum_{j \in 2\ldots n}B_{k+1}^{1'}Q_{k+1}^jB_{k+1}^jN_{k}^{1j}) \label{Wk1na}
\end{align}

\vspace*{-5ex}  \begin{multline}
W_{k+1}^{0} =  -(R_{k+1}^{11})^{-1} (\sum_{j \in 2\ldots n}B_{k+1}^{1'}Q_{k+1}^jA_{k+1}C_{k}^j  \\+ \sum_{j \in 2\ldots n}B_{k+1}^{1'}Q_{k+1}^jB_{k+1}^jN_{k}^{0j}) \\ \label{Wk0na}
\end{multline} \vspace*{-5ex} 

\vspace*{-5ex}  \begin{multline}
w_{k+1} =  -(R_{k+1}^{11})^{-1} (-B_{k+1}^{1'}Q_{k+1}^1 \tilde{x}_{k+1}^i  + B_{k+1}^{1'}l_{k+1} \\+ \sum_{j \in 2\ldots n}B_{k+1}^{1'}Q_{k+1}^jA_{k+1}c_{k}^j + \sum_{j \in 2\ldots n}B_{k+1}^{1'}Q_{k+1}^jB_{k+1}^jn_k^j) + \tilde{u}_{k+1}^{11}\\ \label{wkna}
\end{multline} \vspace*{-5ex}

\vspace*{-5ex}  \begin{align}
T_{k+1}^{i} = -(R_{k+1}^{ii})^{-1}B_{k+1}^{i'}M_{k+1}^{i}  \; ; \; i \in \{2\ldots n\}\label{Tkna}
\end{align}

\vspace*{-5ex}  \begin{align}
t_{k+1}^{i} = -(R_{k+1}^{ii})^{-1}B_{k+1}^{i'}(m_{k+1}^i - Q_{k+1}^i\tilde{x}_{k+1}^{i'}) + \tilde{u}_{k+1}^{ii}  \; ; \; i \in \{2\ldots n\}\label{tkna}
\end{align}

\vspace*{-5ex}  \begin{multline}
-R_{k+1}^{1i}(R_{k+1}^{ii})^{-1}B_{k+1}^{i'}M_{k+1}^{i} + B_{k+1}^{i'}L_{k+1} +  (B_{k+1}^{i'}Q_{k+1}^iB_{k+1}^i \\+ R_{k+1}^{ii})N_{k}^{1i} +\sum_{j \in 2\ldots n \;,\;j \not = i}B_{k+1}^{i'}Q_{k+1}^jB_{k+1}^jN_{k}^{1j} = 0  \; ; \; i \in \{2\ldots n\}\\ \label{Nk1na}
\end{multline} \vspace*{-5ex} 

\vspace*{-5ex}  \begin{multline}
\sum_{j \in 2\ldots n}B_{k+1}^{i'}Q_{k+1}^jA_{k+1}C_{k}^j +  (B_{k+1}^{i'}Q_{k+1}^iB_{k+1}^i + R_{k+1}^{ii})N_{k}^{0i} \\ + \sum_{j \in 2\ldots n \;,\;j \not = i}B_{k+1}^{i'}Q_{k+1}^jB_{k+1}^jN_{k}^{0j} = 0  \; ; \; i \in \{2\ldots n\}\\ \label{Nk0na}
\end{multline} \vspace*{-5ex} 

\vspace*{-5ex}  \begin{multline}
-B_{k+1}^{i'}Q_{k+1}^1\tilde{x}_{k+1}^{i*} + R_{k+1}^{1i}(-(R_{k+1}^{ii})^{-1}B_{k+1}^{i'}(m_{k+1}^i - Q_{k+1}^i\tilde{x}_{k+1}^{i'}) + \tilde{u}_{k+1}^{ii} \\- \tilde{u}_{k+1}^{1i}) + B_{k+1}^{i'}l_{k+1} + \sum_{j \in 2\ldots n}B_{k+1}^{i'}Q_{k+1}^jA_{k+1}c_{k}^j +  (B_{k+1}^{i'}Q_{k+1}^iB_{k+1}^i \\+ R_{k+1}^{ii})n_k^i + \sum_{j \in 2\ldots n \;,\;j \not = i}B_{k+1}^{i'}Q_{k+1}^jB_{k+1}^jn_k^j = 0  \; ; \; i \in \{2\ldots n\}\\ \label{nkna}
\end{multline} \vspace*{-5ex}

\vspace*{-5ex}  \begin{align}
M_{k}^{i} = Q_{k}^i + A_{k+1}^{'}M_{k+1}^{i}\Phi_{k}  \; ; \; M_T^i = Q_T^i\; ; \; i \in \{2\ldots n\} \label{Mikna}
\end{align}

\vspace*{-5ex}  \begin{align}
m_k^i = A_{k+1}^{'}[M_{k+1}^{i}\phi_k + m_{k+1}^i - Q_{k+1}^i\tilde{x}_{k+1}^i] \; ; \; m_T^i = 0\; ; \; i \in \{2\ldots n\} \label{mikna}
\end{align}

\vspace*{-5ex}  \begin{multline}
L_{k} = Q_{k}^1 + A_{k+1}^{'}L_{k+1}\Phi_{k} + \sum_{j \in 2\ldots n}A_{k+1}^{'}Q_{k+1}^jA_{k+1}C_{k}^j \\+ \sum_{j \in 2\ldots n}A_{k+1}^{'}Q_{k+1}^jB_{k+1}^j(N_{k}^{1j}\Phi_{k} + N_{k}^{0j}) \; ; \; L_T = Q_T^1\\\label{Lkna}
\end{multline} \vspace*{-5ex} 

\vspace*{-5ex}  \begin{multline}
l_k = A_{k+1}^{'}(L_{k+1}\phi_k + l_{k+1}) - A_{k+1}^{'}Q_{k+1}^1\tilde{x}_{k+1}^1+ \sum_{j \in 2\ldots n}A_{k+1}^{'}Q_{k+1}^jA_{k+1}c_{k}^j \\+ \sum_{j \in 2\ldots n}A_{k+1}^{'}Q_{k+1}^jB_{k+1}^j(N_{k}^{1j}\phi_k + n_k^j) \; ; \; l_T = 0\\  \label{lkna}
\end{multline} \vspace*{-5ex}

\vspace*{-5ex}  \begin{align}
C_{k+1}^{i} = A_{k+1}C_{k}^i\Phi_{k}^{-1} + B_{k+1}^i(N_{k}^{1i} + N_{k}^{0i}\Phi_{k}^{-1}) \; ; \; C_0^i = 0\; ; \; i \in \{2\ldots n\} \label{Cikna}
\end{align}

\vspace*{-5ex}  \begin{multline}
c_{k+1}^{i} = A_{k+1}(-C_{k}^i\Phi_{k}^{-1}\phi_k + c_{k}^i) \\+ B_{k+1}^i(-N_{k}^{0i}\Phi_{k}^{-1}\phi_k + n_k^i) \; ; \; c_0^i = 0\; ; \;i \in \{2\ldots n\} \\\label{cikna}
\end{multline} \vspace*{-5ex}

\vspace*{-5ex}  \begin{align}
\gamma_{k+1}^{i*}(x_0) = u_{k+1}^{i*} = P_{k+1}^{i}x_k^* + \alpha_{k+1}^i  \label{gammaolsna}
\end{align}

\vspace*{-5ex}  \begin{align}
P_{k+1}^{1} =  W_{k+1}^{1}\Phi_{k} + W_{k+1}^{0} \label{Pk1na}
\end{align}

\vspace*{-5ex}  \begin{align}
\alpha_{k+1}^{1} = W_{k+1}^{1}\phi_k + w_{k+1} \label{alphak1na}
\end{align}

\vspace*{-5ex}  \begin{align}
P_{k+1}^{i} = T_{k+1}^{i}\Phi_{k}   \; ; \; i \in \{2\ldots n\} \label{Pkina}
\end{align}

\vspace*{-5ex}  \begin{align}
\alpha_{k+1}^{i} = T_{k+1}^{i}\phi_k + t_{k+1}^i   \; ; \; i \in \{2\ldots n\} \label{alphakina}
\end{align}\\\\

\label{OLS2i}
\end{theo}

\textsc {Proof:}

Theorem \eqref{OLS_opt} can be applied to the given affine-quadratic game, since all conditions are satisfied for the given state equation \eqref{folsna} and cost functionals \eqref{lolsna}. Furthermore $g_k^i$ is strictly convex in $u_k^i$. This can be seen by applying Corollary \eqref{pdc} to \eqref{conolsna}. Therefore there has to be a unique optimal equilibrium solution.

\vspace*{-5ex}  \begin{multline}
\frac{\partial^2}{\partial u_k^{i^2}} g_k^i(x_k,u_k^1 ,\ldots , u_k^n, x_{k-1}) = B_k^{i'}Q_k^iB_k^i + R_k^{ii} \\\label{conolsna}
\end{multline} \vspace*{-5ex} 

To obtain relations which satisfy this unique solution we have to adapt \eqref{OLS_x} - \eqref{OLS_Hi} to the given state equation and cost functionals. This yields

\vspace*{-5ex}  \begin{multline}
H_k^1 =  \frac{1}{2} (x_k^{'}Q_k^1x_k + \sum_{j \in N}u_k^{j'}R_k^{1j}u_k^j)\;  +
 \frac{1}{2} (\tilde{x}_k^{1'}Q_k^1\tilde{x}_k^1 + \sum_{j \in N}\tilde{u}_k^{1j'}R_k^{1j}\tilde{u}_k^{1j}) \\- \tilde{x}_k^{1'}Q_k^1x_k - \sum_{j\in N}\tilde{u}_k^{1j'}R_k^{1j}u_k^j + \lambda_k^{'}(A_kx_{k-1} + \sum_{j \in N}B_k^ju_k^j+s_k) \\+ \sum_{j \in 2\ldots n}\mu_{k-1}^{j'}A_{k}^{'}[p_{k}^i + Q_{k}^j( x_{k}^* - \tilde{x}_{k}^j )] \\+ \sum_{j \in 2\ldots n}\nu_{k-1}^{j'}(B_k^{j'}Q_k^jx_k + R_k^{jj}u_k^{j} -  B_k^{j'}Q_k^j\tilde{x}_k^{j'} - R_k^{jj}\tilde{u}_k^{jj} + B_k^{j'}p_{k}^j)  \\\label{H1ksna}
\end{multline} \vspace*{-5ex} 

\vspace*{-5ex}  \begin{multline}
H_k^i =   \frac{1}{2} (x_k^{'}Q_k^ix_k + \sum_{j\in N}u_k^{j'}R_k^{ij}u_k^j)\;  +
 \frac{1}{2} (\tilde{x}_k^{i'}Q_k^i\tilde{x}_k^i + \sum_{j\in N}\tilde{u}_k^{ij'}R_k^{ij}\tilde{u}_k^{ij}) \\- \tilde{x}_k^{i'}Q_k^ix_k - \sum_{j\in N}\tilde{u}_k^{ij'}R_k^{ij}u_k^j + p_k^{i'}(A_kx_{k-1} + \sum_{j\in N}B_k^ju_k^j+s_k) \; ; \; i \in \{2\ldots n\}\\\label{Hiksna}
\end{multline} \vspace*{-5ex} 

\vspace*{-5ex}  \begin{multline}
\nonumber \frac{\partial}{\partial u_k^1} \eqref{H1ksna} = 0\; \Rightarrow \;\;\ B_k^{1'}Q_k^1(x_k^*-\tilde{x}_k^i) + R_k^{11}(u_k^{1*}- \tilde{u}_k^{11}) + B_k^{1'}\lambda_k \\+ \sum_{j \in 2\ldots n}B_k^{1'}Q_k^jA_k\mu_{k-1}^{j} + \sum_{j \in 2\ldots n}B_k^{1'}Q_k^jB_k^j\nu_{k-1}^{j}  = 0\\
\end{multline} \vspace*{-5ex} 

\vspace*{-5ex}  \begin{multline}
u_k^{1*} = -(R_k^{11})^{-1}(B_k^{1'}Q_k^1(x_k^*-\tilde{x}_k^i)  + B_k^{1'}\lambda_k \\+ \sum_{j \in 2\ldots n}B_k^{1'}Q_k^jA_k\mu_{k-1}^{j} + \sum_{j \in 2\ldots n}B_k^{1'}Q_k^jB_k^j\nu_{k-1}^{j}) + \tilde{u}_k^{11}\\\label{uk1sna}
\end{multline} \vspace*{-5ex} 

\vspace*{-5ex}  \begin{multline}
\frac{\partial}{\partial u_k^i} \eqref{H1ksna} = 0\; \Rightarrow \;\; B_k^{i'}Q_k^1(x_k^*-\tilde{x}_k^i) + R_k^{1i}(u_k^{i*}- \tilde{u}_k^{1i}) + B_k^{i'}\lambda_k \\+ \sum_{j \in 2\ldots n}B_k^{i'}Q_k^jA_k\mu_{k-1}^{j} +  (B_k^{i'}Q_k^iB_k^i + R_k^{ii})\nu_{k1}^{i} +\sum_{j \in 2\ldots n \;,\;j \not = i}B_k^{i'}Q_k^jB_k^j\nu_{k-1}^{j} = 0 \; ; \; i \in \{2\ldots n\}\\\label{nuksna}
\end{multline} \vspace*{-5ex} 

\vspace*{-5ex}  \begin{multline}
\lambda_{k-1} = A_k^{'}Q_k^1(x_k^*-\tilde{x}_k^1) + A_k^{'}\lambda_k \\+ \sum_{j \in 2\ldots n}A_k^{'}Q_k^jA_k\mu_{k-1}^{j} + \sum_{j \in 2\ldots n}A_k^{'}Q_k^jB_k^j\nu_{k-1}^{j} \; ; \; \lambda_T = 0\\\label{lambdaksna}
\end{multline} \vspace*{-5ex} 

\vspace*{-5ex}  \begin{multline}
\mu_k^{i} =  A_k\mu_{k-1}^{i} + B_k^i\nu_{k-1}^{i}\; ; \; \mu_0^i = 0\; ; \; i \in \{2\ldots n\} \\ \label{muksna}
\end{multline} \vspace*{-5ex} 

\vspace*{-5ex}  \begin{multline}
\nonumber \frac{\partial}{\partial u_i} \eqref{Hiksna} = 0\; \Rightarrow \;\;\ B_k^{i'}Q_k^ix_k^* + R_k^{ii}u_k^{i*} -  B_k^{i'}Q_k^i\tilde{x}_k^{i'} - R_k^{ii}\tilde{u}_k^{ii} + B_k^{i'}p_{k}^i = 0 \; ; \; i \in \{2\ldots n\}\\
\end{multline} \vspace*{-5ex} 

\vspace*{-5ex}  \begin{multline}
u_k^{i*} =  -(R_k^{ii})^{-1}B_k^{i'}(Q_k^i(x_k^*-\tilde{x}_k^{i'}) + p_{k}^i ) + \tilde{u}_k^{ii} \; ; \; i \in \{2\ldots n\}\\\label{ukisna}
\end{multline} \vspace*{-5ex} 

\vspace*{-5ex}  \begin{multline}
p_{k-1}^{i*} = A_{k}^{'}[p_{k}^i + Q_{k}^i( x_{k}^* - \tilde{x}_{k}^i )] \; ; \; p_T^i = 0\; ; \; i \in \{2\ldots n\}\\\label{pkisna}
\end{multline} \vspace*{-5ex} 

\vspace*{-5ex}  \begin{multline}
x_{k}^* = A_kx_{k-1}^* + \sum_{j\in N}B_k^ju_k^{j*}+s_k \; ; \; x_0^* = x_0\\\label{xksna}
\end{multline} \vspace*{-5ex}

In the following induction arguments we will give proof that \eqref{psubsna} - \eqref{musubsna} are valid and the recursive relations for $C^i_{k+1}$, $c^i_{k+1}$, $M_{k}^i$, $m_{k}^i$, $L_{k}$ and $l_{k}$ (stated in the above theorem) are correct.

\vspace*{-5ex}  \begin{multline}
p_{k}^i = (M_{k}^{i}-Q_{k}^i)x_{k}^* + m_k^i  \; ; \; i \in \{2\ldots n\} \\\label{psubsna}
\end{multline} \vspace*{-5ex} 

\vspace*{-5ex}  \begin{multline}
\lambda_{k} = (L_{k}-Q_{k}^1)x_{k}^* + l_k \\\label{lambdasubsna}
\end{multline} \vspace*{-5ex} 

\vspace*{-5ex}  \begin{multline}
\mu_{k}^i = C_{k}x_{k}^* + c_k \; ; \; i \in \{2\ldots n\} \\\label{musubsna}
\end{multline} \vspace*{-5ex}

\textbf{Basis:}

The induction for $p_k^i$ and $\lambda_k$ starts at \textit{k} = \textit{T}. First we make use of the general optimality conditions at stage \textit{T}.

\vspace*{-5ex}  \begin{multline}
\nonumber \eqref{pkisna}_{k=T} \;\; p_{T}^i = 0 \; ; \; i \in \{2\ldots n\}\\
\end{multline} \vspace*{-5ex} 

\vspace*{-5ex}  \begin{multline}
\nonumber \eqref{lambdaksna}_{k=T} \;\; \lambda_{T} = 0\\
\end{multline} \vspace*{-5ex} 

Now we can substitute $p_T^i$ and $\lambda_T$ with functions affinely dependent on $x_T$.

\vspace*{-5ex}  \begin{multline}
\nonumber\eqref{psubsna}_{k=T} \;\; p_{T}^i = (M_{T}^{i}-Q_{T}^i)x_{T}^* + m_T^i \; ; \; i \in \{2\ldots n\}\\
\end{multline} \vspace*{-5ex} 

\vspace*{-5ex}  \begin{multline}
\nonumber (M_{T}^{i}-Q_{T}^i)x_{T}^* + m_T^i  = 0 \; ; \; i \in \{2\ldots n\}\\
\end{multline} \vspace*{-5ex} 

\vspace*{-5ex}  \begin{multline}
\nonumber\eqref{lambdasubsna}_{k=T} \;\; \lambda_{T} = (L_{T}-Q_{T}^1)x_{T}^* + l_T\\
\end{multline} \vspace*{-5ex} 

\vspace*{-5ex}  \begin{multline}
\nonumber (L_{T}-Q_{T}^1)x_{T}^* + l_T  = 0 \; ; \; i \in \{2\ldots n\}\\
\end{multline} \vspace*{-5ex}

Comparing coefficients gives

\vspace*{-5ex}  \begin{multline}
\nonumber\eqref{Mikna}_{k=T} \;\; M_T^{ix} = Q_{T}^i \; ; \; i \in \{2\ldots n\}\\
\end{multline} \vspace*{-5ex} 

\vspace*{-5ex}  \begin{multline}
\nonumber\eqref{mikna}_{k=T} \;\; m_T^{i} = 0 \; ; \; i \in \{2\ldots n\}\\
\end{multline} \vspace*{-5ex}

\vspace*{-5ex}  \begin{multline}
\nonumber\eqref{Lkna}_{k=T} \;\; L_T = Q_{T}^1 \\
\end{multline} \vspace*{-5ex}

\vspace*{-5ex}  \begin{multline}
\nonumber\eqref{lkna}_{k=T} \;\;  l_T = 0 \\
\end{multline} \vspace*{-5ex}

The induction for $\mu_k^i$ starts at \textit{k} = 0. First we make use of the general optimality conditions at stage 0.

\vspace*{-5ex}  \begin{multline}
\nonumber \eqref{muksna}_{k=0} \;\; \mu_{0}^i = 0 \; ; \; i \in \{2\ldots n\}\\
\end{multline} \vspace*{-5ex} 

Now we can substitute the $\mu_k^i$ with functions linearly dependent on $x_0$.

\vspace*{-5ex}  \begin{multline}
\nonumber\eqref{musubsna}_{k=0} \;\; \mu_{0}^i = C_{0}^{i}x_{0}^* + c_0^i \; ; \; i \in \{2\ldots n\}\\
\end{multline} \vspace*{-5ex} 

\vspace*{-5ex}  \begin{multline}
\nonumber C_{0}^{i}x_{0}^* + c_0^i  = 0 \; ; \; i \in \{2\ldots n\}\\
\end{multline} \vspace*{-5ex} 

Comparing coefficients gives

\vspace*{-5ex}  \begin{multline}
\nonumber\eqref{Cikna}_{k=0} \;\; C_0^i = 0 \; ; \; i \in \{2\ldots n\}\\
\end{multline} \vspace*{-5ex}

\vspace*{-5ex}  \begin{multline}
\nonumber\eqref{cikna}_{k=0} \;\; c_0^{i} = 0 \; ; \; i \in \{2\ldots n\}\\
\end{multline} \vspace*{-5ex} \\\\

\textbf{The interwoven inductive steps:}

To prove the statements given in the above theorem, two inductive steps are combined. \\

As an induction hypothesis, the system of equations \eqref{psubsna} and equation \eqref{lambdasubsna} are assumed to be true at stage \textit{l}+2 and the system of equations \eqref{musubsna} is assumed to be true at stage \textit{l}. Now we have to show that these equations are fulfilled at stage \textit{l}+1 and determine the corresponding recursive relations for $C^i_{l+1}$, $c^i_{l+1}$, $M_{l}^i$, $m_{l}^i$, $L_{l}$ and $l_{l}$.

\vspace*{-5ex}  \begin{multline}
\nonumber\eqref{psubsna}_{k=l+2} \;\; p_{l+1}^i = (M_{l+1}^{i}-Q_{l+1}^i)x_{l+1}^* + m_{l+1}^i \; ; \; i \in \{2\ldots n\}\\
\end{multline} \vspace*{-5ex} 

\vspace*{-5ex}  \begin{multline}
\nonumber\eqref{lambdasubsna}_{k=l+2} \;\; \lambda_{l+1} = (L_{l+1}-Q_{l+1}^1)x_{l+1}^* + l_{l+1}\\
\end{multline} \vspace*{-5ex} 

\vspace*{-5ex}  \begin{multline}
\nonumber\eqref{musubsna}_{k=l} \;\; \mu_{l}^i = C_{l}^ix_{l}^* + c_{l}^i \; ; \; i \in \{2\ldots n\}  \\
\end{multline} \vspace*{-5ex}

First these induction hypotheses are used in the general optimality conditions for $p_k^i$, $\lambda_k$ and $\mu^i_k$ at stage \textit{l}+1.

\vspace*{-5ex}  \begin{multline}
\nonumber\eqref{pkisna}_{k=l+1} \;\; p_l^{i*} = A_{l+1}^{'}[p_{l+1}^i + Q_{l+1}^i( x_{l+1}^* - \tilde{x}_{l+1}^i )] \; ; \; i \in \{2\ldots n\}\\
\end{multline} \vspace*{-5ex}

\vspace*{-5ex}  \begin{multline}
\nonumber  p_l^{i*} = A_{l+1}^{'}[M_{l+1}^{i}x_{l+1}^* + m_{l+1}^i - Q_{l+1}^i\tilde{x}_{l+1}^i]  \; ; \; i \in \{2\ldots n\}\\
\end{multline} \vspace*{-5ex}

\vspace*{-5ex}  \begin{multline}
\nonumber\eqref{lambdaksna}_{k=l+1} \;\; \lambda_{l} = A_{l+1}^{'}Q_{l+1}^1(x_{l+1}^*-\tilde{x}_{l+1}^1) + A_{l+1}^{'}\lambda_{l+1} \\+ \sum_{j \in 2\ldots n}A_{l+1}^{'}Q_{l+1}^jA_{l+1}\mu_{l}^{j} + \sum_{j \in 2\ldots n}A_{l+1}^{'}Q_{l+1}^jB_{l+1}^j\nu_{l}^{j} \\
\end{multline} \vspace*{-5ex} 

\vspace*{-5ex}  \begin{multline}
\nonumber  \lambda_{l} = A_{l+1}^{'}(L_{l+1}x_{l+1}^* + l_{l+1}) - A_{l+1}^{'}Q_{l+1}^1\tilde{x}_{l+1}^1\\+ \sum_{j \in 2\ldots n}A_{l+1}^{'}Q_{l+1}^jA_{l+1}\mu_{l}^{j} + \sum_{j \in 2\ldots n}A_{l+1}^{'}Q_{l+1}^jB_{l+1}^j\nu_{l}^{j} \\
\end{multline} \vspace*{-5ex} 

\vspace*{-5ex}  \begin{multline}
\nonumber\eqref{muksna}_{k=l+1} \;\; \mu_{l+1}^{i} =  A_{l+1}\mu_{l}^{i} + B_{l+1}^i\nu_{l}^{i} \\
\end{multline} \vspace*{-5ex} 

\vspace*{-5ex}  \begin{multline}
\nonumber \mu_{l+1}^{i} =  A_{l+1}(C_{l}^ix_{l}^* + c_{l}^i) + B_{l+1}^i\nu_{l}^{i} \\
\end{multline} \vspace*{-5ex}

To complete the inductive step we have to prove that the $p_l^i$ and $\lambda_l$ can be written as affine functions of $x_l$ and that the $\mu_{l+1}^i$ can be written as affine functions of $x_{l+1}$. Therefore an interrelation between $x_l$ and $x_{l+1}$ which does not depend on the controls of the players nor on costate ($p_{l+1}^i$, $\lambda_{l+1}^i$, $\mu_l^i$) or cocontrol ($\nu_l^i$) variables has to be deduced. To do so, first we have to substitute $u_{l+1}^{1*}, \ldots,u_{l+1}^{n*}$ in the equation stated below for the evolution of the optimal state vector $x_{l+1}^*$ by terms which are affine in  $x_l$ and $x_{l+1}$ and furthermore only contain $C^i_{l}$, $c^i_{l}$, $M_{l+1}^i$, $m_{l+1}^i$, $L_{l+1}$, $l_{l+1}$ and matrices and vectors given by the game definition.

\vspace*{-5ex}  \begin{multline}
\nonumber\eqref{xksna}_{k=l+1} \;\; x_{l+1}^* = A_{l+1}x_{l}^* + \sum_{j\in N}B_{l+1}^ju_{l+1}^{j*}+s_{l+1} \\
\end{multline} \vspace*{-5ex}

In the first instance the optimality condition for $u_{l+1}^{i*}$ (defined for $i \in \{2\ldots n\}$) can be rewritten with the help of $\eqref{psubsna}_{k=l+1}$, which are the induction hypotheses for $p_{l+1}^i$. 

\vspace*{-5ex}  \begin{multline}
\nonumber\eqref{ukisna}_{k=l+1} \;\; u_{l+1}^{i*} =   -(R_{l+1}^{ii})^{-1}B_{l+1}^{i'}(Q_{l+1}^i(x_{l+1}^*-\tilde{x}_{l+1}^{i'}) + p_{l+1}^i ) + \tilde{u}_{l+1}^{ii}  \; ; \; i \in \{2\ldots n\}\\
\end{multline} \vspace*{-5ex} 

\vspace*{-5ex}  \begin{multline}
\nonumber u_{l+1}^{i*} =   -(R_{l+1}^{ii})^{-1}B_{l+1}^{i'}(M_{l+1}^{i}x_{l+1}^* + m_{l+1}^i - Q_{l+1}^i\tilde{x}_{l+1}^{i'}) + \tilde{u}_{l+1}^{ii}  \; ; \; i \in \{2\ldots n\}\\
\end{multline} \vspace*{-5ex} 

The structure of the above equations justifies the following substitutions

\vspace*{-5ex}  \begin{multline}
u_{l+1}^{i*} =  T_{l+1}^{i}x_{l+1}^* + t_{l+1}^i   \; ; \; i \in \{2\ldots n\}\\\label{uisubsna}
\end{multline} \vspace*{-5ex} 

\vspace*{-5ex}  \begin{multline}
T_{l+1}^{i}x_{l+1}^* + t_{l+1}^i =   -(R_{l+1}^{ii})^{-1}B_{l+1}^{i'}(M_{l+1}^{i}x_{l+1}^* + \\m_{l+1}^i - Q_{l+1}^i\tilde{x}_{l+1}^{i'}) + \tilde{u}_{l+1}^{ii}  \; ; \; i \in \{2\ldots n\}\\ \label{uilgna}
\end{multline} \vspace*{-5ex} 

By comparing coefficients it follows that

\vspace*{-5ex}  \begin{multline}
\nonumber\eqref{uilgna}_{x_{l+1}^*}=\eqref{Tkna}_{k=l} \;\;  T_{l+1}^{i} = -(R_{l+1}^{ii})^{-1}B_{l+1}^{i'}M_{l+1}^{i}  \; ; \; i \in \{2\ldots n\}\\
\end{multline} \vspace*{-5ex} 

\vspace*{-5ex}  \begin{multline}
\nonumber\eqref{uilgna}_{const.}=\eqref{tkna}_{k=l} \;\; t_{l+1}^{i} = -(R_{l+1}^{ii})^{-1}B_{l+1}^{i'}(m_{l+1}^i - Q_{l+1}^i\tilde{x}_{l+1}^{i'}) + \tilde{u}_{l+1}^{ii}  \; ; \; i \in \{2\ldots n\}\\
\end{multline} \vspace*{-5ex}

The optimality condition for $u_{l+1}^{1*}$ can be rewritten with the help of $\eqref{lambdasubsna}_{k=l+1}$ and $\eqref{musubsna}_{k=l}$, which are the induction hypotheses for $\lambda_{l+1}$ and $\mu_l^i$. 

\vspace*{-5ex}  \begin{multline}
\nonumber\eqref{uk1sna}_{k=l+1} \;\; u_{l+1}^{1*} =  -(R_{l+1}^{11})^{-1} (B_{l+1}^{1'}Q_{l+1}^1 (x_{l+1}^*-\tilde{x}_{l+1}^i)  + B_{l+1}^{1'}\lambda_{l+1} \\+ \sum_{j \in 2\ldots n}B_{l+1}^{1'}Q_{l+1}^jA_{l+1} \mu_{l}^{j} + \sum_{j \in 2\ldots n}B_{l+1}^{1'}Q_{l+1}^jB_{l+1}^j\nu_{l}^{j}) + \tilde{u}_{l+1}^{11} \\
\end{multline} \vspace*{-5ex} 

\vspace*{-5ex}  \begin{multline}
u_{l+1}^{1*} =  -(R_{l+1}^{11})^{-1} (-B_{l+1}^{1'}Q_{l+1}^1 \tilde{x}_{l+1}^i  + B_{l+1}^{1'}(L_{l+1}x_{l+1}^* + l_{l+1}) \\+ \sum_{j \in 2\ldots n}B_{l+1}^{1'}Q_{l+1}^jA_{l+1}(C_{l}^jx_{l}^* + c_{l}^j) + \sum_{j \in 2\ldots n}B_{l+1}^{1'}Q_{l+1}^jB_{l+1}^j\nu_{l}^{j}) + \tilde{u}_{l+1}^{11} \\\label{u1vor}
\end{multline} \vspace*{-5ex}

To make $u_{l+1}^{1*}$ only dependent on $C^i_{l}$, $c^i_{l}$, $M_{l+1}^i$, $m_{l+1}^i$, $L_{l+1}$, $l_{l+1}$ and matrices and vectors given by the game definition, we also have to substitute $\nu_{l}^{j}$ ($j \in \{2 \ldots n\}$) by a term affine in $x_l$ and $x_{l+1}$. For that purpose $\nu_{l}^{j}$ has to be explicated from $\eqref{nuksna}_{l+1}$, because this is the only optimality condition at stage \textit{l}+1 which has not been used sofar. As a start $\lambda_{l+1}$, $\mu_l^i$ and $u_{l+1}^{i*}$ ($i \in \{2\ldots n\}$) are substituted with the help of $\eqref{lambdasubsna}_{k=l+1}$, $\eqref{musubsna}_{k=l}$ and \eqref{uisubsna}.

\vspace*{-5ex}  \begin{multline}
\nonumber\eqref{nuksna}_{k=l+1} \;\; B_{l+1}^{i'}Q_{l+1}^1(x_{l+1}^*-\tilde{x}_{l+1}^{i*}) + R_{l+1}^{1i}(u_{l+1}^{i*} - \tilde{u}_{l+1}^{1i}) + B_{l+1}^{i'}\lambda_{l+1} + \sum_{j \in 2\ldots n}B_{l+1}^{i'}Q_{l+1}^jA_{l+1}\mu_{l}^{j} \\+  (B_{l+1}^{i'}Q_{l+1}^iB_{l+1}^i + R_{l+1}^{ii})\nu_{l}^{i} +\sum_{j \in 2\ldots n \;,\;j \not = i}B_{l+1}^{i'}Q_{l+1}^jB_{l+1}^j\nu_{l}^{j} = 0  \; ; \; i \in \{2\ldots n\}\\
\end{multline} \vspace*{-5ex} 

\vspace*{-5ex}  \begin{multline}
\nonumber -B_{l+1}^{i'}Q_{l+1}^1\tilde{x}_{l+1}^{i*} + R_{l+1}^{1i}(-(R_{l+1}^{ii})^{-1}B_{l+1}^{i'}(M_{l+1}^{i}x_{l+1}^* + m_{l+1}^i - Q_{l+1}^i\tilde{x}_{l+1}^{i'}) \\+ \tilde{u}_{l+1}^{ii} - \tilde{u}_{l+1}^{1i}) + B_{l+1}^{i'}(L_{l+1}x_{l+1}^* + l_{l+1}) + \sum_{j \in 2\ldots n}B_{l+1}^{i'}Q_{l+1}^jA_{l+1}(C_{l}^jx_{l}^* + c_{l}^j) \\+  (B_{l+1}^{i'}Q_{l+1}^iB_{l+1}^i + R_{l+1}^{ii})\nu_{l}^{i} +\sum_{j \in 2\ldots n \;,\;j \not = i}B_{l+1}^{i'}Q_{l+1}^jB_{l+1}^j\nu_{l}^{j} = 0  \; ; \; i \in \{2\ldots n\}\\
\end{multline} \vspace*{-5ex}

The above equations only contain constant expressions and terms linear in $x_l$ or $x_{l+1}$. This fact justifies the following substitutions.

\vspace*{-5ex}  \begin{multline}
\nu_{l}^i = N_{l}^{1i}x_{l+1}^* + N_{l}^{0i}x_{l}^* + n_l^i   \; ; \; i \in \{2\ldots n\}\\\label{nusubsna}
\end{multline} \vspace*{-5ex} 

\vspace*{-5ex}  \begin{multline}
-B_{l+1}^{i'}Q_{l+1}^1\tilde{x}_{l+1}^{i*} + R_{l+1}^{1i}(-(R_{l+1}^{ii})^{-1}B_{l+1}^{i'}(M_{l+1}^{i}x_{l+1}^* + m_{l+1}^i - Q_{l+1}^i\tilde{x}_{l+1}^{i'}) \\+ \tilde{u}_{l+1}^{ii} - \tilde{u}_{l+1}^{1i}) + B_{l+1}^{i'}(L_{l+1}x_{l+1}^* + l_{l+1}) + \sum_{j \in 2\ldots n}B_{l+1}^{i'}Q_{l+1}^jA_{l+1}(C_{l}^jx_{l}^* + c_{l}^j) \\+  (B_{l+1}^{i'}Q_{l+1}^iB_{l+1}^i + R_{l+1}^{ii})(N_{l}^{1i}x_{l+1}^* + N_{l}^{0i}x_{l}^* + n_l^i) \\+ \sum_{j \in 2\ldots n \;,\;j \not = i}B_{l+1}^{i'}Q_{l+1}^jB_{l+1}^j(N_{l}^{1j}x_{l+1}^* + N_{l}^{0j}x_{l+1}^* + n_l^j) = 0  \; ; \; i \in \{2\ldots n\}\\  \label{nulgna}
\end{multline} \vspace*{-5ex} 

Comparing coefficients gives the following three systems of equations that admit unique solutions $N_{l}^{1i}$, $N_{l}^{0i}$ and $n_{l}^{i}$ ($i \in \{2,\ldots,n\}$) per assumption.

\vspace*{-5ex}  \begin{multline}
\nonumber\eqref{nulgna}_{x_{l+1}^*}=\eqref{Nk1na}_{k=l} \;\;  -R_{l+1}^{1i}(R_{l+1}^{ii})^{-1}B_{l+1}^{i'}M_{l+1}^{i} + B_{l+1}^{i'}L_{l+1}\\+  (B_{l+1}^{i'}Q_{l+1}^iB_{l+1}^i + R_{l+1}^{ii})N_{l}^{1i} +\sum_{j \in 2\ldots n \;,\;j \not = i}B_{l+1}^{i'}Q_{l+1}^jB_{l+1}^jN_{l}^{1j} = 0  \; ; \; i \in \{2\ldots n\}\\ 
\end{multline} \vspace*{-5ex} 

\vspace*{-5ex}  \begin{multline}
\nonumber\eqref{nulgna}_{x_{l}^*}=\eqref{Nk0na}_{k=l} \;\;  \sum_{j \in 2\ldots n}B_{l+1}^{i'}Q_{l+1}^jA_{l+1}C_{l}^j \\+  (B_{l+1}^{i'}Q_{l+1}^iB_{l+1}^i + R_{l+1}^{ii})N_{l}^{0i} +\sum_{j \in 2\ldots n \;,\;j \not = i}B_{l+1}^{i'}Q_{l+1}^jB_{l+1}^jN_{l}^{0j} = 0  \; ; \; i \in \{2\ldots n\}\\ 
\end{multline} \vspace*{-5ex} 

\vspace*{-5ex}  \begin{multline}
\nonumber\eqref{nulgna}_{const.}=\eqref{nkna}_{k=l} \;\;  -B_{l+1}^{i'}Q_{l+1}^1\tilde{x}_{l+1}^{i*} + R_{l+1}^{1i}(-(R_{l+1}^{ii})^{-1}B_{l+1}^{i'}(m_{l+1}^i - Q_{l+1}^i\tilde{x}_{l+1}^{i'}) \\+ \tilde{u}_{l+1}^{ii} - \tilde{u}_{l+1}^{1i}) + B_{l+1}^{i'}l_{l+1} + \sum_{j \in 2\ldots n}B_{l+1}^{i'}Q_{l+1}^jA_{l+1}c_{l}^j \\+  (B_{l+1}^{i'}Q_{l+1}^iB_{l+1}^i + R_{l+1}^{ii})n_l^i + \sum_{j \in 2\ldots n \;,\;j \not = i}B_{l+1}^{i'}Q_{l+1}^jB_{l+1}^jn_l^j = 0  \; ; \; i \in \{2\ldots n\}\\ 
\end{multline} \vspace*{-5ex}

Using the elaborated relation for $\nu_{l}^{i}$ ($i \in \{2 \ldots n\}$) in \eqref{u1vor} yields

\vspace*{-5ex}  \begin{multline}
\nonumber u_{l+1}^{1*} =  -(R_{l+1}^{11})^{-1} (-B_{l+1}^{1'}Q_{l+1}^1 \tilde{x}_{l+1}^i  + B_{l+1}^{1'}(L_{l+1}x_{l+1}^* + l_{l+1}) \\+ \sum_{j \in 2\ldots n}B_{l+1}^{1'}Q_{l+1}^jA_{l+1}(C_{l}^jx_{l}^* + c_{l}^j) + \sum_{j \in 2\ldots n}B_{l+1}^{1'}Q_{l+1}^jB_{l+1}^j(N_{l}^{1j}x_{l+1}^* + N_{l}^{0j}x_{l}^* + n_l^j)) + \tilde{u}_{l+1}^{11} \\
\end{multline} \vspace*{-5ex} 

The structure of the above equation justifies the following substitution

\vspace*{-5ex}  \begin{multline}
u_{l+1}^{1*} =  W_{l+1}^{1}x_{l+1}^* + W_{l+1}^{0}x_{l}^* + w_{l+1} \\\label{u1subsna}
\end{multline} \vspace*{-5ex} 

\vspace*{-5ex}  \begin{multline}
W_{l+1}^{1}x_{l+1}^* + W_{l+1}^{0}x_{l}^* + w_{l+1} =  -(R_{l+1}^{11})^{-1} (-B_{l+1}^{1'}Q_{l+1}^1 \tilde{x}_{l+1}^i  \\+ B_{l+1}^{1'}(L_{l+1}x_{l+1}^* + l_{l+1}) + \sum_{j \in 2\ldots n}B_{l+1}^{1'}Q_{l+1}^jA_{l+1}(C_{l}^jx_{l}^* + c_{l}^j) \\+ \sum_{j \in 2\ldots n}B_{l+1}^{1'}Q_{l+1}^jB_{l+1}^j(N_{l}^{1j}x_{l+1}^* + N_{l}^{0j}x_{l}^* + n_l^j)) + \tilde{u}_{l+1}^{11} \\\label{u1lgna}
\end{multline} \vspace*{-5ex} 

By comparing coefficients it follows that

\vspace*{-5ex}  \begin{multline}
\nonumber\eqref{u1lgna}_{x_{l+1}^*}=\eqref{Wk1na}_{k=l} \;\;  W_{l+1}^{1} =   -(R_{l+1}^{11})^{-1} (B_{l+1}^{1'}L_{l+1} + \sum_{j \in 2\ldots n}B_{l+1}^{1'}Q_{l+1}^jB_{l+1}^jN_{l}^{1j})\\
\end{multline} \vspace*{-5ex} 

\vspace*{-5ex}  \begin{multline}
\nonumber\eqref{u1lgna}_{x_{l}^*}=\eqref{Wk0na}_{k=l} \;\;  W_{l+1}^{0} =  -(R_{l+1}^{11})^{-1} (\sum_{j \in 2\ldots n}B_{l+1}^{1'}Q_{l+1}^jA_{l+1}C_{l}^j + \sum_{j \in 2\ldots n}B_{l+1}^{1'}Q_{l+1}^jB_{l+1}^jN_{l}^{0j}) \\
\end{multline} \vspace*{-5ex} 

\vspace*{-5ex}  \begin{multline}
\nonumber\eqref{u1lgna}_{const.}=\eqref{wkna}_{k=l} \;\; w_{l+1} =  -(R_{l+1}^{11})^{-1} (-B_{l+1}^{1'}Q_{l+1}^1 \tilde{x}_{l+1}^i  + B_{l+1}^{1'}l_{l+1} \\+ \sum_{j \in 2\ldots n}B_{l+1}^{1'}Q_{l+1}^jA_{l+1}c_{l}^j + \sum_{j \in 2\ldots n}B_{l+1}^{1'}Q_{l+1}^jB_{l+1}^jn_l^j) + \tilde{u}_{l+1}^{11}\\
\end{multline} \vspace*{-5ex}

Now the control variables can be replaced in the optimal state equation by terms affine in $x_l$ and $x_{l+1}$.

\vspace*{-5ex}  \begin{multline}
\nonumber\eqref{xksna}_{k=l+1} \;\; x_{l+1}^* = A_{l+1}x_{l}^* + B_{l+1}^1u_{l+1}^{1*} + \sum_{j\in \{2 \ldots n\}}B_{l+1}^ju_{l+1}^{j*}+s_{l+1} \\
\end{multline} \vspace*{-5ex} 

\vspace*{-5ex}  \begin{multline}
\nonumber x_{l+1}^* = A_{l+1}x_{l}^* + B_{l+1}^1(W_{l+1}^{1}x_{l+1}^* + W_{l+1}^{0}x_{l}^* + w_{l+1} ) + \sum_{j\in \{2 \ldots n\}}B_{l+1}^j(T_{l+1}^{j}x_{l+1}^* + t_{l+1}^j) + s_{l+1} \\
\end{multline} \vspace*{-5ex} 

Making $x_{l+1}$ explicit gives

\vspace*{-5ex}  \begin{multline}
\nonumber (I - B_{l+1}^1W_{l+1}^{1} - \sum_{j\in \{2 \ldots n\}}B_{l+1}^jT_{l+1}^{j})x_{l+1}^* = A_{l+1}x_{l}^* \\+ B_{l+1}^1(W_{l+1}^{0}x_{l}^* + w_{l+1} ) + \sum_{j\in \{2 \ldots n\}}B_{l+1}^jt_{l+1}^j + s_{l+1} \\
\end{multline} \vspace*{-5ex} 

\vspace*{-5ex}  \begin{multline}
\nonumber x_{l+1}^* = (I - B_{l+1}^1W_{l+1}^{1} - \sum_{j\in \{2 \ldots n\}}B_{l+1}^jT_{l+1}^{j})^{-1}(A_{l+1}x_{l}^* \\+ B_{l+1}^1(W_{l+1}^{0}x_{l}^* + w_{l+1} ) + \sum_{j\in \{2 \ldots n\}}B_{l+1}^jt_{l+1}^j + s_{l+1}) \\
\end{multline} \vspace*{-5ex} 

The structure of the above equation justifies the following substitution

\vspace*{-5ex}  \begin{multline}
\nonumber\eqref{xoptsna}_{k=l} \;\; x_{l+1}^* = \Phi_{l}x_{l}^* + \phi_l\\
\end{multline} \vspace*{-5ex} 

\vspace*{-5ex}  \begin{multline}
 (I- B_{l+1}^1W_{l+1}^{1} - \sum_{j\in \{2 \ldots n\}}B_{l+1}^jT_{l+1}^{j})(\Phi_{l}x_{l}^* + \phi_l) = A_{l+1}x_{l}^* \\+ B_{l+1}^1(W_{l+1}^{0}x_{l}^* + w_{l+1} ) + \sum_{j\in \{2 \ldots n\}}B_{l+1}^jt_{l+1}^j + s_{l+1} \\\label{xlgna}
\end{multline} \vspace*{-5ex} 

Comparing coefficients gives

\vspace*{-5ex}  \begin{multline}
\nonumber\eqref{xlgna}_{x_{l}^*}=\eqref{Phikna}_{k=l} \;\; \Phi_{l} = (I- B_{l+1}^1W_{l+1}^{1} - \sum_{j\in \{2 \ldots n\}}B_{l+1}^jT_{l+1}^{j})^{-1}(A_{l+1} + B_{l+1}^1W_{l+1}^{0})\\
\end{multline} \vspace*{-5ex} 

\vspace*{-5ex}  \begin{multline}
\nonumber\eqref{xlgna}_{const.}=\eqref{phikna}_{k=l} \;\; \phi_{l} = (I- B_{l+1}^1W_{l+1}^{1} - \sum_{j\in \{2 \ldots n\}}B_{l+1}^jT_{l+1}^{j}){-1}\\(B_{l+1}^1w_{l+1} + \sum_{j\in \{2 \ldots n\}}B_{l+1}^jt_{l+1}^j + s_{l+1}) \\
\end{multline} \vspace*{-5ex}

The above relation between $x_{l+1}$ and $x_l$ can be used to finish the inductive step of $p_l^i$, $\lambda_l$ and $\mu_{l+1}^i$.
As a start $\eqref{xoptsna}_{k=l}$ is used to continue the derivation of $p_l^i$.

\vspace*{-5ex}  \begin{multline}
\nonumber  p_l^{i*} = A_{l+1}^{'}[M_{l+1}^{i}x_{l+1}^* + m_{l+1}^i - Q_{l+1}^i\tilde{x}_{l+1}^i]  \; ; \; i \in \{2\ldots n\}\\
\end{multline} \vspace*{-5ex} 

\vspace*{-5ex}  \begin{multline}
\nonumber  p_l^{i*} = A_{l+1}^{'}[M_{l+1}^{i}(\Phi_{l}x_{l}^* + \phi_l) + m_{l+1}^i - Q_{l+1}^i\tilde{x}_{l+1}^i]  \; ; \; i \in \{2\ldots n\}\\
\end{multline} \vspace*{-5ex} 

The structure of the above equations justifies the following substitutions

\vspace*{-5ex}  \begin{multline}
\nonumber\eqref{psubsna}_{k=l+1} \;\; p_{l}^i = (M_{l}^{i}-Q_{l}^2)x_{l}^* + m_{l}^i \; ; \; i \in \{2\ldots n\}\\
\end{multline} \vspace*{-5ex} 

\vspace*{-5ex}  \begin{multline}
(M_{l}^{i}-Q_{l}^2)x_{l}^* + m_{l}^i = A_{l+1}^{'}[M_{l+1}^{i}(\Phi_{l}x_{l}^* + \phi_l) \\+ m_{l+1}^i - Q_{l+1}^i\tilde{x}_{l+1}^i]  \; ; \; i \in \{2\ldots n\}\\\label{plgna}
\end{multline} \vspace*{-5ex} 

By comparing coefficients it follows that

\vspace*{-5ex}  \begin{multline}
\nonumber\eqref{plgna}_{x_{l}^*}=\eqref{Mikna}_{k=l} \;\; M_{l}^{i} = Q_{l}^i + A_{l+1}^{'}M_{l+1}^{i}\Phi_{l}  \; ; \; i \in \{2\ldots n\}\\
\end{multline} \vspace*{-5ex}

\vspace*{-5ex}  \begin{multline}
\nonumber\eqref{plgna}_{const.}=\eqref{mikna}_{k=l} \;\; m_l^i = A_{l+1}^{'}[M_{l+1}^{i}\phi_l + m_{l+1}^i - Q_{l+1}^i\tilde{x}_{l+1}^i] \; ; \; i \in \{2\ldots n\}\\
\end{multline} \vspace*{-5ex}

As a next step $\eqref{xoptsna}_{k=l}$, $\eqref{muksna}_{k=l+1}$ and \eqref{nusubsna} are used to continue the derivation of $\lambda_l$.

\vspace*{-5ex}  \begin{multline}
\nonumber  \lambda_{l} = A_{l+1}^{'}(L_{l+1}x_{l+1}^* + l_{l+1}) - A_{l+1}^{'}Q_{l+1}^1\tilde{x}_{l+1}^1\\+ \sum_{j \in 2\ldots n}A_{l+1}^{'}Q_{l+1}^jA_{l+1}\mu_{l}^{j} + \sum_{j \in 2\ldots n}A_{l+1}^{'}Q_{l+1}^jB_{l+1}^j\nu_{l}^{j} \\
\end{multline} \vspace*{-5ex} 

\vspace*{-5ex}  \begin{multline}
\nonumber  \lambda_{l} = A_{l+1}^{'}(L_{l+1}x_{l+1}^* + l_{l+1}) - A_{l+1}^{'}Q_{l+1}^1\tilde{x}_{l+1}^1\\+ \sum_{j \in 2\ldots n}A_{l+1}^{'}Q_{l+1}^jA_{l+1}(C_{l}^jx_{l}^* + c_{l}^j) + \sum_{j \in 2\ldots n}A_{l+1}^{'}Q_{l+1}^jB_{l+1}^j(N_{l}^{1j}x_{l+1}^* + N_{l}^{0j}x_{l}^* + n_l^j) \\
\end{multline} \vspace*{-5ex}

\vspace*{-5ex}  \begin{multline}
\nonumber  \lambda_{l} = A_{l+1}^{'}(L_{l+1}(\Phi_{l}x_{l}^* + \phi_l) + l_{l+1}) - A_{l+1}^{'}Q_{l+1}^1\tilde{x}_{l+1}^1+ \sum_{j \in 2\ldots n}A_{l+1}^{'}Q_{l+1}^jA_{l+1}(C_{l}^jx_{l}^* + c_{l}^j) \\+ \sum_{j \in 2\ldots n}A_{l+1}^{'}Q_{l+1}^jB_{l+1}^j(N_{l}^{1j}(\Phi_{l}x_{l}^* + \phi_l) + N_{l}^{0j}x_{l}^* + n_l^j) \\
\end{multline} \vspace*{-5ex} 

The structure of the above equation justifies the following substitution

\vspace*{-5ex}  \begin{multline}
\nonumber\eqref{lambdasubsna}_{k=l+1} \;\; \lambda_{l} = (L_{l}-Q_{l}^1)x_{l}^* + l_{l}\\
\end{multline} \vspace*{-5ex} 

\vspace*{-5ex}  \begin{multline}
 (L_{l}-Q_{l}^1)x_{l}^* + l_{l} = A_{l+1}^{'}(L_{l+1}(\Phi_{l}x_{l}^* + \phi_l) + l_{l+1}) - A_{l+1}^{'}Q_{l+1}^1\tilde{x}_{l+1}^1 \\+ \sum_{j \in 2\ldots n}A_{l+1}^{'}Q_{l+1}^jA_{l+1}(C_{l}^jx_{l}^* + c_{l}^j) + \sum_{j \in 2\ldots n}A_{l+1}^{'}Q_{l+1}^jB_{l+1}^j(N_{l}^{1j}(\Phi_{l}x_{l}^* + \phi_l) + N_{l}^{0j}x_{l}^* + n_l^j) \\\label{lambdaglna}
\end{multline} \vspace*{-5ex} 

Comparing coefficients gives

\vspace*{-5ex}  \begin{multline}
\nonumber\eqref{lambdaglna}_{x_{l}^*}=\eqref{Lkna}_{k=l} \;\; L_{l} = Q_{l}^1 + A_{l+1}^{'}L_{l+1}\Phi_{l} + \sum_{j \in 2\ldots n}A_{l+1}^{'}Q_{l+1}^jA_{l+1}C_{l}^j \\+ \sum_{j \in 2\ldots n}A_{l+1}^{'}Q_{l+1}^jB_{l+1}^j(N_{l}^{1j}\Phi_{l} + N_{l}^{0j}) \\
\end{multline} \vspace*{-5ex} 

\vspace*{-5ex}  \begin{multline}
\nonumber\eqref{lambdaglna}_{const.}=\eqref{lkna}_{k=l} \;\; l_l = A_{l+1}^{'}(L_{l+1}\phi_l + l_{l+1}) - A_{l+1}^{'}Q_{l+1}^1\tilde{x}_{l+1}^1\\+ \sum_{j \in 2\ldots n}A_{l+1}^{'}Q_{l+1}^jA_{l+1}c_{l}^j + \sum_{j \in 2\ldots n}A_{l+1}^{'}Q_{l+1}^jB_{l+1}^j(N_{l}^{1j}\phi_l + n_l^j) \\
\end{multline} \vspace*{-5ex}

Eventually $\eqref{xoptsna}_{k=l}$ and \eqref{nusubsna} are utilized to deduce a relation for $\mu_{l+1}$ that is affinely dependent on $x_{l+1}$.

\vspace*{-5ex}  \begin{multline}
\nonumber \mu_{l+1}^{i} =  A_{l+1}(C_{l}^ix_{l}^* + c_{l}^i) + B_{l+1}^i\nu_{l}^{i} \\
\end{multline} \vspace*{-5ex} 

\vspace*{-5ex}  \begin{multline}
\nonumber \mu_{l+1}^{i} =  A_{l+1}(C_{l}^ix_{l}^* + c_{l}^i) + B_{l+1}^i(N_{l}^{1i}x_{l+1}^* + N_{l}^{0i}x_{l}^* + n_l^i) \\
\end{multline} \vspace*{-5ex} 

\vspace*{-5ex}  \begin{multline}
\nonumber \mu_{l+1}^{i} =  A_{l+1}(C_{l}^i\Phi_{l}^{-1}(x_{l+1}^* - \phi_l) + c_{l}^i) + B_{l+1}^i(N_{l}^{1i}x_{l+1}^* + N_{l}^{0i}\Phi_{l}^{-1}(x_{l+1}^* - \phi_l) + n_l^i) \\
\end{multline} \vspace*{-5ex}

The structure of the above equations justifies the following substitutions

\vspace*{-5ex}  \begin{multline}
\nonumber \mu_{l+1}^i =  C_{l+1}^ix_{l+1}^*  + c_{l+1}^i \; ; \; i \in \{2\ldots n\}\\
\end{multline} \vspace*{-5ex} 

\vspace*{-5ex}  \begin{multline}
C_{l+1}^ix_{l+1}^*  + c_{l+1}^i =  A_{l+1}(C_{l}^i\Phi_{l}^{-1}(x_{l+1}^* - \phi_l) + c_{l}^i) \\+ B_{l+1}^i(N_{l}^{1i}x_{l+1}^* + N_{l}^{0i}\Phi_{l}^{-1}(x_{l+1}^* - \phi_l) + n_l^i) \; ; \; i \in \{2\ldots n\}\\\label{mulgna}
\end{multline} \vspace*{-5ex} 

By comparing coefficients it follows that

\vspace*{-5ex}  \begin{multline}
\nonumber\eqref{mulgna}_{x_{l+1}^*}=\eqref{Cikna}_{k=l} \;\; C_{l+1}^{i} = A_{l+1}C_{l}^i\Phi_{l}^{-1} + B_{l+1}^i(N_{l}^{1i} + N_{l}^{0i}\Phi_{l}^{-1}) \; ; \; i \in \{2\ldots n\}\\
\end{multline} \vspace*{-5ex} 

\vspace*{-5ex}  \begin{multline}
\nonumber\eqref{mulgna}_{const}=\eqref{cikna}_{k=l} \;\; c_{l+1}^{i} = A_{l+1}(-C_{l}^i\Phi_{l}^{-1}\phi_l + c_{l}^i) + B_{l+1}^i(-N_{l}^{0i}\Phi_{l}^{-1}\phi_l + n_l^i) \; ; \; i \in \{2\ldots n\}\\
\end{multline} \vspace*{-5ex}

At this point the inductive step and hence the induction argument is completed, but we will try to transform $u_{l+1}^{1*}, \ldots,u_{l+1}^{n*}$ so that their evolution depends affinely on $x_l^*$, the optimal equilibrium value of the state vector of the previous stage.\\

Let us start with $u_{l+1}^{1*}$, where we apply $\eqref{xoptsna}_{k=l}$ to \eqref{u1subsna}

\vspace*{-5ex}  \begin{multline}
\nonumber \eqref{u1subsna} \; \; u_{l+1}^{1*} =  W_{l+1}^{1}x_{l+1}^* + W_{l+1}^{0}x_{l}^* + w_{l+1} \\
\end{multline} \vspace*{-5ex} 

\vspace*{-5ex}  \begin{multline}
\nonumber u_{l+1}^{1*} =  W_{l+1}^{1}(\Phi_{l}x_{l}^* + \phi_l) + W_{l+1}^{0}x_{l}^* + w_{l+1} \\
\end{multline} \vspace*{-5ex} 

The structure of the above equation justifies the following substitution

\vspace*{-5ex}  \begin{multline}
\nonumber u_{l+1}^{1*} = P_{l+1}^{1}x_l + \alpha_{l+1}^1\\
\end{multline} \vspace*{-5ex} 

\vspace*{-5ex}  \begin{multline}
P_{l+1}^{1}x_l + \alpha_{l+1}^1 =  W_{l+1}^{1}(\Phi_{l}x_{l}^* + \phi_l) + W_{l+1}^{0}x_{l}^* + w_{l+1} \\\label{u1evna}
\end{multline} \vspace*{-5ex} 

Comparing coefficients gives

\vspace*{-5ex}  \begin{multline}
\nonumber\eqref{u1evna}_{x_{l}^*}=\eqref{Pk1na}_{k=l} \;\; P_{l+1}^{1} =  W_{l+1}^{1}\Phi_{l} + W_{l+1}^{0}\\
\end{multline} \vspace*{-5ex} 

\vspace*{-5ex}  \begin{multline}
\nonumber\eqref{u1evna}_{const.}=\eqref{alphak1na}_{k=l} \;\; \alpha_{l+1}^{1} = W_{l+1}^{1}\phi_l + w_{l+1}\\
\end{multline} \vspace*{-5ex}

Finally $\eqref{xoptsna}_{k=l}$ is used in \eqref{uisubsna}:

\vspace*{-5ex}  \begin{multline}
\nonumber \eqref{uisubsna} \; \; u_{l+1}^{i*} =  T_{l+1}^{i}x_{l+1}^* + t_{l+1}^i   \; ; \; i \in \{2\ldots n\}\\
\end{multline} \vspace*{-5ex} 

\vspace*{-5ex}  \begin{multline}
\nonumber u_{l+1}^{i*} =  T_{l+1}^{i}(\Phi_{l}x_{l}^* + \phi_l) + t_{l+1}^i   \; ; \; i \in \{2\ldots n\}\\
\end{multline} \vspace*{-5ex} 

The structure of the above equations justifies the following substitutions

\vspace*{-5ex}  \begin{multline}
\nonumber u_{l+1}^{i*} = P_{l+1}^{i}x_l + \alpha_{l+1}^i  \; ; \; i \in \{2\ldots n\}\\
\end{multline} \vspace*{-5ex} 

\vspace*{-5ex}  \begin{multline}
P_{l+1}^{i}x_l + \alpha_{l+1}^i =  T_{l+1}^{i}(\Phi_{l}x_{l}^* + \phi_l) + t_{l+1}^i  \; ; \; i \in \{2\ldots n\}\\\label{uievna}
\end{multline} \vspace*{-5ex} 

By comparing coefficients it follows that

\vspace*{-5ex}  \begin{multline}
\nonumber\eqref{uievna}_{x_{l}^*}=\eqref{Pkina}_{k=l} \;\; P_{l+1}^{i} = T_{l+1}^{i}\Phi_{l}   \; ; \; i \in \{2\ldots n\}\\
\end{multline} \vspace*{-5ex} 

\vspace*{-5ex}  \begin{multline}
\nonumber\eqref{uievna}_{const.}=\eqref{alphakina}_{k=l} \;\; \alpha_{l+1}^{i} = T_{l+1}^{i}\phi_l + t_{l+1}^i   \; ; \; i \in \{2\ldots n\} \;\;\;\;\;\;\;\;\;\;\;\;\;\;\;\;\boxed{}\\
\end{multline} \vspace*{-5ex}

\label{OLS_ss2i}
\clearpage

\subsection{Special case: "Interwoven-inductions-results" for linear-quadratic games with one leader and one follower}

In this subsection we first specialize the results of the previous subsection \eqref{OLS_ss2i} to a linear-quadratic 2-person game in Corollary \eqref{scols2i} and then in Proposition \eqref{bsolsäa} the specialized results are transformed into the terminology used in Corollary 7.1 in Ba\c{s}ar and Olsder (1999, pp. 371-372)\cite{baol} to point out some serious mistakes stated there.

\begin{cor}
A 2-person linear-quadratic dynamic game (cf. Def. \eqref{defspielaq}) admits a \emph{unique open-loop Stackelberg equilibrium solution with one leader and  one follower} if  
\begin{itemize}
\item $Q_k^i \geq 0$ and $R_k^{ii} > 0$ (defined for $k \in K$ , $i \in N$).
\item $(I- B_{k+1}^1W_{k+1}^{1} - B_{k+1}^2T_{k+1}^{2})^{-1}$, $(A_{k+1} + B_{k+1}^1W_{k+1}^{0})^{-1}$ (defined for $k \in K$) exists.
\end{itemize}
If these conditions are satisfied, the unique equilibrium strategies are given by \eqref{gammaolsna2}, where the associated state trajectory $x_{k+1}^*$ is given by \eqref{xoptsna2}.\footnote{For all equations belonging to this corollary and its proof, $k \in \{0 , \ldots ,T-1\}$ if nothing different is stated.}

\vspace*{-5ex}  \begin{align}
f_{k-1}(x_{k-1},u_k^1,\ldots,u_k^n) = A_kx_{k-1} + B_k^1u_k^1 + B_k^2u_k^2 \; ; \; k \in K\label{folsna2}
\end{align}

\vspace*{-5ex}  \begin{align}
L^i(x_0,u^{1},\ldots,u^{n})= \sum_{k = 1}^T g_k^i(x_k,u_k^{1}, u_k^2,x_{k-1}) 
\end{align}

\vspace*{-5ex}  \begin{multline}
 g_k^i(x_k, u_k^{1}, u_k^2,x_{k-1}) = \frac{1}{2} (x_k^{'}Q_k^ix_k + u_k^{i'}u_k^i \\+ u_k^{j'}R_k^{ij}u_k^j) \; ; \; k \in K\; ; \; i,j \in \{1,2\} \; ; \; i \not = j\\\label{golsna2}
\end{multline} \vspace*{-5ex}

\vspace*{-5ex}  \begin{align}
 x_{k+1}^* = \Phi_{k}x_{k}^* \; ; \; x_0^* = x_0\label{xoptsna2}
\end{align}

\vspace*{-5ex}  \begin{align}
\Phi_{k} = (I- B_{k+1}^1W_{k+1}^{1} - B_{k+1}^2T_{k+1}^{2})^{-1}(A_{k+1} + B_{k+1}^1W_{k+1}^{0})\label{Phikna2}
\end{align}

\vspace*{-5ex}  \begin{align}
W_{k+1}^{1} =   -B_{k+1}^{1'}L_{k+1} -  B_{k+1}^{1'}Q_{k+1}^2B_{k+1}^2N_{k}^{12} \label{Wk1na2}
\end{align}

\vspace*{-5ex}  \begin{align}
W_{k+1}^{0} =  -B_{k+1}^{1'}Q_{k+1}^2A_{k+1}C_{k}^2  - B_{k+1}^{1'}Q_{k+1}^2B_{k+1}^2N_{k}^{02}  \label{Wk0na2}
\end{align}

\vspace*{-5ex}  \begin{align}
T_{k+1}^{2} = -B_{k+1}^{2'}M_{k+1}^{2}  \label{Tkna2}
\end{align}

\vspace*{-5ex}  \begin{align}
N_{k}^{12} = -(B_{k+1}^{2'}Q_{k+1}^2B_{k+1}^2 + I)^{-1}(B_{k+1}^{2'}L_{k+1} - R_{k+1}^{12}B_{k+1}^{2'}M_{k+1}^{2})  \label{Nk1na2}
\end{align}

\vspace*{-5ex}  \begin{align}
N_{k}^{02} = -(B_{k+1}^{2'}Q_{k+1}^2B_{k+1}^2 + I)^{-1}(B_{k+1}^{2'}Q_{k+1}^2A_{k+1}C_{k}^2)   \label{Nk0na2}
\end{align}

\vspace*{-5ex}  \begin{align}
M_{k}^{2} = Q_{k}^2 + A_{k+1}^{'}M_{k+1}^{2}\Phi_{k} \; ; \; M_T^2 = Q_T^2 \label{Mikna2}
\end{align}

\vspace*{-5ex}  \begin{multline}
L_{k} = Q_{k}^1 + A_{k+1}^{'}(L_{k+1}\Phi_{k} + Q_{k+1}^2A_{k+1}C_{k}^2 \\+ Q_{k+1}^2B_{k+1}^2(N_{k}^{12}\Phi_{k} + N_{k}^{02})) \; ; \; L_T = Q_T^1 \\ \label{Lkna2}
\end{multline} \vspace*{-5ex}

\vspace*{-5ex}  \begin{align}
C_{k+1}^{2} = A_{k+1}C_{k}^2\Phi_{k}^{-1} + B_{k+1}^2(N_{k}^{12} + N_{k}^{02}\Phi_{k}^{-1}) \; ; \; C_0 = 0\label{Cikna2}
\end{align}

\vspace*{-5ex}  \begin{align}
\gamma_{k+1}^{i*}(x_0) = u_{k+1}^{i*} = P_{k+1}^{i}x_k^* \; ; \; i \in \{1,2\}\label{gammaolsna2}
\end{align}

\vspace*{-5ex}  \begin{align}
P_{k+1}^{1} =  W_{k+1}^{1}\Phi_{k} + W_{k+1}^{0} \label{Pk1na2}
\end{align}

\vspace*{-5ex}  \begin{align}
P_{k+1}^{2} = T_{k+1}^{2}\Phi_{k}   \; ; \; i \in \{2\ldots n\} \label{Pkina2}
\end{align}\\\\

\label{scols2i}
\end{cor}

\textsc {Proof:}

Corollary \eqref{scols2i} is proven in the same way as Theorem \eqref{OLS2i} taking into consideration simplifications resulting from the different number of followers and the modified state equation and cost functionals.\;\;\;\;\;\;\;\;\;\;\;\;\;\;\;\;\boxed{}\\

\begin{rem}
Special attention should be paid to the fact that the assumption about the existence of unique solutions of the systems of equations \eqref{Nk1na}, \eqref{Nk0na} and \eqref{nkna} in Theorem \eqref{OLS2i} is equivalent to the existence of $(B_{k+1}^{2'}Q_{k+1}^2B_{k+1}^2 + I)^{-1}$ in this special case. But the existence of $(B_{k+1}^{2'}Q_{k+1}^2B_{k+1}^2 + I)^{-1}$ is assured because of the assumption made on $Q_{k+1}^2$.
\end{rem}

\begin{pro}
The systems of equations defining the unique equilibrium strategies $\gamma_{k+1}^{i*}(x_0)$ (i $\in \{1,2\}$) and the associated state trajectory $x_{k+1}^*$ in Corollary \eqref{scols2i} can also be written in the following way:\footnote{For all equations belonging to this proposition and its proof, $k \in \{0 , \ldots ,T-1\}$ if nothing different is stated. \eqref{Mikna2b}, \eqref{Lkna2b} and \eqref{Cikna2b} are wrong in Ba\c{s}ar and Olsder.}

\vspace*{-5ex}  \begin{align}
 x_{k+1}^* = \Phi_{k}x_{k}^* \; ; \; x_0^* = x_0\label{xoptsna2b}
\end{align}

\vspace*{-5ex}  \begin{multline}
\Phi_{k} = (I + B_{k+1}^1B_{k+1}^{1'}(I + Q_{k+1}^2B_{k+1}^2B_{k+1}^{2'})^{-1}\Lambda_{k+1}  + B_{k+1}^1B_{k+1}^{1'}Q_{k+1}^2B_{k+1}^2\\(B_{k+1}^{2'}Q_{k+1}^2B_{k+1}^2 + I)^{-1}R_{k+1}^{12}B_{k+1}^{2'}P_{k+1} + B_{k+1}^2B_{k+1}^{2'}P_{k+1})^{-1}\\(A_{k+1} - B_{k+1}^1B_{k+1}^{1'}Q_{k+1}^2(I + B_{k+1}^2B_{k+1}^{2'}Q_{k+1}^2)^{-1}A_{k+1}M_{k}) \\\label{Phikna2b}
\end{multline} \vspace*{-5ex}

\vspace*{-5ex}  \begin{align}
P_{k} = Q_{k}^2 + A_{k+1}^{'}P_{k+1}\Phi_{k} \; ; \; P_T = Q_T^2 \label{Mikna2b}
\end{align}

\vspace*{-5ex}  \begin{align}
\Lambda_{k} = Q_{k}^1 + A_{k+1}^{'}(\Lambda_{k+1}\Phi_{k} + Q_{k+1}^2A_{k+1}M_{k} - Q_{k+1}^2B_{k+1}^2N_{k}) \; ; \; \Lambda_T = Q_T^1 \label{Lkna2b}
\end{align}

\vspace*{-5ex}  \begin{align}
M_{k+1} = (A_{k+1}M_{k} - B_k^2N_k)\Phi_{k}^{-1} \; ; \; M_0 = 0\label{Cikna2b}
\end{align}

\vspace*{-5ex}  \begin{multline}
N_{k+1} = (B_{k+1}^{2'}Q_{k+1}^2B_{k+1}^2 + I)^{-1}(B_{k+1}^{2'}\Lambda_{k+1}\Phi_k \\- R_{k+1}^{12}B_{k+1}^{2'}K_{k+1}^2 + B_{k+1}^{2'}Q_{k+1}^2A_{k+1}M_{k}) \\\label{Nkb}
\end{multline} \vspace*{-5ex}

\vspace*{-5ex}  \begin{align}
\gamma_{k+1}^{i*}(x_0) = u_{k+1}^{i*} = -B_{k+1}^{i'}K_{k+1}^{i}x_k^* \; ; \; i \in \{1,2\}\label{gammaolsna2b}
\end{align}

\vspace*{-5ex}  \begin{multline}
K_{k+1}^{1} =  [I + Q_k^2B_k^2B_k^{2'}]^{-1}\Lambda_{k+1}\Phi_k \\+ Q_k^2B_k^2[I + B_k^{2'}Q_k^2B_k^2]^{-1}R_k^{12}B_k^{2'}P_{k+1}\Phi_k + Q_k^2[I + B_k^2B_k^{2'}Q_k^2]^{-1}A_kM_k \\\label{Pk1na2b}
\end{multline} \vspace*{-5ex}

\vspace*{-5ex}  \begin{align}
K_{k+1}^{2} = P_{k+1}\Phi_{k}   \; ; \; i \in \{2\ldots n\} \label{Pkina2b}
\end{align}\\\\

\label{bsolsäa}
\end{pro}

\textsc {Proof:}

The proof is carried out by renaming the costate matrices and then showing that the relations for the state and costate matrices and for the control vectors of Corollary \eqref{scols2i} can be rewritten in the way stated above.

Let us start by renaming the costate matrices $C_{k}^{2}$, $L_k$ and $M_{k}^{2}$.

\vspace*{-5ex}  \begin{multline}
C_{k}^{2} \;\widehat {=}\; M_k \; ; \;   L_k\; \widehat {=}\; \Lambda_k \; ; \; M_{k}^{2}\;   \widehat {=}\; P_k \; ; \; k \in \{ 0 , \ldots , T\}\\\label{äqusubs2}
\end{multline} \vspace*{-5ex} 

Next we prove that the costate matrices fulfill \eqref{Mikna2b}, \eqref{Cikna2b} and \eqref{Lkna2b} respectively.\\

Taking consideration of the renaming \eqref{Mikna2} gives

\vspace*{-5ex}  \begin{multline}
\nonumber \eqref{Mikna2b} \; \; P_k = Q_{k}^2 + A_{k+1}^{'}P_{k+1}^{2}\Phi_{k} \\
\end{multline} \vspace*{-5ex} 

Now we show (using \eqref{äqusubs2}) that there is a relation between the cocontrol matrices \eqref{Nk1na2}, \eqref{Nk0na2} and \eqref{Nkb}:

\vspace*{-5ex}  \begin{multline}
\nonumber N_{k}^{12}\Phi_{k} + N_{k}^{02} = -(B_{k+1}^{2'}Q_{k+1}^2B_{k+1}^2 + I)^{-1}(B_{k+1}^{2'}\Lambda_{k+1} - R_{k+1}^{12}B_{k+1}^{2'}P_{k+1})\Phi_{k} \\- (B_{k+1}^{2'}Q_{k+1}^2B_{k+1}^2 + I)^{-1}(B_{k+1}^{2'}Q_{k+1}^2A_{k+1}M_{k}) \\
\end{multline} \vspace*{-5ex} 

Making use of \eqref{Pkina2b} gives

\vspace*{-5ex}  \begin{multline}
\nonumber -N_{k}^{12}\Phi_{k} - N_{k}^{02} = (B_{k+1}^{2'}Q_{k+1}^2B_{k+1}^2 + I)^{-1}\\(B_{k+1}^{2'}\Lambda_{k+1}\Phi_k - R_{k+1}^{12}B_{k+1}^{2'}K_{k+1}^2) + (B_{k+1}^{2'}Q_{k+1}^2B_{k+1}^2 + I)^{-1}(B_{k+1}^{2'}Q_{k+1}^2A_{k+1}M_{k}) = N_{k+1} \\
\end{multline} \vspace*{-5ex} 

Therefore the following relation is satisfied

\vspace*{-5ex}  \begin{multline}
\;\;\;\;-\eqref{Nk1na2}\Phi_k \; - \; \eqref{Nk0na2} = \eqref{Nkb} \\\label{Nkr}
\end{multline} \vspace*{-5ex}

Applying \eqref{äqusubs2} and \eqref{Nkr} to \eqref{Cikna2} yields

\vspace*{-5ex}  \begin{multline}
\nonumber M_{k+1} = A_{k+1}M_{k}\Phi_{k}^{-1} + B_{k+1}^2(N_{k}^{12} + N_{k}^{02}\Phi_{k}^{-1}) \\
\end{multline} \vspace*{-5ex} 

\vspace*{-5ex}  \begin{multline}
\nonumber M_{k+1} = (A_{k+1}M_{k} + B_{k+1}^2(N_{k}^{12}\Phi_{k} + N_{k}^{02}))\Phi_{k}^{-1} \\
\end{multline} \vspace*{-5ex} 

\vspace*{-5ex}  \begin{multline}
\nonumber \eqref{Cikna2b} \; \;M_{k+1} = (A_{k+1}M_{k} - B_{k+1}^2N_{k+1})\Phi_{k}^{-1} \\
\end{multline} \vspace*{-5ex}

Making use of \eqref{äqusubs2} and \eqref{Nkr} in \eqref{Lkna2} gives

\vspace*{-5ex}  \begin{multline}
\nonumber \Lambda_{k} = Q_{k}^1 + A_{k+1}^{'}(\Lambda_{k+1}\Phi_{k} + Q_{k+1}^2A_{k+1}M_{k} + Q_{k+1}^2B_{k+1}^2(N_{k}^{12}\Phi_{k} + N_{k}^{02})) \\
\end{multline} \vspace*{-5ex} 

\vspace*{-5ex}  \begin{multline}
\nonumber \eqref{Lkna2b} \; \; \Lambda_{k} = Q_{k}^1 + A_{k+1}^{'}(\Lambda_{k+1}\Phi_{k} + Q_{k+1}^2A_{k+1}M_{k} - Q_{k+1}^2B_{k+1}^2N_{k+1}) \\
\end{multline} \vspace*{-5ex}

Next we show the correctness of the relation of $\Phi_k$ stated in \eqref{Phikna2b}.\\

First use \eqref{Wk1na2}, \eqref{Wk0na2}, \eqref{Tkna2} and \eqref{äqusubs2} in \eqref{Phikna2}:

\vspace*{-5ex}  \begin{multline}
\nonumber\eqref{Phikna2} \; \; \Phi_{k} = (I - B_{k+1}^1W_{k+1}^{1} - B_{k+1}^2T_{k+1}^{2})^{-1}(A_{k+1} + B_{k+1}^1W_{k+1}^{0}) \\
\end{multline} \vspace*{-5ex} 

\vspace*{-5ex}  \begin{multline}
\nonumber \Phi_{k} = (I + B_{k+1}^1(B_{k+1}^{1'}\Lambda_{k+1} +  B_{k+1}^{1'}Q_{k+1}^2B_{k+1}^2N_{k}^{12}) + B_{k+1}^2B_{k+1}^{2'}P_{k+1})^{-1}\\(A_{k+1} - B_{k+1}^1(B_{k+1}^{1'}Q_{k+1}^2A_{k+1}M_{k}  + B_{k+1}^{1'}Q_{k+1}^2B_{k+1}^2N_{k}^{02})) \\
\end{multline} \vspace*{-5ex} 

Substituting $N_k^{12}$ with the help of \eqref{Nk1na2} yields

\vspace*{-5ex}  \begin{multline}
\nonumber \Phi_{k} = (I + B_{k+1}^1(B_{k+1}^{1'}\Lambda_{k+1} -  B_{k+1}^{1'}Q_{k+1}^2B_{k+1}^2(B_{k+1}^{2'}Q_{k+1}^2B_{k+1}^2 + I)^{-1}\\(B_{k+1}^{2'}\Lambda_{k+1} - R_{k+1}^{12}B_{k+1}^{2'}P_{k+1})) + B_{k+1}^2B_{k+1}^{2'}P_{k+1})^{-1}(A_{k+1} - B_{k+1}^1\\(B_{k+1}^{1'}Q_{k+1}^2A_{k+1}M_{k}  - B_{k+1}^{1'}Q_{k+1}^2B_{k+1}^2(B_{k+1}^{2'}Q_{k+1}^2B_{k+1}^2 + I)^{-1}B_{k+1}^{2'}Q_{k+1}^2A_{k+1}M_{k})) \\
\end{multline} \vspace*{-5ex} 

Putting together $\Lambda_{k+1}$ and $M_k$ gives

\vspace*{-5ex}  \begin{multline}
\nonumber \Phi_{k} = (I + B_{k+1}^1B_{k+1}^{1'}(I - Q_{k+1}^2B_{k+1}^2(B_{k+1}^{2'}Q_{k+1}^2B_{k+1}^2 + I)^{-1}B_{k+1}^{2'})\Lambda_{k+1} \\+ B_{k+1}^1B_{k+1}^{1'}Q_{k+1}^2B_{k+1}^2(B_{k+1}^{2'}Q_{k+1}^2B_{k+1}^2 + I)^{-1}R_{k+1}^{12}B_{k+1}^{2'}P_{k+1} + B_{k+1}^2B_{k+1}^{2'}P_{k+1})^{-1}\\(A_{k+1} - B_{k+1}^1B_{k+1}^{1'}Q_{k+1}^2(I - B_{k+1}^2(B_{k+1}^{2'}Q_{k+1}^2B_{k+1}^2 + I)^{-1}B_{k+1}^{2'}Q_{k+1}^2)A_{k+1}M_{k}) \\
\end{multline} \vspace*{-5ex} 

Finally applying Lemma \eqref{maid} 1. and 2. to the particular expressions above leads to

\vspace*{-5ex}  \begin{multline}
\nonumber \eqref{Phikna2b} \; \; \Phi_{k} = (I + B_{k+1}^1B_{k+1}^{1'}(I + Q_{k+1}^2B_{k+1}^2B_{k+1}^{2'})^{-1}\Lambda_{k+1} \\+ B_{k+1}^1B_{k+1}^{1'}Q_{k+1}^2B_{k+1}^2(B_{k+1}^{2'}Q_{k+1}^2B_{k+1}^2 + I)^{-1}R_{k+1}^{12}B_{k+1}^{2'}P_{k+1} \\+ B_{k+1}^2B_{k+1}^{2'}P_{k+1})^{-1}(A_{k+1} - B_{k+1}^1B_{k+1}^{1'}Q_{k+1}^2(I + B_{k+1}^2B_{k+1}^{2'}Q_{k+1}^2)^{-1}A_{k+1}M_{k}) \\
\end{multline} \vspace*{-5ex}

Eventually the correctness of the rewritten equilibrium strategies $\gamma_{k+1}^{i*}$ ($i \in \{1,2\}$, $k \in \{0, \ldots, T\}$), given by \eqref{gammaolsna2b} - \eqref{Pkina2b}, has to be shown.\\

First note that \eqref{gammaolsna2} $\widehat {=}$ \eqref{gammaolsna2b} if and only if $P_{k+1}^i$ = $-B_{k+1}^iK_{k+1}^i$. 

\vspace*{-5ex}  \begin{multline}
\nonumber\eqref{gammaolsna2} \; \; \gamma_{k+1}^{i*}(x_0) = u_{k+1}^{i*} = P_{k+1}^{i}x_k^* \\
\end{multline} \vspace*{-5ex} 

\vspace*{-5ex}  \begin{multline}
\nonumber\eqref{gammaolsna2b} \; \; \gamma_{k+1}^{i*}(x_0) = u_{k+1}^{i*} = -B_{k+1}^{i'}K_{k+1}^{i}x_k^* \\
\end{multline} \vspace*{-5ex} 

Therefore, as a start use \eqref{Wk1na2} and \eqref{Wk0na2} in \eqref{Pk1na2} (also considering \eqref{äqusubs2})

\vspace*{-5ex}  \begin{multline}
\nonumber \eqref{Pk1na2} \; \; P_{k+1}^{1} =  W_{k+1}^{1}\Phi_{k} + W_{k+1}^{0} \\
\end{multline} \vspace*{-5ex} 

\vspace*{-5ex}  \begin{multline}
\nonumber P_{k+1}^{1} =  -(B_{k+1}^{1'}\Lambda_{k+1} +  B_{k+1}^{1'}Q_{k+1}^2B_{k+1}^2N_{k}^{12})\Phi_{k} - (B_{k+1}^{1'}Q_{k+1}^2A_{k+1}M_{k}  + B_{k+1}^{1'}Q_{k+1}^2B_{k+1}^2N_{k}^{02}) \\
\end{multline} \vspace*{-5ex} 

Next substitute $N_{k}^{12}\Phi_{k} + N_{k}^{02}$ with the help of \eqref{Nkr}

\vspace*{-5ex}  \begin{multline}
\nonumber P_{k+1}^{1} =  - B_{k+1}^{1'}\Lambda_{k+1}\Phi_{k} -  B_{k+1}^{1'}Q_{k+1}^2B_{k+1}^2(N_{k}^{12}\Phi_{k} + N_{k}^{02}) - B_{k+1}^{1'}Q_{k+1}^2A_{k+1}M_{k} \\
\end{multline} \vspace*{-5ex} 

\vspace*{-5ex}  \begin{multline}
\nonumber P_{k+1}^{1} =  - B_{k+1}^{1'}\Lambda_{k+1}\Phi_{k} +  B_{k+1}^{1'}Q_{k+1}^2B_{k+1}^2N_{k+1} - B_{k+1}^{1'}Q_{k+1}^2A_{k+1}M_{k} \\
\end{multline} \vspace*{-5ex} 

Now make use of \eqref{Nkb}

\vspace*{-5ex}  \begin{multline}
\nonumber P_{k+1}^{1} =  -B_{k+1}^{1'}\Lambda_{k+1}\Phi_{k} +  B_{k+1}^{1'}Q_{k+1}^2B_{k+1}^2((B_{k+1}^{2'}Q_{k+1}^2B_{k+1}^2 + I)^{-1}(B_{k+1}^{2'}\Lambda_{k+1}\Phi_k \\- R_{k+1}^{12}B_{k+1}^{2'}K_{k+1}^2) + (B_{k+1}^{2'}Q_{k+1}^2B_{k+1}^2 + I)^{-1}(B_{k+1}^{2'}Q_{k+1}^2A_{k+1}M_{k})) - B_{k+1}^{1'}Q_{k+1}^2A_{k+1}M_{k} \\
\end{multline} \vspace*{-5ex} 

Putting together $\Lambda_{k+1}$ and $M_k$ gives

\vspace*{-5ex}  \begin{multline}
\nonumber P_{k+1}^{1} =  -B_{k+1}^{1'}(I - Q_{k+1}^2B_{k+1}^2(B_{k+1}^{2'}Q_{k+1}^2B_{k+1}^2 + I)^{-1}B_{k+1}^{2'})\Lambda_{k+1}\Phi_k \\- B_{k+1}^{1'}Q_{k+1}^2B_{k+1}^2(B_{k+1}^{2'}Q_{k+1}^2B_{k+1}^2 + I)^{-1}R_{k+1}^{12}B_{k+1}^{2'}P_{k+1}\Phi_k \\- B_{k+1}^{1'}Q_{k+1}^2(I - B_{k+1}^2(B_{k+1}^{2'}Q_{k+1}^2B_{k+1}^2 + I)^{-1}B_{k+1}^{2'}Q_{k+1}^2)A_{k+1}M_{k} \\
\end{multline} \vspace*{-5ex} 

Finally applying Lemma \eqref{maid} 1. and 2. to the particular expressions above and using \eqref{Pk1na2b} leads to

\vspace*{-5ex}  \begin{multline}
\nonumber P_{k+1}^{1} =  -B_{k+1}^{1'}[(I + Q_{k+1}^2B_{k+1}^2B_{k+1}^{2'})^{-1}\Lambda_{k+1}\Phi_k + Q_{k+1}^2B_{k+1}^2(B_{k+1}^{2'}Q_{k+1}^2B_{k+1}^2 + I)^{-1}\\R_{k+1}^{12}B_{k+1}^{2'}P_{k+1}\Phi_k + Q_{k+1}^2(I + B_{k+1}^2B_{k+1}^{2'}Q_{k+1}^2)^{-1}A_{k+1}M_{k}] \\
\end{multline} \vspace*{-5ex} 

\vspace*{-5ex}  \begin{multline}
P_{k+1}^{1} =  - B_{k+1}^{1'}K_{k+1}^1 \\
\end{multline} \vspace*{-5ex}

As a last step $P_{k+1}^{2}=-B_{k+1}^{2'}K_{k+1}^2$ can be shown by using \eqref{Tkna2} in \eqref{Pkina2} (also considering \eqref{äqusubs2})

\vspace*{-5ex}  \begin{multline}
P_{k+1}^{2} = T_{k+1}^{2}\Phi_{k}\\
\end{multline} \vspace*{-5ex} 

\vspace*{-5ex}  \begin{multline}
P_{k+1}^{2} = -B_{k+1}^{2'}P_{k+1}\Phi_{k} = -B_{k+1}^{2'}K_{k+1}^2 \;\;\;\;\;\;\;\;\;\;\;\;\;\;\;\;\boxed{}\\\\
\end{multline} \vspace*{-5ex}

%\vspace*{-5ex}  \begin{multline}
%
%\end{multline} \vspace*{-5ex} 
%
%\vspace*{-5ex}  \begin{multline}
%
%\end{multline} \vspace*{-5ex} 

\begin{rem}
The equations given in Proposition \eqref{bsolsäa} are the same as in Ba\c{s}ar and Olsder (1999, p.371) \cite{baol} (except for the three equations already mentionned above) if k is replaced by k', whereas k'=k+1 and $k' \in K$, and the indices of the state vector x and of the matrices related to it (M, P, $\Lambda$, $\Phi$, $Q^i$) are augmented by 1.
\end{rem}

\label{sslqols}

\clearpage
\subsection{The one-induction-results for affine-quadratic games with one leader and arbitrarily many followers}

In the following, the results of Theorem \eqref{OLS_opt} are applied to an affine-quadratic dynamic game with one leader and arbitrarily many followers. Theorem \eqref{OLS1i} presents equilibirum equations that can easily be used for an algorithmic disintegration of the given Stackelberg game.

\begin{theo}
An n-person affine-quadratic dynamic game (cf. Def. \eqref{defspielaq}) admits a \emph{unique open-loop Stackelberg equilibrium solution with one leader and arbitrarily many followers} if  
\begin{itemize}
\item $Q_k^i \geq 0$, $R_k^{ii} > 0$ (defined for $k \in K$ , $i \in N$).
\item $(I - B_{k+1}^1W_{k+1}^{x} - \sum_{j\in 2\ldots n}B_{k+1}^jT_{k+1}^{jx})^{-1}$ (defined for $k \in K$) exists.
\item \eqref{Nkx}, \eqref{Nkmu} and \eqref{nk} admit unique solutions $N_k^{ix}$, $N_k^{ij\mu}$ and $n_k^{i}$ (defined for: $k \in K$ , $i,j \in \{2,\ldots,n\}$).
\end{itemize}
If these conditions are satisfied, the unique equilibrium strategies are given by \eqref{gammaols}, where the associated state trajectory $x_{k+1}^*$ is given by \eqref{xopts}.\footnote {For all equations belonging to this theorem and its proof, $i \in N$ and $k \in \{0 , \ldots ,T-1\}$ if nothing different is stated.}

\begin{align}
f_{k-1}(x_{k-1},u_k^1,\ldots,u_k^n) = A_kx_{k-1} + \sum_{j\in N}B_k^ju_k^j+s_k \; ; \; k \in K\label{fols}
\end{align}

\begin{align}
L^i(x_0,u^{1},\ldots,u^{n})= \sum_{k = 1}^T g_k^i(x_k,u_k^{1},\ldots,u_k^{n},x_{k-1}) \label{lols}
\end{align}

\vspace*{-5ex}  \begin{multline}
 g_k^i(x_k,u_k^1,\ldots,u_k^n,x_{k-1}) = \frac{1}{2} (x_k^{'}Q_k^ix_k + \sum_{j\in N}u_k^{j'}R_k^{ij}u_k^j)\; \\ +
 \frac{1}{2} (\tilde{x}_k^{i'}Q_k^i\tilde{x}_k^i + \sum_{j\in N}\tilde{u}_k^{ij'}R_k^{ij}\tilde{u}_k^{ij}) - \tilde{x}_k^{i'}Q_k^ix_k - \sum_{j\in N}\tilde{u}_k^{ij'}R_k^{ij}u_k^j \; ; \; k \in K\\\label{gols}
\end{multline} \vspace*{-5ex}

\begin{align}
 x_{k+1}^* = \Phi_{k}^xx_{k}^* + \sum_{j \in 2\ldots n }\Phi_{k}^{j\mu}\mu_{k}^j + \phi_k \; ; \; x_0^* = x_0\label{xopts}
\end{align}

\begin{align}
\Phi_{k}^x = (I - B_{k+1}^1W_{k+1}^{x} - \sum_{j\in 2\ldots n}B_{k+1}^jT_{k+1}^{jx})^{-1}A_{k+1} \label{Phikx}
\end{align}

\vspace*{-5ex}  \begin{multline}
\Phi_{k}^{i\mu} = (I - B_{k+1}^1W_{k+1}^{x} - \sum_{j\in 2\ldots n}B_{k+1}^jT_{k+1}^{jx})^{-1}\\(B_{k+1}^1W_{k+1}^{i\mu} + \sum_{j\in 2\ldots n}B_{k+1}^jT_{k+1}^{ji\mu} ) \; ; \; i \in \{2\ldots n\}\\\label{Phikmu}
\end{multline} \vspace*{-5ex} 

\vspace*{-5ex}  \begin{multline}
\phi_{k} = (I - B_{k+1}^1W_{k+1}^{x} - \sum_{j\in 2\ldots n}B_{k+1}^jT_{k+1}^{jx})^{-1}\\(B_{k+1}^1w_{k+1} + \sum_{j\in 2\ldots n}B_{k+1}^jt_{k+1}^j + s_{k+1})  \\\label{phik}
\end{multline} \vspace*{-5ex}

\begin{align}
\mu_{k+1}^i =  \Psi_{k}^{ix}x_{k}^* + \sum_{j \in 2\ldots n }\Psi_{k}^{ij\mu}\mu_{k}^j + \psi_k^i \; ; \; \mu_0^i = 0 \; ; \; i \in \{2\ldots n\} \label{muopts}
\end{align}

\begin{align}
\Psi_{k}^{ix} = B_{k+1}^iN_{k}^{ix}\Phi_{k}^x \; ; \; i \in \{2\ldots n\} \label{Psikx}
\end{align}

\begin{align}
\Psi_{k}^{ii} =  A_{k+1} + B_{k+1}^i(N_{k}^{ix}\Phi_{k}^{i\mu} + N_{k}^{ii\mu})\; ; \; i \in \{2\ldots n\} \label{Psikmui}
\end{align}

\begin{align}
\Psi_{k}^{im} =  B_{k+1}^i(N_{k}^{ix}\Phi_{k}^{m\mu} + N_{k}^{im\mu})\; ; \; i,m \in \{2\ldots n\} \; ; \; m \not = i \label{Psikmum}
\end{align}

\begin{align}
\psi_{k}^{i} =  B_{k+1}^i(N_{k}^{ix}\phi_k + n_k^i) \; ; \; i \in \{2\ldots n\} \label{psik}
\end{align}

\vspace*{-5ex}  \begin{multline}
W_{k+1}^{x} =   -(R_{k+1}^{11})^{-1} (B_{k+1}^{1'}(L_{k+1}^x + \sum_{j \in 2\ldots n}L_{k+1}^{j\mu}B_{k+1}^jN_{k}^{jx}) \\+ \sum_{j \in 2\ldots n}B_{k+1}^{1'}Q_{k+1}^jB_{k+1}^jN_{k}^{jx})\\\label{Wkx}
\end{multline} \vspace*{-5ex} 

\vspace*{-5ex}  \begin{multline}
W_{k+1}^{m\mu} =  -(R_{k+1}^{11})^{-1} (B_{k+1}^{1'}(L_{k+1}^{m\mu}A_{k+1} + \sum_{j \in 2\ldots n}L_{k+1}^{j\mu}B_{k+1}^jN_{k}^{jm\mu}) \\+ B_{k+1}^{1'}Q_{k+1}^mA_{k+1} + \sum_{j \in 2\ldots n}B_{k+1}^{1'}Q_{k+1}^jB_{k+1}^jN_{k}^{jm\mu})  \; ; \; m \in \{2\ldots n\}\\\label{Wkmu}
\end{multline} \vspace*{-5ex} 

\vspace*{-5ex}  \begin{multline}
w_{k+1} =  -(R_{k+1}^{11})^{-1} (-B_{k+1}^{1'}Q_{k+1}^1\tilde{x}_{k+1}^1  + B_{k+1}^{1'}(\sum_{j \in 2\ldots n}L_{k+1}^{j\mu}B_{k+1}^jn_k^j \\+ l_{k+1}) + \sum_{j \in 2\ldots n}B_{k+1}^{1'}Q_{k+1}^jB_{k+1}^jn_k^j) + \tilde{u}_{k+1}^{11} \\\label{wk}
\end{multline} \vspace*{-5ex}

\begin{align}
T_{k+1}^{ix} = -(R_{k+1}^{ii})^{-1}B_{k+1}^{i'}(M_{k+1}^{ix} + \sum_{j \in 2\ldots n}M_{k+1}^{ij\mu}B_{k+1}^jN_{k}^{jx})  \; ; \; i \in \{2\ldots n\} \label{Tkx}
\end{align}

\vspace*{-5ex}  \begin{multline}
T_{k+1}^{im\mu} =  -(R_{k+1}^{ii})^{-1}B_{k+1}^{i'}(M_{k+1}^{im\mu}A_{k+1} + \\\sum_{j \in 2\ldots n}M_{k+1}^{ij\mu}B_{k+1}^jN_{k}^{jm\mu})\; ; \; i,m \in \{2\ldots n\}\\\label{Tkmu}
\end{multline} \vspace*{-5ex} 

\vspace*{-5ex}  \begin{multline}
t_{k+1}^{i} =  -(R_{k+1}^{ii})^{-1}B_{k+1}^{i'}(\sum_{j \in 2\ldots n}M_{k+1}^{ij\mu} B_{k+1}^jn_k^j \\+ m_{k+1}^i - Q_{k+1}^i\tilde{x}_{k+1}^{i'}) + \tilde{u}_{k+1}^{ii}  \; ; \; i \in \{2\ldots n\}\\\label{tk}
\end{multline} \vspace*{-5ex}

\vspace*{-5ex}  \begin{multline}
-R_{k+1}^{1i}(R_{k+1}^{ii})^{-1}B_{k+1}^{i'}M_{k+1}^{ix} + B_{k+1}^{i'}L_{k+1}^x +  (B_{k+1}^{i'}(Q_{k+1}^i+L_{k+1} ^{i\mu})B_{k+1}^i \\+ R_{k+1}^{ii} - R_{k+1}^{1i}(R_{k+1}^{ii})^{-1}B_{k+1}^{i'}M_{k+1}^{ii\mu}B_{k+1}^i)N_{k}^{ix}  + \sum_{j \in 2\ldots n \;,\;j \not = i}(B_{k+1}^{i'}(Q_{k+1}^j+L_{k+1} ^{j\mu})B_{k+1}^j\\ -R_{k+1}^{1i}(R_{k+1}^{ii})^{-1}B_{k+1}^{i'}M_{k+1}^{ij\mu}B_{k+1}^j)N_{k}^{jx} = 0  \; ; \; i \in \{2\ldots n\}\\\label{Nkx}
\end{multline} \vspace*{-5ex} 

\vspace*{-5ex}  \begin{multline}
-R_{k+1}^{1i}(R_{k+1}^{ii})^{-1}B_{k+1}^{i'}M_{k+1}^{im\mu}A_{k+1} + B_{k+1}^{i'}L_{k+1} ^{m\mu}A_{k+1} \\+ B_{k+1}^{i'}Q_{k+1}^mA_{k+1} +  (B_{k+1}^{i'}(Q_{k+1}^i+L_{k+1} ^{i\mu})B_{k+1}^i + R_{k+1}^{ii} \\- R_{k+1}^{1i}(R_{k+1}^{ii})^{-1}B_{k+1}^{i'}M_{k+1}^{ii\mu}B_{k+1}^i)N_{k}^{im\mu} + \sum_{j \in 2\ldots n \;,\;j \not = i}(B_{k+1}^{i'}(Q_{k+1}^j+L_{k+1} ^{j\mu})B_{k+1}^j\\ -R_{k+1}^{1i}(R_{k+1}^{ii})^{-1}B_{k+1}^{i'}M_{k+1}^{ij\mu}B_{k+1}^j)N_{k}^{jm\mu} = 0 \; ; \; i,m \in \{2\ldots n\}\\\label{Nkmu}
\end{multline} \vspace*{-5ex} 

\vspace*{-5ex}  \begin{multline}
-B_{k+1}^{i'}Q_{k+1}^1\tilde{x}_{k+1}^i + R_{k+1}^{1i}(-(R_{k+1}^{ii})^{-1}B_{k+1}^{i'}( m_{k+1}^i - Q_{k+1}^i\tilde{x}_{k+1}^{i'}) + \tilde{u}_{k+1}^{ii} \\- \tilde{u}_{k+1}^{1i}) + B_{k+1}^{i'}l_{k+1}  +  (B_{k+1}^{i'}(Q_{k+1}^i+L_{k+1} ^{i\mu})B_{k+1}^i + R_{k+1}^{ii} - R_{k+1}^{1i}(R_{k+1}^{ii})^{-1}B_{k+1}^{i'}M_{k+1}^{ii\mu}B_{k+1}^i)n_k^i \\+ \sum_{j \in 2\ldots n \;,\;j \not = i}(B_{k+1}^{i'}(Q_{k+1}^j+L_{k+1} ^{j\mu})B_{k+1}^j -R_{k+1}^{1i}(R_{k+1}^{ii})^{-1}B_{k+1}^{i'}M_{k+1}^{ij\mu}B_{k+1}^j)n_k^j = 0  \; ; \; i \in \{2\ldots n\}\\\label{nk}
\end{multline} \vspace*{-5ex}

\begin{align}
M_{k}^{ix} = Q_{k}^i + A_{k+1}^{'}[M_{k+1}^{ix}\Phi_{k}^x + \sum_{j \in 2\ldots n}M_{k+1}^{ij\mu}\Psi_{k}^{jx}]  \; ; \; M_T^{ix} = Q_{T}^i \; ; \; i \in \{2\ldots n\} \label{Mikx}
\end{align}

\begin{align}
M_{k}^{im\mu} = A_{k+1}^{'}[M_{k+1}^{ix}\Phi_{k}^{m\mu} + \sum_{j \in 2\ldots n}M_{k+1}^{ij\mu}\Psi_{k}^{jm\mu}] \; ; \; M_T^{im\mu} = 0 \; ; \; i,m \in \{2\ldots n\} \label{Mikmu}
\end{align}

\vspace*{-5ex}  \begin{multline}
m_k^i = A_{k+1}^{'}[M_{k+1}^{ix}\phi_k + \sum_{j \in 2\ldots n}M_{k+1}^{ij\mu}\psi_k^j + m_{k+1}^i \\- Q_{k+1}^i\tilde{x}_{k+1}^i] \; ; \; m_T^{i} = 0 \; ; \; i \in \{2\ldots n\}\\\label{mik}
\end{multline} \vspace*{-5ex}

\vspace*{-5ex}  \begin{multline}
L_{k}^{x} = Q_{k}^1 + A_{k+1}^{'}[L_{k+1}^{x}\Phi_{k}^x + \sum_{j \in 2\ldots n}L_{k+1}^{j\mu}\Psi_{k}^{jx} \\+ \sum_{j \in 2\ldots n}Q_{k+1}^jB_{k+1}^jN_{k}^{jx}\Phi_{k}^x] \; ; \;L_T^x = Q_{T}^1\\\label{Lkx}
\end{multline} \vspace*{-5ex} 

\vspace*{-5ex}  \begin{multline}
L_{k}^{i\mu} = A_{k+1}^{'}[L_{k+1}^{x}\Phi_{k}^{i\mu} + \sum_{j \in 2\ldots n}L_{k+1}^{j\mu}\Psi_{k}^{ji\mu} + Q_{k+1}^iA_{k+1} \\+ \sum_{j \in 2\ldots n}Q_{k+1}^jB_{k+1}^j(N_{k}^{jx}\Phi_{k}^{i\mu} + N_{k}^{ji\mu})] \; ; \;L_T^{j\mu} = 0 \; ; \; i \in \{2\ldots n\}\\\label{Lkmu}
\end{multline} \vspace*{-5ex} 

\vspace*{-5ex}  \begin{multline}
l_k = A_{k+1}^{'}[L_{k+1}^{x}\phi_k + \sum_{j \in 2\ldots n}L_{k+1}^{j\mu}\psi_k^j + l_{k+1} \\- Q_{k+1}^1\tilde{x}_{k+1}^1 + \sum_{j \in 2\ldots n}Q_{k+1}^jB_{k+1}^j(N_{k}^{jx}\phi_k + n_k^j)] \; ; \;l_T = 0\\\label{lk}
\end{multline} \vspace*{-5ex}

\begin{align}
\gamma_{k+1}^{i*}(x_0) = u_{k+1}^{i*} = P_{k+1}^{ix}x_k^* +  \sum_{j \in 2\ldots n }P_{k+1}^{ij\mu}\mu_k^j + \alpha_{k+1}^i\label{gammaols}
\end{align}

\begin{align}
P_{k+1}^{1x} = W_{k+1}^{x}\Phi_{k}^x \label{Pk1x}
\end{align}

\begin{align}
P_{k+1}^{1i\mu} = W_{k+1}^{x}\Phi_{k}^{i\mu} + W_{k+1}^{i\mu} \; ; \; i \in \{2\ldots n\} \label{Pk1mu}
\end{align}

\begin{align}
\alpha_{k+1}^{1} = W_{k+1}^{x}\phi_k + w_{k+1} \label{alphak1}
\end{align}

\begin{align}
P_{k+1}^{ix} = T_{k+1}^{ix}\Phi_{k}^x   \; ; \; i \in \{2\ldots n\} \label{Pkix}
\end{align}

\begin{align}
P_{k+1}^{im\mu} = T_{k+1}^{ix}\Phi_{k}^{m\mu} + T_{k+1}^{im\mu} \; ; \; i,m \in \{2\ldots n\} \label{Pkimu}
\end{align}

\begin{align}
\alpha_{k+1}^{i} = T_{k+1}^{ix}\phi_k + t_{k+1}^i   \; ; \; i \in \{2\ldots n\}  \label{alphaki}
\end{align}\\\\

\label{OLS1i}
\end{theo}

\textsc {Proof:}\footnote{The first part of the proof, which is the derivation of the general optimality conditions for the game, is the same as in Theorem \eqref{OLS2i}.}

Theorem \eqref{OLS_opt} can be applied to the given affine-quadratic game, since all conditions are satisfied for the given state equation \eqref{fols} and cost functionals \eqref{lols}. Furthermore $g_k^i$ is strictly convex in $u_k^i$ ($i \in N$). This can be seen by applying Corollary \eqref{pdc} to \eqref{conols}. Therefore there has to be a unique optimal equilibrium solution.

\vspace*{-5ex}  \begin{multline}
\frac{\partial^2}{\partial u_k^{i^2}} g_k^i(x_k,u_k^1 ,\ldots , u_k^n, x_{k-1}) = B_k^{i'}Q_k^iB_k^i + R_k^{ii} \\\label{conols}
\end{multline} \vspace*{-5ex} 

To obtain relations which satisfy this unique solution we have to adapt \eqref{OLS_x} - \eqref{OLS_Hi} to the given state equation and cost functionals. This yields

\vspace*{-5ex}  \begin{multline}
H_k^1 =  \frac{1}{2} (x_k^{'}Q_k^1x_k + \sum_{j \in N}u_k^{j'}R_k^{1j}u_k^j)\;  +
 \frac{1}{2} (\tilde{x}_k^{1'}Q_k^1\tilde{x}_k^1 + \sum_{j \in N}\tilde{u}_k^{1j'}R_k^{1j}\tilde{u}_k^{1j}) \\- \tilde{x}_k^{1'}Q_k^1x_k - \sum_{j\in N}\tilde{u}_k^{1j'}R_k^{1j}u_k^j + \lambda_k^{'}(A_kx_{k-1} + \sum_{j \in N}B_k^ju_k^j+s_k) \\+ \sum_{j \in 2\ldots n}\mu_{k-1}^{j'}A_{k}^{'}[p_{k}^i + Q_{k}^j( x_{k}^* - \tilde{x}_{k}^j )] \\+ \sum_{j \in 2\ldots n}\nu_{k-1}^{j'}(B_k^{j'}Q_k^jx_k + R_k^{jj}u_k^{j} -  B_k^{j'}Q_k^j\tilde{x}_k^{j'} - R_k^{jj}\tilde{u}_k^{jj} + B_k^{j'}p_{k}^j)  \\\label{H1ks}
\end{multline} \vspace*{-5ex} 

\vspace*{-5ex}  \begin{multline}
H_k^i =   \frac{1}{2} (x_k^{'}Q_k^ix_k + \sum_{j\in N}u_k^{j'}R_k^{ij}u_k^j)\;  +
 \frac{1}{2} (\tilde{x}_k^{i'}Q_k^i\tilde{x}_k^i + \sum_{j\in N}\tilde{u}_k^{ij'}R_k^{ij}\tilde{u}_k^{ij}) \\- \tilde{x}_k^{i'}Q_k^ix_k - \sum_{j\in N}\tilde{u}_k^{ij'}R_k^{ij}u_k^j + p_k^{i'}(A_kx_{k-1} + \sum_{j\in N}B_k^ju_k^j+s_k) \; ; \; i \in \{2\ldots n\}\\\label{Hiks}
\end{multline} \vspace*{-5ex} 

\vspace*{-5ex}  \begin{multline}
\nonumber \frac{\partial}{\partial u_k^1} \eqref{H1ks} = 0\; \Rightarrow \;\;\ B_k^{1'}Q_k^1(x_k^*-\tilde{x}_k^1) + R_k^{11}(u_k^{1*}- \tilde{u}_k^{11}) + B_k^{1'}\lambda_k \\+ \sum_{j \in 2\ldots n}B_k^{1'}Q_k^jA_k\mu_{k-1}^{j} + \sum_{j \in 2\ldots n}B_k^{1'}Q_k^jB_k^j\nu_{k-1}^{j}  = 0\\
\end{multline} \vspace*{-5ex} 

\vspace*{-5ex}  \begin{multline}
u_k^{1*} = -(R_k^{11})^{-1}(B_k^{1'}Q_k^1(x_k^*-\tilde{x}_k^1)  + B_k^{1'}\lambda_k \\+ \sum_{j \in 2\ldots n}B_k^{1'}Q_k^jA_k\mu_{k-1}^{j} + \sum_{j \in 2\ldots n}B_k^{1'}Q_k^jB_k^j\nu_{k-1}^{j}) + \tilde{u}_k^{11}\\\label{uk1s}
\end{multline} \vspace*{-5ex} 

\vspace*{-5ex}  \begin{multline}
\frac{\partial}{\partial u_k^i} \eqref{H1ks} = 0\; \Rightarrow \;\; B_k^{i'}Q_k^1(x_k^*-\tilde{x}_k^i) + R_k^{1i}(u_k^{i*}- \tilde{u}_k^{1i}) + B_k^{i'}\lambda_k \\+ \sum_{j \in 2\ldots n}B_k^{i'}Q_k^jA_k\mu_{k-1}^{j} +  (B_k^{i'}Q_k^iB_k^i + R_k^{ii})\nu_{k1}^{i} +\sum_{j \in 2\ldots n \;,\;j \not = i}B_k^{i'}Q_k^jB_k^j\nu_{k-1}^{j} = 0 \; ; \; i \in \{2\ldots n\}\\\label{nuks}
\end{multline} \vspace*{-5ex} 

\vspace*{-5ex}  \begin{multline}
\lambda_{k-1} = A_k^{'}Q_k^1(x_k^*-\tilde{x}_k^1) + A_k^{'}\lambda_k \\+ \sum_{j \in 2\ldots n}A_k^{'}Q_k^jA_k\mu_{k-1}^{j} + \sum_{j \in 2\ldots n}A_k^{'}Q_k^jB_k^j\nu_{k-1}^{j} \; ; \; \lambda_T = 0\\\label{lambdaks}
\end{multline} \vspace*{-5ex} 

\vspace*{-5ex}  \begin{multline}
\mu_k^{i} =  A_k\mu_{k-1}^{i} + B_k^i\nu_{k-1}^{i}\; ; \; \mu_0^i = 0\; ; \; i \in \{2\ldots n\}\\\label{muks}
\end{multline} \vspace*{-5ex} 

\vspace*{-5ex}  \begin{multline}
\nonumber \frac{\partial}{\partial u_i} \eqref{Hiks} = 0\; \Rightarrow \;\;\ B_k^{i'}Q_k^ix_k^* + R_k^{ii}u_k^{i*} -  B_k^{i'}Q_k^i\tilde{x}_k^{i'} - R_k^{ii}\tilde{u}_k^{ii} + B_k^{i'}p_{k}^i = 0 \; ; \; i \in \{2\ldots n\}\\
\end{multline} \vspace*{-5ex} 

\vspace*{-5ex}  \begin{multline}
u_k^{i*} =  -(R_k^{ii})^{-1}B_k^{i'}(Q_k^i(x_k^*-\tilde{x}_k^{i'}) + p_{k}^i ) + \tilde{u}_k^{ii}\\\label{ukis}
\end{multline} \vspace*{-5ex} 

\vspace*{-5ex}  \begin{multline}
p_{k-1}^{i*} = A_{k}^{'}[p_{k}^i + Q_{k}^i( x_{k}^* - \tilde{x}_{k}^i )] \; ; \; p_T^i = 0\; ; \; i \in \{2\ldots n\}\\\label{pkis}
\end{multline} \vspace*{-5ex} 

\vspace*{-5ex}  \begin{multline}
x_{k}^* = A_kx_{k-1}^* + \sum_{j\in N}B_k^ju_k^{j*}+s_k \; ; \; x_0^* = x_0\\\label{xks}
\end{multline} \vspace*{-5ex}

In the following induction argument, we will give proof that \eqref{psubs} and \eqref{lambdasubs} are valid and the recursive relations for $M_{k}^{ix}$, $M_{k}^{ij\mu}$, $m_{k}^i$, $L_{k}^x$, $L_{k}^{i\mu}$ and $l_{k}$ ($ i,j \in \{2 \ldots n\}$) (stated in the above theorem) are correct.

\vspace*{-5ex}  \begin{multline}
p_{k}^i = (M_{k}^{ix}-Q_{k}^i)x_{k}^* + \sum_{j \in 2\ldots n}M_k^{ij\mu}\mu_k^j + m_k^i  \; ; \; i \in \{2\ldots n\} \\\label{psubs}
\end{multline} \vspace*{-5ex}

\vspace*{-5ex}  \begin{multline}
\lambda_{k} = (L_{k}^x-Q_{k}^1)x_{k}^* + \sum_{j \in 2\ldots n}L_k ^{j\mu}\mu_k^j + l_k \\\label{lambdasubs}
\end{multline} \vspace*{-5ex}

\textbf{Basis:}

The induction starts at \textit {k} = \textit {T}. First we make use of the general optimality conditions for $p_k^i$ and $\lambda_k$ at stage \textit {T}.

\vspace*{-5ex}  \begin{multline}
\nonumber \eqref{pkis}_{k=T} \;\; p_{T}^i = 0 \; ; \; i \in \{2\ldots n\}\\
\end{multline} \vspace*{-5ex} 

\vspace*{-5ex}  \begin{multline}
\nonumber \eqref{lambdaks}_{k=T} \;\; \lambda_{T} = 0\\
\end{multline} \vspace*{-5ex}

Now we can substitute $p_T^i$ and $\lambda_T$ with functions affinely dependent on ($x_T^* , \mu^2_T , \ldots, \mu^n_T$).

\vspace*{-5ex}  \begin{multline}
\nonumber\eqref{psubs}_{k=T} \;\; p_{T}^i = (M_{T}^{ix}-Q_{T}^i)x_{T}^* + \sum_{j \in 2\ldots n}M_T^{ij\mu}\mu_T^j + m_T^i \; ; \; i \in \{2\ldots n\}\\
\end{multline} \vspace*{-5ex} 

\vspace*{-5ex}  \begin{multline}
\nonumber (M_{T}^{ix}-Q_{T}^i)x_{T}^* + \sum_{j \in 2\ldots n}M_T^{ij\mu}\mu_T^j + m_T^i  = 0 \; ; \; i \in \{2\ldots n\}\\
\end{multline} \vspace*{-5ex} 

\vspace*{-5ex}  \begin{multline}
\nonumber\eqref{lambdasubs}_{k=T} \;\; \lambda_{T} = (L_{T}^x-Q_{T}^1)x_{T}^* + \sum_{j \in 2\ldots n}L_T ^{j\mu}\mu_T^j + l_T\\
\end{multline} \vspace*{-5ex} 

\vspace*{-5ex}  \begin{multline}
\nonumber (L_{T}^x-Q_{T}^1)x_{T}^* + \sum_{j \in 2\ldots n}L_T ^{j\mu}\mu_T^j + l_T = 0\\
\end{multline} \vspace*{-5ex}

Comparing coefficients gives

\vspace*{-5ex}  \begin{multline}
\nonumber\eqref{Mikx}_{k=T} \;\; M_T^{ix} = Q_{T}^i \; ; \; i \in \{2\ldots n\}\\
\end{multline} \vspace*{-5ex} 

\vspace*{-5ex}  \begin{multline}
\nonumber\eqref{Mikmu}_{k=T} \;\; M_T^{ij\mu} = 0 \; ; \; i,j \in \{2\ldots n\}\\
\end{multline} \vspace*{-5ex} 

\vspace*{-5ex}  \begin{multline}
\nonumber\eqref{mik}_{k=T} \;\; m_T^{i} = 0 \; ; \; i \in \{2\ldots n\}\\
\end{multline} \vspace*{-5ex}

\vspace*{-5ex}  \begin{multline}
\nonumber\eqref{Lkx}_{k=T} \;\; L_T^x = Q_{T}^1 \\
\end{multline} \vspace*{-5ex} 

\vspace*{-5ex}  \begin{multline}
\nonumber\eqref{Lkmu}_{k=T} \;\;  L_T^{j\mu} = 0 \; ; \; j \in \{2\ldots n\}\\
\end{multline} \vspace*{-5ex} 

\vspace*{-5ex}  \begin{multline}
\nonumber\eqref{lk}_{k=T} \;\;  l_T = 0 \\
\end{multline} \vspace*{-5ex} \\\\

\textbf{Inductive step:}

As induction hypotheses, the system of equations \eqref{psubs} and equation \eqref{lambdasubs} are assumed to be true at stage \textit {l}+2. Now we have to prove that these equations are fulfilled at stage \textit {l}+1 and determine the corresponding recursive relations for $M_{l}^{ix}$, $M_{l}^{ij\mu}$, $m_{l}^i$, $L_{l}^x$, $L_{l}^{i\mu}$ and $l_{l}$ ($i,j \in \{2, \ldots, n\}$).

\vspace*{-5ex}  \begin{multline}
\nonumber\eqref{psubs}_{k=l+2} \;\; p_{l+1}^i = (M_{l+1}^{ix}-Q_{l+1}^i)x_{l+1}^* + \sum_{j \in 2\ldots n}M_{l+1}^{ij\mu}\mu_{l+1}^j + m_{l+1}^i \; ; \; i \in \{2\ldots n\}\\
\end{multline} \vspace*{-5ex} 

\vspace*{-5ex}  \begin{multline}
\nonumber\eqref{lambdasubs}_{k=l+2} \;\; \lambda_{l+1} = (L_{l+1}^x-Q_{l+1}^1)x_{l+1}^* + \sum_{j \in 2\ldots n}L_{l+1} ^{j\mu}\mu_{l+1}^j + l_{l+1}\\
\end{multline} \vspace*{-5ex}

First the induction hypotheses are used in the general optimality conditions for $p_k^i$ and $\lambda_k$ at stage \textit {l}+1.

\vspace*{-5ex}  \begin{multline}
\nonumber\eqref{pkis}_{k=l+1} \;\; p_l^{i*} = A_{l+1}^{'}[p_{l+1}^i + Q_{l+1}^i( x_{l+1}^* - \tilde{x}_{l+1}^i )] \; ; \; i \in \{2\ldots n\}\\
\end{multline} \vspace*{-5ex} 

\vspace*{-5ex}  \begin{multline}
\nonumber  p_l^{i*} = A_{l+1}^{'}[M_{l+1}^{ix}x_{l+1}^* + \sum_{j \in 2\ldots n}M_{l+1}^{ij\mu}\mu_{l+1}^j + m_{l+1}^i - Q_{l+1}^i\tilde{x}_{l+1}^i]  \; ; \; i \in \{2\ldots n\}\\
\end{multline} \vspace*{-5ex} 

\vspace*{-5ex}  \begin{multline}
\nonumber\eqref{lambdaks}_{k=l+1} \;\; \lambda_{l} = A_{l+1}^{'}Q_{l+1}^1(x_{l+1}^*-\tilde{x}_{l+1}^1) + A_{l+1}^{'}\lambda_{l+1} \\+ \sum_{j \in 2\ldots n}A_{l+1}^{'}Q_{l+1}^jA_{l+1}\mu_{l}^{j} + \sum_{j \in 2\ldots n}A_{l+1}^{'}Q_{l+1}^jB_{l+1}^j\nu_{l}^{j} \\
\end{multline} \vspace*{-5ex} 

\vspace*{-5ex}  \begin{multline}
\nonumber  \lambda_{l} = A_{l+1}^{'}(L_{l+1}^xx_{l+1}^* + \sum_{j \in 2\ldots n}L_{l+1} ^{j\mu}\mu_{l+1}^j + l_{l+1}) -A_{l+1}^{'}Q_{l+1}^1\tilde{x}_{l+1}^1\\+ \sum_{j \in 2\ldots n}A_{l+1}^{'}Q_{l+1}^jA_{l+1}\mu_{l}^{j} + \sum_{j \in 2\ldots n}A_{l+1}^{'}Q_{l+1}^jB_{l+1}^j\nu_{l}^{j} \\
\end{multline} \vspace*{-5ex}

To complete the inductive step, we have to show that the $p_l^i$ and $\lambda_l$ can be written as affine functions of the variables ($x_l^* , \mu^2_l , \ldots, \mu^n_l$). Therefore, interrelations between $x_{l+1}^*$ and ($x_l^* ,\mu_l^2 , \ldots , \mu_l^n$) and  between $\mu_{l+1}^i$ and ($x_l^* ,\mu_l^2 , \ldots , \mu_l^n$) ($i \in \{2, \ldots, n\}$) that do not depend on the controls of the players nor on costate ($p_{l+1}^i$, $\lambda_{l+1}^i$) or cocontrol ($\nu_l^i$) variables have to be deduced. To do so, first we have to substitute $u_{l+1}^{1*}, \ldots,u_{l+1}^{n*}$ in the equation stated below for the evolution of the optimal state vector $x_{l+1}^*$ by terms that are affine in ($x_l^* , x_{l+1}^*,\mu_l^2 , \ldots , \mu_l^n$) and furthermore only contain $M_{l+1}^{ix}$, $M_{l+1}^{ij\mu}$, $m_{l+1}^i$, $L_{l+1}^x$, $L_{l+1}^{i\mu}$, $l_{l+1}$ and matrices and vectors given by the game definition.

\vspace*{-5ex}  \begin{multline}
\nonumber\eqref{xks}_{k=l+1} \;\; x_{l+1}^* = A_{l+1}x_{l}^* + \sum_{j\in N}B_{l+1}^ju_{l+1}^{j*}+s_{l+1} \\
\end{multline} \vspace*{-5ex}

In the first instance the optimality condition for $u_{l+1}^{i*}$ ($i \in \{2\ldots n\}$) can be rewritten with the help of $\eqref{psubs}_{k=l+2}$, which are the induction hypotheses for $p_{l+1}^i$. 

\vspace*{-5ex}  \begin{multline}
\nonumber\eqref{ukis}_{k=l+1} \;\; u_{l+1}^{i*} =   -(R_{l+1}^{ii})^{-1}B_{l+1}^{i'}(Q_{l+1}^i(x_{l+1}^*-\tilde{x}_{l+1}^{i'}) + p_{l+1}^i ) + \tilde{u}_{l+1}^{ii}  \; ; \; i \in \{2\ldots n\}\\
\end{multline} \vspace*{-5ex} 

\vspace*{-5ex}  \begin{multline}
\nonumber u_{l+1}^{i*} =   -(R_{l+1}^{ii})^{-1}B_{l+1}^{i'}(M_{l+1}^{ix}x_{l+1}^* + \sum_{j \in 2\ldots n}M_{l+1}^{ij\mu}\mu_{l+1}^j \\+ m_{l+1}^i - Q_{l+1}^i\tilde{x}_{l+1}^{i'}) + \tilde{u}_{l+1}^{ii}  \; ; \; i \in \{2\ldots n\}\\
\end{multline} \vspace*{-5ex} 

As a next step the stage indices of $\mu^i$ have to be reduced from \textit {l}+1 to \textit {l} with the help of the optimality conditions for $\mu_{l+1}^i$.

\vspace*{-5ex}  \begin{multline}
\nonumber\eqref{muks}_{k=l+1} \;\; \mu_{l+1}^i = A_{l+1}\mu_{l}^{i} + B_{l+1}^i\nu_{l}^{i} \; ; \; i \in \{2\ldots n\}\\
\end{multline} \vspace*{-5ex} 

\vspace*{-5ex}  \begin{multline}
u_{l+1}^{i*} =   -(R_{l+1}^{ii})^{-1}B_{l+1}^{i'}(M_{l+1}^{ix}x_{l+1}^* + \sum_{j \in 2\ldots n}M_{l+1}^{ij\mu}\\(A_{l+1}\mu_{l}^{j} + B_{l+1}^j\nu_{l}^{j}) + m_{l+1}^i - Q_{l+1}^i\tilde{x}_{l+1}^{i'}) + \tilde{u}_{l+1}^{ii}  \; ; \; i \in \{2\ldots n\}\\\label{uinusubs}
\end{multline} \vspace*{-5ex}

Moreover, we have to substitute $\nu_{l}^{j}$ ($j \in \{2 \ldots n\}$). For that purpose, $\nu_{l}^{j}$ has to be explicated from $\eqref{nuks}_{l+1}$. As a start $u_{l+1}^{i*}$ is replaced using \eqref{uinusubs} and the $\nu_l^j$ ($j \in \{2 , \ldots , n\}$) are abstracted in one term.

\vspace*{-5ex}  \begin{multline}
\nonumber\eqref{nuks}_{k=l+1} \;\; B_{l+1}^{i'}Q_{l+1}^1(x_{l+1}^*-\tilde{x}_{l+1}^i) + R_{l+1}^{1i}(u_{l+1}^{i*} - \tilde{u}_{l+1}^{1i}) + B_{l+1}^{i'}\lambda_{l+1} + \sum_{j \in 2\ldots n}B_{l+1}^{i'}Q_{l+1}^jA_{l+1}\mu_{l}^{j} \\+  (B_{l+1}^{i'}Q_{l+1}^iB_{l+1}^i + R_{l+1}^{ii})\nu_{l}^{i} +\sum_{j \in 2\ldots n \;,\;j \not = i}B_{l+1}^{i'}Q_{l+1}^jB_{l+1}^j\nu_{l}^{j} = 0  \; ; \; i \in \{2\ldots n\}\\
\end{multline} \vspace*{-5ex} 

\vspace*{-5ex}  \begin{multline}
\nonumber B_{l+1}^{i'}Q_{l+1}^1(x_{l+1}^*-\tilde{x}_{l+1}^i) + R_{l+1}^{1i}(-(R_{l+1}^{ii})^{-1}B_{l+1}^{i'}\\(M_{l+1}^{ix}x_{l+1}^* + \sum_{j \in 2\ldots n}M_{l+1}^{ij\mu}(A_{l+1}\mu_{l}^{j} + B_{l+1}^j\nu_{l}^{j}) + m_{l+1}^i - Q_{l+1}^i\tilde{x}_{l+1}^{i'}) \\+ \tilde{u}_{l+1}^{ii} - \tilde{u}_{l+1}^{1i}) + B_{l+1}^{i'}\lambda_{l+1} + \sum_{j \in 2\ldots n}B_{l+1}^{i'}Q_{l+1}^jA_{l+1}\mu_{l}^{j} \\+  (B_{l+1}^{i'}Q_{l+1}^iB_{l+1}^i + R_{l+1}^{ii})\nu_{l}^{i} +\sum_{j \in 2\ldots n \;,\;j \not = i}B_{l+1}^{i'}Q_{l+1}^jB_{l+1}^j\nu_{l}^{j} = 0  \; ; \; i \in \{2\ldots n\}\\
\end{multline} \vspace*{-5ex} 

As a next step $\lambda_{l+1}$ is substituted with the help of $\eqref{lambdasubs}_{k=l+1}$.

\vspace*{-5ex}  \begin{multline}
\nonumber  -B_{l+1}^{i'}Q_{l+1}^1\tilde{x}_{l+1}^i + R_{l+1}^{1i}(-(R_{l+1}^{ii})^{-1}B_{l+1}^{i'}(M_{l+1}^{ix}x_{l+1}^* + \sum_{j \in 2\ldots n}M_{l+1}^{ij\mu}A_{l+1}\mu_{l}^{j} + m_{l+1}^i \\- Q_{l+1}^i\tilde{x}_{l+1}^{i'}) + \tilde{u}_{l+1}^{ii}- \tilde{u}_{l+1}^{1i}) + B_{l+1}^{i'}(L_{l+1}^xx_{l+1}^* + \sum_{j \in 2\ldots n}L_{l+1} ^{j\mu}\mu_{l+1}^j + l_{l+1}) \\+ \sum_{j \in 2\ldots n}B_{l+1}^{i'}Q_{l+1}^jA_{l+1}\mu_{l}^{j} +  (B_{l+1}^{i'}Q_{l+1}^iB_{l+1}^i + R_{l+1}^{ii} - R_{l+1}^{1i}(R_{l+1}^{ii})^{-1}B_{l+1}^{i'}M_{l+1}^{ii\mu}B_{l+1}^i)\nu_{l}^{i} \\+\sum_{j \in 2\ldots n \;,\;j \not = i}(B_{l+1}^{i'}Q_{l+1}^jB_{l+1}^j -R_{l+1}^{1i}(R_{l+1}^{ii})^{-1}B_{l+1}^{i'}M_{l+1}^{ij\mu}B_{l+1}^j)\nu_{l}^{j} = 0  \; ; \; i \in \{2\ldots n\}\\
\end{multline} \vspace*{-5ex} 

Eventually the $\mu_{l+1}^i$ are replaced making use of $\eqref{muks}_{k=l+1}$ and the $\nu_l^j$ ($j \in \{2 , \ldots , n\}$) are again abstracted in one term.

\vspace*{-5ex}  \begin{multline}
\nonumber  -B_{l+1}^{i'}Q_{l+1}^1\tilde{x}_{l+1}^i + R_{l+1}^{1i}(-(R_{l+1}^{ii})^{-1}B_{l+1}^{i'}(M_{l+1}^{ix}x_{l+1}^* + \sum_{j \in 2\ldots n}M_{l+1}^{ij\mu}A_{l+1}\mu_{l}^{j} + m_{l+1}^i \\- Q_{l+1}^i\tilde{x}_{l+1}^{i'}) + \tilde{u}_{l+1}^{ii}- \tilde{u}_{l+1}^{1i}) + B_{l+1}^{i'}(L_{l+1}^xx_{l+1}^* + \sum_{j \in 2\ldots n}L_{l+1} ^{j\mu}(A_{l+1}\mu_{l}^{j} + B_{l+1}^j\nu_{l}^{j}) + l_{l+1}) \\+ \sum_{j \in 2\ldots n}B_{l+1}^{i'}Q_{l+1}^jA_{l+1}\mu_{l}^{j} +  (B_{l+1}^{i'}Q_{l+1}^iB_{l+1}^i + R_{l+1}^{ii} - R_{l+1}^{1i}(R_{l+1}^{ii})^{-1}B_{l+1}^{i'}M_{l+1}^{ii\mu}B_{l+1}^i)\nu_{l}^{i} \\+\sum_{j \in 2\ldots n \;,\;j \not = i}(B_{l+1}^{i'}Q_{l+1}^jB_{l+1}^j -R_{l+1}^{1i}(R_{l+1}^{ii})^{-1}B_{l+1}^{i'}M_{l+1}^{ij\mu}B_{l+1}^j)\nu_{l}^{j} = 0  \; ; \; i \in \{2\ldots n\}\\
\end{multline} \vspace*{-5ex} 

\vspace*{-5ex}  \begin{multline}
\nonumber  -B_{l+1}^{i'}Q_{l+1}^1\tilde{x}_{l+1}^i + R_{l+1}^{1i}(-(R_{l+1}^{ii})^{-1}B_{l+1}^{i'}(M_{l+1}^{ix}x_{l+1}^* + \sum_{j \in 2\ldots n}M_{l+1}^{ij\mu}A_{l+1}\mu_{l}^{j} + m_{l+1}^i \\- Q_{l+1}^i\tilde{x}_{l+1}^{i'}) + \tilde{u}_{l+1}^{ii}- \tilde{u}_{l+1}^{1i}) + B_{l+1}^{i'}(L_{l+1}^xx_{l+1}^* + \sum_{j \in 2\ldots n}L_{l+1} ^{j\mu}A_{l+1}\mu_{l}^{j} + l_{l+1}) \\+ \sum_{j \in 2\ldots n}B_{l+1}^{i'}Q_{l+1}^jA_{l+1}\mu_{l}^{j} +  (B_{l+1}^{i'}(Q_{l+1}^i+L_{l+1} ^{i\mu})B_{l+1}^i + R_{l+1}^{ii} - R_{l+1}^{1i}(R_{l+1}^{ii})^{-1}B_{l+1}^{i'}M_{l+1}^{ii\mu}B_{l+1}^i)\nu_{l}^{i} \\+\sum_{j \in 2\ldots n \;,\;j \not = i}(B_{l+1}^{i'}(Q_{l+1}^j+L_{l+1} ^{j\mu})B_{l+1}^j -R_{l+1}^{1i}(R_{l+1}^{ii})^{-1}B_{l+1}^{i'}M_{l+1}^{ij\mu}B_{l+1}^j)\nu_{l}^{j} = 0  \; ; \; i \in \{2\ldots n\}\\
\end{multline} \vspace*{-5ex} 

The above equations only contain constant expressions and terms linear in $x_{l+1}$ or $\mu_{l}^j$ ($j \in \{2 , \ldots , n\}$). This fact justifies the following substitutions:

\vspace*{-5ex}  \begin{multline}
\nu_{l}^i = N_{l}^{ix}x_{l+1}^* + \sum_{j \in 2\ldots n }N_{l}^{ij\mu}\mu_{l}^j + n_l^i   \; ; \; i \in \{2\ldots n\}\\\label{nusubs}
\end{multline} \vspace*{-5ex} 

\vspace*{-5ex}  \begin{multline}
  -B_{l+1}^{i'}Q_{l+1}^1\tilde{x}_{l+1}^i + R_{l+1}^{1i}(-(R_{l+1}^{ii})^{-1}B_{l+1}^{i'}(M_{l+1}^{ix}x_{l+1}^* + \sum_{j \in 2\ldots n}M_{l+1}^{ij\mu}A_{l+1}\mu_{l}^{j} \\+ m_{l+1}^i - Q_{l+1}^i\tilde{x}_{l+1}^{i'}) + \tilde{u}_{l+1}^{ii}- \tilde{u}_{l+1}^{1i}) + B_{l+1}^{i'}(L_{l+1}^xx_{l+1}^* + \sum_{j \in 2\ldots n}L_{l+1} ^{j\mu}A_{l+1}\mu_{l}^{j} + l_{l+1}) \\+ \sum_{j \in 2\ldots n}B_{l+1}^{i'}Q_{l+1}^jA_{l+1}\mu_{l}^{j} +  (B_{l+1}^{i'}(Q_{l+1}^i+L_{l+1} ^{i\mu})B_{l+1}^i + R_{l+1}^{ii} - R_{l+1}^{1i}(R_{l+1}^{ii})^{-1}B_{l+1}^{i'}M_{l+1}^{ii\mu}B_{l+1}^i)\\(N_{l}^{ix}x_{l+1}^* + \sum_{m \in 2\ldots n }N_{l}^{im\mu}\mu_{l}^m + n_l^i) +\sum_{j \in 2\ldots n \;,\;j \not = i}(B_{l+1}^{i'}(Q_{l+1}^j+L_{l+1} ^{j\mu})B_{l+1}^j\\ -R_{l+1}^{1i}(R_{l+1}^{ii})^{-1}B_{l+1}^{i'}M_{l+1}^{ij\mu}B_{l+1}^j)(N_{l}^{jx}x_{l+1}^* + \sum_{m \in 2\ldots n }N_{l}^{jm\mu}\mu_{l}^m + n_l^j) = 0  \; ; \; i \in \{2\ldots n\}\\\label{nulg}
\end{multline} \vspace*{-5ex} 

Comparing coefficients gives the following systems of equations that are assumed to admit unique solutions $N_{l}^{ix}$, $N_{l}^{im\mu}$ and $n_{l}^{i}$ ($i,m \in \{2,\ldots,n\}$).

\vspace*{-5ex}  \begin{multline}
\nonumber\eqref{nulg}_{x_{l+1}^*}=\eqref{Nkx}_{k=l} \;\;  -R_{l+1}^{1i}(R_{l+1}^{ii})^{-1}B_{l+1}^{i'}M_{l+1}^{ix} + B_{l+1}^{i'}L_{l+1}^x +  (B_{l+1}^{i'}(Q_{l+1}^i+L_{l+1} ^{i\mu})B_{l+1}^i \\+ R_{l+1}^{ii} - R_{l+1}^{1i}(R_{l+1}^{ii})^{-1}B_{l+1}^{i'}M_{l+1}^{ii\mu}B_{l+1}^i)N_{l}^{ix}  + \sum_{j \in 2\ldots n \;,\;j \not = i}(B_{l+1}^{i'}(Q_{l+1}^j+L_{l+1} ^{j\mu})B_{l+1}^j\\ -R_{l+1}^{1i}(R_{l+1}^{ii})^{-1}B_{l+1}^{i'}M_{l+1}^{ij\mu}B_{l+1}^j)N_{l}^{jx} = 0  \; ; \; i \in \{2\ldots n\}\\
\end{multline} \vspace*{-5ex} 

\vspace*{-5ex}  \begin{multline}
\nonumber\eqref{nulg}_{\mu_l^m}=\eqref{Nkmu}_{k=l} \;\; -R_{l+1}^{1i}(R_{l+1}^{ii})^{-1}B_{l+1}^{i'}M_{l+1}^{im\mu}A_{l+1} + B_{l+1}^{i'}L_{l+1} ^{m\mu}A_{l+1} \\+ B_{l+1}^{i'}Q_{l+1}^mA_{l+1} +  (B_{l+1}^{i'}(Q_{l+1}^i+L_{l+1} ^{i\mu})B_{l+1}^i + R_{l+1}^{ii} \\- R_{l+1}^{1i}(R_{l+1}^{ii})^{-1}B_{l+1}^{i'}M_{l+1}^{ii\mu}B_{l+1}^i)N_{l}^{im\mu} + \sum_{j \in 2\ldots n \;,\;j \not = i}(B_{l+1}^{i'}(Q_{l+1}^j+L_{l+1} ^{j\mu})B_{l+1}^j\\ -R_{l+1}^{1i}(R_{l+1}^{ii})^{-1}B_{l+1}^{i'}M_{l+1}^{ij\mu}B_{l+1}^j)N_{l}^{jm\mu} = 0 \; ; \; i,m \in \{2\ldots n\}\\
\end{multline} \vspace*{-5ex} 

\vspace*{-5ex}  \begin{multline}
\nonumber\eqref{nulg}_{const.}=\eqref{nk}_{k=l} \;\;  -B_{l+1}^{i'}Q_{l+1}^1\tilde{x}_{l+1}^i + R_{l+1}^{1i}(-(R_{l+1}^{ii})^{-1}B_{l+1}^{i'}( m_{l+1}^i - Q_{l+1}^i\tilde{x}_{l+1}^{i'}) + \tilde{u}_{l+1}^{ii} \\- \tilde{u}_{l+1}^{1i}) + B_{l+1}^{i'}l_{l+1}  +  (B_{l+1}^{i'}(Q_{l+1}^i+L_{l+1} ^{i\mu})B_{l+1}^i + R_{l+1}^{ii} - R_{l+1}^{1i}(R_{l+1}^{ii})^{-1}B_{l+1}^{i'}M_{l+1}^{ii\mu}B_{l+1}^i)n_l^i \\+ \sum_{j \in 2\ldots n \;,\;j \not = i}(B_{l+1}^{i'}(Q_{l+1}^j+L_{l+1} ^{j\mu})B_{l+1}^j -R_{l+1}^{1i}(R_{l+1}^{ii})^{-1}B_{l+1}^{i'}M_{l+1}^{ij\mu}B_{l+1}^j)n_l^j = 0  \; ; \; i \in \{2\ldots n\}\\
\end{multline} \vspace*{-5ex}

Using the elaborated relation for $\nu_{l}^{j}$ ($j \in \{2 \ldots n\}$) in \eqref{uinusubs} yields

\vspace*{-5ex}  \begin{multline}
\nonumber u_{l+1}^{i*} =   -(R_{l+1}^{ii})^{-1}B_{l+1}^{i'}(M_{l+1}^{ix}x_{l+1}^* + \sum_{j \in 2\ldots n}M_{l+1}^{ij\mu}(A_{l+1}\mu_{l}^{j} + B_{l+1}^j(N_{l}^{jx}x_{l+1}^* \\+ \sum_{m \in 2\ldots n }N_{l}^{jm\mu}\mu_{l}^m + n_l^j)) + m_{l+1}^i - Q_{l+1}^i\tilde{x}_{l+1}^{i'}) + \tilde{u}_{l+1}^{ii}  \; ; \; i \in \{2\ldots n\}\\
\end{multline} \vspace*{-5ex} 

The structure of the above equations justifies the following substitutions

\vspace*{-5ex}  \begin{multline}
u_{l+1}^{i*} =  T_{l+1}^{ix}x_{l+1}^* + \sum_{j \in 2\ldots n }T_{l+1}^{ij\mu}\mu_{l}^j + t_{l+1}^i   \; ; \; i \in \{2\ldots n\}\\\label{uisubs}
\end{multline} \vspace*{-5ex} 

\vspace*{-5ex}  \begin{multline}
 T_{l+1}^{ix}x_{l+1}^* + \sum_{j \in 2\ldots n }T_{l+1}^{ij\mu}\mu_{l}^j + t_{l+1}^i =   -(R_{l+1}^{ii})^{-1}B_{l+1}^{i'}(M_{l+1}^{ix}x_{l+1}^* \\+ \sum_{j \in 2\ldots n}M_{l+1}^{ij\mu}(A_{l+1}\mu_{l}^{j} + B_{l+1}^j(N_{l}^{jx}x_{l+1}^* + \sum_{m \in 2\ldots n }N_{l}^{jm\mu}\mu_{l}^m + n_l^j)) \\+ m_{l+1}^i - Q_{l+1}^i\tilde{x}_{l+1}^{i'}) + \tilde{u}_{l+1}^{ii}  \; ; \; i \in \{2\ldots n\}\\\label{uilg}
\end{multline} \vspace*{-5ex} 

By comparing coefficients it follows that

\vspace*{-5ex}  \begin{multline}
\nonumber\eqref{uilg}_{x_{l+1}^*}=\eqref{Tkx}_{k=l} \;\;  T_{l+1}^{ix} = -(R_{l+1}^{ii})^{-1}B_{l+1}^{i'}(M_{l+1}^{ix} + \sum_{j \in 2\ldots n}M_{l+1}^{ij\mu}B_{l+1}^jN_{l}^{jx})  \; ; \; i \in \{2\ldots n\}\\
\end{multline} \vspace*{-5ex} 

\vspace*{-5ex}  \begin{multline}
\nonumber\eqref{uilg}_{\mu_l^m}=\eqref{Tkmu}_{k=l} \;\; T_{l+1}^{im\mu} =  -(R_{l+1}^{ii})^{-1}B_{l+1}^{i'}(M_{l+1}^{im\mu}A_{l+1} \\+ \sum_{j \in 2\ldots n}M_{l+1}^{ij\mu}B_{l+1}^jN_{l}^{jm\mu})\; ; \; i,m \in \{2\ldots n\}\\
\end{multline} \vspace*{-5ex} 

\vspace*{-5ex}  \begin{multline}
\nonumber\eqref{uilg}_{const.}=\eqref{tk}_{k=l} \;\; t_{l+1}^{i} =  -(R_{l+1}^{ii})^{-1}B_{l+1}^{i'}(\sum_{j \in 2\ldots n}M_{l+1}^{ij\mu} B_{l+1}^jn_l^j \\+ m_{l+1}^i - Q_{l+1}^i\tilde{x}_{l+1}^{i'}) + \tilde{u}_{l+1}^{ii}  \; ; \; i \in \{2\ldots n\}\\
\end{multline} \vspace*{-5ex}

The optimality condition for $u_{l+1}^{1*}$ can be rewritten with the help of $\eqref{lambdasubs}_{k=l+1}$ and $\eqref{muks}_{k=l+1}$. These are the induction hypothesis for $\lambda_{l+1}$ and the conditions for the optimal evolution of $\mu^i_{l+1}$. 

\vspace*{-5ex}  \begin{multline}
\nonumber\eqref{uk1s}_{k=l+1} \;\; u_{l+1}^{1*} =  -(R_{l+1}^{11})^{-1} (B_{l+1}^{1'}Q_{l+1}^1 (x_{l+1}^*-\tilde{x}_{l+1}^1)  + B_{l+1}^{1'}\lambda_{l+1} \\+ \sum_{j \in 2\ldots n}B_{l+1}^{1'}Q_{l+1}^jA_{l+1} \mu_{l}^{j} + \sum_{j \in 2\ldots n}B_{l+1}^{1'}Q_{l+1}^jB_{l+1}^j\nu_{l}^{j}) + \tilde{u}_{l+1}^{11} \\
\end{multline} \vspace*{-5ex} 

\vspace*{-5ex}  \begin{multline}
\nonumber u_{l+1}^{1*} =  -(R_{l+1}^{11})^{-1} (-B_{l+1}^{1'}Q_{l+1}^1 \tilde{x}_{l+1}^1  + B_{l+1}^{1'}(L_{l+1}^xx_{l+1}^* + \sum_{j \in 2\ldots n}L_{l+1} ^{j\mu}\mu_{l+1}^j + l_{l+1}) \\+ \sum_{j \in 2\ldots n}B_{l+1}^{1'}Q_{l+1}^jA_{l+1} \mu_{l}^{j} + \sum_{j \in 2\ldots n}B_{l+1}^{1'}Q_{l+1}^jB_{l+1}^j\nu_{l}^{j}) + \tilde{u}_{l+1}^{11} \\
\end{multline} \vspace*{-5ex} 

\vspace*{-5ex}  \begin{multline}
\nonumber u_{l+1}^{1*} =  -(R_{l+1}^{11})^{-1} (-B_{l+1}^{1'}Q_{l+1}^1 \tilde{x}_{l+1}^1  + B_{l+1}^{1'}(L_{l+1}^xx_{l+1}^* + \sum_{j \in 2\ldots n}L_{l+1} ^{j\mu}(A_{l+1}\mu_{l}^{j} + B_{l+1}^j\nu_{l}^{j}) \\+ l_{l+1}) + \sum_{j \in 2\ldots n}B_{l+1}^{1'}Q_{l+1}^jA_{l+1} \mu_{l}^{j} + \sum_{j \in 2\ldots n}B_{l+1}^{1'}Q_{l+1}^jB_{l+1}^j\nu_{l}^{j}) + \tilde{u}_{l+1}^{11} \\
\end{multline} \vspace*{-5ex} 

Now $\nu_l^i$ ($i \in \{2 , \ldots, n\}$) can be replaced by the relations deduced above stated in \eqref{nusubs}. 

\vspace*{-5ex}  \begin{multline}
\nonumber u_{l+1}^{1*} =  -(R_{l+1}^{11})^{-1} (-B_{l+1}^{1'}Q_{l+1}^1\tilde{x}_{l+1}^1  + B_{l+1}^{1'}(L_{l+1}^xx_{l+1}^* + \sum_{j \in 2\ldots n}L_{l+1} ^{j\mu}(A_{l+1}\mu_{l}^{j} \\+ B_{l+1}^j(N_{l}^{jx}x_{l+1}^* + \sum_{m \in 2\ldots n }N_{l}^{jm\mu}\mu_{l}^m + n_l^j)) + l_{l+1}) + \sum_{j \in 2\ldots n}B_{l+1}^{1'}Q_{l+1}^jA_{l+1}\mu_{l}^{j} \\+ \sum_{j \in 2\ldots n}B_{l+1}^{1'}Q_{l+1}^jB_{l+1}^j(N_{l}^{jx}x_{l+1}^* + \sum_{m \in 2\ldots n }N_{l}^{jm\mu}\mu_{l}^m + n_l^j)) + \tilde{u}_{l+1}^{11} \\
\end{multline} \vspace*{-5ex} 

The structure of the above equation justifies the following substitution

\vspace*{-5ex}  \begin{multline}
u_{l+1}^{1*} =  W_{l+1}^{x}x_{l+1}^* + \sum_{j \in 2\ldots n }W_{l+1}^{j\mu}\mu_{l}^j + w_{l+1} \\\label{u1subs}
\end{multline} \vspace*{-5ex} 

\vspace*{-5ex}  \begin{multline}
W_{l+1}^{x}x_{l+1}^* + \sum_{j \in 2\ldots n }W_{l+1}^{j\mu}\mu_{l}^j + w_{l+1} =  -(R_{l+1}^{11})^{-1} (-B_{l+1}^{1'}Q_{l+1}^1\tilde{x}_{l+1}^1  \\+ B_{l+1}^{1'}(L_{l+1}^xx_{l+1}^* + \sum_{j \in 2\ldots n}L_{l+1}^{j\mu}(A_{l+1}\mu_{l}^{j} + B_{l+1}^j(N_{l}^{jx}x_{l+1}^* \\+ \sum_{m \in 2\ldots n }N_{l}^{jm\mu}\mu_{l}^m + n_l^j)) + l_{l+1}) + \sum_{j \in 2\ldots n}B_{l+1}^{1'}Q_{l+1}^jA_{l+1}\mu_{l}^{j} \\+ \sum_{j \in 2\ldots n}B_{l+1}^{1'}Q_{l+1}^jB_{l+1}^j(N_{l}^{jx}x_{l+1}^* + \sum_{m \in 2\ldots n }N_{l}^{jm\mu}\mu_{l}^m + n_l^j)) + \tilde{u}_{l+1}^{11} \\\label{u1lg}
\end{multline} \vspace*{-5ex} 

Comparing coefficients gives

\vspace*{-5ex}  \begin{multline}
\nonumber\eqref{u1lg}_{x_{l+1}^*}=\eqref{Wkx}_{k=l} \;\;  W_{l+1}^{x} =   -(R_{l+1}^{11})^{-1} (B_{l+1}^{1'}(L_{l+1}^x + \sum_{j \in 2\ldots n}L_{l+1}^{j\mu}B_{l+1}^jN_{l}^{jx}) \\+ \sum_{j \in 2\ldots n}B_{l+1}^{1'}Q_{l+1}^jB_{l+1}^jN_{l}^{jx})\\
\end{multline} \vspace*{-5ex} 

\vspace*{-5ex}  \begin{multline}
\nonumber\eqref{u1lg}_{\mu_l^m}=\eqref{Wkmu}_{k=l} \;\; W_{l+1}^{m\mu} =  -(R_{l+1}^{11})^{-1} (B_{l+1}^{1'}(L_{l+1}^{m\mu}A_{l+1} + \sum_{j \in 2\ldots n}L_{l+1}^{j\mu}B_{l+1}^jN_{l}^{jm\mu}) \\+ B_{l+1}^{1'}Q_{l+1}^mA_{l+1} + \sum_{j \in 2\ldots n}B_{l+1}^{1'}Q_{l+1}^jB_{l+1}^jN_{l}^{jm\mu})  \; ; \; m \in \{2\ldots n\}\\
\end{multline} \vspace*{-5ex} 

\vspace*{-5ex}  \begin{multline}
\nonumber\eqref{u1lg}_{const.}=\eqref{wk}_{k=l} \;\; w_{l+1} =  -(R_{l+1}^{11})^{-1} (-B_{l+1}^{1'}Q_{l+1}^1\tilde{x}_{l+1}^1  + B_{l+1}^{1'}(\sum_{j \in 2\ldots n}L_{l+1}^{j\mu}B_{l+1}^jn_l^j \\+ l_{l+1}) + \sum_{j \in 2\ldots n}B_{l+1}^{1'}Q_{l+1}^jB_{l+1}^jn_l^j) + \tilde{u}_{l+1}^{11} \\
\end{multline} \vspace*{-5ex}

At this point it is possible to replace the control variables in the optimal state equation at stage l+1 by terms affine in ($x_{l+1}^* ,\mu_l^2 , \ldots , \mu_l^n$).

\vspace*{-5ex}  \begin{multline}
\nonumber\eqref{xks}_{k=l+1} \;\; x_{l+1}^* = A_{l+1}x_{l}^* + B_{l+1}^1u_{l+1}^{1*} + \sum_{j\in \{2...n\}}B_{l+1}^ju_{l+1}^{j*}+s_{l+1} \\
\end{multline} \vspace*{-5ex} 

\vspace*{-5ex}  \begin{multline}
\nonumber x_{l+1}^* = A_{l+1}x_{l}^* + B_{l+1}^1(W_{l+1}^{x}x_{l+1}^* + \sum_{j \in 2\ldots n }W_{l+1}^{j\mu}\mu_{l}^j + w_{l+1}) \\+ \sum_{j\in 2\ldots n}B_{l+1}^j(T_{l+1}^{jx}x_{l+1}^* + \sum_{m \in 2\ldots n }T_{l+1}^{jm\mu}\mu_{l}^m + t_{l+1}^j) + s_{l+1} \\
\end{multline} \vspace*{-5ex} 

Making $x_{l+1}^*$ explicit yields

\vspace*{-5ex}  \begin{multline}
\nonumber  (I - B_{l+1}^1W_{l+1}^{x} - \sum_{j\in \{2...n\}}B_{l+1}^jT_{l+1}^{jx})x_{l+1}^* = A_{l+1}x_{l}^* + B_{l+1}^1(\sum_{j \in 2\ldots n }W_{l+1}^{j\mu}\mu_{l}^j \\+ w_{l+1}) + \sum_{j\in 2\ldots n}B_{l+1}^j(\sum_{m \in 2\ldots n }T_{l+1}^{jm\mu}\mu_{l}^m + t_{l+1}^j) + s_{l+1} \\
\end{multline} \vspace*{-5ex} 

\vspace*{-5ex}  \begin{multline}
\nonumber  x_{l+1}^* = (I - B_{l+1}^1W_{l+1}^{x} - \sum_{j\in \{2...n\}}B_{l+1}^jT_{l+1}^{jx})^{-1}(A_{l+1}x_{l}^* + B_{l+1}^1(\sum_{j \in 2\ldots n }W_{l+1}^{j\mu}\mu_{l}^j \\+ w_{l+1}) + \sum_{j\in 2\ldots n}B_{l+1}^j(\sum_{m \in 2\ldots n }T_{l+1}^{jm\mu}\mu_{l}^m + t_{l+1}^j) + s_{l+1}) \\
\end{multline} \vspace*{-5ex} 

The structure of the above equation justifies the following substitution

\vspace*{-5ex}  \begin{multline}
\nonumber\eqref{xopts}_{k=l} \;\; x_{l+1}^* = \Phi_{l}^xx_{l}^* + \sum_{j \in 2\ldots n }\Phi_{l}^{j\mu}\mu_{l}^j + \phi_l\\
\end{multline} \vspace*{-5ex} 

\vspace*{-5ex}  \begin{multline}
 (I - B_{l+1}^1W_{l+1}^{x} - \sum_{j\in \{2...n\}}B_{l+1}^jT_{l+1}^{jx})(\Phi_{l}^xx_{l}^* + \sum_{j \in 2\ldots n }\Phi_{l}^{j\mu}\mu_{l}^j + \phi_l) = A_{l+1}x_{l}^* \\+ B_{l+1}^1(\sum_{j \in 2\ldots n }W_{l+1}^{j\mu}\mu_{l}^j + w_{l+1}) + \sum_{j\in 2\ldots n}B_{l+1}^j(\sum_{m \in 2\ldots n }T_{l+1}^{jm\mu}\mu_{l}^m + t_{l+1}^j) + s_{l+1} \\\label{xlg}
\end{multline} \vspace*{-5ex} 

By comparing coefficients it follows that

\vspace*{-5ex}  \begin{multline}
\nonumber\eqref{xlg}_{x_{l}^*}=\eqref{Phikx}_{k=l} \;\; \Phi_{l}^x = (I - B_{l+1}^1W_{l+1}^{x} - \sum_{j\in 2\ldots n}B_{l+1}^jT_{l+1}^{jx})^{-1}A_{l+1}\\
\end{multline} \vspace*{-5ex} 

\vspace*{-5ex}  \begin{multline}
\nonumber\eqref{xlg}_{\mu_{l}^i}=\eqref{Phikmu}_{k=l} \;\; \Phi_{l}^{i\mu} = (I - B_{l+1}^1W_{l+1}^{x} - \sum_{j\in 2\ldots n}B_{l+1}^jT_{l+1}^{jx})^{-1}\\(B_{l+1}^1W_{l+1}^{i\mu} + \sum_{j\in 2\ldots n}B_{l+1}^jT_{l+1}^{ji\mu} ) \; ; \; i \in \{2\ldots n\}\\
\end{multline} \vspace*{-5ex} 

\vspace*{-5ex}  \begin{multline}
\nonumber\eqref{xlg}_{const.}=\eqref{phik}_{k=l} \;\; \phi_{l} = (I - B_{l+1}^1W_{l+1}^{x} - \sum_{j\in 2\ldots n}B_{l+1}^jT_{l+1}^{jx})^{-1}\\(B_{l+1}^1w_{l+1} + \sum_{j\in 2\ldots n}B_{l+1}^jt_{l+1}^j + s_{l+1}) \\
\end{multline} \vspace*{-5ex}

As a next step affine relations between ($x_{l}^*, \mu_{l+1}^i ,\mu_l^2 , \ldots , \mu_l^n$) ($i \in \{2, \ldots, n\}$) are derived. Therefore $\nu_l^i$ is substituted in $\eqref{muks}_{k=l+1}$ by making use of \eqref{nusubs}.

\vspace*{-5ex}  \begin{multline}
\nonumber \eqref{muks}_{k=l+1} \; \; \mu_{l+1}^i =  A_{l+1}\mu_{l}^i + B_{l+1}^i\nu_{l}^i \; ; \; i \in \{2\ldots n\}\\
\end{multline} \vspace*{-5ex} 

\vspace*{-5ex}  \begin{multline}
\nonumber \mu_{l+1}^i =  A_{l+1}\mu_{l}^i + B_{l+1}^i(N_{l}^{ix}x_{l+1}^* + \sum_{j \in 2\ldots n }N_{l}^{ij\mu}\mu_{l}^j + n_l^i) \; ; \; i \in \{2\ldots n\}\\
\end{multline} \vspace*{-5ex} 

Using $\eqref{xopts}_{k=l}$ yields

\vspace*{-5ex}  \begin{multline}
\nonumber \mu_{l+1}^i =  A_{l+1}\mu_{l}^i + B_{l+1}^i(N_{l}^{ix}(\Phi_{l}^xx_{l}^* + \sum_{j \in 2\ldots n }\Phi_{l}^{j\mu}\mu_{l}^j + \phi_l) + \sum_{j \in 2\ldots n }N_{l}^{ij\mu}\mu_{l}^j + n_l^i) \; ; \; i \in \{2\ldots n\}\\
\end{multline} \vspace*{-5ex} 

The structure of the above equations justifies the following substitutions

\vspace*{-5ex}  \begin{multline}
\nonumber\eqref{muopts}_{k=l} \;\; \mu_{l+1}^i =  \Psi_{l}^{ix}x_{l}^* + \sum_{j \in 2\ldots n }\Psi_{l}^{ij\mu}\mu_{l}^j + \psi_l^i \; ; \; i \in \{2\ldots n\}\\
\end{multline} \vspace*{-5ex} 

\vspace*{-5ex}  \begin{multline}
\Psi_{l}^xx_{l}^* + \sum_{j \in 2\ldots n }\Psi_{l}^{j\mu}\mu_{l}^j + \psi_l =  A_{l+1}\mu_{l}^i + B_{l+1}^i(N_{l}^{ix}(\Phi_{l}^xx_{l}^* \\+ \sum_{j \in 2\ldots n }\Phi_{l}^{j\mu}\mu_{l}^j + \phi_l) + \sum_{j \in 2\ldots n }N_{l}^{ij\mu}\mu_{l}^j + n_l^i) \; ; \; i \in \{2\ldots n\}\\\label{mulg}
\end{multline} \vspace*{-5ex} 

Comparing coefficients gives

\vspace*{-5ex}  \begin{multline}
\nonumber\eqref{mulg}_{x_{l}^*}=\eqref{Psikx}_{k=l} \;\; \Psi_{l}^{ix} = B_{l+1}^iN_{l}^{ix}\Phi_{l}^x \; ; \; i \in \{2\ldots n\}\\
\end{multline} \vspace*{-5ex} 

\vspace*{-5ex}  \begin{multline}
\nonumber\eqref{mulg}_{\mu_{l}^i}=\eqref{Psikmui}_{k=l} \;\; \Psi_{l}^{ii} =  A_{l+1} + B_{l+1}^i(N_{l}^{ix}\Phi_{l}^{i\mu} + N_{l}^{ii\mu})\; ; \; i \in \{2\ldots n\}\\
\end{multline} \vspace*{-5ex} 

\vspace*{-5ex}  \begin{multline}
\nonumber\eqref{mulg}_{\mu_{l}^m}=\eqref{Psikmum}_{k=l} \;\; \Psi_{l}^{im} =  B_{l+1}^i(N_{l}^{ix}\Phi_{l}^{m\mu} + N_{l}^{im\mu})\; ; \; i,m \in \{2\ldots n\} \; ; \; m \not = i\\
\end{multline} \vspace*{-5ex} 

\vspace*{-5ex}  \begin{multline}
\nonumber\eqref{mulg}_{const}=\eqref{psik}_{k=l} \;\; \psi_{l}^{i} =  B_{l+1}^i(N_{l}^{ix}\phi_l + n_l^i) \; ; \; i \in \{2\ldots n\}\\
\end{multline} \vspace*{-5ex}

Now it is possible to finish the inductive step of $p_l^i$ and $\lambda_l$. As a start $\eqref{xopts}_{k=l}$ and  $\eqref{muopts}_{k=l}$ are used to continue the derivation of $p_l^i$.

\vspace*{-5ex}  \begin{multline}
\nonumber  p_l^{i*} = A_{l+1}^{'}[M_{l+1}^{ix}x_{l+1}^* + \sum_{j \in 2\ldots n}M_{l+1}^{ij\mu}\mu_{l+1}^j + m_{l+1}^i - Q_{l+1}^i\tilde{x}_{l+1}^i]  \; ; \; i \in \{2\ldots n\}\\
\end{multline} \vspace*{-5ex} 

\vspace*{-5ex}  \begin{multline}
\nonumber  p_l^{i*} = A_{l+1}^{'}[M_{l+1}^{ix}(\Phi_{l}^xx_{l}^* + \sum_{j \in 2\ldots n }\Phi_{l}^{j\mu}\mu_{l}^j + \phi_l) + \sum_{j \in 2\ldots n}M_{l+1}^{ij\mu}(\Psi_{l}^{jx}x_{l}^* \\+ \sum_{m \in 2\ldots n }\Psi_{l}^{jm\mu}\mu_{l}^m + \psi_l^j) + m_{l+1}^i - Q_{l+1}^i\tilde{x}_{l+1}^i]  \; ; \; i \in \{2\ldots n\}\\
\end{multline} \vspace*{-5ex} 

The structure of the above equations justifies the following substitutions

\vspace*{-5ex}  \begin{multline}
\nonumber\eqref{psubs}_{k=l+1} \;\; p_{l}^i = (M_{l}^{ix}-Q_{l}^2)x_{l}^* + \sum_{j \in 2\ldots n}M_{l}^{ij\mu}\mu_{l}^j + m_{l}^i \; ; \; i \in \{2\ldots n\}\\
\end{multline} \vspace*{-5ex} 

\vspace*{-5ex}  \begin{multline}
(M_{l}^{ix}-Q_{l}^i)x_{l}^* + \sum_{j \in 2\ldots n}M_{l}^{ij\mu}\mu_{l}^j + m_{l}^i = A_{l+1}^{'}[M_{l+1}^{ix}(\Phi_{l}^xx_{l}^* + \sum_{j \in 2\ldots n }\Phi_{l}^{j\mu}\mu_{l}^j \\+ \phi_l) + \sum_{j \in 2\ldots n}M_{l+1}^{ij\mu}(\Psi_{l}^{jx}x_{l}^* + \sum_{m \in 2\ldots n }\Psi_{l}^{jm\mu}\mu_{l}^m + \psi_l^j) + m_{l+1}^i - Q_{l+1}^i\tilde{x}_{l+1}^i]  \; ; \; i \in \{2\ldots n\}\\\label{plg}
\end{multline} \vspace*{-5ex} 

By comparing coefficients it follows that

\vspace*{-5ex}  \begin{multline}
\nonumber\eqref{plg}_{x_{l}^*}=\eqref{Mikx}_{k=l} \;\; M_{l}^{ix} = Q_{l}^i + A_{l+1}^{'}[M_{l+1}^{ix}\Phi_{l}^x + \sum_{j \in 2\ldots n}M_{l+1}^{ij\mu}\Psi_{l}^{jx}]  \; ; \; i \in \{2\ldots n\}\\
\end{multline} \vspace*{-5ex} 

\vspace*{-5ex}  \begin{multline}
\nonumber\eqref{plg}_{\mu_{l}^m}=\eqref{Mikmu}_{k=l} \;\; M_{l}^{im\mu} = A_{l+1}^{'}[M_{l+1}^{ix}\Phi_{l}^{m\mu} + \sum_{j \in 2\ldots n}M_{l+1}^{ij\mu}\Psi_{l}^{jm\mu}] \; ; \; i,m \in \{2\ldots n\}\\
\end{multline} \vspace*{-5ex} 

\vspace*{-5ex}  \begin{multline}
\nonumber\eqref{plg}_{const.}=\eqref{mik}_{k=l} \;\; m_l^i = A_{l+1}^{'}[M_{l+1}^{ix}\phi_l + \sum_{j \in 2\ldots n}M_{l+1}^{ij\mu}\psi_l^j \\+ m_{l+1}^i - Q_{l+1}^i\tilde{x}_{l+1}^i] \; ; \; i \in \{2\ldots n\}\\
\end{multline} \vspace*{-5ex}

As a next step $\eqref{xopts}_{k=l}$ is used twice and $\eqref{muopts}_{k=l}$ and \eqref{nusubs} are applied once to continue the derivation of $\lambda_l$.

\vspace*{-5ex}  \begin{multline}
\nonumber \lambda_{l} = A_{l+1}^{'}(L_{l+1}^xx_{l+1}^* + \sum_{j \in 2\ldots n}L_{l+1} ^{j\mu}\mu_{l+1}^j + l_{l+1}) -A_{l+1}^{'}Q_{l+1}^1\tilde{x}_{l+1}^1\\+ \sum_{j \in 2\ldots n}A_{l+1}^{'}Q_{l+1}^jA_{l+1}\mu_{l}^{j} + \sum_{j \in 2\ldots n}A_{l+1}^{'}Q_{l+1}^jB_{l+1}^j\nu_{l}^{j} \\
\end{multline} \vspace*{-5ex} 

\vspace*{-5ex}  \begin{multline}
\nonumber \lambda_{l} = A_{l+1}^{'}(L_{l+1}^x(\Phi_{l}^xx_{l}^* + \sum_{j \in 2\ldots n }\Phi_{l}^{j\mu}\mu_{l}^j + \phi_l) + \sum_{j \in 2\ldots n}L_{l+1} ^{j\mu}(\Psi_{l}^{jx}x_{l}^* + \sum_{m \in 2\ldots n }\Psi_{l}^{jm\mu}\mu_{l}^m + \psi_l^j) \\+ l_{l+1}) - A_{l+1}^{'}Q_{l+1}^1\tilde{x}_{l+1}^1+ \sum_{j \in 2\ldots n}A_{l+1}^{'}Q_{l+1}^jA_{l+1}\mu_{l}^{j} \\+ \sum_{j \in 2\ldots n}A_{l+1}^{'}Q_{l+1}^jB_{l+1}^j(N_{l}^{jx}x_{l+1}^* + \sum_{m \in 2\ldots n }N_{l}^{jm\mu}\mu_{l}^m + n_l^j) \\
\end{multline} \vspace*{-5ex} 

\vspace*{-5ex}  \begin{multline}
\nonumber \lambda_{l} = A_{l+1}^{'}(L_{l+1}^x(\Phi_{l}^xx_{l}^* + \sum_{j \in 2\ldots n }\Phi_{l}^{j\mu}\mu_{l}^j + \phi_l) + \sum_{j \in 2\ldots n}L_{l+1} ^{j\mu}(\Psi_{l}^{jx}x_{l}^* + \sum_{m \in 2\ldots n }\Psi_{l}^{jm\mu}\mu_{l}^m + \psi_l^j) \\+ l_{l+1}) - A_{l+1}^{'}Q_{l+1}^1\tilde{x}_{l+1}^1 + \sum_{j \in 2\ldots n}A_{l+1}^{'}Q_{l+1}^jA_{l+1}\mu_{l}^{j} \\+ \sum_{j \in 2\ldots n}A_{l+1}^{'}Q_{l+1}^jB_{l+1}^j(N_{l}^{jx}(\Phi_{l}^xx_{l}^* + \sum_{j \in 2\ldots n }\Phi_{l}^{j\mu}\mu_{l}^j + \phi_l) + \sum_{m \in 2\ldots n }N_{l}^{jm\mu}\mu_{l}^m + n_l^j) \\
\end{multline} \vspace*{-5ex} 

The structure of the above equation justifies the following substitution

\vspace*{-5ex}  \begin{multline}
\nonumber\eqref{lambdasubs}_{k=l+1} \;\; \lambda_{l} = (L_{l}^x-Q_{l}^1)x_{l}^* + \sum_{j \in 2\ldots n}L_{l} ^{j\mu}\mu_{l}^j + l_{l}\\
\end{multline} \vspace*{-5ex} 

\vspace*{-5ex}  \begin{multline}
 (L_{l}^x-Q_{l}^1)x_{l}^* + \sum_{j \in 2\ldots n}L_{l} ^{j\mu}\mu_{l}^j + l_{l} = A_{l+1}^{'}(L_{l+1}^x(\Phi_{l}^xx_{l}^* + \sum_{j \in 2\ldots n }\Phi_{l}^{j\mu}\mu_{l}^j + \phi_l) \\+ \sum_{j \in 2\ldots n}L_{l+1} ^{j\mu}(\Psi_{l}^{jx}x_{l}^* + \sum_{m \in 2\ldots n }\Psi_{l}^{jm\mu}\mu_{l}^m + \psi_l^j) + l_{l+1}) - A_{l+1}^{'}Q_{l+1}^1\tilde{x}_{l+1}^1+ \sum_{j \in 2\ldots n}A_{l+1}^{'}Q_{l+1}^jA_{l+1}\mu_{l}^{j} \\+ \sum_{j \in 2\ldots n}A_{l+1}^{'}Q_{l+1}^jB_{l+1}^j(N_{l}^{jx}(\Phi_{l}^xx_{l}^* + \sum_{m \in 2\ldots n }\Phi_{l}^{m\mu}\mu_{l}^m + \phi_l) + \sum_{m \in 2\ldots n }N_{l}^{jm\mu}\mu_{l}^m + n_l^j) \\\label{lambdagl}
\end{multline} \vspace*{-5ex} 

Comparing coefficients gives

\vspace*{-5ex}  \begin{multline}
\nonumber\eqref{lambdagl}_{x_{l}^*}=\eqref{Lkx}_{k=l} \;\; L_{l}^{x} = Q_{l}^1 + A_{l+1}^{'}[L_{l+1}^{x}\Phi_{l}^x + \sum_{j \in 2\ldots n}L_{l+1}^{j\mu}\Psi_{l}^{jx} + \sum_{j \in 2\ldots n}Q_{l+1}^jB_{l+1}^jN_{l}^{jx}\Phi_{l}^x] \\
\end{multline} \vspace*{-5ex} 

\vspace*{-5ex}  \begin{multline}
\nonumber\eqref{lambdagl}_{\mu_{l}^i}=\eqref{Lkmu}_{k=l} \;\; L_{l}^{i\mu} = A_{l+1}^{'}[L_{l+1}^{x}\Phi_{l}^{i\mu} + \sum_{j \in 2\ldots n}L_{l+1}^{j\mu}\Psi_{l}^{ji\mu} + Q_{l+1}^iA_{l+1} \\+ \sum_{j \in 2\ldots n}Q_{l+1}^jB_{l+1}^j(N_{l}^{jx}\Phi_{l}^{i\mu} + N_{l}^{ji\mu})] \; ; \; i \in \{2\ldots n\}\\
\end{multline} \vspace*{-5ex} 

\vspace*{-5ex}  \begin{multline}
\nonumber\eqref{lambdagl}_{const.}=\eqref{lk}_{k=l} \;\; l_l = A_{l+1}^{'}[L_{l+1}^{x}\phi_l + \sum_{j \in 2\ldots n}L_{l+1}^{j\mu}\psi_l^j + l_{l+1} \\- Q_{l+1}^1\tilde{x}_{l+1}^1 + \sum_{j \in 2\ldots n}Q_{l+1}^jB_{l+1}^j(N_{l}^{jx}\phi_l + n_l^j)] \\
\end{multline} \vspace*{-5ex}

At this point the inductive step and hence the induction argument is completed. But we will try to transform $u_{l+1}^{1*}, \ldots,u_{l+1}^{n*}$ so that their evolution depends affinely on ($x_l^* , \mu_l^2 , \ldots , \mu_l^n$) and  so therefore their algorithmic computation is straightforward.\\

Let us start with $u_{l+1}^{1*}$ by applying $\eqref{xopts}_{k=l}$ to \eqref{u1subs}

\vspace*{-5ex}  \begin{multline}
\nonumber u_{l+1}^{1*} =  W_{l+1}^{x}x_{l+1}^* + \sum_{j \in 2\ldots n }W_{l+1}^{j\mu}\mu_{l}^j + w_{l+1} \\
\end{multline} \vspace*{-5ex} 

\vspace*{-5ex}  \begin{multline}
\nonumber u_{l+1}^{1*} =  W_{l+1}^{x}(\Phi_{l}^xx_{l}^* + \sum_{j \in 2\ldots n }\Phi_{l}^{j\mu}\mu_{l}^j + \phi_l) + \sum_{j \in 2\ldots n }W_{l+1}^{j\mu}\mu_{l}^j + w_{l+1} \\
\end{multline} \vspace*{-5ex} 

The structure of the above equation justifies the following substitution

\vspace*{-5ex}  \begin{multline}
\nonumber u_{l+1}^{1*} = P_{l+1}^{1x}x_l +  \sum_{j \in 2\ldots n }P_{l+1}^{1j\mu}\mu_l^j + \alpha_{l+1}^1\\
\end{multline} \vspace*{-5ex} 

\vspace*{-5ex}  \begin{multline}
P_{l+1}^{1x}x_l +  \sum_{j \in 2\ldots n }P_{l+1}^{1j\mu}\mu_l^j + \alpha_{l+1}^1 =  W_{l+1}^{x}(\Phi_{l}^xx_{l}^* \\+ \sum_{j \in 2\ldots n }\Phi_{l}^{j\mu}\mu_{l}^j + \phi_l) + \sum_{j \in 2\ldots n }W_{l+1}^{j\mu}\mu_{l}^j + w_{l+1} \\\label{u1ev}
\end{multline} \vspace*{-5ex}

By comparing coefficients it follows that

\vspace*{-5ex}  \begin{multline}
\nonumber\eqref{u1ev}_{x_{l}^*}=\eqref{Pk1x}_{k=l} \;\; P_{l+1}^{1x} = W_{l+1}^{x}\Phi_{l}^x\\
\end{multline} \vspace*{-5ex} 

\vspace*{-5ex}  \begin{multline}
\nonumber\eqref{u1ev}_{\mu_{l}^i}=\eqref{Pk1mu}_{k=l} \;\; P_{l+1}^{1i\mu} = W_{l+1}^{x}\Phi_{l}^{i\mu} + W_{l+1}^{i\mu} \; ; \; i \in \{2\ldots n\}\\
\end{multline} \vspace*{-5ex} 

\vspace*{-5ex}  \begin{multline}
\nonumber\eqref{u1ev}_{const.}=\eqref{alphak1}_{k=l} \;\; \alpha_{l+1}^{1} = W_{l+1}^{x}\phi_l + w_{l+1}\\
\end{multline} \vspace*{-5ex}

Finally $\eqref{xopts}_{k=l}$ is used in \eqref{uisubs}.

\vspace*{-5ex}  \begin{multline}
\nonumber u_{l+1}^{i*} =  T_{l+1}^{ix}x_{l+1}^* + \sum_{j \in 2\ldots n }T_{l+1}^{ij\mu}\mu_{l}^j + t_{l+1}^i   \; ; \; i \in \{2\ldots n\}\\
\end{multline} \vspace*{-5ex} 

\vspace*{-5ex}  \begin{multline}
\nonumber u_{l+1}^{i*} =  T_{l+1}^{ix}(\Phi_{l}^xx_{l}^* + \sum_{j \in 2\ldots n }\Phi_{l}^{j\mu}\mu_{l}^j + \phi_l) + \sum_{j \in 2\ldots n }T_{l+1}^{ij\mu}\mu_{l}^j + t_{l+1}^i   \; ; \; i \in \{2\ldots n\}\\
\end{multline} \vspace*{-5ex} 

The structure of the above equations justifies the following substitutions

\vspace*{-5ex}  \begin{multline}
\nonumber u_{l+1}^{i*} = P_{l+1}^{ix}x_l +  \sum_{j \in 2\ldots n }P_{l+1}^{ij\mu}\mu_l^j + \alpha_{l+1}^i  \; ; \; i \in \{2\ldots n\}\\
\end{multline} \vspace*{-5ex} 

\vspace*{-5ex}  \begin{multline}
P_{l+1}^{ix}x_l +  \sum_{j \in 2\ldots n }P_{l+1}^{ij\mu}\mu_l^j + \alpha_{l+1}^i =  T_{l+1}^{ix}(\Phi_{l}^xx_{l}^* \\+ \sum_{j \in 2\ldots n }\Phi_{l}^{j\mu}\mu_{l}^j + \phi_l) + \sum_{j \in 2\ldots n }T_{l+1}^{ij\mu}\mu_{l}^j + t_{l+1}^i   \; ; \; i \in \{2\ldots n\}\\\label{uiev}
\end{multline} \vspace*{-5ex} 

Comparing coefficients gives

\vspace*{-5ex}  \begin{multline}
\nonumber\eqref{uiev}_{x_{l}^*}=\eqref{Pkix}_{k=l} \;\; P_{l+1}^{ix} = T_{l+1}^{ix}\Phi_{l}^x   \; ; \; i \in \{2\ldots n\}\\
\end{multline} \vspace*{-5ex} 

\vspace*{-5ex}  \begin{multline}
\nonumber\eqref{uiev}_{\mu_{l}^m}=\eqref{Pkimu}_{k=l} \;\; P_{l+1}^{im\mu} = T_{l+1}^{ix}\Phi_{l}^{m\mu} + T_{l+1}^{im\mu} \; ; \; i,m \in \{2\ldots n\}\\
\end{multline} \vspace*{-5ex} 

\vspace*{-5ex}  \begin{multline}
\nonumber\eqref{uiev}_{const.}=\eqref{alphaki}_{k=l} \;\; \alpha_{l+1}^{i} = T_{l+1}^{ix}\phi_l + t_{l+1}^i   \; ; \; i \in \{2\ldots n\} \;\;\;\;\;\;\;\;\;\;\;\;\;\;\;\;\boxed{}\\
\end{multline} \vspace*{-5ex}

\begin{rem}
To solve the Stackelberg game algorithmically, the following order of application of the equations of Theorem \eqref{OLS1i} is advisable ($i,j \in \{2, \ldots, n\}$):
\begin{enumerate}
\item $M_{T}^{ix}$, $M_{T}^{ij\mu}$, $m_{T}^{i}$ 
\item $L_{T}^{x}$, $L_{T}^{i\mu}$, $l_{T}$ 
\item For k running backward from T-1 to 0
\begin{enumerate}
\item $N_k^{ix}$, $N_k^{ij\mu}$ and $n_k^{i}$
\item $T_{k+1}^{ix}$, $T_{k+1}^{ij\mu}$, $t_{k+1}^{i}$ 
\item $W_{k+1}^{x}$, $W_{k+1}^{i\mu}$, $w_{k+1}$ 
\item $\Phi_k^x$, $\Phi_k^{i\mu}$, $\phi_k$
\item $\Psi_k^{ix}$, $\Psi_k^{ij\mu}$, $\psi_k^i$
\item $M_{k}^{ix}$, $M_{k}^{ij\mu}$, $m_{k}^{i}$ 
\item $L_{k}^{x}$, $L_{k}^{i\mu}$, $l_{k}$ 
\end{enumerate}
\item $x_{0}^*$, $\mu_{0}^{i}$
\item For k running forward from 1 to T
\begin{enumerate}
\item $P_{k}^{1x}$, $P_{k}^{1i\mu}$, $\alpha_{k}^1$ 
\item $P_{k}^{ix}$, $P_{k}^{ij\mu}$, $\alpha_{k}^i$ 
\item $u_{k}^{1*}$, $u_{k}^{i*}$
\item $x_{k}^*$, $\mu_{k}^{i}$
\item $g_{k}^1(x_k,u_k^1,\ldots,u_k^n,x_{k-1})$,
  $g_{k}^i(x_k,u_k^1,\ldots,u_k^n,x_{k-1})$ 
\end{enumerate}
\item $L^1(x_0,u^{1},\ldots,u^{n})$, $L^i(x_0,u^{1},\ldots,u^{n})$
\end{enumerate}
\end{rem}

\label{OLS_ss1i}

\clearpage

\subsection{Special case: "One-induction" linear-quadratic games with one leader and one follower}

In this subsection the results of the previous subsection \eqref{OLS_ss1i} are specialized to a linear-quadratic 2-person game to allow comparison with Corollary \eqref{lqfbs2p} and Corollary \eqref{scols2i} and to point out that the number and length of the equations of the game grow rapidly with the number of followers and the consideration of constant terms.

\begin{cor}
A 2-person linear-quadratic dynamic game (cf. Def. \eqref{defspielaq}) admits a \emph{unique open-loop Stackelberg equilibrium solution with one leader and  one follower} if  
\begin{itemize}
\item $Q_k^i \geq 0$, $R_k^{ii} > 0$ (defined for $k \in K$ , $i \in N$).
\item $(I - B_{k+1}^1W_{k+1}^x - B_{k+1}^2T_{k+1}^x)^{-1}$ and $(B_{k+1}^{2'}(Q_{k+1}^2 + L_{k+1}^{\mu})B_{k+1}^2 + I - R_{k+1}^{12}B_{k+1}^{2'}M_{k+1}^{\mu}B_{k+1}^2)^{-1}$ (defined for $k \in K$) exist.
\end{itemize}
If these conditions are satisfied, the unique equilibrium strategies $\gamma_{k+1}^{i*}(x_0)$ are given by \eqref{gammaols2}, where the associated state trajectory $x_{k+1}^*$ is given by \eqref{xopts2}.\footnote {For all equations belonging to this corollary and its proof, $k \in \{0 , \ldots ,T-1\}$ if nothing different is stated.}

\begin{align}
f_{k-1}(x_{k-1},u_k^1,u_k^2) = A_kx_{k-1} + B_k^1u_k^1 + B_k^2u_k^2 \; ; \; k \in K\label{fols2}
\end{align}

\begin{align}
L^i(x_0,u^{1}, u^2)= \sum_{k = 1}^T g_k^i(x_k,u_k^{1}, u_k^2 ,x_{k-1}) 
\end{align}

\vspace*{-5ex}  \begin{multline}
 g_k^i(x_k,u_k^1,\ldots,u_k^n,x_{k-1}) = \frac{1}{2} (x_k^{'}Q_k^ix_k + u_k^{i'}u_k^i \\+ u_k^{j'}R_k^{ij}u_k^j) \; ; \; k \in K\; ; \; i,j \in \{1,2\} \; ; \; i \not = j\\\label{gols2}
\end{multline} \vspace*{-5ex}

\begin{align}
 x_{k+1}^* = \Phi_{k}^xx_{k}^* + \Phi_{k}^{\mu}\mu_{k} \; ; \; x_0^* = x_0 \label{xopts2}
\end{align}

\begin{align}
 \Phi_{k}^x = (I - B_{k+1}^1W_{k+1}^x - B_{k+1}^2T_{k+1}^x)^{-1}A_{k+1}
\end{align}

\begin{align}
  \Phi_{k}^{\mu} = (I - B_{k+1}^1W_{k+1}^x - B_{k+1}^2T_{k+1}^x)^{-1}(B_{k+1}^1W_{k+1}^{\mu} + B_{k+1}^2T_{k+1}^{\mu})
\end{align}

\begin{align}
  \mu_{k+1}  =  \Psi_{k}^xx_{k}^* + \Psi_{k}^{\mu}\mu_{k} \; ; \; \mu_0 = 0  
\end{align}

\begin{align}
  \Psi_{k}^x = B_{k+1}^2N_{k}^x\Phi_{k}^x 
\end{align}

\begin{align}
  \Psi_{k}^{\mu} = A_{k+1} + B_{k+1}^2(N_{k}^x\Phi_{k}^{\mu} + N_{k}^{\mu}) 
\end{align}

\begin{align}
  W_{k+1}^x =   -B_{k+1}^{1'}L_{k+1}^x - B_{k+1}^{1'}L_{k+1}^{\mu}N_k^x - B_{k+1}^{1'}Q_{k+1}^2B_{k+1}^2N_{k}^x
\end{align}

\vspace*{-5ex}  \begin{multline}
  W_{k+1}^{\mu} = - B_{k+1}^{1'}L_{k+1}^{\mu}(A_{k+1} + B_{k+1}^2N_k^{\mu}) \\- B_{k+1}^{1'}Q_{k+1}^2A_{k+1} - B_{k+1}^{1'}Q_{k+1}^2B_{k+1}^2N_{k}^{\mu}\\
\end{multline} \vspace*{-5ex}

\begin{align}
  T_{k+1}^x =  -B_{k+1}^{2'}(M_{k+1}^x +  M_{k+1}^{\mu}B_{k+1}^2N_{k}^x)
\end{align}

\begin{align}
  T_{k+1}^{\mu} = -B_{k+1}^{2'}M_{k+1}^{\mu}(A_{k+1} + B_{k+1}^2N_{k}^{\mu})
\end{align}

\vspace*{-5ex}  \begin{multline}
  N_{k}^x = -(B_{k+1}^{2'}(Q_{k+1}^2 + L_{k+1}^{\mu})B_{k+1}^2 + I - R_{k+1}^{12}B_{k+1}^{2'}M_{k+1}^{\mu}B_{k+1}^2)^{-1}\\(B_{k+1}^{2'}L_{k+1}^x - R_{k+1}^{12}B_{k+1}^{2'}M_{k+1}^x) \\
\end{multline} \vspace*{-5ex} 

\vspace*{-5ex}  \begin{multline}
  N_{k}^{\mu} = -(B_{k+1}^{2'}(Q_{k+1}^2 + L_{k+1}^{\mu})B_{k+1}^2 + I - R_{k+1}^{12}B_{k+1}^{2'}M_{k+1}^{\mu}B_{k+1}^2)^{-1}\\ (B_{k+1}^{2'}(Q_{k+1}^2 + L_{k+1}^{\mu}) - R_{k+1}^{12}B_{k+1}^{2'}M_{k+1}^{\mu})A_{k+1}\\
\end{multline} \vspace*{-5ex}

%\vspace*{-5ex}  \begin{multline}
%  p_{k} = (M_{k}^x-Q_{k}^2)x_{k}^* + M_k^{\mu}\mu_k \\
%\end{multline} \vspace*{-5ex} 

\begin{align}
  M_k^x = Q_k^2 + A_{k+1}^{'}[M_{k+1}^x\Phi_{k}^x + M_{k+1}^{\mu}\Psi_{k}^x] \; ; \; M_T^x = Q_{T}^2 
\end{align}

\begin{align}
  M_k^{\mu} = A_{k+1}^{'}[M_{k+1}^x\Phi_{k}^{\mu} + M_{k+1}^{\mu}\Psi_{k}^{\mu}] \; ; \; M_T^{\mu} = 0
\end{align}

%\vspace*{-5ex}  \begin{multline}
%  \lambda_{k} = (L_{k}^x-Q_{k}^1)x_{k}^* + L_k ^{\mu}\mu_k  \\
%\end{multline} \vspace*{-5ex} 

\vspace*{-5ex}  \begin{align}
L_k^x = Q_k^1 + A_{k+1}^{'}L_{k+1}^x\Psi_{k}^x + A_{k+1}^{'}L_{k+1}^{\mu}\Psi_{k}^x + A_{k+1}^{'}Q_{k+1}^2B_{k+1}^2N_{k}^x\Psi_{k}^x \; ; \; L_T^x = Q_{T}^1 
\end{align}

\vspace*{-5ex}  \begin{multline}
L_k^{\mu} = A_{k+1}^{'}L_{k+1}^x\Psi_{k}^{\mu} + A_{k+1}^{'}L_{k+1}^{\mu}\Psi_{k}^{\mu} \\+ A_{k+1}^{'}Q_{k+1}^2A_{k+1}\mu_{k} + A_{k+1}^{'}Q_{k+1}^2B_{k+1}^2(N_{k}^x\Psi_{k}^{\mu} + N_{k}^{\mu}) \; ; \; L_T^{\mu} = 0 \\
\end{multline} \vspace*{-5ex}

\vspace*{-5ex}  \begin{align}
\gamma_{k+1}^{i*}(x_0) = u_{k+1}^{i*} = P_{k+1}^{ix}x_{k}^* + P_{k+1}^{i\mu}\mu_{k} \; ; \; i \in \{1,2\} \label{gammaols2}
\end{align}

\vspace*{-5ex}  \begin{align}
  P_{k+1}^{1x} = W_{k+1}^x\Phi_{k}^x
\end{align}

\vspace*{-5ex}  \begin{align}
  P_{k+1}^{1\mu} =  W_{k+1}^x\Phi_{k}^{\mu} + W_{k+1}^{\mu} 
\end{align}

\vspace*{-5ex}  \begin{align}
  P_{k+1}^{2x} = T_{k+1}^x\Phi_{k}^x
\end{align}

\vspace*{-5ex}  \begin{align}
  P_{k+1}^{2\mu} = T_{k+1}^x\Phi_{k}^{\mu} + T_{k+1}^{\mu} 
\end{align}\\\\

\label{scols1i}
\end{cor}

\textsc {Proof:}

Corollary \eqref{scols1i} is proven in the same way as Theorem \eqref{OLS1i} taking into consideration simplifications resulting from the different number of followers and the modified state equation and cost functionals.\;\;\;\;\;\;\;\boxed{}\\

\begin{rem}
Special attention should be paid to the fact that the assumption about the existence of unique solutions of the systems of equations \eqref{Nkx}, \eqref{Nkmu} and \eqref{nk} in Theorem \eqref{OLS1i} is equivalent to assuming the existence of $(B_{k+1}^{2'}(Q_{k+1}^2 + L_{k+1}^{\mu})B_{k+1}^2 + I - R_{k+1}^{12}B_{k+1}^{2'}M_{k+1}^{\mu}B_{k+1}^2)^{-1}$ in this special case.
\end{rem}

\label{OLS_1i}

\clearpage

\setcounter{chapter}{4}
\setcounter{fussall}{\thefootnote}
\chapter {Conclusion}
\setcounter{footnote}{\thefussall}

Beside some corrections of results for \textit{n}-person discrete-time affine-quadratic dynamic games of prespecified fixed duration with open-loop and feedback information patterns in the literature, extensions were presented for the open-loop and feedback Stackelberg equilibrium solutions of \textit{n}-person discrete-time affine-quadratic dynamic games of prespecified fixed duration concerning the number of followers and the possibility of algorithmic disintegration. These extensions enable a better modeling and numerical solution of real-world applications characterized by a hierarchical structure of the players' interactions.\\

Extending the number of leaders from 1 to \textit{n} for open-loop and feedback Stackelberg discrete-time affine-quadratic dynamic games of prespecified fixed duration and finding an interrelation between the assumptions for the unique existence of Nash and Stackelberg equilibrium solutions  for affine-quadratic dynamic games with open-loop and feedback information patterns and the matrices defining the affine-quadratic dynamic games are challenging tasks for future research.

\clearpage
% Anhang
%
%\backmatter

\bibliographystyle{acm}
\bibliography{DP_lit}
%\appendix
%\input{DP_lIT}
%\input{dp_erkl}

\end{document}